\documentclass{article}

\usepackage{amsfonts, amsmath, amssymb, amsthm, amscd}
\usepackage{ascmac}

\usepackage{verbatim}

\usepackage{graphics}
\usepackage[all,2cell]{xypic}
\usepackage{mathptmx}

\usepackage{helvet}
\usepackage{courier}        
\usepackage{type1cm} 

\objectmargin+{1mm}
\labelmargin+{0.8mm}
\SelectTips{cm}{12}

\theoremstyle{plain}  
\newtheorem{thm}{Theorem}[subsection]

\newtheorem{cor}[thm]{Corollary}

\theoremstyle{definition}

\newtheorem{para}[thm]{}

\theoremstyle{remark}

\DeclareMathOperator{\cA}{\mathcal{A}}
\DeclareMathOperator{\cB}{\mathcal{B}}
\DeclareMathOperator{\cC}{\mathcal{C}}
\DeclareMathOperator{\calD}{\mathcal{D}}
\DeclareMathOperator{\cE}{\mathcal{E}}
\DeclareMathOperator{\cF}{\mathcal{F}}
\DeclareMathOperator{\cG}{\mathcal{G}}
\DeclareMathOperator{\calH}{\mathcal{H}}
\DeclareMathOperator{\cI}{\mathcal{I}}

\DeclareMathOperator{\calL}{\mathcal{L}}
\DeclareMathOperator{\cM}{\mathcal{M}}
\DeclareMathOperator{\cN}{\mathcal{N}}
\DeclareMathOperator{\cO}{\mathcal{O}}
\DeclareMathOperator{\cP}{\mathcal{P}}
\DeclareMathOperator{\cQ}{\mathcal{Q}}
\DeclareMathOperator{\calR}{\mathcal{R}}
\DeclareMathOperator{\cS}{\mathcal{S}}
\DeclareMathOperator{\cT}{\mathcal{T}}
\DeclareMathOperator{\cU}{\mathcal{U}}

\DeclareMathOperator{\cW}{\mathcal{W}}
\DeclareMathOperator{\cX}{\mathcal{X}}
\DeclareMathOperator{\cY}{\mathcal{Y}}
\DeclareMathOperator{\cZ}{\mathcal{Z}}

\DeclareMathOperator{\bC}{\mathbf{C}}

\DeclareMathOperator{\bE}{\mathbf{E}}
\DeclareMathOperator{\bF}{\mathbf{F}}
\DeclareMathOperator{\bG}{\mathbf{G}}

\DeclareMathOperator{\bX}{\mathbf{X}}
\DeclareMathOperator{\bY}{\mathbf{Y}}
\DeclareMathOperator{\bZ}{\mathbf{Z}}

\DeclareMathOperator{\bbD}{\mathbb{D}}

\DeclareMathOperator{\bbK}{\mathbb{K}}

\DeclareMathOperator{\bbZ}{\mathbb{Z}}

\DeclareMathOperator{\fc}{\mathfrak{c}}

\DeclareMathOperator{\fd}{\mathfrak{d}}

\DeclareMathOperator{\fe}{\mathfrak{e}}

\DeclareMathOperator{\fG}{\mathfrak{G}}

\DeclareMathOperator{\fH}{\mathfrak{H}}

\DeclareMathOperator{\fraki}{\mathfrak{i}}

\DeclareMathOperator{\fj}{\mathfrak{j}}

\DeclareMathOperator{\fL}{\mathfrak{L}}

\DeclareMathOperator{\fM}{\mathfrak{M}}

\DeclareMathOperator{\fQ}{\mathfrak{Q}}

\DeclareMathOperator{\fr}{\mathfrak{r}}

\DeclareMathOperator{\fs}{\mathfrak{s}}

\DeclareMathOperator{\fu}{\mathfrak{u}}

\DeclareMathOperator{\fv}{\mathfrak{v}}

\DeclareMathOperator{\fX}{\mathfrak{X}}

\UseAllTwocells

\def\sn{\smallskip\noindent}
\def\mn{\medskip\noindent}

\def\enumidef{\renewcommand{\labelenumi}{$\mathrm{(\arabic{enumi})}$}}

\newcommand{\Cat}{\operatorname{\bf Cat}}
\newcommand{\cf}{\textrm{cf.}\;}
\newcommand{\Ch}{\operatorname{\bf Ch}}

\newcommand{\Comp}{\operatorname{\bf Comp}}
\newcommand{\CompEx}{\operatorname{\bf CompEx}}
\newcommand{\Cone}{\operatorname{Cone}}
\newcommand{\colim}{\operatorname{colim}}
\newcommand{\coker}{\operatorname{Coker}}

\newcommand{\Cyl}{\operatorname{Cyl}}
\newcommand{\CW}{\operatorname{CW}}

\newcommand{\deq}{\operatorname{deq}}

\newcommand{\dom}{\operatorname{dom}}

\newcommand{\Ex}{\operatorname{Ex}}
\newcommand{\ExCat}{\operatorname{\bf ExCat}}

\newcommand{\frob}{\operatorname{frob}}

\newcommand{\heq}{\operatorname{heq}}
\newcommand{\hFib}{\operatorname{Hfib}}
\newcommand{\Ho}{\operatorname{Ho}}

\newcommand{\hocolim}{\operatorname{hocolim}}
\newcommand{\HOM}{\mathcal{H}\!\!om}
\newcommand{\Hom}{\operatorname{Hom}}

\newcommand{\id}{\operatorname{id}}
\newcommand{\im}{\operatorname{Im}}
\newcommand{\inn}[2]{{#1}\ \operatorname{ in }\ {#2}}

\newcommand{\Isom}{\operatorname{isom}}
\newcommand{\isoto}{\overset{\scriptstyle{\sim}}{\to}}

\newcommand{\Ker}{\operatorname{Ker}}

\newcommand{\length}{\operatorname{length}}
\newcommand{\Lex}{\operatorname{\bf Lex}}

\newcommand{\linf}{\leftarrowtail}

\newcommand{\lv}{\operatorname{lv}}

\newcommand{\Mor}{\operatorname{Mor}}

\newcommand{\NC}{\operatorname{NC}}

\newcommand{\nor}{\operatorname{nor}}
\newcommand{\nullclass}{\operatorname{null}}

\newcommand{\Ob}{\operatorname{Ob}}
\newcommand{\onto}[1]{\stackrel{#1}{\to}}
\newcommand{\op}{\operatorname{op}}

\newcommand{\Pat}{\operatorname{Pat}}

\newcommand{\prenull}{\operatorname{prenull}}

\newcommand{\qfis}{\operatorname{q^fis}}
\newcommand{\qis}{\operatorname{qis}}

\newcommand{\ran}{\operatorname{ran}}
\newcommand{\rdef}{\twoheadrightarrow}
\newcommand{\RelCat}{\operatorname{\bf RelCat}}
\newcommand{\RelComp}{\operatorname{\bf RelComp}}
\newcommand{\RelEx}{\operatorname{\bf RelEx}}

\newcommand{\rinc}{\hookrightarrow}
\newcommand{\rinf}{\rightarrowtail}
\newcommand{\rlto}[2]{\overset{#1}{\underset{#2}{\rightleftarrows}}}

\newcommand{\seminull}{\operatorname{semi-null}}

\newcommand{\st}{\operatorname{st}}
\newcommand{\standard}{\operatorname{sta}}

\newcommand{\stnor}{\operatorname{sn}}

\newcommand{\thi}{\operatorname{thi}}
\newcommand{\Tot}{\operatorname{Tot}}

\newcommand{\Tri}{\operatorname{Tri}}

\newcommand{\TriCat}{\operatorname{\bf TriCat}}

\title{A d\'evissage theorem of non-connective $K$-theory}
\date{}

\author{Satoshi Mochizuki}

\begin{document}

\maketitle

\tableofcontents

\section*{Introduction}

The purpose of this article is to show a version of d\'evissage theorem 
of non-connective $K$-theory. 
As a mathematical jargon, 
the term `d\'evissage' is introduced by Grothendieck. 
It is a technique for verifying assertions about coherent sheaves 
on noetherian schemes \cite[Th\'eor\`eme 3.1.2]{DG61} 
and in the mid-80's he returned back to this theme in the context of 
`tame topology' \cite[\S 5]{Gro84}. 
This argument is imported into calculation of lower algebraic $K$-groups 
and abstracted by Heller \cite{BHS64}, \cite{Hel65}. 
In the renowned paper \cite{Qui73}, 
Quillen provided a general d\'evissage theorem for higher $K$-theory 
of abelian categories. 
Namely for an essentially small abelian category $\cA$ 
and a topologizing subcategory (\ref{df:toplogizing subcat}) $\cB$ of $\cA$ 
such that the inclusion functor $\cB\rinc \cA$ satisfies Quillen's  
d\'evissage condition (\ref{df:classical Devissage condition}), 
the inclusion functor 
induces a homotopy equivalence $K(\cB)\to K(\cA)$ on $K$-theory. 
It might be well-known that 
for an exact category which satisfies a suitable Jordan-H\"older 
theorem, 
the proof of d\'evissage theorem in \cite{Qui73} with 
appropriate refinements still works well. 
See for example \cite[p.175, l.5]{Mit94}. 
Notice that 
there exists a d\'evissage theorem for $K$-theory of exact $\infty$-categories 
by Barwick in \cite{Bar13}. 
These approach to a d\'evissage theorem crucially depends upon 
(a generalization of) Quillen's $Q$-construction. 
Especially in the situation above, we analyze the homotopy fiber of 
$Q\cB \rinc Q\cA$ by utilizing a (generalized) Jordan-H\"older theorem. 
We will illustrate this phenomena in the context of appropriate 
Waldhausen $K$-theory in Appendix.

In \cite[p.188]{Wal84}, Waldhausen wrote:

\sn
`{\it
As to a general attack on the spaces $\operatorname{M}(\ast,p,n)$, 
the first {\rm(}and perhaps main{\rm)} step should be the search for a d\'evissage theorem.}'

\sn
In this scheme, he actually developed the cell filtration theorem 
in \cite[1.7.1]{Wal85} 
(for precise statement for non-connective $K$-theory, see 
Corollary~\ref{cor:Waldhausen's cell filtration theorem}). 
Even though its statement is apparently different from 
Quillen's d\'evissage theorem, 
in this article, we will regard it as a version of 
d\'evissage theorem. 
In line with this thinking, 
there also exists similar flavoured approaches for d\'evissage theorem by Yao 
in \cite{Yao95} and by Blumberg and Mandell in \cite{BM08} and by 
Barwick \cite{Bar15} 
and recently Raptis \cite{Rap18}. 
These investigations are close to theorem of heart originated with 
Neeman in series of papers 
\cite{Nee98}, \cite{Nee99}, \cite{Nee01a} 
(for recent development about theorem of heart into $K$-theory 
of infinite categories, 
see \cite{Sos17}, \cite{Fon18}, \cite{Rap18}, \cite{Bar15}). 
The theorem says $K$-theory of a triangulated category 
is equivalent to $K$-theory of its heart 
(for our version, see Corollary~\ref{cor:Theorem of heart}). 

Recall that 
in the celebrated paper written by 
Thomason (and Trobaugh) \cite[1.11.1]{TT90}, 
Thomason proposed an open problem;

\sn 
`{\it Find a general result for Waldhausen categories 
that specialises to Quillen's d\'evissage theorem when applied to 
the category of bounded complexes in an abelian category.}'

\sn
The d\'evissage theorem~\ref{thm:divissage for complicial exact categories} 
which we will establish in this article contains 
Waldhausen's cell filtration 
theorem~\ref{cor:Waldhausen's cell filtration theorem}, 
theorem of heart~\ref{cor:Theorem of heart} 
and Quillen's d\'evissage theorem~\ref{cor:Quillen's devissage theorem} 
as special cases. 
Namely, we will give an affirmative answer to 
Waldhausen's and Thomason's problems above in some sense. 
We will turn in the notions of cell structures 
(\ref{df:cell str}) and 
d\'evissage spaces 
(\ref{df:Derived devissage condition}) 
in Section~\ref{sec:devissage spaces} and our 
d\'evissage theorem states a structure of non-connective $\bbK$-theory 
of d\'evissage spaces in terms of non-connective $\bbK$-theory of 
heart of cell structures.  
The specific feature in our proof is `motivic' in the sense that 
properties of $K$-theory which we will utilize to prove the theorem 
is only categorical homotopy invariance, localization and cocontinuity 
(see Corollary~\ref{cor:devissage for localization theory}). 
On the other hands, the analogue of the d\'evissage theorem for 
$K$-theory does not hold for 
Hochschild homology theory (see \cite[1.11]{Kel99}). 
In this point of view, we could say that d\'evissage theorem is 
not `motivic' over dg-categories. 
To overcome this dilemma, 
the notion of d\'evissage spaces should not be 
expressed by the language of dg-categories. 
From Section~\ref{sec:complicial exact cat} to 
Section~\ref{sec:hom th of rel comp ex cat} 
are devoted to the foundation of 
our model of stable $(\infty,1)$-categories which we will play on to 
give a description of d\'evissage spaces. 

In a standpoint of geometry over categories 
(for example, topoi theory or non-commutative geometry of abelian categories or dg-categories), 
a d\'evissage condition is related to nilpotent immersions 
(see \ref{df:classical Devissage condition}). 
In the future work, the author plans to clarify the meaning of 
a d\'evissage condition in a perspective of noncommutative motive theory. 

\section{Complicial exact categories}
\label{sec:complicial exact cat}

The conception of complicial biWaldhausen categories 
(closed under the formations of the canonical homotopy push-outs and 
the canonical homotopy pull-backs) 
is 
introduced by Thomason in \cite{TT90}, 
which is 
a special class of Waldhausen categories 
whose underlying categories are full subcategories of 
categories of chain complexes on abelian categories. 
The notion is abstracted and further studied 
by Schlichting in the survey article 
of algebraic $K$-theory \cite{Sch11} under the name of 
{\it complicial exact categories} ({\it with weak equivalences}), 
whose underlying exact categories are equipped with 
monoidal actions of the symmetric monoidal category of 
bounded chain complexes of finitely generated free abelian groups 
(there exists a similar concept which is called 
{\it $C(k)$-model categories} and discussed in \cite[\S 2.1]{Toe11} and 
there exists a correlative notion which is designated as 
cat\'egories d\'erivables and established in \cite{Cis10}. 
See also \ref{para:derivable Waldhausen category}). 
In the series of papers \cite{Moc10} and \cite{Moc13b} 
under the name of {\it bicomplicial categories} 
(or {\it bicomplicial pairs}), 
we discuss the possibility of simplifying Schlichting's axioms to establish the theory of complicial exact categories with weak equivalences in \cite{Sch11}. 
In this section, following the papers \cite{Moc10} and \cite{Moc13b}, 
we recall the notion of complicial exact categories 
(with slightly different conventions) from Ibid. and study them further. 
Now we give a guide for the structure of this section. 
In the first subsection \ref{subsec:locally exact 2-cat}, 
we start by recalling the notion of exact categories. 
In the next subsection \ref{subsec:complicial structures}, 
we will introduce the notion of complicial objects 
in a locally exact $2$-category. 
Practical examples of complicial objects are 
normal ordinary complicial exact categories and 
ordinary complicial structure will be explained in 
subsection \ref{subsec:ord comp str}. 
In the final subsection \ref{subsec:Frob exact str}, 
we will discuss Frobenius complicial exact structure.

\subsection{Locally exact $2$-categories}
\label{subsec:locally exact 2-cat}

In this subsection, we recall the notion of exact categories 
in the sense of Quillen. 
In particular we study $2$-category $\underline{\ExCat}$ of 
small exact categories and the $2$-category ${\underline{\ExCat}}^{\cI}$ of 
$\cI$-diagrams in $\underline{\ExCat}$ for a small category $\cI$. 
The Hom categories in the both categories are equipped with the natural 
exact structures but the both categories are not 
enriched over the category of small exact categories for the reason that 
the composition functors $\underline{\HOM}_{\ExCat}(x,y)\times
\underline{\HOM}_{\ExCat}(y,z)\to\underline{\HOM}_{\ExCat}(x,z)$ 
are not exact. 
Thus as a substitutive concept, we introduce the notion of 
{\it locally exact $2$-categories} 
(see \ref{df:locally exact 2-categories}) which contains 
$\underline{\ExCat}$ and ${\underline{\ExCat}}^{\cI}$ 
for a small category $\cI$ 
(see \ref{nt:diagram cat over locally exact 2-cat}) as typical examples. 
In the next subsection, this notion will make our treatments of complicial structures simplify.

\begin{para}
\label{nt:exact categories}
{\bf (Exact categories).}\ \ 
Basically, for the conventions of {\it exact categories}, 
we follows the notations in \cite{Qui73}. 
Recall that a functor between exact categories 
$f\colon\cE\to\cF$ {\it reflects exactness} 
if for a sequence $x\to y\to  z$ in $\cE$ 
such that $fx\to fy\to fz$ is an admissible exact sequence 
in $\cF$, 
$x\to y\to z$ is an admissible exact sequence 
in $\cE$. 
For an exact category $\cE$, 
we say that its full subcategory $\cF$ is 
an {\it exact subcategory} 
if it is an exact category and 
the inclusion functor 
$\cF\rinc \cE$ is exact 
and say that $\cF$ is a {\it strict exact subcategory} if it 
is an exact subcategory and 
moreover the inclusion functor reflects exactness. 
We say that $\cF$ is an {\it extension closed} ({\it full}) 
{\it subcategory} of $\cE$ or 
{\it closed under extensions} in $\cE$ if for any admissible exact 
sequence $x\rinf y\rdef z$ in $\cE$, $x$ and $z$ are isomorphic to objects 
in $\cF$ respectively, then $y$ is isomorphic to an object in $\cF$. 
\end{para}

We will sometimes use the following lemmata.

\begin{para}
\label{lem:kel90}
{\bf Lemma.}\ (\cf \cite[Step 1 in the proof of A.1]{Kel90}.)\ \  
{\it
Let $\cE$ be an exact category and let 
\begin{equation}
\label{eq:comm square}
\xymatrix{
x \ar[r]^i \ar[d]_j & y \ar[d]^q\\
z \ar[r]_p & w 
} 
\end{equation}
be a commutative square in $\cE$. 
If the square $\mathrm{(\ref{eq:comm square})}$ 
is a push-out {\rm(}resp. pull-back{\rm)} 
and the morphism $i$ {\rm(}resp. $q${\rm)} is an admissible 
monomorphism {\rm(}resp. epimorphism{\rm)}, 
then the sequence 
$x\overset{\tiny{\begin{pmatrix}i\\ -j\end{pmatrix}}}{\rinf} y\oplus z 
\overset{\tiny{\begin{pmatrix}q & p\end{pmatrix}}}{\rdef} w $ 
is an admissible exact sequence in $\cE$.
}
\qed
\end{para}

\begin{para}
\label{lem:char of admissible square}
{\bf Lemma.}\ \ 
{\it
Let $\cE$ be an exact category and let 
$$
\xymatrix{
x \ar@{>->}[r] \ar@{>->}[d] & y \ar@{>->}[d]\\
z \ar@{>->}[r] & w
}
$$
be a commutative square of admissible monomorphisms in $\cE$. 
Then the following conditions are equivalent.

\sn
$\mathrm{(1)}$ 
The induced morphism $z\sqcup_x y\to w$ is an admissible monomorphism. 

\sn
$\mathrm{(2)}$ 
The induced morphism $z/x \to w/y$ is an admissible monomorphism. 

\sn
$\mathrm{(3)}$ 
The induced morphism $y/x\to w/z$ is an admissible monomorphism. 
}
\end{para}

\begin{proof}
Assertion that $\mathrm{(1)}$ implies $\mathrm{(2)}$ and $\mathrm{(3)}$ 
is proven in \cite[1.1.1]{Wal85} and 
assertion that $\mathrm{(2)}$ or $\mathrm{(3)}$ implies $\mathrm{(1)}$ 
is proven in \cite[1.7.4]{TT90}. 
\end{proof}

\begin{para}
\label{df:2-categories}
{\bf ($2$-categories).}\ \ 
In this article, 
a {\it $2$-category} means a category enriched over 
the category of small categories. 
For a $2$-category $\cC$ and a pair of objects 
$x$ and $y$, we write $\underline{\HOM}_{\cC}(x,y)$ for 
the Hom category from $x$ to $y$. 
We call an object in $\underline{\HOM}_{\cC}(x,y)$ an {\it $1$-morphism} 
({\it from $x$ to $y$}) and call a morphism $a\colon f\to g$ 
in $\underline{\HOM}_{\cC}(x,y)$ a {\it $2$-morphism} ({\it from $f$ to $g$}). 
The compositions of $2$-morphisms in the category 
$\underline{\HOM}_{\cC}(x,y)$ is called 
the {\it vertical compositions} and denoted by the letter $\cdot$ 
and the composition of $2$-morphisms 
comes from the composition functors 
$\underline{\HOM}_{\cC}(x,y)\times\underline{\HOM}_{\cC}(y,z)\to
\underline{\HOM}_{\cC}(x,z)$ is called the {\it horizontal compositions} 
and denoted by $\ast$. 
For example let $x$, $y$ and $z$ be a triple of objects in $\cC$ and 
$f_1,\ f_2, f_3\colon x\to y$ and $g_1,\ g_2,\ g_3\colon y\to z$ 
be $1$-morphisms in $\cC$ and $\alpha\colon f_1\to f_2$, 
$\beta\colon f_2\to f_3$, $\gamma\colon g_1\to g_2$ and 
$\delta\colon g_2\to g_3$ be $2$-morphisms in $\cC$, 
then we have the equality 
\begin{equation}
\label{eq:Dedekind identity}
(\delta\ast\beta)\cdot (\gamma\ast\alpha)=
(\delta\cdot\gamma)\ast (\beta\cdot\alpha).
\end{equation}

For $1$-morphisms $f_1,\ f_2\colon x\to y$ and $g_1,\ g_2\colon y\to z$ 
and $2$-morphisms $a\colon f_1\to f_2$ and $b\colon g_1\to g_2$ 
in a $2$-category, we denote $b\ast a$ by $b\ast f_1$ 
if $f_1=f_2$ and $a=\id_{f_1}$ and 
by $g_1\ast a$ if $g_1=g_2$ and $b=\id_{g_1}$.

Let $x$ and $y$ be a pair of objects in a $2$-category $\cC$. 
We say that an $1$-morphism $f\colon x\to y$ is an {\it equivalence} 
if there exists an $1$-morphism $g\colon y\to x$ and 
a pair of $2$-isomorphisms 
\begin{equation}
\label{eq:pair of two isomorphisms}
\alpha\colon fg\isoto \id_y \text{\ \ and\ \ } 
\beta\colon \id_x\isoto gf
\end{equation}
We say that a quadruple $(f,g,\alpha,\beta)$ consisting of a pair of 
$1$-morphisms $f\colon x\to y$ and $g\colon y\to x$ and 
a pair of $2$-isomorphisms as in 
$\mathrm{(\ref{eq:pair of two isomorphisms})}$ 
in a $2$-category $\cC$ 
is an {\it adjoint equivalence} from 
$x$ to $y$ if the following equalities hold:
\begin{equation}
\label{eq:equivalence identity 1}
(\alpha\ast f)\cdot(f\ast \beta)=\id_f,
\end{equation}
\begin{equation}
\label{eq:equivalence identity 2}
(g\ast\alpha)\cdot(\beta\ast g)=\id_g.
\end{equation}
We can show that for any equivalence $f\colon x\to y$, 
there exists an $1$-morphism $g\colon y\to x$ and 
a pair of $2$-isomorphisms as in 
$\mathrm{(\ref{eq:pair of two isomorphisms})}$ 
in a $2$-category $\cC$ such that the quadruple 
$(f,g,\alpha,\beta)$ is an adjoint equivalence. 

An exemplar of $2$-categories is 
the $2$-category of small categories $\underline{\Cat}$ 
whose objects are small categories, whose $1$-morphisms are 
functors and whose $2$-morphisms are natural transformations. 
The horizontal compositions and the vertical compositions 
of natural transformations are given in the following way. 
For functors between categories, 
$f_1$, $f_2$, $f_3\colon \cX \to \cY$ and 
$g_1$, $g_2\colon \cY\to \cZ$ 
and natural transformations 
$a\colon f_1\to f_2$, $b\colon f_2 \to f_3$ and $c\colon g_1\to g_2$, 
we define the vertical composition $b$ and $a$, 
$b\cdot a\colon f_1\to f_3$ and 
the horizontal composition $c$ and $a$, 
$c\ast a\colon g_1f_1\to g_2f_2$ by the formulas
\begin{equation}
\label{eq:vertical comp df}
(b\cdot a)(x):=b(x)a(x),
\end{equation}
\begin{equation}
\label{eq:horizontal comp df}
(c\ast a)(x):=c(f_2(x))g_1(a(x))=g_2(a(x))c(f_1(x))
\end{equation}
for any object $x$ in $\cX$. 
\end{para}

\begin{para}
\label{df:locally exact 2-categories}
{\bf Definition (Locally exact $2$-categories).}\ \ 
A {\it local exact $2$-category} is a $2$-category $\cC$ such that 
for any pair of objects $x$ and $y$ in $\cC$, the Hom category 
$\underline{\HOM}_{\cC}(x,y)$ is equipped with an exact category structure 
which subjects to the condition that 
for any triple of objects $x$, $y$, $z$ in $\cC$ and 
a $1$-morphisms $f\colon x\to y$, the induced functors 
$\underline{\HOM}_{\cC}(z,f)\colon\underline{\HOM}_{\cC}(z,x)\to\underline{\HOM}_{\cC}(z,y)$ and $\underline{\HOM}_{\cC}(f,z)\colon
\underline{\HOM}_{\cC}(y,z)\to\underline{\HOM}_{\cC}(x,z)$ 
are exact functors. 
\end{para}

In \ref{nt:ExCat} and \ref{nt:diagram cat over locally exact 2-cat}, 
we will provide representatives of locally exact $2$-categories. 
First we consider general remarks of diagram categories 
over small exact categories. 

\begin{para}
\label{nt:functor category}
{\bf (Functors categories).}\ \ 
Let $\cE$ be an exact category and let $\cI$ be a category. 
We say that a sequence 
$x\onto{f}y\onto{g}z$ of morphisms in $\cE^{\cI}$ 
the functor category from $\cI$ to $\cE$ 
is a {\it level admissible exact sequence} if 
$x(i)\onto{f(i)}y(i)\onto{g(i)}z(i)$ 
is an admissible exact sequence 
in $\cE$ for any object $i$ in $\cI$. 
We assume that 
$\cI$ is essentially small. 
Then $\cE^{\cI}$ with level admissible exact sequences is an exact category. 
Let $f\colon \cE\to \cE'$ be an exact functor between exact categories. 
We define $f^{\cI}\colon\cE^{\cI}\to\cE'^{\cI}$ to be a functor 
by sending an $\cI$-diagram $x$ in $\cE$ to $fx$ the composition with $f$. 
Then $f^{\cI}$ be an exact functor. 
Let $\theta\colon f\to f'$ be a natural transformation 
between exact functors $f$, $f'\colon\cE\to\cE'$. 
We define $\theta^{\cI}\colon f^{\cI}\to f'^{\cI}$ to 
be a natural transformation by setting
\begin{equation}
\label{eq:diagram natural transformation df}
\theta^{\cI}(x)(i):=\theta(x(i))
\end{equation}
for any objects $x$ in $\cE$ and $i$ in $\cI$. 
Moreover we denote the full subcategory of $\cE^{\cI}$ consisting of 
exact functors $\cI\to\cE$ by $\Ex(\cI,\cE)$. 
Then $\Ex(\cI,\cE)$ is a strict exact subcategory of $\cE^{\cI}$. 
Notice that for any level admissible monomorphism $f \onto{i} g$ in 
$\Ex(\cI,\cE)$ and for any morphism $x\onto{a} y$ in $\cE$, 
the square below is admissible in the sense that 
the induced morphism $f(y)\sqcup_{f(x)}g(x)\rinf g(y)$ 
is an admissible monomorphism in $\cE$. 
(See \cite[\S A.2]{Moc10}.)
$$\xymatrix{
f(x) \ar[r]^{i_x} \ar[d]_{f(a)} & g(x) \ar[d]^{g(a)}\\
f(y) \ar[r]_{i_y} & g(y).
}$$
\end{para}

\begin{para}
\label{nt:ExCat}
{\bf ($\ExCat$).}\ \ 
We write $\ExCat$ for the category of 
small exact categories and exact functors 
and if we regard $\ExCat$ as 
a $2$-category where the class of $2$-morphisms 
is the class of all natural transformations, 
then we denote it by $\underline{\ExCat}$. 
Notice that for any pair of small exact categories 
$\cE$ and $\cF$, the Hom category $\HOM(\cE,\cF)$ 
with level admissible exact sequences is an exact category 
by \ref{nt:functor category} and $\underline{\ExCat}$ 
is a locally exact $2$-category. 
\end{para}

\begin{para}
\label{nt:diagram cat over locally exact 2-cat}
{\bf (Diagrams over locally exact $2$-categories).}\ \ 
Let $\cI$ be a small category and let 
$\cC$ be a $2$-category and 
$\cE$ and $\cF$ 
be a pair of $\cI$-diagrams in $\cC$, 
namely $1$-functors $\cI\to\cC$ 
and $f$ and $g$ be 
a pair of natural transformations from $\cE$ to $\cF$. 
A {\it modification} from $f$ to $g$ 
is a family of $2$-morphisms 
$\theta=\{\theta_i\colon f_i\to g_i\}_{i\in\Ob\cI}$ in $\cC$ 
indexed by the set of objects of $\cI$ subjects to the conditions 
that for any morphism $a\colon i\to j$ in $\cI$, 
we have the equality
\begin{equation}
\label{eq:modification cond}
\theta_j\ast\cE_a=\cF_a\ast \theta_i.
\end{equation}
We denote this situation by $\theta\colon f\to g$. 

Moreover 
let $\cG$ be an $\cI$-diagram in $\cC$ 
and let $h\colon\cE\to\cF$ and $f'$, $g'\colon \cF\to \cG$ 
be natural transformations and let 
$\alpha \colon f\to g$, $\beta\colon g\to h$ and $\gamma\colon f'\to g'$ 
be modifications. 
Then we define $\beta\cdot\alpha\colon f\to h$ and 
$\gamma\ast\alpha\colon f'f\to g'g$ to be modifications by setting 
for any object $i$ of $\cI$, 
\begin{equation}
\label{eq:vertical comp of modifi}
{(\beta\cdot\alpha)}_i:=\beta_i\cdot\alpha_i,
\end{equation}
\begin{equation}
\label{eq:horizontal comp of modifi}
{(\alpha'\ast\alpha)}_i=\alpha'_i\ast\alpha_i.
\end{equation}

We call $\beta\cdot\alpha$ and $\alpha'\ast\alpha$ 
the {\it vertical composition of $\beta$ and $\alpha$} 
and the {\it horizontal composition $\alpha'$ and $\alpha$} respectively. 
We denote the $2$-category of $\cI$-diagrams in $\cC$ by $\cC^{\cI}$ 
whose objects are $\cI$-diagrams in $\cC$, 
whose $1$-morphisms are natural transformations and 
whose $2$-morphisms are modifications. 

Assume that $\cC$ is locally exact and let  
$\cE$ and $\cF$ be a pair of $\cI$-diagrams in $\cC$. 
A sequence $f\onto{\alpha}g\onto{\beta}h$ 
in $\underline{\HOM}_{\cC^{\cI}}(\cE,\cF)$ the category of 
natural transformations from $\cE$ to $\cF$ and 
modifications is a {\it level admissible exact sequence} if 
for any object $i$ in $\cI$, 
the sequence $f_i\onto{\alpha_i}g_i\onto{\beta_i}h_i$ 
is a level admissible exact sequence of exact functors 
from $\cE_i$ to $\cF_i$. 
We can show that $\underline{\HOM}_{\cC^{\cI}}(\cE,\cF)$ with 
the set of level admissible exact sequences is an exact category 
and $\cC^{\cI}$ is a locally exact $2$-category. 
In particular for a small category $\cI$, the $2$-category 
$\underline{\ExCat}^{\cI}$ 
of $\cI$-diagrams over $\underline{\ExCat}$ is 
a locally exact $2$-category. 
\end{para}

Next we define the suitable notion of morphisms between 
locally exact $2$-categories. 

\begin{para}
\label{df:locally exact 2-functor}
{\bf Definition (Locally exact $2$-functors).}\ \ 
Let $\cC$ and $\calD$ be locally exact $2$-categories. 
A $2$-functor $f\colon \cC\to \calD$ is {\it locally exact} 
if for any pair of objects $x$ and $y$ in $\cC$, 
the functor 
$f\colon\underline{\HOM}_{\cC}(x,y)\to\underline{\HOM}_{\calD}(fx,fy)$ 
is exact. 
\end{para}

\begin{para}
\label{ex:filtered colimit in ExCat}
{\bf (Filtered colimit in $\ExCat$).}\ \ 
Let $\cI$ be a small category. 
Assume that $\cI$ is {\it filtered}, namely $\cI$ satisfies the 
following three conditions:

\sn
$\bullet$ $\cI$ is a non-empty category. 

\sn
$\bullet$ For any pair of objects $i$ and $j$ in $\cI$, 
there exists an object $k$ in $\cI$ and a pair of morphisms 
$i\to k$ and $j\to k$. 

\sn
$\bullet$ For any pair of morphisms $a$, $b\colon i\to j$ in $\cI$, 
there exists a morphism $c\colon j\to k$ such that $ca=cb$.

Then we have the {\it colimit} $2$-functor 
$\displaystyle{\underset{\cI}{\colim}\colon\underline{\ExCat}^{\cI}\to\underline{\ExCat}}$. 
We briefly recall the construction of this functor. 
Let $\cE\colon\cI\to\ExCat$ be an $\cI$-diagram in $\ExCat$. 
Then we have the equalities
\begin{equation}
\label{eq:colim ob}
\Ob\underset{\cI}{\colim}\cE=\underset{\cI}{\colim}\Ob\cE,
\end{equation}
\begin{equation}
\label{eq:colim mor}
\Mor\underset{\cI}{\colim}\cE=\underset{\cI}{\colim}\Mor\cE.
\end{equation}
Namely for example, for the set of objects, 
we have the equality
\begin{equation}
\label{eq:colim eq 2}
\Ob\underset{\cI}{\colim}\cE=\Ob\underset{\cI}{\colim}\cE:=\bigsqcup_{i\in\Ob\cI}\Ob\cE_i/\sim
\end{equation}
where the equivalence relation $\sim$ is defined as follows. 
For any pair of objects $x\in\Ob\cE_i$ and $y\in\Ob\cE_j$, 
we say that $x$ and $y$ are equivalent if 
there exists a pair of morphisms $a\colon i\to k$ and $b\colon j\to k$ 
such that $\cE_a(x)=\cE_b(y)$ and 
we denote this situation by $x\sim y$. 
Then we can show that the relation $\sim$ is an equivalence relation on 
$\displaystyle{\bigsqcup_{i\in\Ob\cI}\Ob\cE_i}$. 
We say that 
a sequence $x\to y\to z$ in $\displaystyle{\underset{\cI}{\colim}\cE}$ 
is an {\it admissible exact sequence} if 
it represented by an admissible exact sequence in some $\cE_i$. 
We can show that $\displaystyle{\underset{\cI}{\colim}\cE}$ 
with the set of admissible exact sequences is an exact category and 
we can show that $2$-functor $\displaystyle{\underset{\cI}{\colim}\colon 
\underline{\ExCat}^{\cI}\to \underline{\ExCat}}$ is 
locally exact $2$-functor. 
\end{para}

\begin{para}
\label{df:coproduct of exact categories}
{\bf Definition (Coproduct of exact categories).}\ \ 
Let $\{\cE_i\}_{i\in\cI}$ be 
a family of pointed categories indexed by a set $\cI$. 
We denote the full subcategory of $\displaystyle{\prod_{i\in\cI}\cE_i}$ 
consisting of 
those objects $x={(x_i)}_{i\in\cI}$ such that 
$\#\{i\in\cI;x_i\neq 0 \}<+\infty$ by $\displaystyle{\bigvee_{i\in\cI}\cE_i}$. 
We assume that $\cE_i$ is an exact category for all $i\in\cI$. 
We say that a sequence 
${(x_i)}_{i\in\cI}\to{(y_i)}_{i\in\cI}\to{(z_i)}_{i\in\cI}$ in 
$\displaystyle{\prod_{i\in\cI}\cE_i}$ is 
a {\it level admissible exact sequence} if 
all $i\in \cI$, a sequence $x_i\to y_i\to z_i$ 
is an admissible exact sequence in $\cE_i$. 
$\displaystyle{\prod_{i\in\cI}\cE_i}$ with level admissible exact sequences 
is an exact category and 
$\displaystyle{\bigvee_{i\in\cI}\cE_i}$ 
is a strict exact subcategory of $\displaystyle{\prod_{i\in\cI}\cE_i}$. 
\end{para}

\begin{para}
\label{df:Idempotent completion}
{\bf (Idempotent completion functor).}\ \ 
One of a typical example for a locally exact $2$-functor 
$\ExCat\to\ExCat$ is the idempotent completion functor. 
We recall the definition and fundamental properties from \cite{Kar68} and \cite[\S A]{TT90}. 
An additive category $\cA$ is 
{\it idempotent complete} if any idempotent 
$e\colon x \to x$ with $e^2=e$, 
arises from a splitting of $x$, 
$x \isoto \im(e) \oplus \Ker(e)$. 
For an additive category $\cA$, 
its {\it idempotent completion} $\widehat{\cA}$ 
is a category whose objects are pair $(x,e)$ 
consisting of an object $x$ in $\cA$ 
and an idempotent endomorphism $e$ of $x$. 
A morphism $a\colon (x,e) \to (x',e')$ in $\widehat{\cA}$ 
is a morphism $a\colon x\to x'$ in $\cA$ 
subjects to the condition that $ae=e'a=a$. 
Then for a pair of objects $(x,e)$ and $(y,e')$ in $\widehat{\cA}$, 
$\displaystyle{(x,e)\oplus (y,e')\isoto (x\oplus y,\begin{pmatrix}e & 0\\ 0 &e' \end{pmatrix})}$
In particular $\widehat{\cA}$ is an additive category. 
There is the fully faithful additive functor 
$
i_{\cA}\colon \cA \to \widehat{\cA},\ 
x\mapsto (x,\id_x)
$ 
satisfying the following universal property:

\sn
For any idempotent complete additive category $\cB$ 
and any additive functor $f\colon \cA \to \cB$, 
up to natural equivalence, 
$f$ factor in a unique way through $i_{\cA}$.

Let $\cE$ be an exact category. 
Then we can make its idempotent completion $\widehat{\cE}$ 
into an exact category by declaring that a composable sequence 
in $\widehat{\cE}$ is an admissible exact sequence 
if and only if it is a direct summand of 
an admissible exact sequence 
in $\cE$. 
Then the functor $i_{\cE}\colon \cE\to \widehat{\cE}$ 
is exact and reflects exactness.

Now we define the idempotent completion $2$-functor 
$\widehat{(-)}\colon\underline{\ExCat} \to \underline{\ExCat}$. 
For an exact functor $f\colon\cE \to \cF$, 
we define 
$\widehat{f}\colon\widehat{\cE} \to \widehat{\cF}$ to be an exact functor 
by sending an object  $(x,e)$ in $\widehat{\cE}$ to $(f(x),f(e))$. 
For a natural transformation 
$\theta\colon f \to g$ between 
exact functors $f$, $g\colon \cE \to \cF$, 
we define $\widehat{\theta}\colon f \to g$ 
to be a natural transformation by 
setting for an object $(x,e)$ in $\widehat{\cE}$, 
\begin{equation}
\label{eq:widehat theta df}
\widehat{\theta}(x,e):=g(e)\theta(x)
\end{equation}
We can check that 
this association is $2$-functorial and this $2$-functor 
is locally exact. 
\end{para}

\begin{para}
\label{df:diagram functors}
{\bf (Diagram functor).}\ \ 
Let $\cI$ be a small category. 
Then we define ${(-)}^{\cI}\colon\ExCat\to \ExCat$ 
to be a $2$-functor by sending 
an exact category $\cE$ to $\cE^{\cI}$ the category of 
$\cI$-diagrams in $\cE$ with level admissible exact structure 
(See \ref{nt:functor category}) and sending an exact functor 
$f\colon\cE\to\cE'$ to 
$f^{\cI}\colon\cE^{\cI}\to\cE'^{\cI}$ 
and sending a natural transformation $\theta\colon f\to f'$ 
between exact functors $f$, $f'\colon\cE\to\cE'$ to 
$\theta^{\cI}\colon f^{\cI}\to f'^{\cI}$. 
Then we can show that ${(-)}^{\cI}$ 
is locally exact and we call it the 
{\it $\cI$-diagram association $2$-functor}. 
\end{para}

\begin{para}
\label{df:toplogizing subcat}
{\bf (Topologizing subcategories).}\ \ 
Let $\cE$ be an exact category and let $\calD$ 
be a non-empty full subcategory of $\cE$. 
We say that $\calD$ is a {\it topologizing subcategory} of $\cE$ if 
$\calD$ is closed under finite direct sums and closed under 
admissible sub- and quotient objects. 
The last condition means that 
for an admissible exact sequence $x\rinf y \rdef z$ in $\cE$, 
if $y$ is in $\calD$, 
then $x$ and $z$ are also in $\calD$. 
Thus in this case, $\calD$ contains the zero object and 
closed under isomorphisms in $\cE$. 
Namely for any object $x$ in $\cE$ which is isomorphic to an object in $\calD$ 
is also in $\calD$. 
The naming of the term `topologizing' comes from 
noncommutative geometry of abelian categories 
by Rosenberg (see \cite[Lecture 2 1.1]{Ros08}). 
By \cite[5,3]{Moc13a}, a topologizing subcategory of $\cE$ 
naturally becomes a strict exact subcategory of $\cE$. 
\end{para}

\begin{para}
\label{df:Serre subcategory}
{\bf (Serre subcategories, Serre radical, semi-Serre subcategories).}\ \ 
Let $\cE$ be an exact category. 
We say that a full subcategory $\calD$ of $\cE$ is a {\it Serre subcategory} 
if it is an extensional closed topologizing subcategory of $\cE$. 
For any full subcategory $\calD$ of $\cE$, 
we write ${}^S\!\!\!\sqrt{\calD}$ 
for intersection of all Serre subcategories 
which contain $\calD$ and call it the {\it Serre radical of $\calD$} 
({\it in $\cE$}). 

We say that a full subcategory $\calD$ of $\cE$ is 
a {\it semi-Serre subcategory} if it is closed under 
finite direct sums and for an admissible sequence 
$x\rinf y\rdef z$ in $\cE$, 
if $x$ and $y$ are in $\calD$, then 
$z$ is also in $\calD$ and if $y$ and $z$ are in $\calD$, 
then $x$ is also in $\calD$. 
\end{para}

\begin{para}
\label{lem:semi-Serre subcategory}
{\bf Lemma.}\ \ 
{\it
Let $\cE$ be an exact category and let $\calD$ be a 
semi-Serre subcategory. 
Then $\calD$ is a strict exact subcategory of $\cE$.
}
\end{para}

\begin{proof}
We will show that\\
$\mathrm{(i)}$ 
for a pair of admissible monomorphisms $x\rinf y$ and $y\rinf z$ 
with $x$, $y$ and $z$ are in $\calD$, $\coker(x\rinf z)$ is in $\calD$,\\
$\mathrm{(ii)}$ 
and in the pushout diagram in $\cE$ below
$$
\xymatrix{
x\ar@{>->}[r] \ar[d] & y\ar[d]\\
z\ar@{>->}[r] & w,
}
$$
if $x$, $y$ and $z$ are in $\calD$, then $w$ and $\coker (z\rinf w)$ 
are also in $\calD$.

\begin{proof}[Proof of $\mathrm{(i)}$]
In the admissible exact sequence $x\rinf z\rdef \coker(x\rinf z)$ in $\cE$, 
since $x$ and $z$ are in $\calD$, $\coker(x\rinf z)$ are also in $\calD$.
\end{proof}
\begin{proof}[Proof of $\mathrm{(ii)}$]
By \ref{lem:kel90}, the sequence $x\rinf z\oplus y \rdef w$ 
is an admissible exact 
sequence in $\cE$ with $x$ and $z\oplus y$ are in $\calD$. 
Thus $w$ is in $\calD$ and considering the 
admissible exact sequence $z\rinf w\rdef \coker(z\rinf w)$ with 
$z$ and $w$ are in $\calD$, 
it turns out that $\coker(z\rinf w)$ is in $\calD$. 
\end{proof}
We can show 
similar statements for admissible epimorphisms and 
hence $\calD$ is a strict exact subcategory of $\cE$. 
\end{proof}

\subsection{Complicial structures}
\label{subsec:complicial structures}

\begin{para}
\label{df:complicial exact category}
{\bf (Complicial objects).}\ \ 
Let $\cC$ be a $2$-category and $x$ an object in $\cC$. 
A {\it semi-commutative unitary magmad on $x$} is 
a quadruple $(C,\iota,r,\sigma)$ consisting of 
an $1$-morphism $C\colon x\to x$, 
a pair of $2$-morphisms 
$\iota\colon\id_{\cC}\to C$ and $r\colon CC\to C$ and 
a $2$-isomorphism 
$\sigma\colon CC\isoto CC$ which satisfy the following conditions:

\sn
{\bf (Unitary law).}\ 
$r\cdot (C\ast\iota)=r\cdot (\iota\ast C)=\id_{C}$.

\sn
{\bf (Semi-commutativity law).}\ 
$\sigma \cdot (C\ast\iota)=\iota\ast C$.

\sn
{\bf (Involution law).}\ 
$\sigma\cdot\sigma=\id_{CC}$.
$$
\xymatrix{
C \ar[r]^{C\ast\iota} \ar[rd]_{\id_C} & 
CC \ar[d]^r & 
C \ar[l]_{\iota\ast C} \ar[ld]^{\id_C}\\ 
& C, &\\
} 
\xymatrix{
& C \ar[rd]^{C\ast\iota} \ar[ld]_{\iota\ast C}&\\
CC \ar[rr]^{\sim}_{\sigma} & & CC.
}
$$

Moreover assume that $\cC$ is a locally exact $2$-category 
(see Definition~\ref{df:locally exact 2-categories}). 
A semi-commutative unitary magmad $(C,\iota,r,\sigma)$ on an 
object $x$ in $\cC$ is {\it exact} if 
the $2$-morphism $\iota\colon\id_{x}\to C$ 
is an admissible monomorphism 
in the Hom category $\underline{\HOM}_{\cC}(x,x)$ from $x$ to $x$. 
When we consider an exact semi-commutative unitary magmad, we 
fix an admissible exact sequence 
\begin{equation}
\label{eq:id-> C-> T}
\id_x \overset{\iota}{\rinf} C \overset{\pi}{\rdef} T. 
\end{equation}
We call $T$ the {\it suspension} ({\it $1$-morphism associated with $C$}).

An exact semi-commutative unitary magmad on an 
object $x$ is {\it stable} if 
the suspension $1$-morphism $T\colon x\to x$ 
is an equivalence. 
When we consider a stable exact semi-commutative unitary magmad, 
we fix an adjoint equivalence $(T,T^{-1},\alpha,\beta)$ 
from $x$ to $x$. 
Namely 
$T^{-1}\colon x\to x$ is an $1$-morphism and 
a pair of $2$-isomorphisms 
$\alpha\colon TT^{-1}\isoto \id_{x}$ and 
$\beta\colon \id_{x}\isoto T^{-1}T$ which satisfy the equalities:
\begin{equation}
\label{eq:T-eq}
(\alpha\ast T)\cdot(T\ast \beta)=\id_T,  
\end{equation}
\begin{equation}
\label{eq:T-inverse-eq}
(T^{-1}\ast\alpha)\cdot (\beta\ast T^{-1})=\id_{T^{-1}}.
\end{equation}

In the paper \cite{Moc10}, 
we consider $\cC=\underline{\ExCat}$ the locally exact 
$2$-category of small exact categories 
(see \ref{nt:ExCat}) 
and 
we call an exact semi-commutative unitary magmad 
on an exact category $\cE$ 
a complicial structure on $\cE$. 
But in this paper 
we only use stable complicial structures and for simplicity 
we call stable semi-commutative unitary magmads on an object $x$ 
in a locally exact $2$-category 
{\it complicial structures} on $x$. 
Thus more precisely, a complicial structure on $x$ is a ninefold of 
$(C,\iota,r,\sigma,T,T^{-1},\pi,\alpha,\beta)$ consisting of 
a triple of $1$-morphisms 
$C$, $T$ and $T^{-1}\colon x\to x$ and 
$2$-morphisms $\iota\colon\id_{x}\to C$, $r\colon CC\to C$, 
$\pi\colon C\to T$ and 
a triple of $2$-isomorphisms 
$\sigma\colon CC\isoto CC$, 
$\alpha\colon TT^{-1}\isoto \id_{x}$ and $\beta\colon \id_{x}\isoto T^{-1}T$ 
such that  
a sequence $\id_{x}\overset{\iota}{\rinf}C\overset{\pi}{\rdef}T$ 
is exact in the Hom category $\underline{\HOM}_{\cC}(x,x)$ 
from $x$ to $x$ and 
satisfy the equalities 
$r\cdot (C\ast\iota)=r\cdot (\iota\ast C)=\id_{C}$, 
$\sigma\cdot\sigma=\id_{CC}$, $\sigma \cdot (C\ast\iota)=\iota\ast C$, 
$(T^{-1}\ast \alpha) \cdot (\beta\ast T^{-1})=\id_{T^{-1}}$ and 
$(\alpha\ast T)\cdot(T\ast \beta)=\id_{T}$. 
We call an object equipped with a complicial structure 
a {\it complicial object}. 
In particular we call a complicial object in $\underline{\ExCat}$ 
the $2$-category of small exact categories 
a {\it complicial exact category}. 
When we denote a complicial structure, we often omit 
$T$, $T^{-1}$, $\pi$, $\alpha$ and $\beta$ in the notations. 
\end{para}

\begin{para}
\label{ex:standard complicial structure}
{\bf Example (Standard complicial structure).}\ \ 
Here we will illustrate an archetype of complicial exact categories. 
Let $\cE$ be an exact category. 
Now we give the complicial structure on $\Ch_{b}(\cE)$ 
the category of bounded chain complexes over $\cE$ 
 as follows. 
Here we use the homological notations for chain complexes. 
Namely boundary morphisms are of degree $-1$. 

The functor $C\colon \Ch_b(\cE) \to \Ch_b(\cE)$ 
is given by sending $x$ to $Cx:=\Cone\id_x$ the canonical mapping cone 
of the identity morphism of $x$. 
Namely the degree $n$ part of $Cx$ is ${(Cx)}_n=x_{n-1}\oplus x_n$ and 
the degree $n$ boundary morphism $d_n^{Cx}\colon {(Cx)}_n\to {(Cx)}_{n-1}$ 
is given by 
$\displaystyle{d^{Cx}_n=
\begin{pmatrix} 
-d^x_{n-1} & 0\\
-\id_{x_{n-1}} & d_n^x
\end{pmatrix}}$. 
For a chain morphism $f\colon x\to y$, we define 
$Cf\colon Cx \to Cy$ to be a chain morphism by setting 
$(Cf)_n:=\begin{pmatrix}f_{n-1} & 0\\ 0 & f_n\end{pmatrix}$ 
for an integer $n$. 
For any complex $x$, 
we define 
$\iota_x\colon x \to C(x)$, 
$r_x\colon CC(x) \to C(x)$ 
and $\sigma_x\colon CC(x) \to CC(x)$ 
to be chain morphisms 
by setting 
\begin{equation}
\label{eq:standard iota}
{(\iota_x)}_n=
\begin{pmatrix}
0\\
\id_{x_n}
\end{pmatrix},
\end{equation}
\begin{equation}
\label{eq:standard r}
{(r_x)}_n= 
\begin{pmatrix}
0 &\!\! \id_{x_{n-1}} &\!\! \id_{x_{n-1}} &\!\! 0\\
0 &\!\! 0 &\!\! 0 &\!\! \id_{x_n}
\end{pmatrix}
\ \text{and} 
\end{equation}
\begin{equation}
\label{eq:standard sigma} 
{(\sigma_x)}_n=
\begin{pmatrix}
\!\!-\id_{x_{n-2}} &\!\! 0 &\!\! 0 &\!\! 0\\
\!\!0 &\!\! 0 &\!\! \id_{x_{n-1}} &\!\! 0\\
\!\!0 &\!\! \id_{x_{n-1}} &\!\! 0 &\!\! 0\\
\!\!0 &\!\! 0 &\!\! 0 &\!\! \id_{x_n}
\end{pmatrix}.
\end{equation}
Then $\Ch_{b}(\cE)$ together with the quadruple $(C,\iota,r,\sigma)$ 
forms a complicial exact category. 

In this example, we give the functor $T$, 
$T^{-1}\colon\Ch_b(\cE)\to\Ch_b(\cE)$ 
and a natural transformation $\pi\colon C\to T$ and natural equivalences 
$\alpha\colon TT^{-1}\isoto\id_{\Ch_b(\cE)} $ and 
$\beta\colon\id_{\Ch_b(\cE)}\isoto T^{-1}T$ 
by the following formulas for any complex $x$ and integer $n$;
\begin{equation}
\label{eq:standard T}
{(Tx)}_n=x_{n-1}, \ \ d_n^{Tx}=-d^x_{n-1},
\end{equation}
\begin{equation}
\label{eq:standard T-1}
{(T^{-1}x)}_n=x_{n+1}, \ \ d_n^{T^{-1}x}=-d^x_{n+1},
\end{equation}
\begin{equation}
\label{eq:standard pi}
{(\pi_x)}_n=\begin{pmatrix}\id_{x_{n-1}} & 0 \end{pmatrix},
\end{equation}
\begin{equation}
\label{eq:standard alpha beta}
\alpha=\beta=\id_{\Ch_b(\cE)}.
\end{equation}
\end{para}

\begin{para}
\label{df:mor zeta}
{\bf ($\zeta\colon TC\to CC$).}\ \ 
Let $x$ be a complicial object in a locally 
exact $2$-category $\cC$. 
For the exact sequence of $1$-morphisms
in $\underline{\HOM}_{\cC}(x,x)$
\begin{equation}
\label{eq:fund split seq}
C\overset{\iota\ast C}{\rinf} CC\overset{\pi\ast C}{\rdef} TC,
\end{equation}
since we have the retraction $r\colon CC\to C$ of $\iota\ast C$, 
there exists the section $\zeta\colon TC\to CC$ 
such that we have equalities
\begin{equation}
\label{eq:zeta 1}
r\cdot\zeta=0,
\end{equation} 
\begin{equation}
\label{eq:zeta 2}
(\pi\ast C) \cdot \zeta=\id_{TC}, 
\end{equation} 
\begin{equation}
\label{eq:zeta 3}
(\iota\ast C)\cdot r + \zeta\cdot (\pi\ast C)=\id_{CC}.
\end{equation} 
In Example~\ref{ex:typicla split seq}, 
we will give an explicit construction of $\zeta$.

In the standard example~\ref{ex:standard complicial structure}, 
for any chain complex $x$ over an exact category and any integer $n$, 
the morphism ${(\zeta_x)}_n\colon {(TCx)}_n\to {(TCx)}_n$ is given by 
\begin{equation}
\label{eq:standard zeta}
{(\zeta_x)}_n=\begin{pmatrix}
\id_{x_{n-2}} & 0\\
0 & \id_{x_{n-1}}\\
0 & -\id_{x_{n-1}}\\
0 & 0
\end{pmatrix}.
\end{equation}
\end{para}

\begin{para}
\label{df:functor P}
{\bf (Path $1$-morphism).}\ \ 
Let $x$ be a complicial object in a locally exact $2$-category $\cC$. 
Then we define 
$P\colon x\to x$ to be an $1$-morphism and 
$j\colon T^{-1}\to P$ and $q\colon P\to \id_{\cC}$ 
to be $2$-morphisms by setting 
$P:=CT^{-1}$ and $j:=\iota\ast T^{-1}$ and 
$q:=\alpha\cdot( \pi \ast T^{-1} )$. 
Then we have the admissible exact sequence 
$T^{-1}\overset{i}{\rinf} P \overset{q}{\rdef} \id_{\cC}$ in 
$\underline{\HOM}_{\cC}(x,x)$. 
We say that $P$ is the {\it path $1$-morphism} ({\it associated with $C$}). 

In the standard example~\ref{ex:standard complicial structure}, 
for any chain complex $x$ over an exact category, 
the complex $Px$ 
and the chain morphisms 
$j_x\colon T^{-1}x\to Px$ and $q_x\colon Px\to x$ are given by 
\begin{equation}
\label{eq:standard P}
{(Px)}_n:=x_n\oplus x_{n+1},\ \ d_n^{Px}:=
\begin{pmatrix}d_n^x & 0\\ -\id_{x_n} & -d^x_{n+1}\end{pmatrix},
\end{equation}
\begin{equation}
\label{eq:standard j}
{(j_x)}_n:=\begin{pmatrix}0\\ \id_{x_{n+1}} \end{pmatrix},
\end{equation}
\begin{equation}
\label{eq:standard q}
{(q_x)}_n:=\begin{pmatrix}\id_{x_{n}} & 0 \end{pmatrix}
\end{equation}
for any integer $n$. 

We will illustrate the $1$-morphism $P$ is a dual notion of the 
$1$-morphism $C$ 
in some sense. 
In Lemma-Definition~\ref{lemdf:dual comp structure}, 
we will show that 
for a complicial structure $(C,\iota,r,\sigma)$ on an exact category $\cE$ 
which satisfies the certain conditions, 
we will define a natural transformations $s\colon P\to PP$ and 
a natural equivalence $\tau\colon PP\isoto PP$ 
such that the quadruple $(P,q,s,\tau)$ is a 
complicial structure on $\cE^{\op}$ the opposite category of $\cE$. 
\end{para}

\begin{para}
\label{df:odtproduct}
{\bf ($\odot$-products).}\ (\cf \cite[2.20]{Moc10}.)\ \ 
Let $x_i$ ($i=1$, $2$, $3$) be objects in a $2$-category $\cC$ and let 
$\cO_i\colon x_i\to x_i$ and 
$f\colon x_1\to x_2$ and $g\colon x_2\to x_3$ be $1$-morphisms 
and let $c\colon \cO_2f\to f\cO_1$ and $d\colon \cO_3g\to g\cO_2$ 
be $2$-morphisms. 
Then we define $d\odot_{\cO}c$ to be a $2$-morphism
from $\cO_3gf$ to $gf\cO_1$ by setting 
\begin{equation}
\label{eq:odtproduct df}
d\odot_{\cO}c:=(g\ast c) \cdot (d\ast f).
\end{equation}
\end{para}

\begin{para}
\label{df:complicial exact functor}
{\bf (Complicial $1$-morphisms).}\ \ 
A {\it complicial $1$-morphism} between complicial 
objects 
$x \to x'$ in a locally exact $2$-category $\cC$ 
is a pair of an $1$-morphism 
$f\colon x \to x'$ and 
a $2$-morphism $c\colon C^{x'}f\isoto fC^{x}$ 
which satisfies the equality 
\begin{equation}
\label{eq:comp ex func cond}
c\cdot(\iota^{x'} \ast f)=f\ast\iota^{x} 
\end{equation}
where $(C^x,\iota^x,r^x,\sigma^x)$ and 
$(C^{x'},\iota^{x'},r^{x'},\sigma^{x'})$ 
are complicial structures on 
$x$ and $x'$ respectively. 
We often omit $c$ in the notation. 
For complicial $1$-morphisms 
$x \onto{(f,c)} x' \onto{(g,d)} x''$, 
we define composition of $(g,d)$ and $(f,c)$ by 
$(g,d)(f,c):=(gf,d\odot_C c)$. 

We say that a complicial $1$-morphism $(f,c)\colon x\to x'$ 
is {\it normal} if it satisfies the following two conditions
\begin{equation}
\label{eq:normal comp ex func r}
(f\ast r^{x})\cdot (c\ast C^{x})\cdot (C^{x'}\ast c)
=c\cdot (r^{x'}\ast f),
\end{equation}
\begin{equation}
\label{eq:normal comp ex func sigma}
(f\ast\sigma^{x})\cdot (c\ast C^{x})\cdot (C^{x'}\ast c) =
(c\ast C^{x})\cdot (C^{x'}\ast c)\cdot (\sigma^{x'}\ast f).
\end{equation}
$$\xymatrix{
C^{x'}C^{x'}f \ar[r]^{C^{x'}\ast c}_{\sim} \ar[d]_{r^{x'}\ast f} & 
C^{x'}fC^{x} \ar[r]^{c\ast C^{x}}_{\sim} & 
fC^{x}C^{x} \ar[d]^{f\ast r^{x}}\\
C^{x'}f \ar[rr]^{\sim}_{c} & & fC^{x},
}\ \ 
\xymatrix{
C^{x'}C^{x'}f \ar[r]^{C^{x'}\ast c}_{\sim} \ar[d]_{\sigma^{x'}\ast f}^{\wr} & 
C^{x'}fC^{x} \ar[r]^{c\ast C^{x}}_{\sim} & 
fC^{x}C^{x} \ar[d]^{f\ast \sigma^{x}}_{\wr}\\
C^{x'}C^{x'}f \ar[r]_{C^{x'}\ast c}^{\sim} & 
C^{x'}fC^{x} \ar[r]_{c\ast C^{x}}^{\sim} & fC^{x}C^{x}.
}
$$

We say that a complicial $1$-morphism $(f,c)\colon x\to x'$ is 
{\it strictly normal} if 
it is normal and moreover it satisfies the following two equalities:
\begin{equation}
\label{eq:strict normal 1}
C^{x'}f=fC^{x},
\end{equation}
\begin{equation}
\label{eq:strict normal 2}
c=\id_{fC^{x}}.
\end{equation}

We can show that 
for a pair of composable normal (resp. strictly normal) 
complicial $1$-morphisms 
$x\onto{(f,c)}x'\onto{(g,d)}x''$, 
the compositions $(gf,d\odot_C c)\colon x\to x''$ is also a 
normal (resp. strictly normal) complicial $1$-morphism. 
In particular we call complicial $1$-morphisms in 
$\underline{\ExCat}$ {\it complicial exact functors}. 
\end{para}

\begin{para}
\label{df:normal complicial structure}
{\bf (Normal complicial structure).}\ \ 
Let $x$ be an object in a locally exact $2$-category $\cC$. 
We say that a complicial structure $(C,\iota,r,\sigma)$  
on $x$ is {\it normal} 
(resp. {\it strictly normal}) if the pair $(C,\sigma)$ 
is a normal (resp. strictly normal) 
complicial $1$-morphism $x\to x$. 
We say that a complicial object is {\it normal} 
(resp. {\it strictly normal}) 
if its complicial structure is normal (resp. strictly normal). 
We can show that the standard complicial structure 
on the category of bounded chain complexes on an exact category 
(see \ref{ex:standard complicial structure}) is strictly normal. 
\end{para}

\begin{para}
\label{df:comp nat transform}
{\bf (Complicial $2$-morphisms).}\ \ 
A {\it complicial $2$-morphisms} 
between complicial $1$-morphisms 
$$(f,c),\ (g,d)\colon x \to x'$$
from $(f,c)$ to $(g,d)$ 
in a locally exact $2$-category $\cC$ 
is a $2$-morphism $\phi\colon f \to g$ 
which subjects to 
the condition that 
\begin{equation}
\label{eq:comp nat trans cond}
d\cdot (C^{x'}\ast\phi)=(\phi\ast C^{x})\cdot c. 
\end{equation}
In particular we call complicial $2$-morphisms in $\underline{\ExCat}$ 
{\it complicial natural transformations}. 
We can check that for any triple of complicial 
$2$-morphisms 
$\varphi\colon (f,c)\to (g,d)$, 
$\varphi'\colon (g,d)\to (h,e)$ and 
$\psi\colon (f',c')\to (g',d')$ between 
complicial $1$-morphisms 
$(f,c)$, $(g,d)$, $(h,e)\colon x\to x'$ and 
$(f',c')$, $(g',d')\colon x'\to x''$ between complicial objects, 
the vertical and the horizontal compositions $\varphi'\cdot\varphi$ and 
$\psi\ast\varphi$ 
are again complicial $2$-morphisms 
$(f,c)\to(h,e)$ and $(f'f,c'\odot_Cc)\to (g'g,d'\odot_Cd)$ respectively. 
\end{para}

\begin{para}
\label{df:Comp(-)}
{\bf Definition ($\Comp(-)$).}\ \ 
Let $\cC$ be a locally exact $2$-category. 
We write $\Comp(\cC)$ for the category of 
complicial objects and complicial $1$-morphisms in $\cC$ 
and if we regard $\Comp(\cC)$ as a $2$-category 
whose $2$-morphisms are complicial $2$-morphisms with 
the usual horizontal and vertical compositions, 
then we denote it by $\underline{\Comp}(\cC)$. 
We write $\Comp_{\nor}(\cC)$ and  $\Comp_{\stnor}(\cC)$ 
(resp. $\underline{\Comp}_{\nor}(\cC)$ and 
$\underline{\Comp}_{\stnor}(\cC)$)
for the ($2$-)subcategory of $\Comp(\cC)$ 
(resp. $\underline{\Comp}(\cC)$) 
consisting of 
complicial objects and normal complicial $1$-morphisms 
and 
complicial objects and 
strictly normal complicial $1$-morphisms respectively 
(and complicial $2$-morphisms).

Let $x$ and $y$ be a pair of complicial objects in $\cC$. 
For a pair of complicial $1$-morphisms 
$(f,c)$, $(g,d)\colon x\to y$ and a pair of complicial 
$2$-morphisms $\phi$, $\psi\colon (f,c)\to(g,d)$, 
we consider the addition $\phi+\psi\colon f\to g$ 
in $\Hom_{\cC}(f,g)$. 
We can show that $\phi+\psi$ 
can be regarded as a complicial $2$-morphism 
$\phi+\psi\colon (f,c)\to(g,d)$ and 
by this addition, 
$\underline{\HOM}_{\Comp(\cC)}(x,y)$ is enriched by 
the category of abelian groups. 
We say that a sequence of complicial $1$-morphisms from $x$ to $y$, 
$(f,c)\onto{u}(g,d)\onto{v}(h,e)$ in 
$\underline{\HOM}_{\Comp(\cC)}(x,y)$ 
is an admissible exact sequence if the sequence 
$f\onto{u}g\onto{v}h$ is an admissible exact sequence 
in $\underline{\HOM}_{\cC}(x,y)$. 
Then we can show that 
$\underline{\HOM}_{\Comp(\cC)}(x,y)$ 
with the class of admissible exact sequences is an exact category. 
For example, push-out of admissible monomorphism 
$(f,c)\rinc (g,d)$ by a morphism $(f,c)\to(h,e)$ in 
$\underline{\HOM}_{\Comp(\cC)}(x,y)$ is given by 
$(h\sqcup_f g,b\cdot (e\sqcup_c d)\cdot a)$ where 
$h\sqcup_f g$ and 
$e\sqcup_c d\colon C^yh\sqcup_{C^yf}C^yg\isoto hC^x\sqcup_{fC^x}gC^x$ 
are cofiber products 
in the category $\underline{\HOM}_{\cC}(x,y)$ and 
$a\colon C^y(h\sqcup_f g)\isoto C^yh\sqcup_{C^yf}C^yg$ and 
$b\colon hC^x\sqcup_{fC^x}gC^x\isoto (h\sqcup_f g)C^x$ 
are canonical isomorphisms. 
We can show that by this exact structures, 
$\underline{\Comp}(\cC)$ and 
$\underline{\Comp}_{\nor}(\cC)$ are locally exact $2$-category.

For a locally exact $2$-category $\cC$ and a small category $\cI$, 
inspection shows the following equality
\begin{equation}
\label{eq:CompsnCI}
\underline{\Comp}_{\stnor}(\cC^{\cI})={(\underline{\Comp}_{\stnor}(\cC))}^{\cI}.
\end{equation}
\end{para}

We can show the following lemma.

\begin{para}
\label{lem:functoriality of Comp}
{\bf Lemma.}\ \ 
{\it
Let $\cC$ and $\cC'$ be a pair of locally exact $2$-categories and 
let $f\colon\cC\to\cC'$ be a locally exact $2$-functor. 
Then $f$ induces a $2$-functors 
$\underline{\Comp}(f)\colon \underline{\Comp}(\cC)
\to\underline{\Comp}(\cC')$ and 
$\underline{\Comp}_{\stnor}(f)\colon\underline{\Comp}_{\stnor}(\cC)\to
\underline{\Comp}_{\stnor}(\cC')$.
}\qed
\end{para}

By the equality $\mathrm{(\ref{eq:CompsnCI})}$ 
and the lemma above, we obtain 
the following:

\begin{para}
\label{cor:colimit in Comp}
{\bf Corollary.}\ \ 
{\it
Let $\cC$ be a locally exact $2$-category and let $\cI$ 
be a small category. 
Assume that $\cC$ is closed under $\cI$-index exact colimits. 
Namely there exists a locally exact $2$-functor 
$\displaystyle{\underset{\cI}{\colim}\colon \cC^{\cI}\to\cC}$. 
Then $\Comp_{\stnor}(\cC)$ is also closed under $\cI$-index exact colimits. 
That is, there exists a locally exact $2$-functor 
$\displaystyle{\underset{\cI}{\colim}\colon{\underline{\Comp}_{\stnor}(\cC)}^{\cI}\to \underline{\Comp}_{\stnor}(\cC)}$. 
}\qed
\end{para}

To calculate general filtered colimits in 
$\underline{\Comp}_{\nor}(\ExCat)$, 
the following strification lemma is useful.

\begin{para}
\label{lem:strification}
{\bf Lemma (Strification).}\ \ 
{\it
Let $\cI$ be a small filtered category and let 
$\cE\colon\cI\to\Comp_{\nor}(\ExCat)$ be a diagram of 
small complicial exact categories. 
Then there exists an $\cI$-diagram 
$\cE'\colon\cI\to\Comp_{\stnor}(\ExCat)$ 
and the natural equivalence $\Theta\colon\cE\isoto \cE'$. 
}
\end{para}

\begin{proof}
For an object $i$ in $\cI$, we define 
$\cE_i'$ to be a category whose objects 
are pairs $(a\colon j\to i,x)$ consisting of a morphism $a\colon j\to i$ 
in $\cI$ and an object $x$ in $\cE_j$ and whose 
morphisms $u\colon (a\colon j\to i,x)\to (a'\colon j'\to i,x')$ 
is a morphism $u\colon\cE_a(x)\to\cE'_{a'}(x)$ in $\cE_i$. 
For a morphism $a\colon i\to j$ in $\cI$, 
we define $\cE_a'\colon \cE_i'\to\cE_j'$ to be 
a functor by sending an object 
$(\alpha\colon k\to i,x)$ in $\cE_i'$ to an object 
$(a\alpha\colon k\to j,x)$ in $\cE_j'$ and 
a morphism $g\colon (\alpha\colon k\to i,x)\to (\beta\colon l\to i,y)$ 
in $\cE_i'$ to a morphism 
$\cE_a(g)\colon (a\alpha\colon k\to j,x)\to (a\beta\colon l\to j,y)$ 
in $\cE_j'$. 
For an object $i$ in $\cI$, 
we define $\Theta_i\colon\cE_i'\to \cE_i$ and 
$\Psi_i\colon \cE_i\to \cE_i'$ to be 
a pair of functors by sending an object 
$(\alpha\colon j\to i,x)$ 
in $\cE_i'$ to 
$\cE_a(x)$ in $\cE_i$ and 
an object $y$ in $\cE_i$ to an object $(\id_i\colon i\to i,y)$ 
in $\cE_i'$ respectively. 
Moreover for an object $(\alpha\colon j\to i,x)$ in $\cE_i'$, 
we set $u(\alpha\colon j\to i,x):=
\id_{\cE_a(x)}\colon(\alpha\colon j\to i,x)\isoto 
(\id_i\colon i\to i,\cE_a(x))$. 
Then we can show that $\Theta_i\Psi_i=\id_{\cE_i}$ and 
$u\colon\id_{\cE_i'}\isoto\Psi_i\Theta_i$ 
and thus $\Theta_i\colon\cE'_i\isoto\cE_i$ 
is an equivalence of categories. 
We can make $\cE'$ into an exact category such that the 
functor $\Theta_i\colon\cE_i'\to\cE_i$ is an exact functor. 
We define $C_i\colon\cE_i'\to\cE_i'$, 
$\iota_i'\colon\id_{\cE_i'}\to C_i'$, 
$r_i'\colon C_i'C_i'\to C_i'$ and 
$\sigma_i'\colon C_i'C_i'\isoto C_i'C_i'$ 
to be a functor and 
natural transformations respectively by sending 
a morphism $g\colon(a\colon j\to i,x)\to (b\colon k\to i,y)$ 
to a morphism 
$(c_b)(C_ig)(c_a^{-1})\colon(a\colon j\to i,C_jx)\to (b\colon k\to i,C_ky)$ 
and by setting 
${\iota_i'}_{(a\colon j\to i,x)}:=(c_a) (\iota_{\cE_a(x)})\colon 
(a\colon j\to i,x)\to (a\colon j\to i,C_jx)$,  
${r_i'}_{(a\colon j\to i,x)}:=({c_a}_x)({r_j}_{{c_a}_x})
({C_i\ast {c_a}_x}^{-1})
({{c_a}_x\ast C_jx}^{-1})$ and 
${\sigma_i'}_{(a\colon j\to i,x)}:={\sigma_i}_x \colon 
(a\colon j\to i,x)\isoto(a\colon j\to i,x)$ 
for an object $(a\colon j\to i,x)$ 
in $\cE_i$ respectively 
where $c_a\colon C_i\cE_a\isoto \cE_aCj$ is a part of complicial $1$-functor 
$(\cE_a,c_a)\colon \cE_i\to\cE_j$ associated to a morphism $a\colon i\to j$ 
and $C_i\colon\cE_i\to\cE_i$, $\iota_i\colon\id_{\cE_i}\to C_i$, 
$r_i\colon C_iC_i\to C_i$ and $\sigma_i\colon C_iC_i\isoto C_iC_i$ 
are parts of complicial structure on $\cE_i$. 
Then we can show that 
the quadruple 
$(C_i',\iota_i',r_i',\sigma_i')$ is a 
complicial structure on $\cE_i'$ and the equivalence 
$\Theta_i\colon\cE_i'\isoto\cE_i$ is a complicial $1$-morphism 
and $u\colon\id_{\cE_i'}\isoto\Psi_i\Theta_i$ is a complicial $2$-morphism. 
Moreover we can regard $\cE'$ is a $\cI$-diagram on 
$\Comp_{\stnor}(\ExCat)$. 
\end{proof}

\begin{para}
\label{df:comp str on idempotent completion}
{\bf (Complicial structure on idempotent completion).}\ \ 
By Lemma~\ref{lem:functoriality of Comp}, 
the idempotent completion $2$-functor 
$\widehat{(-)}\colon\underline{\ExCat}\to\underline{\ExCat}$ 
induces a $2$-functor $\underline{\Comp}(\widehat{(-)})\colon \underline{\Comp}(\underline{\ExCat})\to 
\underline{\Comp}(\underline{\ExCat})$. 
In particular for a complicial exact category $\cC$, 
$\widehat{\cC}$ the idempotent completion of $\cC$ has 
a natural complicial structure induced from the functor 
$\underline{\Comp}(\widehat{(-)})$ 
and the canonical exact functor $i_{\cC}\colon\cC\to\widehat{\cC}$, 
$x\mapsto (x,\id_x)$ is a strictly normal complicial exact functor.
\end{para}

\begin{para}
\label{df:level complicial structure}
{\bf Definition (Level complicial structure).}\ \ 
Let $\cI$ be a small category. 
Then the $\cI$-diagram association $2$-functor 
${(-)}^{\cI}\colon\underline{\ExCat}\to\underline{\ExCat}$ 
(see \ref{df:diagram functors}) 
induces a $2$-functor $\underline{\Comp}({(-)}^{\cI})\colon
\underline{\Comp}(\underline{\ExCat})\to 
\underline{\Comp}(\underline{\ExCat})$ 
by Lemma~\ref{lem:functoriality of Comp}. 
In particular for a complicial exact category $\cC$, 
$\cC^{\cI}$ the category of $\cI$-diagrams in $\cC$ 
has the complicial structure induced 
from the functor $\underline{\Comp}({(-)}^{\cI})$ which we call a 
{\it level complicial structure} ({\it on $\cC$}). 
\end{para}

\begin{para}
\label{lemdf:cF}
{\bf Lemma-Definition.}\ 
(\cf \cite[2.21, 2.33, 2.35]{Moc10}.)\ \ 
{\it
Let $\cC$ be a locally exact $2$-category. 

\sn
$\mathrm{(1)}$ 
Let $(f,c)\colon x \to x'$ be a complicial $1$-morphism 
in $\cC$. 
Then there exists a $2$-isomorphism 
$c_F\colon Ff\isoto fF$ for any $F\in\{T,T^{-1},P\}$ 
which characterized by the following equalities:
\begin{equation}
\label{eq:cT}
c_T\cdot (\pi\ast f)= (f\ast \pi)\cdot c,
\end{equation}
\begin{equation}
\label{eq:cT-1}
c_{T^{-1}}:=(\beta^{-1}\ast fT^{-1})\cdot 
(T^{-1}\ast c_T^{-1} \ast T^{-1})\cdot (T^{-1}f\ast\alpha^{-1}), 
\end{equation}
\begin{equation}
\label{eq:cP}
c_P:=(c\ast T^{-1})\cdot (C\ast c_{T^{-1}}).
\end{equation}
We also have the following equalities:
\begin{equation}
\label{eq:cp compati1}
c_P\cdot (j_{\cC'}\ast f)=(f\ast j_{\cC})\cdot c_{T^{-1}}, 
\end{equation}
\begin{equation}
\label{eq:cp compati2}
(f\ast q_{\cC})\cdot c_P= q_{\cC'}\ast f.
\end{equation}

\sn
$\mathrm{(2)}$ 
Let $(f,c)$, $(g,d) \colon x \to x'$ be complicial $1$-morphisms 
between complicial objects 
and let $\phi\colon (f,c)\to (g,d)$ 
be a complicial $2$-morphism in $\cC$. 
Then for any $F\in\{T,T^{-1},P \}$, we have the equality
\begin{equation}
\label{eq:nat and cF}
d_F\cdot (F\ast \phi)=(\phi\ast F)\cdot c_F.
\end{equation}

\sn
$\mathrm{(3)}$ 
Let $(f,c)\colon x\to x'$ and $(g,d)\colon x'\to x''$ be 
complicial $1$-morphisms between complicial objects in $\cC$. 
Then for any $F\in\{T,T^{-1},P \}$, we have the equality
\begin{equation}
\label{eq:comp and cF}
{(d\odot_C c)}_F=d_F\odot_F c_F.
\end{equation}

\sn
$\mathrm{(4)}$ 
For an $1$-morphism $F\colon x\to x$ in $F\in\{T,T^{-1},P\}$, 
by applying $\mathrm{(1)}$ to the complicial $1$-morphism 
$(C,\sigma)\colon x\to x$, 
there exists a $2$-isomorphism $\sigma_F\colon FC\isoto CF$. 
Then the pair $(F,\sigma_F^{-1})$ is a complicial $1$-morphism 
from $x$ to $x$. 
Namely we have the equalities:
\begin{equation}
\label{eq:can comp func 1}
\sigma_T^{-1}\cdot (\iota\ast T)=T\ast\iota
\end{equation}
\begin{equation}
\label{eq:can comp func 2}
\sigma_{T^{-1}}^{-1}\cdot (\iota\ast T^{-1})=T^{-1}\ast\iota
\end{equation}
\begin{equation}
\label{eq:can comp func 3}
\sigma_P^{-1}\cdot (\iota\ast P)=P\ast\iota
\end{equation}
}\qed
\end{para}

\begin{para}
\label{lemdf:tauFG}
{\bf Lemma-Definition.}\ 
(\cf \cite[2.52, 2.54]{Moc10}.)\ \ 
Let $\cC$ be a locally exact $2$-category and let $x$ 
be a complicial object. 
For any pair $(F,G)$ of $1$-morphisms from $x$ to $x$, where 
$F$, $G\in\{C,T,T^{-1},P\}$ 
by applying \ref{lemdf:cF} $\mathrm{(1)}$ for a complicial $1$-morphism 
$(F,\sigma^{-1}_F)\colon x\to x$, we define 
$\tau^{F,G}$ to be a $2$-isomorphism $FG\isoto GF$ by the formula
\begin{equation}
\label{eq:tauFGdf}
\tau^{F,G}:={(\sigma_G^{-1})}_F,
\end{equation}
where by convention, we set $\sigma_C^{-1}=\sigma^{-1}(=\sigma)$. 
Then we have the equality
\begin{equation}
\label{eq:tauFG inv}
{(\tau^{F,G})}^{-1}=\tau^{G,F}
\end{equation}
for any $F$ and $G$ in $\{C,T,T^{-1},P\}$. 
For simplicity we set $\tau^F:=\tau^{F,F}$. 

In the standard example~\ref{ex:standard complicial structure}, 
the natural equivalence $\tau^{F,G}\colon FG\isoto GF$ is given 
by the following formulas for any chain complex $x$ on an exact category 
and an integer $n$; 
\begin{equation}
\label{eq:tabular}
{\footnotesize{\tiny{
\begin{tabular}{|l||c|c|c|c|} \hline
$G\ \backslash\ F$ & $C$ & $T$ & $T^{-1}$ & $P$\\
\hline\hline
$C$ &$\begin{pmatrix}-\id_{x_{n-2}} & 0 & 0 & 0\\ 0 & 0 & \id_{x_{n-1}} & 0\\ 0 &\id_{x_{n-1}} & 0 & 0\\ 0 & 0 & 0 &\id_{x_{n}} \end{pmatrix}$ &
$\begin{pmatrix}-\id_{x_{n-2}} & 0\\ 0 &\id_{x_{n-1}} \end{pmatrix}$ & 
$\begin{pmatrix}-\id_{x_{n}} & 0\\ 0 &\id_{x_{n+1}} \end{pmatrix}$ & 
$\begin{pmatrix}\id_{x_{n-1}} & 0 & 0 & 0\\ 0 & 0 & \id_{x_{n}} & 0\\ 0 &-\id_{x_{n}} & 0 & 0\\ 0 & 0 & 0 &\id_{x_{n+1}} \end{pmatrix}$\\
\hline
$T$ & $\begin{pmatrix}-\id_{x_{n-2}} & 0\\ 0 &\id_{x_{n-1}} \end{pmatrix}$&
$-\id_{x_{n-2}}$ &$-\id_{x_{n}}$ & 
$\begin{pmatrix}\id_{x_{n-1}} & 0\\ 0 &-\id_{x_{n}} \end{pmatrix}$\\
\hline
$T^{-1}$ & $\begin{pmatrix}-\id_{x_{n}} & 0\\ 0 &\id_{x_{n+1}} \end{pmatrix}$&
$-\id_{x_{n}}$ &$-\id_{x_{n+2}}$ & 
$\begin{pmatrix}\id_{x_{n+1}} & 0\\ 0 &-\id_{x_{n+2}} \end{pmatrix}$\\
\hline
$P$ &$\begin{pmatrix}\id_{x_{n-1}} & 0 & 0 & 0\\ 0 & 0 & -\id_{x_{n}} & 0\\ 0 &\id_{x_{n}} & 0 & 0\\ 0 & 0 & 0 &\id_{x_{n+1}} \end{pmatrix}$ &
$\begin{pmatrix}\id_{x_{n-1}} & 0\\ 0 &-\id_{x_{n}} \end{pmatrix}$ & 
$\begin{pmatrix}\id_{x_{n+1}} & 0\\ 0 &-\id_{x_{n+2}} \end{pmatrix}$& 
$\begin{pmatrix}\id_{x_{n}} & 0 & 0 & 0\\ 0 & 0 & \id_{x_{n+1}} & 0\\ 0 &\id_{x_{n+1}} & 0 & 0\\ 0 & 0 & 0 &-\id_{x_{n+2}} \end{pmatrix}$\\
\hline
\end{tabular}
}}}
\end{equation}
\qed
\end{para}

\begin{para}
\label{lem:comp func and alpha beta}
{\bf Lemma.}\ \ 
{\it
Let $\cC$ be a locally exact $2$-category and 
let $(f,c)\colon x\to x'$ be a complicial $1$-morphism between 
complicial objects in $\cC$. 
Then we have the equalities:
\begin{equation}
\label{eq:c vs alpha}
\alpha\ast f=(f\ast \alpha)\cdot (c_T\ast T^{-1})\cdot(T\ast c_{T^{-1}}),
\end{equation}
\begin{equation}
\label{eq:c vs beta}
f\ast \beta=(c_{T^{-1}}\ast T)\cdot (T^{-1}\ast c_T)\cdot (\beta\ast f).
\end{equation}
$$\xymatrix{
TT^{-1}f \ar[r]_{\!\!\!\!\sim}^{T\ast c_{T^{-1}}} \ar[d]^{\wr}_{\alpha\ast f} & TfT^{-1} \ar[d]_{\wr}^{c_T\ast T^{-1}}\\
f & fTT^{-1}\ar[l]_{\sim}^{f\ast \alpha},
}\ \ \ 
\xymatrix{
f \ar[r]_{\!\!\!\!\!\!\!\!\sim}^{f\ast\beta} \ar[d]^{\wr}_{\beta\ast f} & fT^{-1}T\\
T^{-1}Tf \ar[r]^{\sim}_{T^{-1}\ast c_T} & 
T^{-1}fT \ar[u]^{\wr}_{c_{T^{-1}}\ast T}.
}$$
}
\end{para}

\begin{proof}
The equality $\mathrm{(\ref{eq:c vs alpha})}$ 
follows from the left commutative diagram below 
and the equality $\mathrm{(\ref{eq:c vs beta})}$ 
follows from the right commutative diagram 
below:
$${\footnotesize{
\xymatrix{
TT^{-1}f \ar[r]_{\!\!\!\!\!\!\!\!\sim}^{TT^{-1}f\ast \alpha^{-1}} \ar[d]^{\wr}_{\alpha\ast f} & 
TT^{-1}fTT^{-1} \ar[r]_{\sim}^{TT^{-1}\ast c_T^{-1}\ast T^{-1}} 
\ar[d]_{\wr}^{\alpha\ast fTT^{-1}} & 
TT^{-1}TfT^{-1} \ar[d]_{\wr}^{T\ast \beta^{-1}\ast fT^{-1}} \\
f & 
fTT^{-1} \ar[r]^{\sim}_{c_T^{-1}\ast T^{-1}} \ar[l]_{\sim}^{f\ast \alpha^{-1}} & 
TfT^{-1},
}\ \ 
\xymatrix{
T^{-1}Tf \ar[r]_{\sim}^{T^{-1}\ast c_T} \ar[rd]_{\id_{TT^{-1}f}} & 
T^{-1}fT \ar[r]_{\!\!\!\!\!\!\!\!\sim}^{T^{-1}f\ast\alpha^{-1}\ast T} 
\ar[d]_{\wr}^{T^{-1}\ast c_T^{-1}} &
T^{-1}fTT^{-1}T \ar[d]_{\wr}^{T^{-1}\ast c_T^{-1}\ast T^{-1}T} \\
& 
T^{-1}Tf \ar[r]_{\!\!\!\!\!\!\!\!\sim}^{T^{-1}Tf\ast\beta} & T^{-1}TfT^{-1}T \\
& f \ar[r]_{f\ast \beta}^{\!\!\!\!\sim} \ar[u]^{\beta\ast f}_{\wr} & fTT^{-1}.
\ar[u]^{\wr}_{\beta\ast fT^{-1}T} 
}}}
$$
\end{proof}

\begin{para}
\label{df:comp exact subcat}
{\bf (Complicial exact subcategories).}\ \ 
Let $\cE$ be a complicial exact category and let $\cA$ be a strict 
exact subcategory of $\cE$. 
If for any object $x$ in $\cA$, 
$Cx$ and $T^{\pm 1}x$ are also in $\cA$, 
then the restriction of the complicial structure on $\cE$ to $\cA$ is 
a complicial structure on $\cA$. 
We say that $\cA$ with this complicial structure is a {\it complicial 
exact subcategory of $\cE$}. 
\end{para}

Recall the definition of coproduct of exact categories 
from \ref{df:coproduct of exact categories}. 

\begin{para}
\label{ex:coproduct of complicial exact categories}
{\bf Example (Coproduct of complicial exact categories).}\ \ 
Let $\{\cC_i\}_{i\in\cI}$ be 
a family of complicial exact categories 
indexed by a set $\cI$. 
Then we define 
$\displaystyle{C\colon \prod_{i\in\cI}\cC_i\to \prod_{i\in\cI}\cC_i}$ 
and $\iota\colon\id_{\prod_{i\in\cI}\cC_i}\to C$, $r\colon CC\to C$ and $\sigma\colon CC\isoto CC$ to 
be a functor and natural transformations and a natural equivalence 
by setting component-wisely. 
Then the quadruple $(C,\iota,r,\sigma)$ 
is a complicial structure on $\displaystyle{\prod_{i\in\cI}\cC_i}$ and 
$\displaystyle{\bigvee_{i\in\cI}\cC_i }$ is a complicial exact subcategory of 
$\displaystyle{\prod_{i\in\cI}\cC_i}$. 

\end{para}

\begin{para}
\label{df:comp topol subcat}
{\bf Definition (Complicial topologizing subcategories).}\ \ 
Let $\cE$ be a complicial exact category and 
let $\cA$ be a full subcategory of $\cE$. 
We say that $\cA$ is a {\it complicial topologizing subcategory of $\cE$} 
if $\cA$ is topologizing (see \ref{df:toplogizing subcat}) and 
if for any object $x$ in $\cA$, $Cx$ and $T^{-1}x$ are also in $\cA$. 
Then since $\cA$ is closed under admissible quotient objects, 
for any object $x$ in $\cA$, $Tx$ is in $\cA$. 
Thus $\cA$ is a complicial exact subcategory of $\cE$. 
\end{para}

\subsection{Ordinary complicial structures}
\label{subsec:ord comp str}

In this subsection, 
we introduce the specific class of complicial structures, namely 
ordinary complicial structures. 
We will associate an ordinary complicial structure on an exact category 
with its dual complicial structure in 
Lemma-Definition~\ref{lemdf:dual comp structure}. 
We start by defining commutativity of semi-commutative unital magmads 
(recall the definition of semi-commutative 
unital magmads from \ref{df:complicial exact category}).

\begin{para}
\label{df:commutative comp str}
{\bf Definition (Commutative complicial structure).}\ \ 
We say that 
a semi-commutative unital magmad $(C,\iota,r,\sigma)$ on an object $x$ 
in a locally exact $2$-category category $\cC$ is 
{\it commutative} if we have an equality 
\begin{equation}
\label{eq:comm law}
r\cdot \sigma =r.
\end{equation}
$$\xymatrix{
CC \ar[rr]_{\sim}^{\sigma} \ar[dr]_r & & CC \ar[dl]^r\\
& C. 
}$$
We say that a complicial object is {\it commutative} 
if its complicial structure is commutative. 

The category of bounded chain complexes $\Ch_b(\cE)$ over an exact category 
$\cE$ equipped with the standard complicial structure 
in \ref{ex:standard complicial structure} is a commutative complicial 
exact category.
\end{para}

Recall the definition of the $2$-morphisms 
$\zeta\colon TC\to TT$, 
$\tau^{T^{\pm}}\colon T^{\pm}T^{\pm}\isoto T^{\pm}T^{\pm}$ 
from \ref{df:mor zeta} and 
Lemma-Definition~\ref{lemdf:tauFG} respectively. 

\begin{para}
\label{df:mor l}
{\bf ($l\colon TC\to TC$).}\ \ 
Let $x$ be a commutative complicial object in a locally exact $2$-category. 
Then there exists an involution $l$ of $TC$, 
which is characterized by the equality 
\begin{equation}
\label{eq:l char}
\sigma\cdot \zeta =\zeta \cdot l.
\end{equation}
$$\xymatrix{
TC \ar[d]_l \ar[r]^{\zeta} & CC \ar[r]^r \ar[d]_{\sigma} & C \ar@{=}[d]\\
TC \ar[r]_{\zeta} & CC \ar[r]_{r} & C.
}$$

In the standard example~\ref{ex:standard complicial structure}, 
the chain morphism 
$l\colon TC\to TC$ is given by 
\begin{equation}
\label{eq:standard l}
l=-\id_{TC}.
\end{equation}
\end{para}

\begin{para}
\label{mot:ordinary comp exact cat}
Let $\cC$ be a complicial exact category. 
We wish to define $s\colon P\to PP$ to be a 
natural transformation which makes 
the quadruple 
$(P^{\op},q^{\op},s^{\op},{\tau^P}^{\op})$ a complicial 
structure on $\cC^{\op}$ the opposite category of $\cC$. 
There are two kinds of candidates of definition of the natural transformation 
$s$. Namely we set 
\begin{equation}
\label{eq:s1 df}
s_1:=(C\ast \tau^{C,T^{-1}}\ast T^{-1})\cdot (\sigma \ast T^{-1}T^{-1})\cdot 
(\zeta\ast T^{-1}T^{-1})\cdot (\tau^{C,T}\ast T^{-1}T^{-1})\cdot 
(C\ast \tau^{T^{-1},T}\ast T^{-1})\cdot (P\ast \alpha^{-1}),
\end{equation}
\begin{equation}
\label{eq:s2 df}
s_2:=(C\ast \tau^{C,T^{-1}}\ast T^{-1})\cdot (CC\ast \tau^{T^{-1}})\cdot 
(\zeta\ast T^{-1}T^{-1})\cdot (TC\ast \tau^{T^{-1}})\cdot 
(T\ast \tau^{T^{-1},C}\ast T^{-1})\cdot (\alpha^{-1} \ast P).
\end{equation}
Then the equation 
\begin{equation}
\label{eq:piTT zetaTT=id}
(\pi\ast T^{-1}T^{-1})\cdot (\zeta\ast T^{-1}T^{-1})=\id_{T^{-1}T^{-1}}
\end{equation}
and the commutative diagrams $\mathrm{(\ref{eq:counitary})}$ 
below show the equalities
\begin{equation}
\label{eq:s1 and s2 satisfies counitary law}
(P\ast q)\cdot s_1=(q\ast P)\cdot s_2=\id_P.
\end{equation}
\begin{equation}
\label{eq:counitary}
\xymatrix{
CT^{-1}CT^{-1} \ar[d]_{CT\ast\pi\ast T^{-1}} & 
CCT^{-1}T^{-1} \ar[l]^{\sim}_{C\ast \tau^{C,T^{-1}}\ast T^{-1}} 
\ar[d]_{C\ast\pi\ast T^{-1}T^{-1}}
& 
CCT^{-1}T^{-1} \ar[r]_{\sim}^{CC\ast \tau^{T^{-1}}}
\ar[l]^{\sim}_{\sigma\ast T^{-1}T^{-1}}  
\ar[d]_{\pi\ast CT^{-1}T^{-1}}
&
CCT^{-1}T^{-1} \ar[r]_{\sim}^{C\ast \tau^{T^{-1},C}\ast T^{-1}} 
\ar[d]^{\pi\ast CT^{-1}T^{-1}}
& 
CT^{-1}CT^{-1}
\ar[d]^{\pi\ast CT^{-1}}
\\
CT^{-1}TT^{-1} \ar[r]^{\sim}_{C\ast\tau^{T^{-1},T}\ast T^{-1}} & 
CTT^{-1}T^{-1} \ar[r]^{\sim}_{\tau^{C,T}\ast T^{-1}T^{-1}} & 
TCT^{-1}T^{-1} &
TCT^{-1}T^{-1} \ar[l]_{\sim}^{TC\ast \tau^{T^{-1}}}& 
TT^{-1}CT^{-1} \ar[l]_{\sim}^{T\ast \tau^{T^{-1},C}\ast T^{-1}}.\\
}
\end{equation}
We propose a sufficient condition which provides an equality $s_1=s_2$. 
\end{para}

\begin{para}
\label{df:ordinary comp exact cat}
{\bf Definition (Ordinary complicial structure).}\ \ 
A complicial structure on an object in a locally exact $2$-category 
is {\it ordinary} if 
it is commutative and if the following equality holds 
\begin{equation}
\label{eq:l and taut-1}
l\ast \tau^{T^{-1}}=\id_{TCT^{-1}T^{-1}}.
\end{equation}
A complicial structure on an object is {\it strictly ordinary} if 
it is commutative and we have the equalities 
$\mathrm{(\ref{eq:standard l})}$ and 
\begin{equation}
\label{eq:ordinary 1}
\tau^T=-\id_{TT}, \text{ and}
\end{equation}
\begin{equation}
\label{eq:ordinary 2}
\tau^{T^{-1}}=-\id_{T^{-1}T^{-1}}.
\end{equation}
We can show that a strictly ordinary complicial structure is ordinary. 
We say that a complicial object is {\it ordinary} 
(resp. {\it strictly ordinary}) 
if its complicial structure is ordinary (resp. strictly ordinary). 
For example, the standard complicial structure 
on $\Ch_b(\cE)$ the category of bounded chain complexes over 
an exact category $\cE$ (see \ref{ex:standard complicial structure}) is 
strictly ordinary. 
\end{para}

\begin{para}
\label{lemdf:dual comp structure}
{\bf Lemma-Definition.}\ 
(\cf \cite[2.60]{Moc10}.)\ \ 
Let $\cC$ be an ordinary complicial exact category. Then

\sn
$\mathrm{(1)}$ 
We have the equality
\begin{equation}
\label{eq:s1=s2}
s_1=s_2
\end{equation}
of natural transformations $P\to PP$.

\sn
$\mathrm{(2)}$ 
We define $s\colon P\to PP$ to be a natural transformation by setting 
\begin{equation}
\label{eq:s df}
s:=s_1(=s_2). 
\end{equation}
Then the quadruple $(P,q,s,\tau^P)$ is a commutative complicial structure on 
$\cC^{\op}$ the opposite category of $\cC$ and we call it 
the {\it dual complicial structure of $(C,\iota,r,\sigma)$}.

\sn
$\mathrm{(3)}$ 
We have the equality 
\begin{equation}
\label{eq:s and r}
s\cdot (r\ast T^{-1})=(r\ast T^{-1}P)\cdot (C\ast s).
\end{equation}
\end{para}

\begin{proof}
$\mathrm{(1)}$ 
By applying Lemma~\ref{lem:comp func and alpha beta} to the complicial functor 
$(P,\sigma_P^{-1})\colon \cC\to \cC$, 
we obtain the equality 
\begin{equation}
\label{eq:Pandalpha}
(\tau^{P,T}\ast T^{-1})\cdot(P\ast\alpha^{-1})=
(T\ast \tau^{T^{-1},P})\cdot (\alpha^{-1}\ast P).
\end{equation}
Namely we have the equality
\begin{equation}
\label{eq:multline}
(\tau^{C,T}\ast T^{-1}T^{-1})\cdot (c\ast \tau^{T^{-1},T}\ast T^{-1})\cdot (P\ast \alpha^{-1}) =
(TC\ast \tau^{T^{-1}})\cdot (T\ast \tau^{T^{-1},C}\ast T^{-1})
\cdot (\alpha^{-1}\ast P).
\end{equation} 
On the other hands, by the equality $\mathrm{(\ref{eq:l and taut-1})}$, 
we have the equalities
\begin{multline} 
(CC\ast\tau^{T^{-1}})\cdot (\zeta\ast T^{-1}T^{-1})
=(\zeta\ast T^{-1}T^{-1})\cdot (TC\ast \tau^{T^{-1}})\\
=(\zeta\ast T^{-1}T^{-1}) \cdot (l\ast T^{-1}T^{-1})
=(\sigma\ast T^{-1}T^{-1})\cdot (\zeta\ast T^{-1}T^{-1}).
\end{multline}
Thus we obtain the equality $\mathrm{(\ref{eq:s1=s2})}$. 

\sn
$\mathrm{(2)}$ 
We have the commutative diagram below:
$$\xymatrix{
P \ar[d]^{\wr}_{P\ast \alpha^{-1}} \ar@{=}[rrrr] & & & & 
P \ar[d]_{\wr}^{\alpha^{-1}\ast P}\\
CT^{-1}TT^{-1} \ar[r]_{\sim}^{C\ast\tau^{T^{-1},T}\ast T^{-1}} & 
CTT^{-1}T^{-1} \ar[r]_{\sim}^{\tau^{C,T}\ast T^{-1}T^{-1}} & 
TCT^{-1}T^{-1} \ar[d]_{\zeta\ast T^{-1}T^{-1}} &
TCT^{-1}T^{-1} \ar[l]^{\sim}_{TC\ast \tau^{T^{-1}}}& 
TT^{-1}CT^{-1} \ar[l]^{\sim}_{T\ast \tau^{T^{-1},C}\ast T^{-1}}\\
CT^{-1}CT^{-1} & 
CCT^{-1}T^{-1} \ar[l]_{\sim}^{C\ast \tau^{C,T^{-1}}\ast T^{-1}} & 
CCT^{-1}T^{-1} \ar[r]^{\sim}_{CC\ast \tau^{T^{-1}}}
\ar[l]_{\sim}^{\sigma\ast T^{-1}T^{-1}}  &
CCT^{-1}T^{-1} \ar[r]^{\sim}_{C\ast \tau^{T^{-1},C}\ast T^{-1}} & 
CT^{-1}CT^{-1} \\
PP \ar@{=}[u] \ar[rrrr]_{\tau^P}^{\sim} & & & & PP \ar@{=}[u]
}$$
The diagram implies the equality $\tau^P\cdot s =\tau^P$. 
We apply the equality $\mathrm{(\ref{eq:cp compati2})}$ to 
$(f,c)=(P,\sigma_P^{-1})$, we obtain the equality
\begin{equation}
\label{eq:semi commutativity op}
q\ast P=(P\ast q)\cdot \tau^P.
\end{equation}
Thus it turns out that 
the quadruple $(P,q,s,\tau^P)$ 
is a commutative complicial structure on the opposite 
category $\cC^{\op}$ of $\cC$ by Lemma-Definition~\ref{lemdf:cF} and Lemma-Definition~\ref{lemdf:tauFG}.

\sn
$\mathrm{(3)}$ 
By the virtue of the equality $\mathrm{(\ref{eq:diagonalizable 4})}$ which will be proven 
in Lemma~\ref{lem:fundament equalities in comp exact cat}, 
we have the following commutative diagram and it implies 
the equality $\mathrm{(\ref{eq:s and r})}$. 
Notice that there is no circular reasoning in our argument. 
$${\footnotesize{\xymatrix{
CP \ar[r]_{\sim}^{CP\ast \alpha^{-1}} \ar[d]^{r\ast T^{-1}} & 
CPTT^{-1} \ar[r]_{\sim}^{CC\ast \tau^{T^{-1},T}\ast T^{-1}} 
\ar[d]^{r\ast T^{-1}TT^{-1}} & 
CCTT^{-1}T^{-1} \ar[r]_{\sim}^{C\ast \tau^{C,T}\ast T^{-1}T^{-1}} 
\ar[d]^{r\ast TT^{-1}T^{-1}} & 
CTCT^{-1}T^{-1} \ar[r]^{C\ast\zeta\ast T^{-1}T^{-1}}  & 
CCCT^{-1}T^{-1} \ar[r]_{\sim}^{CCC\ast \tau^{T^{-1}}} \ar[d]_{r\ast PT^{-1}} & 
CCCT^{-1}T^{-1} \ar[r]_{\sim}^{CC\ast \tau^{C,T^{-1}}\ast T^{-1}} \ar[d]_{r\ast PT^{-1}} & 
CPP \ar[d]_{r\ast T^{-1}P}\\
P \ar[r]^{\sim}_{P\ast \alpha^{-1}} & 
PTT^{-1} \ar[r]^{\sim}_{C\ast \tau^{T^{-1},T}\ast T^{-1}} & 
CTT^{-1}T^{-1} \ar[r]^{\sim}_{\tau^{C,T}\ast T^{-1}T^{-1}} & 
TCT^{-1}T^{-1} \ar[r]_{\zeta\ast T^{-1}T^{-1}} & 
CCT^{-1}T^{-1} \ar[r]^{\sim}_{CC\ast \tau^{T^{-1}}} & 
CCT^{-1}T^{-1} \ar[r]^{\sim}_{C\ast \tau^{C,T^{-1}}\ast T^{-1}} & 
PP.
}}}$$
\end{proof}

\subsection{Frobenius exact structures}
\label{subsec:Frob exact str}

In this subsection, 
we recall Frobenius exact structure on complicial exact categories. 
We start by recalling Frobenius categories. 

\begin{para}
\label{df:Frobenius cat}
{\bf (Frobenius categories).}\ \ 
A {\it Frobenius category} is an exact category such that 
it has enough injective objects and projective objects and 
the class of injective objects and projective objects coincided. 
The naming `Frobenius' comes from the following Frobenius theorem: 

{\it A category of modules over a group ring of a finite group over 
a commutative field is a Frobenius category.}

We give a typical example of Frobenius categories. 
Let $\cA$ be an additive category and let 
$\Ch_b(\cA)$ be the category of bounded chain complexes on $\cA$. 
We regard $\Ch_b(\cA)$ as an exact category by declaring 
degree-wised split short exact sequences to be admissible exact sequences. 
Then $\Ch_b(\cA)$ is a Frobenius category and a chain complex in $\Ch_b(\cA)$ 
is projective-injective if and only if it is contractible, namely 
chain homotopic to the zero complex. 
\end{para}

\begin{para}
\label{df:Frobenius structure}
{\bf (Frobenius exact structure).}\ \ 
(See \cite[2.26, 2.27, 2.30, 2.40, 2.46]{Moc10}, \cite[6.5, B.16]{Sch11}, \cite[3.6]{Moc13b}.)\ 
Let $\cE$ be a complicial exact category. 
An admissible monomorphism $x\overset{i}{\rinf}y$ 
in $\cE$ is {\it Frobenius} if 
for any object $u$ of $\cE$ and any morphism $x\onto{f}Cu$, 
there exists a morphism $y\onto{g}Cu$ such that $f=gi$. 
Similarly an admissible epimorphism $x\overset{p}{\rdef}y$ 
in $\cE$ is {\it Frobenius} if 
for any object $u$ of $\cE$ and any morphism $Cu\onto{f}y$, 
there exists a morphism $Cu\onto{g}x$ such that $f=pg$. 

In an admissible exact sequence 
\begin{equation}
\label{eq:adm exact seq}
x\overset{i}{\rinf} y \overset{p}{\rdef} z
\end{equation}
in $\cE$, $i$ is Frobenius if and only if $p$ is Frobenius. 
In this case we call the sequence $\mathrm{(\ref{eq:adm exact seq})}$ 
a {\it Frobenius admissible exact sequence}. 
We denote $\cE$ with the Frobenius exact structure by $\cE_{\frob}$ 
and we can show that $\cE_{\frob}$ is a Frobenius category. 
Moreover, 
we can show that $\cE_{\frob}$ together with the original complicial structure is a complicial exact category. 
\end{para}

\section{Mapping cone and homotopy commutative diagrams on complicial exact categories}
\label{sec:Mapping cone}

The purpose of this section is to establish a theory 
of higher homotopical structures on complicial exact categories. 
Many notions on a category of chain complexes on an additive category 
are generalized to those on a complicial exact category. 
In the first subsection \ref{subsec:C-hom and P-hom}, 
we will generalize the notion of chain homotopies and 
in the next subsection \ref{subsec:homotopy commutative squares and diagrams}, 
we will introduce the category of homotopy commutative diagrams on 
a normal ordinary complicial exact category. 
In the last subsection \ref{subsec:mapp cyl and mapp cone}, 
we will study mapping cone on 
a complicial exact category.

\subsection{$C$-homotopy and $P$-homotopy}
\label{subsec:C-hom and P-hom}

We can generalize the notion of chain homotopies in 
a category of chain complexes on an additive category to 
a complicial exact category. 

\begin{para}
\label{df:C-homotopic}
{\bf ($C$-homotopy, $P$-homotopy).}\ \ 
For a complicial exact category $\cC$, 
morphisms $f$, $g\colon x\to y$ in $\cC$ are ({\it $C$-}){\it homotopic} 
if there exists a morphism $H\colon Cx \to y$ such that 
$f-g=H\iota_x$. 
We denote this situation by $H\colon f\Rightarrow_C g$ and 
we say that $H$ is a {\it $C$-homotopy from $f$ to $g$}. 
A morphism $f\colon x\to y$ in $\cC$ is 
a ({\it $C$-}){\it homotopy equivalence} 
if there exists a morphism $g \colon y \to x$ such that 
$gf$ and $fg$ are homotopic to $\id_x$ and $\id_y$ respectively. 
Then we say that $x$ and $y$ are ({\it $C$-}){\it homotopy equivalent}. 

Similarly we can define the notion of {\it $P$-homotopy}. 
Namely a pair of morphisms $f,\ g\colon x\to y$ in $\cC$ 
is {\it $P$-homotopic} 
if there exists a morphism $H\colon x\to Py$ such that $q_yH=f-g$. 
In this situation, we denote by $H\colon f\Rightarrow_P g$ and 
call $H$ a {\it $P$-homotopy} ({\it from $f$ to $g$}). 

We denote the set of $C$-homotopies from $f$ to $g$ by $\Hom_C(f,g)$ 
and the set of $P$-homotopies from $f$ to $g$ by $\Hom_P(f,g)$.

Assume that $\cC$ is ordinary 
(see Definition~\ref{df:ordinary comp exact cat}). 
Recall the definition of the natural transformation 
$s\colon P\to PP$ from Lemma-Definition~\ref{lemdf:dual comp structure}. 
Then 
first we define $s'\colon C\to PC$ to be a natural transformation 
by setting
\begin{equation}
\label{eq:sdash def}
s':=(PC\ast \beta^{-1})\cdot (s\ast T)\cdot (C\ast \beta).
\end{equation}
Second we define $A_{f,g}\colon\Hom_C(f,g)\to\Hom_P(f,g)$ and 
$B_{f,g}\colon\Hom_P(f,g)\to\Hom_C(f,g)$ to be maps by sending 
a $C$-homotopy $H\colon f\Rightarrow_C g$ to a $P$-homotopy
from $f$ to $g$
\begin{equation}
A_{f,g}(H):=PH\cdot s'_x\cdot \iota_x
\end{equation}
and sending a $P$-homotopy $H'\colon f\Rightarrow_P g$ 
to a $C$-homotopy from $f$ to $g$
\begin{equation}
B_{f,g}(H'):=q_y\cdot (r\ast T^{-1}_y)\cdot CH'.
\end{equation}
\end{para}

\begin{para}
\label{ex:C-homotopy}
{\bf Example.}\ \ 
Let $\Ch_b(\cE)$ be the category of bounded chain complexes on 
an exact category $\cE$ equipped with the standard complicial 
structure in \ref{ex:standard complicial structure} and 
let $f$, $g\colon x\to y$ be a pair of 
chain morphisms in $\Ch_b(\cE)$. 
Recall that a {\it chain homotopy from $f$ to $g$} is 
a family of morphisms $\{h_n\colon x_n\to y_{n+1} \}_{n\in\bbZ}$ in $\cE$ indexed by the set of integers 
such that it satisfies the equality
\begin{equation}
\label{eq:chain homotopy}
d^y_{n+1}h_n+h_{n-1}d_n^x=f_n-g_n
\end{equation}
for any integer $n$. 
Then we define $H\colon Cx\to y$ and $H'\colon x\to Py$ 
to be chain morphisms by setting 
\begin{equation}
\label{eq:Hn df}
H_n:=\begin{pmatrix}-h_{n-1}& f_n-g_n \end{pmatrix} 
\end{equation}
and 
\begin{equation}
\label{eq:H'n df}
{H'}_n:=\begin{pmatrix}f_n-g_n\\ -h_n \end{pmatrix}
\end{equation}
for any integer $n$. 
We can show that $H$ is a $C$-homotopy from $f$ to $g$ and 
$H'$ is a $P$-homotopy form $f$ to $g$. 
Conversely if we give a $C$-homotopy or a $P$-homotopy from 
$f$ to $g$, 
it provides a chain homotopy from $f$ to $g$. 
Thus the notion of chain homotopies, $C$-homotopies and 
$P$-homotopies are equivalent in the standard example of 
the complicial exact category $\Ch_b(\cE)$. 
\end{para}

In the general situation, the relationship of $C$-homotopies 
and $P$-homotopies is summed up with the following. 

\begin{para}
\label{lem:C hom equiv to P hom}
{\bf Lemma.}\ \ 
{\it
Let $f$, $g\colon x\to y$ be a pair of morphisms 
in a complicial exact category $\cC$. Then

\sn
$\mathrm{(1)}$ 
$f\Rightarrow_C g$ if and only if $f\Rightarrow_P g$. 

\sn
$\mathrm{(2)}$ 
Moreover we assume that $\cC$ is ordinary. 
Then for the maps $A_{f,g}\colon \Hom_C(f,g)\to \Hom_P(f,g)$ and 
$B_{f,g}\colon \Hom_P(f,g)\to \Hom_C(f,g)$, we have equalities 
$B_{f,g}A_{f,g}=\id_{\Hom_C(f,g)}$ and $A_{f,g}B_{f,g}=\id_{\Hom_P(f,g)}$.
}
\end{para}

\begin{proof}
$\mathrm{(1)}$ 
If there exists a morphism $H\colon Cx\to y$ such that 
$H\iota_x=f-g$, then since $q_y$ is a Frobenius admissible 
epimorphism, there exists a morphism $H'\colon Cx\to Py$ 
such that $q_yH'=H$. 
Then $q_yH'\iota_x=H\iota_x=f-g$. 
Thus $f\Rightarrow_P g$.
$$\xymatrix{
x \ar[r]^{H''} \ar[d]_{\iota_x} & Py \ar[d]^{q_y}\\
Cx \ar[r]_H \ar@{-->}[ru] & y.
}$$
Conversely if there is a morphism $H''\colon x\to Py$ such that 
$q_yH''=f-g$, then since $\iota_x$ is a Frobenius admissible morphism, 
there exists a morphism $H'\colon Cx\to Py$ such that 
$H'\iota_x=H''$. 
Then $q_yH'\iota_x=q_yH''=f-g$. 
Thus $f\Rightarrow_C g$. 

\sn
$\mathrm{(2)}$ 
Let $H\colon f\Rightarrow_{C}g$ and $H'\colon f\Rightarrow_P g$ be 
a $C$-homotopy and a $P$-homotopy from $f$ to $g$ respectively. 
Then we have the commutative diagrams below by virtue of Lemma~\ref{lem:s'}. 
$${\footnotesize{\xymatrix{
Cx \ar[r]^{C\ast\iota_x} \ar[rd] & 
CCx \ar[r]^{C\ast s'_x} \ar[d]_{r_x} & 
CPCx \ar[d]_{r\ast T^{-1}Cx} \ar[r]^{CPH} & 
CPy \ar[d]^{r\ast T^{-1}_y}\\
& Cx \ar[r]^{s'_x} \ar[rd]_{\id_{Cx}} & 
Pcx \ar[r]^{PH} \ar[d]_{q\ast C_x} & 
Py \ar[d]^{q_y}\\
& & Cx \ar[r]_H & y,
}
\xymatrix{
& x \ar[r]^{H'}\ar[ld]_{\iota_x} \ar[d]_{\iota_x} & 
Py \ar[d]_{\iota\ast P_y} \ar[rd]^{\id_{Py}}\\
Cx \ar[r]_{\id_{Cx}} \ar[rd]_{s'_x} & Cx \ar[r]^{CH'} & 
CPy \ar[r]^{r\ast T^{-1}_y} &
Py\\
& PCx \ar[r]_{PCH'}\ar[u]_{q\ast C_x} & 
PCPy \ar[u]_{q\ast CP_y} \ar[r]_{P\ast r\ast T^{-1}_y} & 
PPy \ar[u]_{q\ast P_y}.
}}}$$
The left diagram implies $B_{f,g}A_{f,g}(H)=H$ and 
the right diagram shows $A_{f,g}B_{f,g}(H')=H'$. 
We complete the proof. 
\end{proof}

\begin{para}
\label{lem:s'}
{\bf Lemma.}\ \ 
{\it
Let $\cC$ be an ordinary complicial exact category. 
Then we have the following equalities.
\begin{equation}
\label{eq:s' 1}
(q\ast C)\cdot s'=\id_C,
\end{equation}
\begin{equation}
\label{eq:s' 2}
(r\ast T^{-1}C)\cdot (C\ast s')=s'\cdot r.
\end{equation}
}
\end{para}

\begin{proof}
By Lemma-Definition~\ref{lemdf:dual comp structure}, 
we have the following commutative diagrams and they imply 
the equalities $\mathrm{(\ref{eq:s' 1})}$ and $\mathrm{(\ref{eq:s' 2})}$.
$${\footnotesize{\xymatrix{
C \ar[rrr]^{C\ast s'} \ar[rd]^{CC\ast\beta} \ar[ddd]_r & & & 
CPC \ar[ld]_{CPC\ast \beta} \ar[ddd]^{r\ast T^{-1}C}\\
& C \ar[r]^{C\ast s\ast T} \ar[d]_{r\ast T^{-1}T} & 
CPPT \ar[d]^{r\ast T^{-1}PT} & \\
&  PT \ar[r]_{s\ast T} \ar[ld]_{C\ast\beta^{-1}}& PPT \ar[rd]^{P\ast \beta^{-1}}&\\
C \ar[rrr]_{s'} & & & PC,
}
\xymatrix{
C \ar[r]^{s'} \ar[d]^{\wr}_{C\ast\beta} & 
PC \ar[r]^{q\ast C} \ar[d]_{\wr}^{PC\ast\beta} & 
C \ar[d]^{C\ast \beta}_{\wr}\\
PT \ar[r]_{s\ast T} \ar@/_2pc/[rr]_{\id_{PT}}& PPT\ar[r]_{q\ast T} & PT.
}}}$$
\end{proof}

\begin{para}
\label{df:C-contractible}
{\bf ($C$-contractible).}\ \ 
(\cf \cite[2.7, 2.24, 2.44]{Moc10}, \cite[B.16]{Sch11}.)\ 
Let $\cC$ be a complicial exact category. 
We say that an object $x$ in $\cC$ is {\it $C$-contractible} 
if it is $C$-homotopy equivalent to the zero object. 
We can show that $x$ is $C$-contractible 
if and only if $x$ is a direct summand of $Cu$ for some object $u$ in $\cC$ 
if and only if $x$ is a retraction of $Cu$ for some object $u$ in $\cC$ 
and if and only if $x$ is a projective-injective object in $\cC_{\frob}$. 
\end{para}

We can prove the following lemma.

\begin{para}
\label{lem:left and right homotopy inverse}
{\bf Lemma.}\ \ 
{\it 
A morphism in a complicial exact category which admits both 
left and right $C$-homotopy inverses is a $C$-homotopy equivalence. 
More precisely, 
let $f\colon x\to y$ and $g$, $h\colon y\to x$ be morphisms 
in a complicial exact category and let 
$H\colon hf\Rightarrow_C \id_x$ and $K\colon fg\Rightarrow_C \id_y$ 
be $C$-homotopies. 
Then $-hK+HCg$ is a $C$-homotopy from $h$ to $g$ and 
$-fhK+fHCg+K$ is a $C$-homotopy from $fh$ to $\id_y$. 
In particular $f$ is a $C$-homotopy equivalence.    
}
\qed
\end{para}

\begin{para}
\label{lem:split Frob exact}
{\bf Lemma.}\ \ 
{\it
Let $\cC$ be a complicial exact category and $f\colon x\to y$ 
a morphism in $\cC$. 
Then\\
$\mathrm{(1)}$ 
If $f$ is a Frobenius admissible monomorphism, then $C(f)$ is 
a split monomorphism. 

\sn
$\mathrm{(2)}$ 
If $f$ is a Frobenius admissible epimorphism, then $C(f)$ is 
a split epimorphism.
}
\end{para}

\begin{proof}
$\mathrm{(1)}$ 
Since 
$C(f)$ is a Frobenius admissible monomorphism and $C(x)$ is an injective 
object in $\cC_{\frob}$ by \ref{df:Frobenius structure} 
and by \ref{df:C-contractible} respectively, 
$C(f)$ is a split monomorphism. 
A proof of $\mathrm{(2)}$ is similar. 
\end{proof}

\subsection{Homotopy commutative squares and diagrams}
\label{subsec:homotopy commutative squares and diagrams}

In this subsection we introduce concepts of homotopy commutative squares and 
commutative diagrams in complicial exact categories. 
We will prove that the category of homotopy commutative diagrams on a 
normal ordinary 
complicial exact category becomes a complicial exact category 
(see \ref{lemdf:additive str on Fbh}). 
We start by defining a notion of homotopy commutative squares.

\begin{para}
\label{df:homotopy commutative square}
{\bf (Homotopy commutative squares).}\ \ 
Let $\cC$ be a complicial exact category and let 
$[f\colon x\to x']$ and 
$[g\colon y \to y']$ be 
a pair of objects in $\cC^{[1]}$ the morphisms category of $\cC$. 
A ({\it $C$-}){\it homotopy commutative square} 
({\it from $[f\colon x\to x']$ to 
$[g\colon y\to y']$}) is a triple ${(a,b,H)}_C$ consisting of 
morphisms $a\colon x\to y $, $b\colon x'\to y'$ and 
$H\colon Cx\to y' $ in $\cC$ such that 
$H\iota_x=ga-bf$. 
Namely $H$ is a $C$-homotopy from $ga$ to $bf$. 
We often write $(a,b,H)$ for ${(a,b,H)}_C$. 
Similarly a {\it $P$-homotopy commutative square} 
({\it from $[f\colon x\to x']$ to 
$[g\colon y\to y']$}) is a triple ${(a,b,H)}_P$ consisting of 
morphisms $a\colon x\to y $, $b\colon x'\to y'$ and 
$H\colon x\to Py' $ in $\cC$ such that 
$q_yH=ga-bf$. 

Notice that for a pair of morphisms $a$, $b\colon x\to y$ 
and for a $C$-homotopy 
$H\colon a\Rightarrow_C b$, 
we can regard it as a $C$-homotopy commutative 
square $[x\onto{\id_x}x]\onto{(a,b,H)}[y\onto{\id_y}y]$. 
For simplicity We denote this situation by $x\onto{(a,b,H)}y$. 

Let $[f\colon x\to x']$, $[g\colon y\to y']$, 
$[f'\colon x'\to x'']$, $[g'\colon y'\to y'']$ 
and 
$[h\colon z\to z']$ be objects in $\cC^{[1]}$ 
and let ${(a,b,H)}_C$, ${(a',b',H')}_C$ 
and ${(b,c,K)}_C$ 
be $C$-homotopy commutative squares 
from $[f\colon x\to x']$ to $[g\colon y\to y']$,  
from $[g\colon y\to y']$ to $[h\colon z\to z']$ 
and from 
$[f'\colon x' \to x'']$ to $[g'\colon y'\to y'']$ 
respectively. 
Then we define ${(a',b',H')}_C{(a,b,H)}_C$ to be a homotopy commutative square 
from $[f\colon x\to x']$ to $[h\colon z\to z']$ by setting
\begin{equation}
\label{eq:comp of hcs}
{(a',b',H')}_C{(a,b,H)}_C:={(a'a,b'b,H'\star_C H)}_C
\end{equation}
where $H'\star_C H$ is a $C$-homotopy from $ha'a$ to $b'bf$ 
given by the formula 
\begin{equation}
\label{eq:homotopy star}
H'\star_C H:=b'H+H'Ca.
\end{equation}
We also define $K\bullet_C H$ to be a $C$-homotopy $cf'f\Rightarrow_C g'ga$ 
by the formula
\begin{equation}
\label{eq: homotopy bullet}
K\bullet_C H:=KCf+ g'H.
\end{equation}
We often write $(a,b,H)$, $H'\star H$ and $K\bullet H$ for 
${(a,b,H)}_C$, 
$H'\star_C H$ and 
$K\bullet_C H$ respectively. 

Next let ${(a,b,H)}_P$, ${(a',b',H')}_P$ and ${(b,c,K)}_P$ 
be $P$-homotopy commutative squares from 
from $[f\colon x\to x']$ to $[g\colon y\to y']$,  
from $[g\colon y\to y']$ to $[h\colon z\to z']$ 
and from 
$[f'\colon x' \to x'']$ to $[g'\colon y'\to y'']$ 
respectively. 
Then we define ${(a',b',H')}_P{(a,b,H)}_P$ to be a homotopy commutative square 
from $[f\colon x\to x']$ to $[h\colon z\to z']$ by setting
\begin{equation}
\label{eq:comp of hcs P}
{(a',b',H')}_P{(a,b,H)}_P:={(a'a,b'b,H'\star_P H)}_P
\end{equation}
where $H'\star_P H$ is a $P$-homotopy from $ha'a$ to $b'bf$ 
given by the formula 
\begin{equation}
\label{eq:homotopy star P}
H'\star_P H:=Pb'H+H'a.
\end{equation}
We also define $K\bullet_P H$ to be a $P$-homotopy $cf'f\Rightarrow_P g'ga$ 
by the formula
\begin{equation}
\label{eq:homotopy bullet P}
K\bullet_P H:=Kf+ Pg'H.
\end{equation}

We define $\cC_h^{[1]}$ to be a category whose objects are 
morphisms in $\cC$ and whose morphisms are 
homotopy commutative squares and compositions of morphisms 
are give by the formula $\mathrm{(\ref{eq:comp of hcs})}$ and 
we define $\cC^{[1]}\to \cC_h^{[1]}$ 
to be a functor by sending an object $[f\colon x\to x']$ to $[f\colon x\to x']$ 
and a morphism $(a,b)\colon[f\colon x\to x']\to [g\colon y\to y'] $ to 
$(a,b,0)\colon[f\colon x\to x']\to [g\colon y\to y']$.  
By this functor, we regard $\cC^{[1]}$ as a subcategory of $\cC_h^{[1]}$. 
\end{para}

Recall the definitions of the maps $A$ and $B$ from 
\ref{df:C-homotopic}. 
We can show the following lemma.

\begin{para}
\label{lem:comp and A and B}
{\bf Lemma.}\ \ 
{\it
Let $\cC$ be an ordinary complicial exact category. 
Let $F$ be a letter in $F\in\{C,P\}$ and let 
$[f\colon x\to x']$, $[g\colon y\to y']$, 
$[f'\colon x'\to x'']$, $[g'\colon y'\to y'']$ 
and 
$[h\colon z\to z']$ be objects in $\cC^{[1]}$ 
and let ${(a,b,H)}_F$, ${(a',b',H')}_F$ 
and ${(b,c,K)}_F$ 
be $F$-homotopy commutative squares 
from $[f\colon x\to x']$ to $[g\colon y\to y']$,  
from $[g\colon y\to y']$ to $[h\colon z\to z']$ 
and from 
$[f'\colon x' \to x'']$ to $[g'\colon y'\to y'']$ 
respectively. Then 

\sn
$\mathrm{(1)}$ 
If $F=C$, then we have the following equalities
\begin{equation}
\label{eq:star and A}
A_{ha'a,b'bf}(H'\star_C H)=A_{ha',b'g}(H')\star_P A_{ga,bf}(H),
\end{equation}
\begin{equation}
\label{eq:bullet and A}
A_{g'ga,cf'f}(K\bullet_C H)=A_{g'b,cf'}(K)\bullet_P A_{ga,bf}(H).
\end{equation}

\sn
$\mathrm{(2)}$ 
If $F=P$, then we have the following equalities
\begin{equation}
\label{eq:star and B}
B_{ha'a,b'bf}(H'\star_P H)=B_{ha',b'g}(H')\star_C B_{ga,bf}(H),
\end{equation}
\begin{equation}
\label{eq:bullet and B}
B_{g'ga,cf'f}(K\bullet_P H)=B_{g'b,cf'}(K)\bullet_C B_{ga,bf}(H).
\end{equation}
}\qed
\end{para}

\begin{para}
\label{df:CC-homotopy}
{\bf Definition ($CC$-homotopy).}\ \ 
Let $\cC$ be a complicial exact category and let $f\colon x\to y$ be 
a morphism in $\cC$ and let $H$, $H'\colon f\Rightarrow_C g$ be a 
pair of $C$-homotopies. 
A {\it $CC$-homotopy} ({\it from $H$ to $H'$}) is 
a morphism $S\colon CCx\to y$ such that 
$L(C\ast\iota_x-\iota\ast Cx)=H-H'$. 
We denote this situation by $L\colon H\Rightarrow_{CC} H'$. 
We say that a pair of $C$-homotopies $H$, $H'\colon f\Rightarrow_{C}g$ 
are {\it $CC$-homotopic} if there exists a $CC$-homotopy $H\Rightarrow_{CC}H'$. As lemma below, the $CC$-homotopic relation between $C$-homotopies is an equivalence relation.

Let $(a,b,H)$, $(a,b,H')\colon [f\colon x\to x']\to [g\colon y\to y']$, 
$(d,e,K)$, $(d,e,K')\colon [g\colon y\to y']\to [h\colon z\to z']$ and 
$(b,c,L)$, $(b,c,L')\colon [f'\colon x'\to x'']\to [g'\colon y'\to y'']$ be 
homotopy commutative squares and let 
$S\colon H\Rightarrow_{CC} H'$, $T\colon K\Rightarrow_{CC}K'$, 
$U\colon L\Rightarrow_{CC} L'$ be a triple of 
$CC$-homotopies. 
Then we define $T\star_{CC}S\colon K\star_{CC} H\Rightarrow_{CC}K'\star_{CC}H'$ and $U\bullet_{CC}S\colon L\bullet_{CC}H\Rightarrow_{CC}L'\bullet_{CC}H'$ 
to be $CC$-homotopies by setting
\begin{equation}
\label{eq:star CC df}
T\star_{CC}S:=b'S+TCCa, 
\end{equation}
\begin{equation}
\label{eq:bullet CC df}
U\bullet_{CC}S:=g'S+UCCf.
\end{equation}
\end{para}

\begin{para}
\label{ex:CC-homotopy}
{\bf Example.}\ \ 
Let $\Ch_b(\cE)$ be the category of bounded chain complexes 
on an exact category $\cE$ equipped with 
the standard complicial structure 
as in \ref{ex:standard complicial structure} 
and let $f$, $g\colon x\to y$ be a pair of chain morphisms 
in $\Ch_b(\cE)$ and let 
$H$, $L\colon Cx\to y$ be a pair of $C$-homotopies 
from $f$ to $g$. 
That is, as in \ref{ex:C-homotopy}, there are 
families $\{h_n\colon x_n\to y_{n+1}\}_{n\in\bbZ}$ 
and $\{l_n\colon x_n\to y_{n+1}\}_{n\in\bbZ}$ 
of morphisms in $\cE$ indexed by set of integers such that 
$H_n=\begin{pmatrix}-h_{n-1}& f_n-g_n\end{pmatrix}$ and 
$L_n=\begin{pmatrix}-l_{n-1}& f_n-g_n\end{pmatrix}$, 
$f_n-g_n=d_{n+1}^yh_n+h_{n-1}d_n^x=d_{n+1}^yl_n+l_{n-1}d_n^x$ for each integer $n$. 
In this case, for a $CC$-homotopy $S\colon CCx\to y$ from 
$H$ to $L$, 
there exists a family of morphisms 
$\{s_n\colon x_n\to y_{n+2} \}_{n\in\bbZ}$ in $\cE$ 
indexed by the set of integers 
such that 
\begin{equation}
\label{eq:Sn and sn}
S_n=\begin{pmatrix}s_{n-2} & 0 & -h_{n-1}+l_{n-1} & 0 \end{pmatrix},
\end{equation}
\begin{equation}
\label{eq:Hn-Ln and sn}
H_n-L_n=d^y_{n+2}s_n-s_{n-1}d_n^x
\end{equation}
for any integer $n$. 
Conversely the family of morphisms 
$\{s_n\colon x_n\to y_{n+2} \}_{n\in\bbZ}$ in $\cE$ 
which satisfies the equalities 
$\mathrm{(\ref{eq:Hn-Ln and sn})}$ 
for each integer $n$ 
gives a $CC$-homotopy $S\colon CCx\to y$ from $H$ to $L$ 
by setting as in 
$\mathrm{(\ref{eq:Sn and sn})}$. 
\end{para}

\begin{para}
\label{rem:2-structure of ChbE}
{\bf Remark.}\ \ 
Let $a$, $b$, $c\colon x\to y$ and 
let $c$, $d$, $e\colon y\to z$ be 
a sextuple of morphisms in 
a complicial exact category and let 
$H\colon a\Rightarrow_C b$, 
$K\colon c\Rightarrow_C d$, 
$L\colon b\Rightarrow_C c$ and 
$M\colon d\Rightarrow_C e$ be 
a quadruple of $C$-homotopies. 
Then as in \ref{lem:CC fundamental results} 
$\mathrm{(4)}$ below, in general 
$(M\star_C L)\bullet_C(K\star_C H)$ and 
$(M\bullet_C K)\star_C(L\bullet_C H)$ are not equal but they are 
only $C$-homotopic to each other. 
Thus to regard $\Ch_b(\cE)$ the category of 
bounded chain complexes on an exact category $\cE$ 
as a $2$-category by declaring the class of chain homotopies 
to be the class of $2$-morphisms, 
we need some arrangement of the class of chain homotopies. 
Namely we adopt the class of $CC$-homotopic classes of 
chain homotopies as the class of $2$-morphisms. 
Then the compositions $\bullet_C$ and $\star_C$ 
give the vertical and horizontal compositions of $2$-morphisms 
respectively and we can make $\Ch_b(\cE)$ into a $2$-category. 
If we regard $\Ch_b(\cE)$ as a $2$-category in this way, 
we denote it by $\underline{\Ch_b(\cE)}$. 
This is a motivation of defining $\underline{F_{b,h}\cC}$ in 
\ref{df:homotopy commutative diagrams} below. 
Since the exact structure and the complicial structure on 
$\Ch_b(\cE)$ only involve $1$-categorical structures and 
$\underline{\Ch_b(\cE)}$ 
does not change a matter of objects and $1$-morphisms 
from $\Ch_b(\cE)$, 
$\underline{\Ch_b(\cE)}$ naturally becomes a complicial 
exact category. 
\end{para}

We can show the following lemma.

\begin{para}
\label{lem:CC fundamental results}
{\bf Lemma.}\ \ 
{\it
Let $\cC$ be a complicial exact category. 
Then

\sn
$\mathrm{(1)}$ 
The $CC$-homotopic relation between $C$-homotopies is an equivalence 
relation. More precisely, 
let $f$, $g\colon x\to y$ be a pair of morphisms in $\cC$ and let 
$H$, $L$, $M\colon f\Rightarrow_C g$ be a triple of $C$-homotopies 
from $f$ to $g$ and let 
$S\colon H\Rightarrow_{CC}L$ and $T\colon L\Rightarrow_{CC}M$ be 
a pair of $CC$-homotopies. Then

\sn
$\mathrm{(a)}$ 
$0\colon CCx \to y$ is a $CC$-homotopy from $H$ to $H$.

\sn
$\mathrm{(b)}$ 
$-S\colon CCx\to y$ is a $CC$-homotopy from $L$ to $H$.

\sn
$\mathrm{(c)}$ 
$S+T\colon CCx\to y$ is a $CC$-homotopy from $H$ to $M$.

\sn
$\mathrm{(2)}$ 
The $C$-homotopic relation between morphisms 
and the $CC$-homotopic relation between $C$-homotopies are 
compatible with additions of morphisms. 
Namely for $i=0$, $1$, let 
$f_i$, $g_i\colon x\to y$ be a morphisms and let 
$H_i$, $L_i\colon f_i\Rightarrow_Cg_i$ be $C$-homotopies and 
let $S_i\colon H_i\Rightarrow_{CC}L_i$ be $CC$-homotopies. 
Then $H_0+H_1$ and $L_0+L_1$ are $C$-homotopies from $f_0+f_1$ to $g_0+g_1$ 
and $S_0+S_1$ is a $CC$-homotopy from $H_0+H_1$ to $L_0+L_1$. 

\sn
$\mathrm{(3)}$ 
$C$-homotopic relation between morphisms and $CC$-homotopic relation between 
$C$-homotopies are compatible with the compositions of morphisms. 
That is, let $f\colon x\to y$, $g$, $g'\colon y\to z$, $h\colon z\to w$ be 
a quadruple of morphisms in $\cC$ and let 
$H$, $H'\colon g\Rightarrow_C g'$ be a pair of $C$-homotopies and 
let $S\colon H\Rightarrow_{CC} H'$ be a $CC$-homotopies from $H$ to $H'$. 
Then $HCf$, $H'Cf\colon gf\Rightarrow_C g'f$ and $hH$, 
$hH'\colon hg\Rightarrow_C hg'$ are $C$-homotopies and 
$SCCf\colon Hcf\Rightarrow_{CC}H'Cf$ and $hS\colon hH\Rightarrow_{CC}hH'$ 
are $CC$-homotopies. 

\sn
$\mathrm{(4)}$ 
Let $(a,b,H)\colon [u\colon x\to x'] \to [v\colon y\to y']$, 
$(d,e,K)\colon [v\colon y\to y'] \to [w\colon z\to z']$, 
$(b,c,L)\colon [u'\colon x'\to x''] \to [v'\colon y'\to y'']$ and 
$(e,f,M)\colon [v'\colon y'\to y''] \to [w'\colon z'\to z'']$ 
be a quadruple of homotopy commutative squares. 
Then the morphism $MCH\colon CCx\to z''$ is a $CC$-homotopy from 
$(M\bullet_C K)\star_C(L\bullet_C H)$ to $(M\star_C L)\bullet_C(K\star_C H)$.

\sn
$\mathrm{(5)}$ 
Let $(f,d)\colon\cC\to \cC'$ be a complicial exact functor 
between complicial exact categories $\cC$ and $\cC'$ and let 
$u$, $v\colon x\to y$ be a morphism in $\cC$ and 
let $H$, $K\colon u\Rightarrow_C v$ be a $C$-homotopy from 
$u$ to $v$ and let $S\colon H\Rightarrow_{CC}K$ be 
a $CC$-homotopy from $H$ to $K$. Then 
$fH\cdot d_x$ and $fK\cdot d_x$ are $C$-homotopies from $fu$ to $fv$ and 
$fS\cdot (d\ast Cx)\cdot (C'\ast d_x)$ is a $CC$-homotopy from $fH\cdot d_x$ to $fK\cdot d_x$. 

\sn
$\mathrm{(6)}$ 
Let $f\colon x\to y$ be a morphism 
in a complicial exact category $\cC$ and 
let $H\colon f\Rightarrow_C f$ be a $C$-homotopy 
from $f$ to $f$. Namely 
we have an equality 
$H\iota_x=0$. 
Assume $x$ is $C$-contractible. 
That is, there exists a morphism $U\colon Cx\to x$ such that 
$U\iota_x=\id_x$. Then $HCU\colon CCx\to y$ is a 
$CC$-homotopy from $H$ to $0$.
}\qed
\end{para}

\begin{para}
\label{df:filtered objects}
{\bf Definition (Filtered objects).}\ \ 
Let $\cX$ be a category with a specific zero object $0$. 
Let $\bbZ$ be the linear ordered set of all integers with 
the usual linear order. 
We regard $\bbZ$ as a category in the usual way. 
We denote $\cX^{\bbZ}$ the functor category from $\bbZ$ to $\cX$ by $F\cX$. 
We call an object $F\cX$ a {\it filtered objects} in $\cX$. 
For an object $x$ in $F\cX$, 
we write $x_n$ and $i^x_n$ for $x(n)$ and $x(n\leq n+1)$ respectively and 
illustrate $x$ as follows:
$$\cdots \onto{i^x_{n-1}} x_n \onto{i^x_n} x_{n+1} \onto{i^x_{n+1}} 
x_{n+2} \onto{i^x_{n+2}}\cdots.$$
Let $a\leq b$ be a pair of integers and $x$ an object in $F\cX$. 
We say that {\it $x$ has amplitude contained in $[a,b]$} 
if for any $k<a$, $x_k=0$ and for any $b\leq k$, 
$x_k=x_b$ and $i^x_k=\id_{x_b}$. 
In this case we write $x_{\infty}$ for $x_b=x_{b+1}=\cdots$. 
We set for $x\neq 0$, 
\begin{equation}
\label{eq:dimension of filtered object}
\dim x:=\min\{k\in\bbZ;x_n=x_{\infty}\text{ for any $n\geq k$.}\}
\end{equation}
and call it the {\it dimension of $x$}. 
By convention, we set $\dim x=-1$ for $x=0$. 
Similarly for any morphism $f\colon x\to y$ in $F\cX$ between 
objects which have amplitude contained in $[a,b]$, 
we denote $f_b=f_{b+1}=\cdots$ by $f_{\infty}$. 
We write $F_{[a,b]}\cX$ 
for the full subcategory of $F\cX$ consisting of those objects 
having amplitude contained in $[a,b]$. 
We also set $\displaystyle{F_b\cX:=\underset{\substack{a<b\\ (a,b)\in\bbZ^2}}{\bigcup} F_{[a,b]}\cX}$, $\displaystyle{F_{b,\leq b}\cX:=\underset{\substack{a<b\\ a\in\bbZ}}{\bigcup} F_{[a,b]}\cX}$ and $\displaystyle{F_{b,\geq a}\cX:=\underset{\substack{a<b\\ b\in\bbZ}}{\bigcup} F_{[a,b]}\cX}$ 
and we call an object in $F_b\cX$ a {\it bounded filtered object} ({\it in $\cX$}).

Let $u\colon x\to y$ be a morphism in $\cX$. 
We write $\fj(x)$ and $\fj(f)\colon \fj(x)\to \fj(y)$ 
for the object and the morphism in $F_{[0,0]}\cX$ such that 
${\fj(x)}_{\infty}=x$ and ${\fj(f)}_{\infty}=f$. 
We define ${(-)}_{\infty}\colon F_b\cX\to \cX$ and 
$\fj\colon\cX\to F_b\cX$ 
to be functors by sending an object $x$ in 
$F_b\cX$ to 
$x_{\infty}$ in $\cX$ and an object $y$ in $\cX$ to $\fj(y)$ 
respectively. 
Obviously we have the equality ${(-)}_{\infty}\fj=\id_{\cX}$. 

If $\cX$ is an exact category, then 
$F\cX$ with the level admissible exact sequences 
(see \ref{nt:functor category}) is an exact category and 
$F_b\cX$ and $F_{[a,b]}\cX$ are strict exact subcategories of $F\cX$. 
We also denote the full subcategory of $F_b\cX$ consisting of those objects 
$x$ such that $i^x_n$ is an admissible monomorphism for any integer $n$ 
by $F^{\rinf}_b\cC$ and 
we set $F_{[a,b]}^{\rinf}\cX=F_{[a,b]}\cX\cap F_b^{\rinf}\cX$ and so on.
By \ref{lem:FrinfbC} below, it turns out that 
$F_b^{\rinf}\cX$ and $F_{[a,b]}^{\rinf}\cX$ 
are strict exact subcategories of $F_b\cX$.

Moreover if $\cX$ is a complicial exact category, then 
$F\cX$ with the level complicial structure 
(see \ref{df:level complicial structure}) is a complicial exact category and 
$F_b\cX$, $F_{[a,b]}\cX$, $F_b^{\rinf}\cX$ and $F_{[a,b]}^{\rinf}\cX$ are complicial exact subcategories of $F\cX$. 
\end{para}

\begin{para}
\label{lem:FrinfbC}
{\bf Lemma.}\ \ 
{\it
Let $\cC$ be an exact category and let $a<b$ be a pair of integers. 
Then $F_b^{\rinf}\cC$ and $F_{[a,b]}^{\rinf}\cC$ 
are strict exact subcategory of $F_b\cC$. 
}
\end{para}

\begin{proof}
We will only give a proof for $F_b^{\rinf}\cC$ and 
a proof for $F_{[a,b]}^{\rinf}\cC$ is similar. 
Let $x\to y\to z$ be a sequence in $F_b^{\rinf}\cC$ such that 
it is an admissible exact sequence in $F_b\cC$. 
Then by \ref{lem:char of admissible square}, 
for each integer $n$, the induced 
morphism $x_n\sqcup_{x_{n-1}}y_{n-1}\to y_n$ is an admissible monomorphism in $\cC$. 
Thus by \cite[1.1.4.]{Wal85}, the class of admissible monomorphisms 
in $F_b^{\rinf}\cC$ is closed under compositions and co-base change 
by arbitrary morphisms. 

Next by applying the similar argument above to $\cC^{\op}$ 
the opposite category of $\cC$, it turns out that 
the class of admissible epimorphisms 
in $F_b^{\rinf}\cC$ is closed under compositions and base change 
by arbitrary morphisms. 
Thus $F_b^{\rinf}\cC$ is a strict exact subcategory of $F_b\cC$. 
\end{proof}

The notion of homotopy commutative diagrams below is 
a generalization of that of homotopy commutative squares 
\ref{df:homotopy commutative square}. 

\begin{para}
\label{df:homotopy commutative diagrams}
{\bf Definition (Homotopy commutative diagrams).}\ \ 
Let $\cC$ be a complicial exact category 
and $x$ and $y$ objects in $F_b\cC$. 
A {\it homotopy commutative diagram} $(f,H)$ from $x$ to $y$ consisting of 
a family of morphisms $f=\{f_n\colon x_n\to y_n\}_{n\in\bbZ}$ in $\cC$ 
indexed by $\bbZ$ and a family of $C$-homotopy 
$\{H_n\colon i_n^yf_n \Rightarrow_C f_{n+1}i^x_n \}_{n\in\bbZ}$ 
indexed by $\bbZ$ such that for sufficiently large $m$, 
$f_m=f_{m+1}=\cdots(:=f_{\infty})$ and $H_m=H_{m+1}=\cdots=0$. 
We illustrate by $(f,H)\colon x\to y$ in this situation. 
We say that a pair of homotopy commutative diagrams 
$(f,H)$, $(g,K)\colon x\to y$ are {\it $CC$-homotopic} 
if $f=g$ and if for any integer $n$, $H_n$ and $K_n$ are 
$CC$-homotopic. 
By \ref{lem:CC fundamental results}, 
we can show that the $CC$-homotopic relation between 
homotopy commutative diagrams is an equivalence relation.

Next we define the compositions of homotopy commutative diagrams. 
Let $(f,H)\colon x\to y$ and 
$(g,K)\colon y\to z$ be 
homotopy commutative diagrams between objects in $F_b\cC$. 
We define $(g,K)(f,H)\colon x\to z$ 
to be a homotopy commutative diagram by the formula
\begin{equation}
\label{eq:comp of hcd}
(g,K)(f,H):=(gf,K\star H)
\end{equation}
where $gf:=\{g_nf_n\colon x_n\to z_n \}_{n\in\bbZ}$ and 
$K\star H:=\{K_n\star_C H_n\colon i_n^zg_nf_n 
\Rightarrow_C g_{n+1}f_{n+1}i_n^x \}_{n\in\bbZ}$. 
(For the definition of the operation $\star$, 
see \ref{df:homotopy commutative square}.) 
We can show that the compositions of 
homotopy commutative diagrams 
is associative.

We define $F_{b,h}\cC$ to be a category 
whose class of objects is $\Ob F_b\cC$ 
and whose morphisms are homotopy commutative diagrams and 
whose compositions are given by compositions 
of homotopy commutative diagrams. 
We also define ${\underline{F_{b,h}\cC}}$ to be a category 
by setting $\Ob{\underline{F_{b,h}\cC}}=\Ob F_b\cC$ and 
for any pair of objects $x$ and $y$, 
\begin{equation}
\label{eq:underlineFbhC Hom df}
\Hom_{{\underline{F_{b,h}\cC}}}(x,y):=\Hom_{F_{b,h}\cC}(x,y)/
(CC\text{-homotopic relation}).
\end{equation}
Then by \ref{lem:CC fundamental results} again, 
compositions of morphisms in $F_{b,h}\cC$ induces the compositions of 
morphisms in 
${\underline{F_{b,h}\cC}}$. 

Similarly for any pair of integers $a\leq b$, we define the category 
$F_{[a,b],h}\cC$, $F_{b,\geq a,h}\cC$, $F_{b,\leq b,h}\cC$, 
$\underline{F_{[a,b],h}\cC}$, $\underline{F_{b,\geq a,h}\cC}$ and 
$\underline{F_{b,\leq b,h}\cC}$. 
For any integer $n$, 
We can regard $F_{[n,n+1]}\cC$ and $F_{[n,n+1],h}\cC$ as 
$\cC^{[1]}$ and $\cC^{[1]}_h$.

\end{para}

\begin{para}
\label{lemdf:additive str on Fbh} 
{\bf Lemma-Definition.}\ \ 
Let $\cC$ be a complicial exact category 
and let $a\leq b$ be a pair of integers. Then

\sn
$\mathrm{(1)}$ 
Let $(f,H)$, $(g,K)\colon x\to y$ be 
a pair of homotopy commutative diagrams. 
We define the addition of homotopy commutative diagrams by the formula
\begin{equation}
\label{eq:add of hcd}
(f,H)+(g,K):=(f+g,H+K).
\end{equation}
Then $F_{b,h}\cC$, $F_{[a,b],h}\cC$, 
$F_{b,\leq b,h}\cC$ and $F_{b,\geq a,h}\cC$ 
naturally become additive categories. 
Here for any objects $x$ and $y$ in $F_{b,h}\cC$, 
$x\oplus y$ is given by ${(x\oplus y)}_n=x_n\oplus y_n$ 
and $\displaystyle{i_n^{x\oplus y}=
\begin{pmatrix}i_n^x & 0\\ 0 &i_n^y\end{pmatrix}}$ for any integer $n$.

\sn
$\mathrm{(2)}$ 
We say that a sequence 
\begin{equation}
\label{eq:Frob exact in Fbh}
x\overset{(f,H)}{\rinf} y \overset{(g,K)}{\rdef} z
\end{equation} 
in $F_{b,h}\cC$ is a {\it level Frobenius admissible exact sequence} 
if we have the equality $(g,K)(f,H)=0$ 
and if for any integer $n$, the sequence 
$x_n\overset{f_n}{\rinf} y_n \overset{g_n}{\rdef} z_n$ is 
a Frobenius admissible exact sequence in $\cC$. 
In this case we call $(f,H)$ a 
{\it level Frobenius admissible monomorphism} 
and call $(g,K)$ a {\it level Frobenius admissible epimorphism}. 
Then 

\sn
$\mathrm{(a)}$ 
For a homotopy commutative diagram $(f,H)\colon x\to y$ in $F_{b,h}\cC$, 
$(f,H)$ is a level Frobenius admissible monomorphism 
if and only if for any integer $n$, 
the morphism $f_n\colon x_n\to y_n$ is a Frobenius admissible monomorphism 
in $\cC$. 
Moreover if we assume that $\cC$ is ordinary 
(see Definition~\ref{df:ordinary comp exact cat}), 
then 
similarly $(f,H)$ is a level Frobenius admissible epimorphism if and only if 
for any integer $n$, 
the morphism $f_n\colon x_n\to y_n$ 
is a Frobenius admissible epimorphism in $\cC$.

\sn
$\mathrm{(b)}$ 
If $\cC$ is ordinary, then 
$F_{b,h}\cC$, $F_{[a,b],h}\cC$, $F_{b,\leq b,h}\cC$ 
and $F_{b,\geq a,h}\cC$ 
endowed with the class of level Frobenius admissible exact sequences 
are exact categories.

\sn
$\mathrm{(3)}$ 
We define $C\colon F_{b,h}\cC\to F_{b,h}\cC$ 
to be a functor by setting 
${Cx}_n:=C(x_n)$ and $i_n^{Cx}:=C(i_n^x)$ and 
${C(f,H)}_n:=(f_n,CH_n\sigma_{x_n})$ 
for an object $x$ in $F_{b,h}\cC$ and a homotopy 
commutative diagram $(f,H)\colon x\to y$ and an integer $n$. 
Moreover we assume $\cC$ is normal and 
we also define  
$\iota\colon \id_{F_{b,h}}\to C $, 
$r\colon CC\to C$ and $\sigma\colon CC\isoto CC$ 
to be natural transformations by setting component-wisely. 
Then the quadruple $(C,\iota,r,\sigma)$ 
is a normal ordinary complicial structure on $F_{b,h}\cC$ and 
$F_{[a,b],h}\cC$, $F_{b,\leq b,h}\cC$ and $F_{b,\geq a,h}\cC$ 
are complicial exact subcategories of $F_{b,h}\cC$. 
We call this complicial structure the {\it level complicial structure} 
({\it on $F_{b,h}\cC$}) (resp. $F_{[a,b],h}\cC$, 
$F_{b,\leq b,h}\cC$, $F_{b,\geq a,h}\cC$).
\end{para}

In the proof of \ref{lemdf:additive str on Fbh}, 
we will use \ref{lem:homotopy commutative diagram induces cyl and cone seq}. 
Notice that there is no circular reasoning in our argument.

\begin{proof}
$\mathrm{(1)}$ 
We can show that $F_{b,h}\cC$ with the addition defined by 
$\mathrm{(\ref{eq:add of hcd})}$ 
is enriched over the category of abelian groups. 
We will check the universal property of direct sums. 
Let $x$, $y$, $z$ and $w$ be objects in $F_{b,h}\cC$. 
For any homotopy commutative diagrams $(f,H)\colon x\to z$, 
$(f',H')\colon y\to z$, $(g,K)\colon w\to x$ and $(g',K')\colon w\to y$, 
there exists unique homotopy commutative diagrams 
$\displaystyle{{\footnotesize{(\begin{pmatrix}f & f' \end{pmatrix},\begin{pmatrix}H & H' \end{pmatrix})}}\colon x\oplus y \to z}$ and 
$\displaystyle{{\footnotesize{\left(\begin{pmatrix}g \\ g' \end{pmatrix},\begin{pmatrix}K \\ K' \end{pmatrix} \right )}}\colon w\to x\oplus y}$ 
which defined by the formula 
$${\begin{pmatrix}f & f'\end{pmatrix}}_n=
\begin{pmatrix}f_n & f'_n \end{pmatrix},\ \ \ 
{\begin{pmatrix}H & H' \end{pmatrix}}_n=
\begin{pmatrix}H_n & H'_n \end{pmatrix},
 $$
$${\begin{pmatrix}g \\ g'\end{pmatrix}}_n=
\begin{pmatrix}g_n \\ g'_n \end{pmatrix},\ \ \ 
{\begin{pmatrix}K \\ K' \end{pmatrix}}_n= 
\begin{pmatrix}K_n \\ K'_n \end{pmatrix} $$
and make the diagrams below commutative. 
$$
\xymatrix{
x \ar[r]^{{\footnotesize{\begin{pmatrix}\id_x\\0\end{pmatrix}}}} \ar[rd]_{(f,H)} & 
x\oplus y 
\ar@{-->}[d] & 
y \ar[l]_{{\footnotesize{\begin{pmatrix}0 \\ \id_y\end{pmatrix}}}} 
\ar[ld]^{(f',H')}\\
& z,
}
\xymatrix{
& w \ar[ld]_{(g,K)} \ar[rd]^{(g',K')} \ar@{-->}[d]\\
x & x\oplus y \ar[r]_{{\footnotesize{\begin{pmatrix}0 &\id_y\end{pmatrix}}}} 
\ar[l]^{{\footnotesize{\begin{pmatrix}\id_x & 0\end{pmatrix}}}} & y.
}
$$
Hence $F_{b,h}\cC$ is an additive category.

\sn
$\mathrm{(2)}$ 
$\mathrm{(a)}$
Let $(f,H)\colon x\to y$ be a homotopy commutative diagram in $F_{b,h}\cC$. 
Assume that for any integer $n$, $f_n\colon x_n\to y_n$ is 
a Frobenius admissible monomorphism in $\cC$. 
Then we set $z_n:=y_n/x_n$ and let $g_n\colon y_n\rdef z_n$ be 
the canonical quotient morphism. 
Then by Lemma~\ref{lem:homotopy commutative diagram induces cyl and cone seq}, 
there exists a morphism $i_n^z\colon z_n\to z_{n+1}$ 
and a $C$-homotopy 
$K_n\colon i_n^zg_n\Rightarrow_C g_{n+1}i_n^y$ such that 
$K_n\star_C H_n=0$. 
Thus the system $\{z_n,i_n^z \}_{n\in\bbZ}$ 
forms an object in $F_{b,h}\cC$ 
and the pair $(g,K)$ forms a homotopy commutative diagram $y\to z$ 
and the sequence $x\overset{(f,H)}{\rinf} y \overset{(g,K)}{\rdef} z$ 
is a level Frobenius admissible exact sequence in $F_{b,h}\cC$.

For level Frobenius admissible epimorphisms, 
we shall apply the argument 
above to $\cC^{\op}$ the opposite category of $\cC$ and 
the dual complicial structure $(P,q,s,\tau^P)$ of $(C,\iota,r,\sigma)$ 
and utilize Lemma~\ref{lem:C hom equiv to P hom} 
and Lemma~\ref{lem:comp and A and B}. 

\sn
$\mathrm{(b)}$ 
What we need to prove are the following assertions:

\sn
$\mathrm{(i)}$ 
For any level Frobenius admissible exact sequence of the form 
$\mathrm{(\ref{eq:Frob exact in Fbh})}$, 
the commutative square 
\begin{equation}  
\label{eq:biCartesian of hcd}
\xymatrix{
x \ar[r]^{(f,H)} \ar[d] & y \ar[d]^{(g,K)}\\
0 \ar[r] & z
}
\end{equation}
is a biCartesian square in $F_{b,h}\cC$. 

\sn
$\mathrm{(ii)}$ 
The class of level Frobenius admissible monomorphisms in $F_{b,h}\cC$ 
is closed under compositions, and push-out along arbitrary morphisms. 
The class of level Frobenius admissible epimorphisms in $F_{b,h}\cC$ 
is closed under compositions, and pull-back along arbitrary morphisms. 

\sn
$\mathrm{(iii)}$ 
Split short exact sequences are level Frobenius admissible exact sequences.

\mn
Assertion $\mathrm{(iii)}$ is clear.

\begin{proof}[Proof of $\mathrm{(i)}$]
We only prove that the commutative square 
$\mathrm{(\ref{eq:biCartesian of hcd})}$ is 
a push-out diagram. 
A proof of the pull-back case is similar. 
For any object $u$ and any morphism $(h,L)\colon y\to u$ 
in $F_{b,h}\cC$ such that 
\begin{equation}
\label{eq:hlfh=0}
(h,L)(f,H)=0,
\end{equation}
we will show that 
there exists a unique morphism $(\psi,M)\colon z\to u$ such that 
$(h,L)=(\varphi,M)(g,K)$. 
First notice that the equality $\mathrm{(\ref{eq:hlfh=0})}$ implies 
the following equalities
\begin{equation}
\label{eq:hf is 0}
hf=0,
\end{equation}
\begin{equation}
\label{eq:h and L}
h_{n+1}H_n+L_nCf_n=0,
\end{equation}
\begin{equation}
\label{eq:g and K}
g_{n+1}H_n+K_nCf_n=0
\end{equation}
for any integer $n$. 
Moreover, since $K_n$ and $L_n$ are $C$-homotopies 
$i_n^zg_n\Rightarrow_C g_{n+1}y_n^y$ and 
$i_n^uh_n\Rightarrow_C h_{n+1}i_n^y$ respectively, 
we have equalities 
\begin{equation}
\label{eq:K_n homotopy}
i_n^zg_n=K_n\iota_{y_n}+g_{n+1}i_n^y,
\end{equation}
\begin{equation}
\label{eq:L_n homotopy}
L_n\iota_{y_n}=i_n^uh_n-h_{n+1}i_n^y
\end{equation}
for any integer $n$.
Let us fix an integer $n$. 
First notice that there exists a morphism 
$\varphi_n\colon z_n\to u_n$ in $\cC$ such that 
$\varphi_ng_n=h_n$ by the universality 
of the left commutative square below
$$\xymatrix{
x_n \ar[r]^{f_n} \ar[d] & y_n \ar[d]^{g_n} \ar@/^2pc/[rdd]^{h_n}\\
0 \ar[r] \ar@/_2pc/[rrd] & z_n \ar@{-->}[rd]_{\varphi_n}\\
& & u_n,
}\ \ \ 
\xymatrix{
Cx_n \ar[r]^{Cf_n} \ar[d] & Cy_n \ar[d]^{Cg_n} \ar@/^2pc/[rdd]^{L_n-\varphi_{n+1}K_n}\\
0 \ar[r] \ar@/_2pc/[rrd] & Cz_n \ar@{-->}[rd]_{M_n}\\
& & u_{n+1}.}
$$
Next since we have the equalities
$$
(L_n-\varphi_{n+1}K_n)Cf_n=L_nCf_n-\varphi_{n+1}K_nCf_n
=-h_{n+1}H_n+\varphi_{n+1}g_{n+1}H_n=0,
$$
there exists a morphism $M_n\colon Cz_n\to u_{n+1}$ such that 
\begin{equation}
\label{eq:L=M star K}
L_n-\varphi_{n+1}K_n=M_nCg_n
\end{equation}
by the universal property of the right commutative square above. 
Next we will show that the equality
\begin{equation}
\label{eq:i and varphi}
i_n^u\varphi_n-M_n\iota_{z_n}=\varphi_{n+1}i_n^z.
\end{equation}
We have the equalities
\begin{multline*}
(i_n^u\varphi_n-M_n\iota_{z_n})g_n=i_n^u\varphi_ng_n-M_n\iota_{z_n}g_n=
i_n^uh_n-M_n\iota_{z_n}g_n\\
= i_n^uh_n+(\varphi_{n+1}K_n-L_n)\iota_{y_n}-M_nCg_n=
\varphi_{n+1}K_n\iota_{y_n}+i_n^uh_n-L_n\iota_{y_n}\\
=\varphi_{n+1}K_n\iota_{y_n}+h_{n+1}i_n^y=
\varphi_{n+1}K_n\iota_{y_n}+\varphi_{n+1}g_{n+1}i_n^y=
\varphi_{n+1}(K_n\iota_{y_n}+g_{n+1}i_n^y)=\varphi_{n+1}i_n^zg_n. 
\end{multline*}
Since $g_n$ is a Frobenius admissible epimorphism, a fortiori, 
an epimorphism, we obtain the equality $\mathrm{(\ref{eq:i and varphi})}$. 
The equality $\mathrm{(\ref{eq:i and varphi})}$ 
means the pair $(\varphi,M)$ is a homotopy commutative diagram 
from $z$ to $u$ and the equality $\mathrm{(\ref{eq:L=M star K})}$ shows 
we have the equality $L=M\star K$. 
By construction the homotopy commutative diagram $(\varphi,M)$ such that 
$(h,L)=(\varphi,M)(g,K)$ is unique. 
Hence we complete the proof. 
\end{proof}

\begin{proof}[Proof of $\mathrm{(ii)}$] 
First notice that for a pair of 
level Frobenius admissible monomorphisms 
$(f,H)\colon x\rinf y$ and $(g,K)\colon y\rinf z$, 
the composition $(gf,K\star H)$ is 
a level Frobenius admissible monomorphism 
by $\mathrm{(a)}$. 

Second, let $(f,H)\colon x\to y$ be a 
level Frobenius admissible monomorphism and 
let $(g,K)\colon x\to z$ be an arbitrary homotopy commutative diagram. 
We will show that there exists a pair of homotopy commutative diagrams 
$(f',H')\colon z\to u $ and $(g',K')\colon y\to u$ 
which is a push out of $(f,H)$ along $(g,K)$ 
and $(g',K')$ is a level Frobenius admissible monomorphism. 
\begin{equation}
\label{eq:coCartesian of hom comm dia}
\xymatrix{
x \ar@{>->}[r]^{(f,H)} \ar[d]_{(g,K)} & y \ar@{-->}[d]^{(g',K')}\\
z \ar@{-->}[r]_{(f',H')} & u.
}
\end{equation}
Then it turns out that the homotopy commutative diagram 
$\displaystyle{\footnotesize{\left (\begin{pmatrix}f\\ g\end{pmatrix},\begin{pmatrix}H\\ K\end{pmatrix}\right )}\colon x\to y\oplus z}$ is 
a level Frobenius admissible monomorphism 
by \ref{lem:kel90} and $\mathrm{(a)}$. 
Thus by the proof of $\mathrm{(a)}$, 
there exists an object $u$ in $F_{b,h}\cC$ 
and a level Frobenius admissible epimorphism 
$\displaystyle{\footnotesize{\left (\begin{pmatrix}g'& -f' \end{pmatrix},\begin{pmatrix}K' & -H' \end{pmatrix} \right )}\colon y\oplus z\to u}$ such that 
the sequence $x\overset{{\footnotesize{\left (\begin{pmatrix}f\\ g\end{pmatrix},\begin{pmatrix}H\\ K\end{pmatrix}\right )}}}{\rinf} y\oplus z \overset{{\footnotesize{\left (\begin{pmatrix}g'& -f' \end{pmatrix},\begin{pmatrix}K' & -H' \end{pmatrix} \right )}}}{\rdef} u$ 
is a level Frobenius admissible exact sequence in $F_{b,h}\cC$. 
Hence the commutative square of the form $\mathrm{(\ref{eq:coCartesian of hom comm dia})}$ is 
a coCartesian square in $F_{b,h}\cC$ and 
it turns out that 
the homotopy commutative diagram $(f',H')\colon z\to u$ is  
a level Frobenius admissible monomorphism 
by \ref{lem:kel90} and $\mathrm{(a)}$ again. 

Finally by applying the previous argument to $\cC^{\op}$ 
the opposite category of $\cC$ and the 
dual complicial structure $(P,q,s,\tau^P)$ 
of $(C,\iota,r,\sigma)$ and 
by utilizing Lemma~\ref{lem:C hom equiv to P hom} and 
Lemma~\ref{lem:comp and A and B}, 
we obtain the results for level Frobenius admissible epimorphisms. 
We complete the proof.
\end{proof}

\sn
$\mathrm{(3)}$ 
Almost all verifications are straightforwards. 
What we use normality assumption are following: 

\sn
For any homotopy commutative diagram $(f,H)\colon x\to y$ 
and any integer $n$, 
the following diagrams are commutative 
and they follow from equalities 
$\mathrm{(\ref{eq:normal comp ex func r})}$ and 
$\mathrm{(\ref{eq:normal comp ex func sigma})}$. 
$${\footnotesize{
\xymatrix{
[CCx_n\onto{CCi_n^x}CCx_{n+1}] \ar[r]^{(r_{x_n},r_{x_{n+1}},0)} 
\ar[d]_{(CCf_n,CCf_{n+1},CCH_nC\sigma_{x_n}\sigma_{Cx_{n}})} & 
[Cx_n\onto{Ci_n^x}Cx_{n+1}] \ar[d]^{(Cf_n,Cf_{n+1},CH_n\sigma_{x_n})}\\
[CCy_n\onto{CCi^y_n}CCy_{n+1}] \ar[r]_{(r_{y_n},r_{y_{n+1}},0)} & [Cy_n\onto{Ci^y_n}Cy_{n+1}],
}}}
$$
$${\footnotesize{ 
\xymatrix{
[CCx_n\onto{CCi_n^x}CCx_{n+1}] \ar[r]^{(\sigma_{x_n},\sigma_{x_{n+1}},0)} 
\ar[d]_{(CCf_n,CCf_{n+1},CCH_nC\sigma_{x_n}\sigma_{Cx_{n}})} & 
[CCx_n\onto{CCi_n^x}CCx_{n+1}] 
\ar[d]^{(Cf_n,Cf_{n+1},CCH_nC\sigma_{x_n}\sigma_{Cx_{n}})}\\
[CCy_n\onto{CCi^y_n}CCy_{n+1}] \ar[r]_{(\sigma_{y_n},\sigma_{y_{n+1}},0)} & 
[CCy_n\onto{CCi^y_n}CCy_{n+1}].
}}}
$$
\end{proof}

By \ref{lem:CC fundamental results} and \ref{lemdf:additive str on Fbh} 
(and observations in \ref{rem:2-structure of ChbE}), 
we obtain the following result.

\begin{para}
\label{cor:underine FbhC}
{\bf Corollary.}\ \ 
{\it
Let $\cC$ be a normal ordinary complicial exact category. 
Then we can naturally 
make $\underline{F_{b,h}\cC}$ into a complicial exact category 
such that the canonical functor $F_{b,h}\cC\to\underline{F_{b,h}\cC}$ 
is a strictly normal 
complicial exact functor. 
Moreover for any pair of integers $a<b$, 
similar statements hold for $\underline{F_{[a,b],h}\cC}$ and so on. 
}\qed
\end{para}

We can show the following. 

\begin{para}
\label{lemdf:functoriality of Fbh} 
{\bf Lemma-Definition.}\ \ 
Let $(f,d)\colon\cC\to \cC'$ be a complicial exact functor 
between complicial exact categories and let 
$(\phi,H)\colon x\to x'$ be a homotopy commutative diagram in $F_{b,h}\cC$. 
Then we define $f(x)$ and $f(\phi,H)$ to be a filtered object on $\cC'$ 
and a homotopy commutative diagram in $F_{b,h}\cC'$ by setting 
${(f(x))}_n:=f(x_n)$, $i_n^{fx}:=f(i_n^x)$ and 
${(f(\phi,H))}_n:=(f(\phi_n),f(H_n)\cdot d_{x_n})$ for any integer $n$. 
We also define 
$F_{b,h}\colon F_{b,h}\cC\to F_{b,h}\cC'$ to be a functor by sending a 
bounded filtered object $x$ on $\cC$ to $f(x)$ and 
a homotopy commutative diagram $(\phi,H)\colon x\to y$ in 
$F_{b,h}\cC$ to $f(\phi,H)$. 
The functor $F_{b,h}f$ sends a level Frobenius admissible exact sequence 
in $F_{b,h}\cC$ 
to a level Frobenius admissible exact sequence in $F_{b,h}\cC'$. 
In particular if $\cC$ and $\cC'$ are ordinary, then 
$F_{b,h}f$ is an exact functor.

Moreover if we assume $(f,d)$ is normal, 
then $d$ induces a natural equivalence 
$d^{\lv}\colon C'^{\lv}F_{b,h}f\isoto F_{b,h}fC^{\lv}$ 
defined 
by the formula ${(d_x^{\lv})}:=(\{d_{x_n}\}_{n\in\bbZ},0)$ 
for any object $x$ in $F_{b,h}\cC$ 
where $C^{\lv}$ and $C'^{\lv}$ are 
the endofunctors on $F_{b,h}\cC$ and $F_{b,h}\cC'$ respectively 
defined by level-wisely. 
In particular if $\cC$ and $\cC'$ are normal ordinary, 
then the pair 
$(F_{b,h}f,d^{\lv})$ is 
a complicial exact functor from $F_{b,h}\cC$ to $F_{b,h}\cC'$. 
\qed
\end{para}

For a complicial exact category $\cC$, we will define several operations 
on $F_{b,h}\cC$ or $\underline{F_{b,h}\cC}$. 

\begin{para}
\label{df:truncations on FbhC}
{\bf Definition (Truncation functors).}\ \ 
Let $\cC$ be a complicial exact category and let $n$ be an integer and 
let $x$ be a bounded filtered object on $\cC$ and let 
$(f,H)\colon x\to y$ be a homotopy commutative diagram. 
We define $\sigma_{\leq n}x$ and $\sigma_{\geq n}x$ 
to be filtered objects on $\cC$ by setting
\begin{equation}
\label{eq:sigmaleqn obj df}
{(\sigma_{\leq n}x)}_k:=
\begin{cases}
x_k & \text{if $k\leq n$}\\
x_n & \text{if $k\geq n$}
\end{cases},\ \ 
i^{\sigma_{\leq n}x}_k:=
\begin{cases}
i^x_k & \text{if $k\leq n-1$}\\
\id_{x_n} & \text{if $k\geq n$}
\end{cases}
\end{equation}
\begin{equation}
\label{eq:sigmageq obj df}
{(\sigma_{\geq n}x)}_k:=
\begin{cases}
x_k & \text{if $k\geq n$}\\
0 & \text{if $k< n$}
\end{cases},\ \ 
i^{\sigma_{\geq n}x}_k:=
\begin{cases}
i^x_k & \text{if $k\geq n$}\\
0 & \text{if $k\leq n-1$}.
\end{cases}
\end{equation}
We also define 
$\sigma_{\leq n}(f,H)\colon \sigma_{\leq n}x\to\sigma_{\leq n}y$ and 
$\sigma_{\geq n}(f,H)\colon \sigma_{\geq n}x\to\sigma_{\geq n}y$ 
to be homotopy commutative diagrams by setting
\begin{equation}
\label{eq:sigman mor df}
{(\sigma_{\leq n}(f,H))}_k:=
\begin{cases}
(f_k,H_k) & \text{if $k< n$}\\
(f_n,0) & \text{if $k\geq n$}
\end{cases},\ \ 
{(\sigma_{\geq n}(f,H))}_k:=
\begin{cases}
(f_n,H_n) & \text{if $k\geq n$}\\
(0,0) & \text{if $k\leq n-1$}.
\end{cases}
\end{equation}
The associations $\sigma_{\leq n}$, $\sigma_{\geq n}\colon 
F_{b,h}\cC\to F_{b,h}\cC$ are functors which send 
a level Frobenius admissible exact sequence to 
a level Frobenius admissible exact sequence. 

For any object $x$ in $F_{b,h}\cC$, 
there exists a pair of functorial 
natural morphism $\sigma_{\geq n}x\to x$ and $\sigma_{\leq n}x\to x$ and 
both $\sigma_{\geq n}\colon F_{b,h}\cC\to F_{b,\geq n,h}\cC$ 
and $\sigma_{\leq n}\colon F_{b,h}\cC\to F_{b,\leq n,h}\cC $ are right adjoint functors of the 
inclusion functors $F_{b,\geq n,h}\cC\to F_{b,h}\cC$ and 
$F_{b,\leq n,h}\cC\to F_{b,h}\cC$ respectively.

Moreover if $\cC$ is normal, then we have the equalities 
\begin{equation}
\label{eq:truncation strictly normal}
\sigma_{\leq n}C^{\lv}=C^{\lv}\sigma_{\leq n},\ \ 
\sigma_{\geq n}C^{\lv}=C^{\lv}\sigma_{\geq n}.
\end{equation}
In particular if $\cC$ is normal ordinary, then 
$\sigma_{\leq n}$ and $\sigma_{\geq n}$ are strictly normal 
complicial exact functors $F_{b,h}\cC\to F_{b,h}\cC$. 
\end{para}

\begin{para}
\label{df:degree shift}
{\bf Definition (Degree shift).}\ \ 
Let $\cC$ be a complicial exact category and let $k$ 
be an integer and let 
$(f,H)\colon x\to y$ be 
a homotopy commutative diagram in $F_{b,h}\cC$. 
We define $x[k]$ and $(f[k],H[k])\colon x[k]\to y[k]$ to be an object 
and a morphism in $F_{b,h}\cC$ by setting 
${x[k]}_n=x_{k+n}$, $i^{x[k]}_n=i^x_{k+n} $, 
${f[k]}_n=f_{k+n}$ and ${H[k]}_n=H_{k+n}$. 
The association $(-)[k]\colon F_{b,h}\cC\to F_{b,h}\cC$, $x\mapsto x[k]$ 
gives a functor which sends a level Frobenius admissible exact 
sequence to a level Frobenius admissible exact sequence. 
In particular if $\cC$ is normal ordinary, then 
the functor $(-)[k]$ is a 
strictly normal complicial exact functor. 
\end{para}

\begin{para}
\label{df:Cone functor}
{\bf Definition (Cone functor).}\ \ 
Let $k$ be an integer and let $\cC$ be a complicial exact category. 
We define $\fc_k\colon F_{b,h}\cC\to F_{b,h}\cC$ 
to be a functor by setting 
\begin{equation}
\label{eq: fck df}
{\scriptstyle{
{(\fc_k(x))}_n:=
\begin{cases}
x_n & \text{if $n\leq k-1$}\\
Cx_n & \text{if $n\geq k$}
\end{cases},\ \ 
i_n^{\fc_k(x)}:=
\begin{cases}
i_n^x & \text{if $n\leq k-2$}\\
\iota_{x_n}i_k^x & \text{if $n=k-1$}\\
Ci_n^x & \text{if $n\geq k$}
\end{cases},\ \ 
{\fc_k(f,H)}_n:=
\begin{cases}
(f_n,H_n) & \text{if $n\leq k-2$}\\
(f_{k-1},0\star H_{k-1}) & \text{if $n=k-1$}\\
(Cf_n,CHn\cdot\sigma_{x_n}) & \text{if $n\geq k$}
\end{cases}}}
\end{equation}
for an object $x$ and a morphism 
$(f,H)\colon x\to y$ in $F_{b,h}\cC$.

For a pair of integers $a<b$, we also denote the restrictions 
of $\fc_k$ to $F_{[a,b],h}\cC$ and $\underline{F_{[a,b],h}\cC}$ 
and so on by the same letters $\fc_k$. 

If $\cC$ is normal ordinary, then we define 
$d^{\fc_k}\colon C^{\lv}\fc_k\isoto \fc_k C^{\lv}$ to be 
a natural equivalence by setting for any object $x$ in $F_{b,h}\cC$, 
\begin{equation}
\label{eq: dck df}
{(d^{\fc_k})}_n=
\begin{cases}
\id_{Cx_n} & \text{if $n\leq k-1$},\\
\sigma_{x_n} & \text{if $n\geq k$}.
\end{cases}
\end{equation}
Then the pair $(\fc_k,d^{\fc_k})$ is a complicial exact functor 
$F_{b,h}\cC\to F_{b,h}\cC$. 
We call $\fc_k$ {\it the $k$th cone functor}.
\end{para}

\begin{para}
\label{df:skip functor}
{\bf Definition (skip functor).}\ \ 
Let $k$ be an integer and let $\cC$ be a complicial exact category. 
We define $\fs_k\colon\underline{F_{b,h}\cC}\to\underline{F_{b,h}\cC}$ 
to be a functor by setting
\begin{equation}
\label{eq:skip df}
{\scriptstyle{
{(\fs_k(x))}_n:=
\begin{cases}
x_n & \text{if $n\leq k-1$}\\
x_{n+1} & \text{if $n\geq k$}
\end{cases},\ \ 
i_n^{\fs_k(x)}:=
\begin{cases}
i_n^x & \text{if $n\leq k-2$}\\
i_k^xi_{k-1}^x & \text{if $n=k-1$}\\
i_{k+1}^x & \text{if $n\geq k$}
\end{cases},\ \ 
{(\fs_k(f,H))}_n:=
\begin{cases}
(f_n,H_n) & \text{if $n\leq k-2$}\\
(f_{k-1},H_k\star H_{k-1}) & \text{if $n=k-1$}\\
(f_{n+1},H_{n+1}) & \text{if $n\geq k$}
\end{cases}}}
\end{equation}
for an object $x$ and a morphism $(f,H)\colon x\to y$ in 
$\underline{F_{b,h}\cC}$. 

Notice that $\fs_k$ is not $1$-functorial on $F_{b,h}\cC$, 
but on $\underline{F_{b,h}\cC}$. 
If $\cC$ is normal ordinary, then $\fs_k$ 
is a strictly normal complicial exact. 
We call $\fs_k$ {\it the $k$th skip functor}.
\end{para}

Recall the definition of $F_{[a,b]}^{\rinf}\cC$ from 
Definition~\ref{df:filtered objects}

\begin{para}
\label{df:rab}
{\bf Definition ($\fr_{[a,b]}$).}\ \ 
Let $\cC$ be a normal ordinary complicial exact category and let 
$a<b$ be a pair of integers. 
We will define $\fr_{[a,b]}\colon F_{[a,b],h}\cC\to F_{[a,b]}^{\rinf}\cC$ 
to be a complicial functor by setting in the following way. 
For a pair of objects $x$ and $y$ in $F_{[a,b],h}\cC$ 
and a pair of integers 
$a\leq s\leq t\leq b$ and a family of morphisms 
$\{f_i\colon x_i\to y_i \}_{i\in[s,t]}$ in $\cC$ indexed by the integers 
$i\in[s,t]$, 
we set 
$\displaystyle{Cx_{[s,t]}:=Cx_t\oplus Cx_{t-1}\oplus \cdots \oplus Cx_s}$ and 
$\displaystyle{Cf_{[s,t]}:=
{\scriptscriptstyle{
\begin{pmatrix}
Cf_t&&&&\\
&Cf_{t-1}&&& \text{{\huge{0}}}\\
&&\ddots&\\
&&&\ddots&\\
\text{{\huge{0}}}&&&& Cf_s
\end{pmatrix}
}}\colon Cx_{[s,t]}\to Cy_{[s,t]}}$ and 
we define $\fr_{[a,b]}x$ to be an object in $F^{\rinf}_{[a,b]}\cC$ by setting
\begin{equation}
\label{eq:rab df}
{\footnotesize{
{(\fr_{[a,b]}x)}_{a+k}:=
\begin{cases}
0 & \text{if $k<0$}\\
x_a & \text{if $k=0$}\\
x_{a+k}\oplus Cx_{[a,a+k-1]} & \text{if $1\leq k< b-a$}\\
x_b\oplus Cx_{[a,b-1]} & \text{if $b-a\leq k$}\\
\end{cases}\ \ \ 
i^{\fr_{[a,b]}x}_{a+k}:=
\begin{cases}
0 & \text{if $k<0$}\\
\begin{pmatrix}
\xi_{1,i^x_{a+k}} & 0\\
0 & \id_{Cx_{[a,a+k-1]}}
\end{pmatrix} 
& \text{if $0\leq k<b-a$}\\
\id_{x_b\oplus Cx_{[a,b-1]}} & \text{if $b-a\leq k$}.
\end{cases}}}
\end{equation} 
For a homotopy commutative diagram $(f,H)\colon x\to y$ in $F_{[a,b],h}\cC$, 
we will inductively construct a morphism 
$\fr_{[a,b]}(f,H)\colon \fr_{[a,b]}x\to \fr_{[a,b]}y$ 
in $F_{[a,b]}^{\rinf}\cC$. 
First we set
\begin{equation}
\label{eq:rabfHa+k  leq 1 df}
{\fr_{[a,b]}(f,H)}_{a+k}:=
\begin{cases}
0 & \text{if $k<0$}\\
f_a & \text{if $k=0$}\\
\Cyl(f_a,f_{a+1},H_a) & \text{if $k=1$}.
\end{cases}
\end{equation}
For $2\leq k\leq b-a$, we will set
\begin{equation}
\label{eq:rabfHa+k 2 leq df}
{\fr_{[a,b]}(f,H)}_{a+k}:=
\begin{pmatrix}
\Cyl(f_{a+k-1},f_{a+k},H_{a+k-1}) & 0_{2,k-1}\\
0_{k-1,2} & Cf_{[a,a+k-2]}
\end{pmatrix}+
\begin{pmatrix}
0_{k,2} & {\calH(f,H)}_k\\
0_{1,2} & 0_{1,k-1}
\end{pmatrix}
\end{equation}
where $O_{s,t}$ is the $s\times t$th zero matrix and 
${\calH(f,H)}_k\colon Cx_{[a,a+k-2]}\to y_{a+k+1}\oplus Cy_{[a+1,a+k]}$ 
to be a morphism in $\cC$ which will be defined by inductively. 
By convention we set ${\calH(f,H)}_1=0$. 
If we define the morphisms 
${\fr_{[a,b]}(f,H)}_{a+k}$ and 
${\fr_{[a,b]}(f,H)}_{a+k+1}$ appropriately, 
then we need to have an equality
\begin{equation}
\label{eq:iyrab=rab ix}
i^y_{a+k}{\fr_{[a,b]}(f,H)}_{a+k}={\fr_{[a,b]}(f,H)}_{a+k+1}i^x_{a+k}.
\end{equation}
Thus inspection shows that we need to have an equality
\begin{equation}
\label{eq:recurrence formula}
{\calH(f,H)}_{k+1}:=
\begin{pmatrix}
\begin{matrix}
-\xi_{1,i^y_{a+k}}H_{a+k-1}\\
0_{k-1,1}
\end{matrix}
&
\begin{pmatrix}
\xi_{1,i^y_{a+k}} & 0_{1,k-1}\\
0_{k-1,1} &\id_{Cy_{[a+1,a+k-1]}}
\end{pmatrix}
{\calH(f,H)}_k
\end{pmatrix}.
\end{equation}
Hence we inductively define ${\calH(f,H)}_k$ by using the recurrence formula 
$\mathrm{(\ref{eq:recurrence formula})}$. 
For $k>b-a$, we set ${\fr_{[a,b]}(f,H)}_{a+k}:={\fr_{[a,b]}(f,H)}_b$. 
Let $(g,K)\colon y\to z$ be another homotopy commutative diagram 
in $F_{b,h}\cC$. 
Then we have an equality
\begin{equation}
\label{eq:rab composition}
{\fr_{[a,b]}(g,K)}_{a+k}{\fr_{[a,b]}(f,H)}_{a+k}=
{\fr_{[a,b]}(gf,K\star H)}_{a+k}
\end{equation}
for $k\leq 1$. Thus by induction on $k$ and 
uniqueness of construction of ${\calH(f,H)}_k$, 
we obtain the equality 
\begin{equation}
\label{eq:rab composition 2}
\fr_{[a,b]}(g,K)\fr_{[a,b]}(f,H)=\fr_{[a,b]}(gf,K\star H)
\end{equation}
as a morphism from $\fr_{[a,b]}x$ to $\fr_{[a,b]}z$. 

The natural equivalence 
$\sigma\colon CC\isoto CC$ induces a natural equivalence 
\begin{equation}
\label{eq:drab df}
d^{\fr_{[a,b]}}\colon C^{\lv}\fr_{[a,b]}\isoto \fr_{[a,b]}C^{\lv}.
\end{equation}
and the pair $(\fr_{[a,b]},d^{\fr_{[a,b]}})$ 
is a complicial exact functor $F_{[a,b],h}\cC\to F_{[a,b]}^{\rinf}\cC$. 
We write $\fraki_{[a,b]} $ and ${\fraki'}_{[a,b]}$ 
for the inclusion functors 
$F_{[a,b]}^{\rinf}\cC\rinc F_{[a,b],h}\cC$ and $F_{[a,b]}^{\rinf}\cC\rinc F_{[a,b]}\cC $ respectively. 
Then 
there are complicial natural transformations 
$\id_{F^{\rinf}_{[a,b]}\cC}\to \fr_{[a,b]}\fraki_{[a,b]}$ and 
$\id_{F_{[a,b]}\cC} \to {\fraki'}_{[a,b]}\fr_{[a,b]}|_{F_{[a,b]}\cC}$. 
\end{para}

\subsection{Mapping cylinder and mapping cone}
\label{subsec:mapp cyl and mapp cone}

We can also generalize the notion of mapping cone and mapping cylinder 
of chain morphisms in 
a category of chain complexes on an additive category to 
a complicial exact category. 

\begin{para}
\label{df:mapping cone and cylinders}
{\bf (Mapping cylinder and mapping cone).}\ \ 
Let $\cC$ be a complicial exact category. 
Like as categories of chain complexes, we will define a mapping cone and 
a mapping cylinder for morphisms in $\cC$ as follows. 
First we write $\dom$ and $\ran$ for the exact functors 
$\cC^{[1]}\to\cC$ induced from 
the functors $[0]\to [1]$ defined by sending $0$ to $0$ and $1$ respectively. 
We define $\epsilon\colon\dom\to\ran$ to be a natural transformation 
by setting for any object $f\colon x\to y$ in $\cC^{[1]}$, 
$\epsilon(f):=f$. 
We define 
$\Cone,\ \Cyl\colon \cC^{[1]} \to \cC$ to be exact functors 
called the {\it mapping cone functor} and 
the {\it mapping cylinder functor} respectively 
by setting  
$\Cyl:=\ran\oplus C\dom$ and 
$\Cone:=\ran\sqcup_{\dom}C\dom$ 
where $\Cone$ is defined the following push-out diagram: 
$$
{\xymatrix{
\dom \ar@{>->}[r]^{\iota\ast\dom} \ar[d]_{\epsilon} 
& 
C\dom \ar[d]^{\mu} \ar@{->>}[r]^{\pi\ast\dom} & 
T\dom \ar[d]^{\id_{T\dom}}& \\
\ran \ar@{>->}[r]_{\kappa} & \Cone \ar@{->>}[r]_{\psi} & T\dom 
}}$$
where $\psi$ is induced from the universal property of $\Cone$. 
We have the Frobenius admissible exact sequence
$$\dom \overset{\xi_1}{\rinf} \Cyl \overset{\eta}{\rdef} \Cone$$ 
where $\xi_1:=
\footnotesize{
\begin{pmatrix}
\epsilon\\
-\iota\ast\dom
\end{pmatrix}
}$ 
and $\eta=
\footnotesize{
\begin{pmatrix}
\kappa & \mu
\end{pmatrix}
}$. 
Moreover we define 
$\xi_2\colon \ran \to \Cyl$, 
$\xi_3\colon C\dom \to\Cyl$ and 
$\upsilon\colon \Cyl \to \ran$ to be natural transformations by setting 
$\xi_2:=
\footnotesize{
\begin{pmatrix}
\id_{\ran}\\
0
\end{pmatrix}
}$, 
$
\xi_3:=
\footnotesize{
\begin{pmatrix}
0\\
\id_{C\dom}
\end{pmatrix}}$
and 
$\upsilon:=
\footnotesize{
\begin{pmatrix}
\id_{\ran} & 0
\end{pmatrix}
}$. 
Then we have the following commutative diagram:
$${\xymatrix{
\dom \ar@{>->}[r]^{\xi_1} \ar[rd]_{\epsilon} & \Cyl \ar[d]_{\upsilon} & 
\ran \ar@{>->}[l]_{\xi_2} \ar[ld]^{\id_{\ran}}\\
& \ran  .
}}$$ 

Since the functor $\Cone$ is defined by the push-out diagram, 
we can choose the cofiber products as follows. 
For an object $x$ in $\cC$, since the diagrams below 
are push-outs, we shall assume that 
$\Cone \id_x=Cx$, $\kappa_{\id_x}=\iota_x$, 
$\mu_{\id_x}=\id_{Cx}$, $\Cone(x\to 0)=Tx$, $\mu_{x\to 0}=\pi_x$, 
$\Cone(0\to x)=x$ and $\kappa_{0\to x}=\id_x$.   
$$\xymatrix{ 
x \ar@{>->}[r]^{\iota_x} \ar[d]_{\id_x} & Cx \ar[d]^{\id_{Cx}}\\
x \ar@{>->}[r]_{\iota_x} & Cx,
}\ 
\xymatrix{
x \ar@{>->}[r]^{\iota_x} \ar[d] & Cx \ar[d]^{\pi_x}\\
0 \ar@{>->}[r] & Tx,
}\ 
\xymatrix{
0 \ar@{>->}[r]^{\iota_0} \ar[d] & C(0)=0 \ar[d]\\
x \ar@{>->}[r]_{\id_x} & x.
}$$
Moreover we define $s$, $t$, $\Delta\colon \cC \to \cC^{[1]}$ 
to be exact functors by sending an object $x$ to 
$[x\to 0]$, $[0\to x]$ and $[x\onto{\id_x} x]$ respectively. 
Then we shall assume that $\Cone s=T$, $\Cone t=\id_{\cC}$ 
and $\Cone \Delta=C$. 

For a homotopy commutative square 
$(a,b,H)\colon[f\colon x\to x']\to [g\colon y\to y']$ in $\cC$, 
we define a morphism $\Cyl(a,b,H)\colon \Cyl f\to \Cyl g$ 
in $\cC$ by setting 
$\displaystyle{\Cyl(a,b,H):=\begin{pmatrix}b & -H\\ 0 & Ca\end{pmatrix}}$. 
Then the left square below is commutative and 
there exists a morphism $\Cone(a,b,H)\colon \Cone f \to \Cone g$ in $\cC$ which makes the right square below commutative: 
$$\xymatrix{
x \ar@{>->}[r]^{{\xi_1}_f} \ar[d]_a & \Cyl f \ar[d]_{\Cyl(a,b,H)} \ar@{->>}[r]^{\eta_f} & \Cone f 
\ar@{-->}[d]^{\Cone(a,b,H)}\\
x' \ar@{>->}[r]_{{\xi_1}_g}  & \Cyl g \ar@{->>}[r]_{\eta_g} & \Cone g. 
}$$
Then we can regard $\Cyl$ and $\Cone$ as additive functors on 
$\cC^{[1]}_h$ and 
$\xi_1\colon \dom \to \Cyl$, $\xi_2\colon \ran \to\Cyl$ 
and $\kappa\colon \ran \to\Cone$ 
as 
natural transformations between 
functors on $\cC^{[1]}_h$. 
Since $\Cyl$ and $\Cone$ are defined by push-out diagrams in 
$\HOM_{\underline{\ExCat}}(\cC^{[1]},\cC)$ the (large) exact category of 
exact functors from $\cC^{[1]}$ to $\cC$, 
$\Cyl$ and $\Cone$ are exact functors from $\cC^{[1]}$ 
to $\cC$. 
\end{para} 

\begin{para}
\label{df:rC,Cone}
{\bf Definition ($r^{C,\Cone}$, $r^{\Cone,C}$).}\ \ 
Let $\cC$ be a complicial exact category. 
We define two natural transformations 
$r^{C,\Cone}\colon C\Cone \to C\ran$ and $r^{\Cone,C}\colon \Cone C\to C\ran$ 
between functors from 
$\cC^{[1]}$ to $\cC$ by 
the universal properties of $C\Cone$ and $\Cone C$ 
as in the commutative diagrams below. 
\begin{equation}
\label{eq:CCone df}
\xymatrix{
C\dom \ar@{>->}[r]^{C\ast\iota} \ar[d]_{C\ast\epsilon} & 
CC\dom \ar[r]^r \ar[d]^{C\ast\mu} & C\dom \ar[dd]^{C\ast\epsilon}\\
C\ran \ar@{>->}[r]^{C\ast\kappa} \ar[rrd]_{\id_{C\dom}} & C\Cone 
\ar[rd]^{r^{C,\Cone}}\\
& & C\ran,
}
\end{equation}
\begin{equation}
\label{eq:ConeC df}
\xymatrix{
C\dom \ar@{>->}[r]^{\iota\ast C} \ar[d]_{\epsilon\ast C} & 
CC\dom \ar[r]^r \ar[d]^{\mu\ast C} & C\dom \ar[dd]^{\epsilon \ast C}\\
C\ran \ar@{>->}[r]^{\kappa \ast C} \ar[rrd]_{\id_{C\dom}} & \Cone C
\ar[rd]^{r^{\Cone,C}}\\
& & C\ran.
}
\end{equation}
\end{para}

\begin{para}
\label{lem:Cone(f,id)}
{\bf Lemma.}\ \ 
{\it
Let $\cC$ be a complicial exact category and let 
$f\colon x\to y$ be a morphism. Then we have an equality
\begin{equation}
\label{eq:Cone(f,id)}
\Cone(f,\id_y)=r^{C,\Cone}_f\cdot \iota_{\Cone f}.
\end{equation}
}
\end{para}

\begin{proof}
Notice that we have equalities
\begin{equation}
\label{eq:rfiotaConefkappaf}
r^{C,\Cone}_f\cdot\iota_{\Cone f}\cdot \kappa_f=
r^{C,\Cone}_f\cdot C\ast\kappa_f \cdot\iota_y =\iota_y,
\end{equation}
\begin{equation}
\label{eq:rfiotaConefmuf}
r^{C,\Cone}_f\cdot\iota_{\Cone f}\cdot \mu_f=
r^{C,\Cone}_f\cdot C\ast\mu_f \cdot\iota_{Cx} =Cf\cdot r_x\cdot \iota_{Cx}=Cf.
\end{equation}
Thus it turns out that the morphism 
$r^{C,\Cone}_f\iota_{\Cone f}\colon\Cone f \to Cy$ 
makes the diagram below commutative. 
In particular we obtain the equality $\mathrm{(\ref{eq:Cone(f,id)})}$. 
$$
\xymatrix{
x \ar@{>->}[r]^{\tiny{\begin{pmatrix}-f\\ \iota_x\end{pmatrix}}} \ar[d]_f &
y\oplus Cx \ar@{->>}[r]^{\tiny{\begin{pmatrix}\kappa_f & \mu_f\end{pmatrix}}}
\ar[d]_{\tiny{\begin{pmatrix}\id_y & 0 \\ 0 & Cf \end{pmatrix}}} & 
\Cone f \ar[d]^{r^{C,\Cone}_f\cdot\iota_{\Cone f}}\\
y \ar@{>->}[r]_{\tiny{\begin{pmatrix}-\id_y\\ \iota_y \end{pmatrix}}} & 
y\oplus Cy \ar@{->>}[r]_{\tiny{\begin{pmatrix}\iota_y & \id_{Cy}\end{pmatrix}}} & Cy.
}
$$
\end{proof}

\begin{para}
\label{df:ThetaX}
{\bf Definition ($\Theta_{\fX}$).}\ \ 
Let $\cC$ be a complicial exact category and let $\fX$ be a commutative 
diagram $\mathrm{(\ref{eq:diagram fx})}$.
\begin{equation}
\label{eq:diagram fx}
\xymatrix{
x \ar[r]^a \ar[d]_f & x' \ar[d]^{f'}\\
y \ar[r]_b \ar[ru]^h & y'.
}
\end{equation}
Then we define $\Theta_{\fX}\colon C\Cone f \to \Cone f'$ to be a 
$C$-homotopy $\Cone (a,b)\Rightarrow_C 0$ by setting
\begin{equation}
\label{eq:ThetafX df} 
\Theta_{\fX}:=\Cone(h,b)\cdot r^{C,\Cone}_f.
\end{equation}
$$
\xymatrix{
\Cone f \ar@{>->}[r]^{\iota_{\Cone f}} \ar[d]_{\Cone (a,b)}
\ar[rd]^{\!\!\!\!\!\!\Cone(f,\id_y)} & C\Cone f \ar[d]^{r^{C,\Cone}_f}\\
\Cone f' & Cy \ar[l]^{\Cone (h,b)}.
}
$$
In particular if $x'=y$, $a=f$, $y'=0$, $f'=b=0$ and $h=\id_y$, then 
we write $\Theta_f$ for $\Theta_{\fX}$. 
\end{para}

\begin{para}
\label{ex:typicla split seq}
{\bf Example.}\ \ 
Let $\cC$ be a complicial exact category. 
Then the exact sequence of exact endofunctors on $\cC_h^{[1]}$, 
\begin{equation}
\label{eq:typ split seq}
[0\to C]\onto{(0,\id_C)}[C\onto{\id_C}C]\onto{(\id_C,0)}[C\to 0]
\end{equation}
is not split on $\cC^{[1]}$, but split on $\cC^{[1]}_h$, 
where a retraction and a section are given by 
$(0,\id_C,-r)\colon[C\onto{\id_C}C]\to [0\to C]$ and 
$(\id_C,0,r)\colon[C\to 0]\to [C\onto{\id_C}C]$ respectively. 
By the additive functor $\Cone\colon \cC^{[1]}_h\to \cC$, the sequence 
$\mathrm{(\ref{eq:typ split seq})}$ 
induces a split exact sequence $\mathrm{(\ref{eq:fund split seq})}$, 
and we have the equalities 
\begin{equation}
\label{eq:r as cone}
r=\Cone(0,\id_C,-r), 
\end{equation}
\begin{equation}
\label{eq:zeta as cone}
\zeta=\Cone(\id_C,0,r),
\end{equation}
where $\zeta\colon TC\to CC$ is defined in \ref{df:mor zeta}. 
\end{para}

\begin{para}
\label{lem:Cone(0,kf,-muf)}
{\bf Lemma.}\ \ 
{\it
Let $f\colon x\to y$ be a morphism 
in a complicial exact category $\cC$. 
Then we have an equality}
\begin{equation}
\label{eq:Cone(0,kf,-muf)}
\Cone(0,\kappa_f,\mu_f)=\id_{\Cone f}.
\end{equation}
\end{para}

\begin{proof}
For a $C$-homotopy commutative square 
$(0,\kappa_f,-\mu_f)\colon[x\onto{f}y]\to [0\to \Cone f]$, 
we have the equality 
\begin{equation}
\label{eq:Cyl(0,kf,-muf)}
\Cyl(0,\kappa_f,-\mu_f)=\begin{pmatrix}\kappa_f & \mu_f \end{pmatrix}
\colon y\oplus Cx \to \Cone f.
\end{equation}
Thus by the following commutative diagram, we obtain the 
equality $\mathrm{(\ref{eq:Cone(0,kf,-muf)})}$. 
$$
\xymatrix{
x \ar[r] \ar@{>->}[d]_{{\tiny{\begin{pmatrix}f & -\iota_x \end{pmatrix}}}} & 0 
\ar[d]\\
y\oplus Cx \ar[r]_{\!\!\!\!\!\!\!{\tiny{\begin{pmatrix}\kappa_f & \mu_f\end{pmatrix}}}} 
\ar@{->>}[d]_{{\tiny{\begin{pmatrix}\kappa_f & \mu_f\end{pmatrix}}}} 
& \Cone f\oplus 0 \ar[d]^{\id_{\Cone f}}\\
\Cone f \ar[r]_{\id_{\Cone f}} & \Cone f.
}
$$
\end{proof}

\begin{para}
\label{lem:fundament equalities in comp exact cat}
{\bf Lemma.}\ \ 
{\it In a complicial exact category, 
we have the following equalities:
\begin{equation}
\label{eq:diagonalizable 1}
(r\ast C)\cdot (C\ast \zeta)\cdot (C\ast\pi\ast C)\cdot (C\ast \sigma)\cdot 
(\zeta\ast C)=0,
\end{equation}
\begin{equation}
\label{eq:diagonalizable 2}
\zeta\cdot (\pi\ast C)\cdot \sigma \cdot (r\ast C)=
(r\ast C)\cdot (C\ast \zeta)\cdot (C\ast \pi \ast C)\cdot (C\ast \sigma),
\end{equation}
\begin{equation}
\label{eq:diagonalizable 3}
(\zeta\ast C)\cdot (T\ast \zeta)\cdot (T\ast \pi\ast C)\cdot (T\ast \sigma)
=(C\ast \zeta)\cdot (C\ast \pi \ast C)\cdot (C\ast \sigma)\cdot (\zeta\ast C),
\end{equation}
\begin{equation}
\label{eq:diagonalizable 4}
\zeta\cdot \tau^{C,T}\cdot (r\ast T)=(r\ast C)\cdot (C\ast \zeta)\cdot (C\ast \tau^{C,T}).
\end{equation}
}
\end{para}

\begin{proof}
First we will prove the equality $\mathrm{(\ref{eq:diagonalizable 1})}$. 
By using the equalities $\mathrm{(\ref{eq:r as cone})}$ and 
$\mathrm{(\ref{eq:zeta as cone})}$, 
we have the equalities
\begin{multline*}
(r\ast C)\cdot (C\ast \zeta)\cdot (C\ast\pi\ast C)\cdot (C\ast \sigma)\cdot 
(\zeta\ast C)
\\=
\Cone(0,\id_{CC},-r\ast C)\cdot\Cone(\zeta,\zeta)\cdot
\Cone(\pi\ast C,\pi\ast C)\cdot
\Cone(\sigma,\sigma)\cdot\Cone(\id_{CC},0,r\ast C)\\ 
=\Cone(0,0,0)=0.
\end{multline*}
Thus we obtain the equality $\mathrm{(\ref{eq:diagonalizable 1})}$ 
and it turns out that 
we have the following commutative diagram
$$\xymatrix{
TCC \ar@{>->}[r]^{\zeta\ast C} \ar[d]^{\wr}_{T\ast\sigma} & 
CCC \ar@{->>}[r]^{r\ast C} \ar[d]_{\wr}^{C\ast \sigma} &
CC \ar[d]_{\wr}^\sigma\\
TCC \ar@{->>}[d]_{T\ast \pi \ast C} & 
CCC \ar@{->>}[d]^{C\ast\pi\ast C} & 
CC \ar@{->>}[d]^{\pi\ast C}\\
TTC \ar@{>->}[d]_{T\ast \zeta} & CTC \ar@{>->}[d]^{C\ast\zeta} & 
TC \ar@{>->}[d]^{\zeta}\\
TCC \ar@{>->}[r]_{\zeta\ast C} & CCC \ar@{->>}[r]_{r\ast C} & CC.
}$$
Namely we obtain the equalities 
$\mathrm{(\ref{eq:diagonalizable 2})}$ and 
$\mathrm{(\ref{eq:diagonalizable 3})}$. 
Finally, the following commutative diagram and the fact that 
the morphism $CC\ast\pi\colon CCC\to CCT$ is an epimorphism imply the equality 
$\mathrm{(\ref{eq:diagonalizable 4})}$.
$$\xymatrix{
CCT \ar@{->>}[rrr]^{r\ast T} \ar[ddd]^{\wr}_{C\ast\tau^{C,T}} & & & 
CT \ar[ddd]_{\wr}^{\tau^{C,T}}\\
& CCC \ar@{->>}[ul]_{CC\ast\pi} 
\ar@{->>}[r]^{r\ast C} \ar[d]^{\wr}_{C\ast\sigma} 
& CC \ar@{->>}[ru]^{C\ast\pi} \ar[d]_{\wr}^{\sigma} & \\
& CCC \ar@{->>}[ld]_{C\ast\pi\ast C} & CC \ar@{->>}[rd]^{\pi\ast C} & \\
CTC \ar@{>->}[d]_{C\ast\zeta} & & & TC \ar@{>->}[d]^{\zeta}\\
CCC \ar@{->>}[rrr]_{r\ast C} & & & CC.
}$$ 
\end{proof}

\begin{para}
\label{ex:cone of composition}
{\bf Example.}\ \ 
Let $x\onto{f} y\onto{g}z$ to be a pair of morphisms in a 
complicial exact categories. 
Then applying $\Cone$ to the 
admissible exact sequence 
$$[x\onto{f}y]\overset{{\tiny{
\left (\begin{pmatrix}\id_x\\ -f\end{pmatrix},
\begin{pmatrix}g\\ -\id_y\end{pmatrix} \right )}}}{\rinf}[x\oplus y
\overset{\tiny{\begin{pmatrix}gf & 0\\ 0 & \id_y\end{pmatrix}}}{\to}z\oplus y]\overset{\tiny{
\left(\begin{pmatrix}f &\id_y\end{pmatrix},
\begin{pmatrix}\id_z & g\end{pmatrix}
\right)
}}{\rdef}[y\onto{g}z]$$
in $\cC^{[1]}$, we obtain the admissible exact sequence
\begin{equation}
\Cone f\overset{\tiny{\begin{pmatrix}\Cone(\id_x,g)\\ -\Cone(f,\id_y)\end{pmatrix}}}{\rinf}\Cone gf\oplus Cy \overset{\tiny{
\begin{pmatrix}\Cone(f,\id_z) & \Cone(\id_y,g)\end{pmatrix}}}{\rdef}\Cone g
\end{equation}
in $\cC$. 
The morphism 
$\Cone(\kappa_{f},\kappa_{gf})\colon \Cone g\to \Cone\Cone(\id_x,g)$ 
is a $C$-homotopy equivalence. (See \cite[3.34]{Moc10}.) 
\end{para}

\begin{para}
\label{ex:fundamental triangle of cone}
{\bf Example.}\ \ 
Let $f\colon x\to y$ be a morphism in a complicial exact category. 
Then there are sequence of 
(homotopy) commutative squares 
$[x\onto{f}y] \onto{(\iota_x,\kappa_f)}[Cx\onto{\mu_f}\Cone f] 
\onto{(0,\id_{\Cone f},-\mu_fr_x)}[0\to \Cone f]$ 
and by applying $\Cone$ to these (homotopy) squares, we obtain the 
equality
\begin{equation}
\label{eq:fundamental triangle of cone}
\Cone(0,\id_{\Cone f},-\mu_fr_x)\Cone(\iota_x,\kappa_f)=\id_{\Cone f}.
\end{equation}
$$
\xymatrix{
\Cone f \ar[rr]^{\id_{\Cone f}} \ar[rd]_{\Cone(\iota_x,\kappa_f)}& &\Cone f\\
& \Cone\mu_f. \ar[ru]_{\Cone(0,\id_{\Cone f},-\mu_fr_x)}
}
$$
\end{para}

We can show the following.

\begin{para}
\label{lemdf:comp functor induce on C_h^[1]}
{\bf Lemma-Definition.}\ \ 
Let $\cC$ and $\cC'$ be complicial exact categories and 
let $(f,d)\colon \cC\to \cC'$ be a complicial exact functor. 
Then $f$ induces a functor $f_h^{[1]}\colon \cC_h^{[1]}\to {\cC'}_h^{[1]}$ 
by Lemma-Definition~\ref{lemdf:functoriality of Fbh}. 
Namely $f^{[1]}_h$ sends an object $[a\colon x\to y]$ 
in $\cC_h^{[1]}$ to an object 
$[fa\colon fx\to fy]$ in ${\cC'}_h^{[1]}$ 
and sends a homotopy commutative square 
$(\psi,\phi,H)\colon [a\colon x\to y] \to [a'\colon x'\to y']$ 
in $\cC_h^{[1]}$ to 
a homotopy commutative square 
$(f\psi,f\phi,fH\cdot c_x)\colon [fa\colon fx\to fy]\to [fa'\colon fx' \to fy']$ in ${\cC'}_h^{[1]}$. 
We define $d^{\Cyl}\colon \Cyl f^{[1]}_h\isoto f\Cyl$ 
to be a natural equivalence by setting 
$d^{\Cyl}:=\begin{pmatrix}\id_f& 0\\ 0 & d\end{pmatrix}$. 
The natural equivalence $d^{\Cyl}$ induces a natural equivalence 
$d^{\Cone}\colon \Cone f_h^{[1]}\isoto f\Cone$ 
characterized by the following equalities:
\begin{equation}
\label{eq:kappa and ca }
d^{\Cone}\cdot (\kappa\ast f)=f\ast\kappa,
\end{equation}
\begin{equation}
\label{eq:mu and ca}
(f\ast \mu)\cdot c=d^{\Cone}\cdot (\mu\ast f).
\end{equation}
$$
\xymatrix{
f\ran \ar@{>->}[r]^{\kappa\ast f} \ar@{>->}[rd]_{f\ast \kappa} & 
\Cone f_h^{[1]} \ar[d]_{\wr}^{d^{\Cone}} & Cf\dom \ar[l]_{\ \ \ \ \mu\ast f} \ar[d]_{\wr}^c\\
& f\Cone & fC\dom. \ar[l]^{\ \ \ \ f\ast \mu}
}$$

We write $C^{\lv}$ and $\iota^{\lv}$ 
for the endofunctor on $\cC^{[1]}_h$ 
and the natural transformation $\id_{\cC^{[1]}_h}\to C^{\lv}$ 
defined by level-wisely. 
By applying the argument above to the complicial 
functor $(C,\sigma)\colon\cC\to\cC$, 
we obtain natural equivalences 
$\sigma^{\Cyl}\colon C\Cyl\isoto \Cyl C^{\lv}$ and 
$\sigma^{\Cone}\colon C\Cone \isoto \Cone C^{\lv}$ 
which satisfy the following equalities:
\begin{equation}
\label{eq:sigmacyl}
\Cyl\ast\iota^{\lv}=\sigma^{\Cyl}\cdot (\iota\ast \Cyl),
\end{equation}
\begin{equation}
\label{eq:sigmacone}
\Cone\ast\iota^{\lv}=\sigma^{\Cone}\cdot (\iota\ast \Cone).
\end{equation}
In particular, $(\Cyl,\sigma^{\Cyl}) $, $(\Cone,\sigma^{\Cone})\colon\cC^{[1]}\to\cC$ are complicial exact functors. 
\qed
\end{para}

Recall the definitions of $\psi_f$ and $\Theta_f$ 
from \ref{df:mapping cone and cylinders} 
and \ref{df:ThetaX} respectively and 
recall the definition of strictly ordinary 
from \ref{df:ordinary comp exact cat}. 

\begin{para}
\label{lem:Cone(psif,0,Thetaf)}
{\bf Lemma.}\ \ 
{\it
Let $\cC$ be a strictly ordinary complicial exact category 
and let $f\colon x\to y$ be a morphism in $\cC$. 
Then a homotopy commutative square 
$(\psi_f,0,\Theta_f)\colon[C\Cone f\to 0] \to [Tf\colon Tx\to Ty]$ 
induces an equality
\begin{equation}
\label{eq:Cone(psi,0,Thetaf)}
\Cone(\psi_f,0,\Theta_f)=-\sigma_T^{\Cone}.
\end{equation}
In particular the left square in the commutative diagram 
$\mathrm{(\ref{eq:Conepsif diagram})}$ below is a push-out diagram.
\begin{equation}
\label{eq:Conepsif diagram}
\xymatrix{
\Cone f \ar@{>->}[r]^{\iota_{\Cone f}} \ar[d]_{\psi_f} & 
C\Cone f \ar@{->>}[r]^{\pi_{\Cone f}} 
\ar[d]_{\tiny{\begin{pmatrix}-\Theta_f\\ C\psi_f\end{pmatrix}}} &
T\Cone f \ar[d]_{\wr}^{-\sigma_T^{\Cone}}\\
Tx \ar@{>->}[r]_{\!\!\!\!\!\!\!\!\!\!\!\tiny{\begin{pmatrix}-Tf \\ \iota_{Tx}\end{pmatrix}}} & 
Ty\oplus CTx \ar@{->>}[r]_{\tiny{\begin{pmatrix}\kappa_{Tf} & \mu_{Tf}\end{pmatrix}}} & \Cone Tf.
}
\end{equation}
Thus we can take 
\begin{equation}
\label{eq:Conepsif}
\Cone \psi_f=Ty\oplus CTx,
\end{equation}
\begin{equation}
\label{eq:kappapsif}
\kappa_{\psi_f}=\begin{pmatrix}-Tf\\ \iota_{Tx}\end{pmatrix}\ \ \text{and},
\end{equation}
\begin{equation}
\label{eq:mupsif}
\mu_{\psi_f}=\begin{pmatrix}-\Theta_f\\ C\psi_f \end{pmatrix}.
\end{equation}
}
\end{para}

\begin{proof}
First we assume that $y=x$ and $f=\id_x$. 
In this case, $\Cone(\psi_f,0,\Theta_f)=\Cone(\pi_x,0,\pi_xr_x)=\Cone(\pi_x,\pi_x)\cdot\zeta_x$. 
Since $\cC$ is strictly ordinary we have an equality 
$\zeta=-\sigma\cdot\zeta$. Thus we have equalities 
$(C\ast\pi) \cdot \pi=-(C\ast\pi) \cdot \sigma\cdot\zeta =
-\tau^{T,C}\cdot(\pi\ast C)\cdot\zeta=-\tau^{T,C}$. 
For general case, 
the equality $\mathrm{(\ref{eq:Cone(psi,0,Thetaf)})}$ 
follows from the commutative diagram below.
$$
\xymatrix{
\Cone f \ar@{>->}[rr]^{\iota\ast\Cone f} \ar[ddd]_{\psi_f} & & 
C\Cone f \ar@{->>}[rr]^{\pi\ast\Cone f} & & 
T\Cone f \ar[ddd]^{-\sigma_T^{\Cone}}_{\wr}\\
& Cx \ar@{>->}[r]^{\iota\ast Cx} \ar@{->>}[d]_{\pi_x} \ar[lu]^{\mu_f} & 
CCx \ar@{->>}[r]^{\pi_{Cx}} \ar[u]_{C\ast\mu_f} 
\ar[d]_{\tiny{\begin{pmatrix}-\pi_xr_x\\ C\ast\pi_x\end{pmatrix}}} & 
TCx \ar[ru]_{T\ast\mu_f} \ar[d]_{\wr}^{-\tau^{T,C}}\\
& Tx \ar@{>->}[r]_{\tiny{\begin{pmatrix}-\id_{Tx}\\ \iota\ast Tx\end{pmatrix}}}
\ar[ld]_{\id_{Tx}} & 
Tx\oplus CTx \ar[d]^{\tiny{\begin{pmatrix}Tf & 0\\ 0 &\id_{CTx}\end{pmatrix}}} 
\ar@{->>}[r]_{\tiny{\begin{pmatrix}\iota\ast Tx & \id_{CTx}\end{pmatrix}}} & 
CTx \ar[rd]^{\mu\ast Tf}\\
Tx \ar@{>->}[rr]_{\tiny{\begin{pmatrix}-Tf \\ \iota\ast Tx\end{pmatrix}}} & &  
Ty\oplus CTx \ar@{->>}[rr]_{\tiny{\begin{pmatrix}\kappa\ast Tf & \mu\ast T f\end{pmatrix}}} & & \Cone Tf.
}
$$
\end{proof}

\begin{para}
\label{lem:homotopy commutative square}
{\bf Lemma.}\ \ 
{\it
Let $\cC$ be a complicial exact category. Then\\
$\mathrm{(1)}$ 
Let $(a,b,H)\colon[x\onto{f_0}y]\to[x'\onto{g_0}y']$ 
be a homotopy commutative square and let 
$V\colon f_0\Rightarrow_C f_1$ and $W\colon g_0\Rightarrow_C g_1$ be 
$C$-homotopies. Then\\
$\mathrm{(i)}$ 
A triple $(a,b,-W Ca+bV+H)$ 
is a homotopy commutative square from $[x\onto{f_1}y]$ to 
$[x'\onto{g_1}y']$.\\
$\mathrm{(ii)}$ 
We define a morphism $\Cyl(V,W,H)\colon C(\Cyl a)\to \Cyl b$ by setting 
\begin{equation}
\label{eq:cyl psi theta H}
\Cyl(V,W,H):=\begin{pmatrix}W & (-W Ca+bV)r_x\\ 
0 & C(V)\sigma \end{pmatrix}.
\end{equation}
Then $\Cyl(V,W,H)$ is a $C$-homotopy from $\Cyl(f_0,g_0,H)$ to 
$\Cyl(f_1,g_1,-W Ca+bV+H)$ and makes the left square below commutative.\\
$\mathrm{(iii)}$ 
There exists a morphism $\Cone(V,W,H)\colon C\Cone a\to \Cone b$ 
which makes the right square below commutative 
and it is a $C$-homotopy from $\Cone(f_0,g_0,H)$ to 
$\Cone(f_1,g_1,-W Ca+bV +H)$.
$$\xymatrix{
Cx \ar@{>->}[r]^{C\ast{\xi_1}_a} \ar[d]_{V} & C \Cyl a \ar@{->>}[r]^{C\ast\eta_a} 
\ar[d]_{\Cyl(V,W,H)} & C\Cone a \ar@{-->}[d]^{\Cone(V,W,H)}\\
y \ar@{>->}[r]_{{\xi_1}_b} & \Cyl b \ar@{->>}[r]_{\eta_b} & \Cone b.
}$$
$\mathrm{(2)}$ 
Let $(a,b,H)\colon[x\onto{f}x']\to[y\onto{g}y']$ be a homotopy 
commutative square in $\cC$ and 
let $f'\colon x'\to x$ and $g\colon y' \to y$ be $C$-homotopy inverse 
morphisms of $f$ and $g$ respectively. 
Namely there exists $C$-homotopies 
$K\colon f'f\Rightarrow_C\id_x$, 
$K'\colon ff'\Rightarrow_C \id_{x'}$, 
$L\colon g'g\Rightarrow_C \id_y$ and 
$L'\colon gg'\Rightarrow_C \id_{y'}$. 
We set $H':=-g'BK'-g'HC(f')+LC(af')$ and $U:=-LC(a)+bK+H'\star H$. 
Then\\
$\mathrm{(i)}$ 
The triple $(b,a,H')$ is a homotopy commutative square 
$[x'\onto{f'}x]\to [y'\onto{g'}y]$.\\
$\mathrm{(ii)}$ 
The morphism $\Cyl(f,g,H)\colon\Cyl a \to \Cyl b$ 
is a $C$-homotopy equivalence.\\
$\mathrm{(iii)}$ 
The morphism $\Cone(f,g,H)\colon\Cone a \to \Cone b$ 
is a $C$-homotopy equivalence. 
}
\end{para}

\begin{proof}
A proof of $\mathrm{(1)}$ is given in \cite[3.11]{Moc10} 
and 
$\mathrm{(2)}$ 
$\mathrm{(i)}$ is straightforward. 

\sn
$\mathrm{(ii)}$ 
By $\mathrm{(1)}$, 
$\Cyl(K,L,H'\star H)$ is a $C$-homotopy 
$\Cyl(f'f,g'g,H'\star H)\Rightarrow_C \Cyl(\id_x,\id_y,U) $ and 
$\Cyl(\id_x,\id_y,U)=
\begin{pmatrix}
\id_y & -U\\ 0 & \id_{Cx}
\end{pmatrix}$ is an isomorphism. 
Thus ${\Cyl(\id_x,\id_y,U)}^{-1}\Cyl(f',g',H')$ 
is a left $C$-homotopy inverse of the morphism $\Cyl(f,g,H)$. 
Similarly it turns out that $\Cyl(f,g,H)$ admits a right $C$-homotopy 
inverse. 
Thus $\Cyl(f,g,H)$ is a $C$-homotopy equivalence by 
Lemma~\ref{lem:left and right homotopy inverse}. 
A proof of $\mathrm{(iii)}$ is similar.
\end{proof}

\begin{para}
\label{lem:homotopy commutative diagram induces cyl and cone seq}
{\bf Lemma.}\ \ 
{\it
Let $\cC$ be a complicial exact category and 
$x\overset{f}{\rinf}y\overset{g}{\rdef}z$ and 
$x'\overset{f'}{\rinf}y'\overset{g'}{\rdef}z'$ be 
Frobenius admissible exact sequences in $\cC$ and 
$a\colon x\to x'$ and $b\colon y\to y'$ morphisms in $\cC$ 
and $H\colon bf\Rightarrow_C f'a$ a $C$-homotopy. 
Then there exists a morphism $c\colon z\to z'$ and 
a retraction 
$\rho\colon Cy\to Cx$ of $Cf\colon Cx\to Cy$ 
such that $-g'H\rho\colon Cy\to z'$ 
is a $C$-homotopy from $cg$ to $g'b$. 
$$\xymatrix{
x \ar@{>->}[r]^f \ar[d]_a & y \ar@{->>}[r]^g \ar[d]_b & z \ar@{-->}[d]^c\\
x' \ar@{>->}[r]_{f'} & y' \ar@{->>}[r]_{g'}& z'.
}$$
Moreover sequences 
$\Cyl a \overset{\tiny{\Cyl(f,f',H)}}{\rinf} \Cyl b 
\overset{\tiny{\Cyl(g,g',-g'H\rho)}}{\rdef} \Cyl c $ and 
$\Cone a \overset{\tiny{\Cone(f,f',H)}}{\rinf} \Cone b 
\overset{\tiny{\Cone(g,g',-g'H\rho)}}{\rdef} \Cone c$ 
are Frobenius admissible sequences and we have the equality
\begin{equation}
\label{eq:-g7Hrho star H=0 }
(-g'H\rho)\star_C H=0.
\end{equation}
}
\end{para}

\begin{proof}
There exists a morphism $\rho\colon Cy\to Cx$ 
which is a retraction of $Cf$ by Lemma~\ref{lem:split Frob exact}. 
Since $\displaystyle{\begin{pmatrix}\id_{y'} & -H\rho\\ 0 & \id_y 
\end{pmatrix}}$ and 
$\displaystyle{\begin{pmatrix}f' & 0\\ 0 &Cf\end{pmatrix}}$ 
are Frobenius admissible monomorphisms, the composition 
$$\displaystyle{\begin{pmatrix}f' & -H\\ 0 & Cf \end{pmatrix}=
\begin{pmatrix}\id_{y'} & -H\rho\\ 0 & \id_{Cy} \end{pmatrix}  
\begin{pmatrix}f' & 0\\ 0 & Cf \end{pmatrix}}$$ 
is also a Frobenius admissible monomorphism. 
Similarly it turns out that the composition 
$$\displaystyle{\begin{pmatrix}g' & g'H\rho\\ 0 & Cg \end{pmatrix}=
\begin{pmatrix}g' & 0\\ 0 & Cg \end{pmatrix}
\begin{pmatrix}\id_{y'} & H\rho\\ 0& \id_{Cy}\end{pmatrix}}$$
is a Frobenius admissible epimorphism. 
In the commutative diagram below, 
by the universality of cokernel, there exists 
the dotted morphism 
$\displaystyle{{\footnotesize{\begin{pmatrix}c\\ u \end{pmatrix}}}\colon z\to z'}$ 
in the diagram below which makes 
the right square below commutative: 
$$\xymatrix{
x \ar@{>->}[r]^f \ar[d]_{\tiny{\begin{pmatrix}a \\ -\iota_x\end{pmatrix}}} & 
y \ar@{->>}[r]^g \ar[d]_{\tiny{\begin{pmatrix}b\\ -\iota_y\end{pmatrix}}} & 
z \ar@{-->}[d]^{\tiny{\begin{pmatrix}c\\ u\end{pmatrix}}}\\
x'\oplus Cx \ar@{>->}[r]^{\tiny{\begin{pmatrix}f' & -H\\ 0 & Cf \end{pmatrix}}} 
\ar@{>->}[rd]_{\tiny{\begin{pmatrix}f' & 0\\ 0 &Cf\end{pmatrix}}} & 
y'\oplus Cy \ar@{->>}[r]^{\tiny{\begin{pmatrix}g' & g'H\rho\\ 0 & Cg \end{pmatrix}}} 
\ar[d]^{\wr}_{\varphi} & 
z'\oplus Cz\\
& y'\oplus Cy \ar@{->>}[ru]_{\tiny{\begin{pmatrix}g' & 0\\ 0& Cg\end{pmatrix}}}
}$$
where the morphism $\varphi$ is $\displaystyle{\begin{pmatrix}\id_{y'} & H\rho\\ 0 & \id_{Cy}\end{pmatrix}}$. 
By the commutativity of the right square above, 
we have the equalities $u g=-Cg\iota_y=-\iota_zg$ and 
it turns out that $u=-\iota_z$ by surjectivity of $g$. 
By applying the $3\times 3$-lemma to the commutative diagram below, 
we obtain a Frobenius admissible exact sequence 
$\Cone a \overset{\tiny{\Cone(f,f',H)}}{\rinf} \Cone b 
\overset{\tiny{\Cone (g,g',g'H\rho)}}{\rdef} \Cone c$. 
$$\xymatrix{
x \ar@{>->}[r]^f \ar@{>->}[d]_{{\xi_1}_a} & y \ar@{->>}[r]^g \ar@{>->}[d]_{{\xi_1}_b} & z \ar@{>->}[d]^{{\xi_1}_c}\\
\Cyl a \ar@{>->}[r]^{\Cyl(f,f',H)} \ar@{->>}[d]_{\eta_a} & 
\Cyl b \ar@{->>}[r]^{\Cyl(g,g',g'H\rho)} \ar@{->>}[d]_{\eta_b} & \Cyl c 
\ar@{->>}[d]^{\eta_c}\\
\Cone a \ar@{>->}[r]_{\Cone(f,f',H)} & \Cone b \ar@{->>}[r]_{\Cone(g,g',g'H\rho)} &\Cone c.
}$$
Finally we have equalities 
$(-g'H\rho)\star_C H=g'H-g'H\rho Cf=g'H-g'H=0$.

\end{proof}

The proof above and \ref{lemdf:comp functor induce on C_h^[1]} 
show the following assertion.

\begin{para}
\label{cor:exactness of Cyl and Cone}
{\bf Corollary.}\ \ 
{\it
Let $\cC$ be a complicial exact category. 
Then the functors $\Cyl$ and $\Cone$ send a level Frobenius admissible 
exact sequence in $\cC_h^{[1]}$ to a Frobenius admissible exact sequence 
in $\cC$. 
In particular if $\cC$ is ordinary, then $\Cyl$ and $\Cone$ 
are exact functors. 

Moreover if $\cC$ is normal, then the pairs 
$(\Cyl,\sigma^{\Cyl})$ and $(\Cone,\sigma^{\Cone})$ are complicial exact functors from $\cC^{[1]}_h$ to $\cC$. 
}\qed
\end{para}

\begin{para}
\label{cor:ConeCone}
{\bf Corollary.}\ \ 
{\it
Let $\cC$ be a complicial exact category and let 
$[f\colon x\to y] \onto{(a,b,H)} [f'\colon x'\to y']$ 
be a homotopy commutative square in $\cC$. 
Then we have the canonical isomorphism
\begin{equation}
\label{eq:ConeCone}
\Cone\Cone(a,b,H)\isoto\Cone\Cone(f,f',-H).
\end{equation}
}
\end{para}

\begin{proof}
If $\cC$ is normal and ordinary, it is just a consequence of 
\ref{lemdf:comp functor induce on C_h^[1]} and 
\ref{cor:exactness of Cyl and Cone}. 
For general $\cC$, it follows from the fact that 
both $\Cone\Cone(a,b,H)$ and $\Cone\Cone(f,f',-H)$ 
are naturally isomorphic to $z$ in the 
commutative diagram below. 
$$
\xymatrix{
x \ar@{>->}[r]^{\tiny{\begin{pmatrix}a\\ -\iota_x\end{pmatrix}}} 
\ar@{>->}[d]_{\tiny{\begin{pmatrix}f\\ -\iota_x\end{pmatrix}}} & 
x'\oplus Cx \ar@{->>}[r]^{\tiny{\begin{pmatrix}\kappa_a & \mu_a\end{pmatrix}}} 
\ar@{>->}[d]_{\tiny{\begin{pmatrix}f' & H\\ -\iota_{x'} & 0\\ 0 & Cf\\ 
0 & -C\ast\iota_x\end{pmatrix}}} &
\Cone a \ar@{>->}[d]^{\tiny{\begin{pmatrix}\Cone(f,f',H)\\ -\iota_{\Cone a}\end{pmatrix}}}\\
y\oplus Cx \ar@{>->}[r]_{\!\!\!\!\!\!\!\!\!\!\!\!\!\!\!\!\!\!\!\!\!\!\!\!\tiny{\begin{pmatrix}b & -H\\ 0 & Ca\\ -\iota_{y} & 0\\ 
0 & -C\ast\iota_x\end{pmatrix}}} 
\ar@{->>}[d]_{\tiny{\begin{pmatrix}\kappa_f &\mu_f\end{pmatrix}}} &
y'\oplus Cx'\oplus Cy\oplus CCx 
\ar@{->>}[r]^{\tiny{\begin{pmatrix}\kappa_b & 0 & \mu_b & 0\\ 
0 & C\ast\kappa_a & 0 & C\ast\mu_a\end{pmatrix}}} 
\ar@{->>}[d]^{\tiny{\begin{pmatrix}\kappa_{f'} & \mu_{f'} & 0 & 0\\ 
0 & 0 & C\ast\kappa_f & C\ast\mu_f\end{pmatrix}}} & 
\Cone b\oplus C\Cone a \ar@{->>}[d]\\
\Cone f \ar@{>->}[r]_{\!\!\!\!\!\!\!\!\!\!\!\!\!\!\!\!\!\!\!\!\!\!\!\!\tiny{\begin{pmatrix}\Cone(a,b,H)\\ -\iota_{\Cone f}\end{pmatrix}}} & 
\Cone f'\oplus C\Cone f \ar@{->>}[r] &
z.
}
$$
\end{proof}

\begin{para}
\label{ex:Conekappaf}
{\bf Example.}\ \ 
Let $\cC$ be a complicial exact category and 
let $f\colon x\to y$ be 
a morphism in $\cC$. 
By applying \ref{cor:ConeCone} to the commutative square 
$$\footnotesize{\xymatrix{
0 \ar[r] \ar[d] & x \ar[d]^f\\
y \ar[r]_{\id_y} & y,
}}$$
we obtain the canonical isomorphism 
\begin{equation}
\label{eq:conekappaf}
\Cone \kappa_f\isoto \Cone \Cone(0,f).
\end{equation}
A morphism $(\id_x,0)\colon [\Cone(0,f)\colon x\to Cy]\to [x\to 0]$ in 
$\cC^{[1]}$ induces a $C$-homotopy equivalence
\begin{equation}
\label{eq:Cone idx0}
\Cone(\id_x,0)\colon\Cone\Cone(0,f)\to Tx
\end{equation}
by \ref{lem:homotopy commutative square} $\mathrm{(2)}$.
\end{para}

\begin{para}
\label{nt:Pat,Fib}
{\bf (Homotopy fiber functor, path space functor).}\ \ 
Let $\cC$ be a complicial exact category. 
Then there are dual notions of the mapping cylinder and 
the mapping cone functors. 
We define 
$\hFib,\ \Pat\colon\cC^{[1]} \to \cC$ to be functors which 
called the {\it homotopy fiber functor} and 
the {\it mapping path space functor} respectively 
in the following way. 
\begin{equation}
\Pat:=\dom\oplus P\ran
\end{equation}
\begin{equation}
\hFib:=\dom\times_{\ran}P\ran
\end{equation}
where $\hFib$ is defined by the following pull back diagram: 
\begin{equation}
\label{eq:Hfib pullback}
\xymatrix{
T^{-1} \ar@{>->}[r]^{\lambda} \ar[d]_{\id_{T^{-1}}} 
& 
\hFib \ar[d]^{\nu} \ar@{->>}[r]^{\upsilon} \ar@{}[dr]|\bigstar& 
\dom \ar[d]^{\epsilon}\\
T^{-1} \ar@{>->}[r]_{j} & P\ran \ar@{->>}[r]_{q_{\ran}} & \ran 
}
\end{equation}
where $\lambda$ is induced from the universal property of $\hFib$. 
Then there is an admissible exact sequence
\begin{equation}
\label{eq:canonical hfib pat sequence}
\hFib \overset{\rho}{\rinf} \Pat \overset{t_2}{\rdef} \ran
\end{equation}
where 
\begin{equation}
\label{eq:t2, rho df}
t_2:=
\begin{pmatrix}
\epsilon & q_{\ran}
\end{pmatrix}
\ \ \text{and}\ \ \rho=
\begin{pmatrix}
\upsilon\\ 
-\nu
\end{pmatrix}
.
\end{equation}

\sn
Moreover we define 
$t_1\colon \Pat \to \dom$ and 
$\chi\colon\dom \to \Pat$ to be natural transformations by 
setting
\begin{equation}
\label{eq:t1 and chi df}
t_1:=
\begin{pmatrix}
\id_{\dom} & 0
\end{pmatrix}
\ \ \text{and  }
\chi:=
\begin{pmatrix}
\id_{\dom}\\ 
0
\end{pmatrix} 
.
\end{equation}
Then we have the following commutative diagram:
\begin{equation}
\label{eq:canonical nat between Pat}
\xymatrix{
& \dom \ar[ld]_{\id_{\dom}} \ar[rd]^{\epsilon} \ar[d]^{\chi} & \\ 
\dom & \Pat \ar@{->>}[r]_{t_2} \ar@{->>}[l]^{t_1} & \ran.
}
\end{equation}
Since the functor $\hFib$ is obtained by the pull-back diagram, 
we can choose the fiber products in the following way. 
For an object $x$ in $\cC$, 
since the diagrams below are pull-backs, 
we shall assume that 
$\hFib\id_x=Px$, $\nu_{\id_{x}}=\id_{Px}$, 
$\upsilon_{\id_x}=q_x$, $\hFib(0\to x)=T^{-1}x$, 
$\nu_{0\to x}=\iota_{T^{-1}x}$, 
$\hFib(x\to 0)=x$ and $\upsilon_{x\to 0}=\id_x$.
$$\xymatrix{ 
Px \ar@{->>}[r]^{q_x} \ar[d]_{\id_{Px}} & x \ar[d]^{\id_{x}}\\
Px \ar@{->>}[r]_{q_x} & x,
}\ 
\xymatrix{
T^{-1}x \ar@{->>}[r] \ar[d]_{\iota_{T^{-1}x}} & 0 \ar[d]\\
Px \ar@{->>}[r]_{q_x} & x,
}\ 
\xymatrix{
x \ar@{->>}[r]^{\id_x} \ar[d] & x \ar[d]\\
P(0)=0 \ar@{->>}[r]_{q_0} & 0.
}
$$
Recall the definitions $s$, $t$ and $\Delta\colon\cC\to\cC^{[1]}$ 
from \ref{df:mapping cone and cylinders}. 
Then we shall assume that $\hFib s=\id_{\cC}$, 
$\hFib t=T^{-1}$ and $\hFib\Delta =P$.
\end{para}

\begin{para}
\label{lem:Cone =THfib}
{\bf Lemma.}\ \ 
{\it
Let $\cC$ be a complicial exact category. 
Then there is a canonical natural equivalence between functors 
$\cC^{[1]}\to\cC$. 
\begin{equation}
\label{eq:Cone=THfib}
\Cone \isoto T\hFib.
\end{equation}
}
\end{para}

\begin{proof}
Let $f\colon x\to y$ be a morphism in $\cC$. 
By applying the operation $\Cone$ to the 
Frobenius exact sequence in $\cC^{[1]}$, 
\begin{equation}
{\footnotesize{
\begin{bmatrix}
\xymatrix{
T^{-1}x \ar[d]_{\scriptstyle{T^{-1}f}}\\
T^{-1}y
}
\end{bmatrix} 
\begin{matrix}
\xymatrix{
\overset{\scriptstyle{\begin{pmatrix}\id_{T^{-1}x}\\ T^{-1}f\end{pmatrix}}}{\to} \ar@{}[d]\\
\underset{\id_{T^{-1}y}}{\to}\\
}
\end{matrix} 
\begin{bmatrix}
\xymatrix{
T^{-1}x\oplus T^{-1}y \ar[d]_{\scriptstyle{\begin{pmatrix}0 & 
\id_{T^{-1}y} \end{pmatrix}}}\\
T^{-1}y
}
\end{bmatrix} 
\begin{matrix}
\xymatrix{
\overset{\scriptstyle{\begin{pmatrix}-T^{-1}f &\id_{T^{-1}y} \end{pmatrix}}}{\to} \ar@{}[d]\\
\to\\
}
\end{matrix} 
\begin{bmatrix}
\xymatrix{
T^{-1}y \ar[d] \\
0
}
\end{bmatrix}
,}}
\end{equation}
we obtain the Frobenius admissible exact sequence in $\cC$
\begin{equation}
\label{eq:Cone is Hfib}
\Cone T^{-1}f \rinf x\oplus Py 
\overset{{\footnotesize{\begin{pmatrix}-f & q_y \end{pmatrix}}}}{\rdef} y.
\end{equation}
Comparing the sequence $\mathrm{(\ref{eq:canonical hfib pat sequence})}$, 
we obtain the canonical isomorphism 
$\Cone T^{-1}\isoto \hFib$ as the functors $\cC^{[1]}\to \cC$. 
\end{proof}

\begin{para}
\label{df:homotopy pushout pullback}
{\bf Definition (Homotopy push-out, homotopy pull-back).}\ \ 
Let $\cC$ be a complicial exact category and let 
$x \overset{f}{\leftarrow} y \overset{g}{\to}z$ be 
a pair of morphisms in $\cC$. 
Then we set
\begin{equation}
\label{eq:homotopy pushout df} 
x\underset{y}{\sqcup_h}z:=
\Cone(y\overset{\tiny{\begin{pmatrix}f\\ g\end{pmatrix}}}{\to}x\oplus z).
\end{equation} 
Then there exists a Frobenius admissible exact sequence
\begin{equation}
\label{eq:homotopy pushout seq}
y\overset{\tiny{\begin{pmatrix}\iota_y\\ f\\ g\end{pmatrix}}}{\rinf} 
Cy\oplus x\oplus z 
\overset{\tiny{\begin{pmatrix}\psi_{\tiny{\begin{pmatrix}f\\ g\end{pmatrix}}} & i_g & -i_f\end{pmatrix}}}{\rdef} x\underset{y}{\sqcup_h}z
\end{equation}
and homotopy commutative diagram
\begin{equation}
\label{eq:homotopy pushout diagram}
\xymatrix{
y \ar[r]^g \ar[d]_f & z \ar[d]^{i_f}\\
x\ar[r]_{i_g} & x\underset{y}{\sqcup_h}z.
}
\end{equation}
We call it the {\it canonical homotopy push out} ({\it of the 
diagram $x \overset{f}{\leftarrow} y \overset{g}{\to}z$}) and it has 
the following universal property:

\sn
For a homotopy commutative square 
$[f\colon y\to x] \onto{(g,g',H)} [f'\colon z\to u]$, there 
exists a unique morphism $a\colon x\underset{y}{\sqcup_h}z \to u$ such that 
$f'=ai_f$, $g'=ai_g$ and $H=\psi_{\tiny{\begin{pmatrix}f\\ g\end{pmatrix}}}a$.

Dually there is a notion of homotopy pull-back. 
Let $x\onto{f} y\overset{g}{\leftarrow} z$ be a pair of morphism in $\cC$. 
Then we set 
\begin{equation}
\label{eq:homotopy pullback df} 
x\underset{y}{\times_h}z:=
\hFib(x\oplus z\overset{\tiny{\begin{pmatrix}f & g\end{pmatrix}}}{\to}y).
\end{equation} 
Then there exists a Frobenius admissible exact sequence
\begin{equation}
\label{eq:homotopy pullback seq}
x\underset{y}{\times_h}z
\overset{\tiny{\begin{pmatrix}u_{f,g}\\ p_g\\ -p_f\end{pmatrix}}}{\rinf} 
Py\oplus x\oplus z 
\overset{\tiny{\begin{pmatrix}q_y & f & g\end{pmatrix}}}{\rdef} y
\end{equation}
and $P$-homotopy commutative diagram
\begin{equation}
\label{eq:homotopy pullback diagram}
\xymatrix{
x\underset{y}{\times_h}z \ar[r]^{p_f} \ar[d]_{p_g} & z \ar[d]^{g}\\
x\ar[r]_{f} & y.
}
\end{equation}
We call it the {\it canonical homotopy pull-back} ({\it of the 
diagram $x \onto{f} y \overset{g}{\leftarrow}z$}) and 
it has the following universal property:

\sn
For a $P$-homotopy commutative square 
$[g'\colon u\to x] \onto{(f,f',H)} [g\colon z\to y]$, there 
exists a unique morphism $a\colon u \to x\underset{y}{\times_h}z$ such that 
$f'=p_fa$, $g'=p_ga$ and $H=u_{f,g}a$.
\end{para}

\section{Homotopy theory of relative complicial exact categories}
\label{sec:hom th of rel comp ex cat}

As pointed out in \cite[Remark 3.6]{Bar15}, 
in the view of $(\infty,1)$-category theory, 
a complicial exact category with weak equivalences is actually 
nothing more than a model for some stable $(\infty,1)$-category. 
In this section, we will further develop such a theory under 
the name of relative complicial exact categories. 
In the first subsection \ref{subsec:null class}, 
we will study the notion of (thick) null classes which 
is a variant of the notion of (thick) 
triangulated subcategories of triangulated categories. 
In the next subsection \ref{subsec:rel compl exact cat}, we will interpret 
the fundamental terminologies about 
relative (complicial) exact categories. 
A typical example of relative complicial exact categories 
is a category of chain complexes on a relative exact category 
which will be detailed in subsection \ref{subsec:chain complexes on cec}. 
In the next subsection 
\ref{subsec:Homotopy category of relative categories}, 
we will study homotopy categories of relative complicial exact categories 
and introduce the notion of derived equivalences. 
In \ref{cor:complification}, we will show the (large) relative category of 
relative exact categories with derived equivalences is 
categorical homotopy equivalent to the (large) relative category of 
relative complicial exact categories with derived equivalences. 
This stands for that a relative exact category is also a model 
for some stable infinity category. 
In the final subsection \ref{subsec:homology theory on complicial exact categories}, 
following \cite[\S 7]{Wal85}, 
we will study the notion of homology theory on complicial exact categories. 
In particular we will consider a homology theory induced from 
a $t$-structure on a homotopy category 
\ref{ex:homoogy theory associated with t-structure}.

\subsection{Null classes}
\label{subsec:null class}

\begin{para}
\label{df:w-closure, w-trivial}
{\bf ($w$-closure, $w$-trivial objects).}\ \ 
Let $\cC$ be a category and let $w$ be a class of morphisms in $\cC$ 
and let $\cS$ be a class of objects in $\cC$. 
We say that an object $x$ is {\it $w$-equivalent to an object $y$} 
if there exists a zig-zag sequence of finite morphisms in $w$ which 
connects $x$ and $y$. 
We denote the class of all objects which are $w$-equivalent to 
some objects in $\cS$ by $\cS_{w,\cC}$ or simply $\cS_w$
and call it the {\it $w$-closure of $\cS$} ({\it in $\cC$}). 
We regard $\cS_{w,\cC}$ as the full subcategory of $\cC$. 

Assume that $\cC$ admits a zero object and all isomorphisms 
between zero objects are in $w$ and $w$ is closed under finite 
compositions. 
Then we say that an object $x$ is {\it $w$-trivial} 
if the canonical morphism $0\to x$ is in $w$. 
By assumption, the definition of $w$-triviality 
does not depend upon the choice of a zero object. 
We denote the class of all $w$-trivial objects in $\cC$ by $\cC^w$. 
We regard $\cC^w$ as the full subcategory of $\cC$. 
\end{para}

\begin{para}
\label{df:heq closure}
{\bf Definition ($C$-homotopy closure).}\ \ 
Let $\cC$ be a complicial exact category and 
let $\cA$ be a non-empty full subcategory of $\cC$ and 
$\heq$ be a class of all $C$-homotopy equivalences in $\cC$. 
Then by the notation in \ref{df:w-closure, w-trivial}, 
let $\cA_{\heq}$ denote the full subcategory of $\cC$ 
consisting of those objects $x$ which are $C$-homotopy 
equivalent to objects in $\cA$. 
We call $\cA_{\heq}$ the 
{\it $C$-homotopy closure of $\cA$} ({\it in $\cC$}). 
Since $C$-homotopy equivalence relation is transitive, 
an object in $\cC$ which is $C$-homotopy equivalent to 
an object in $\cA_{\heq}$ is also in $\cA_{\heq}$. 
\end{para}

\begin{para}
\label{df:null class}
{\bf Definition (Prenull classes, semi-null classes, null classes).}\ \ 
Let $\cC$ be a complicial exact category and 
let $\cN$ be a full subcategory of $\cC$. 
We say that $\cN$ is a {\it null class} of $\cC$ 
if $\cN$ contains all $C$-contractible objects and 
if for any admissible exact sequence $x\rinf y\rdef z$ 
in $\cC$, 
if two of $x$, $y$ and $z$ are in $\cN$, 
then the third one is also in $\cN$. 
We say that $\cN$ is a {\it prenull class} if 
$\cN$ contains a zero object and 
if for any morphism $f\colon x\to y$ in $\cN$, 
the objects $\Cone f$ and $x\oplus y$ are also in $\cN$. 
Since $Tx=\Cone(x\to 0)$, in this case, $Tx$ is in $\cN$ for any 
object $x$ in $\cN$. 
We say that $\cN$ is a {\it semi-null class} if 
$\cN$ is a prenull class a and for any object $x$ in $\cN$, 
$T^{-1}x$ is in $\cN$. 

We say that a null class (resp. prenull class, semi-null class) $\cN$ 
is {\it thick} if 
it is closed under retractions. 
Namely for any object $x$ in $\cC$ if 
there exists an object $y$ in $\cN$ and a pair of 
morphisms $f\colon x\to y$ and $g\colon y\to x$ 
such that $gf=\id_x$, then $x$ is also in $\cN$. 
\end{para}

\begin{para}
\label{lem:semi null characterization}
{\bf Lemma.}\ \ 
{\it
Let $\cC$ be a complicial exact category and let $\cS$ be 
a full subcategory of $\cC$. Then the following conditions are equivalent.
\begin{enumerate}
\enumidef
\item
$\cS$ is closed under the operations $\Cone$ and $\hFib$.  
\item
$\cS$ is closed under the operations $\Cone$ and $T^{-1}$. 
\item
$\cS$ is closed under the operations $\hFib$ and $T$. 
\end{enumerate}
In particular, $\cS$ is a semi-null class of $\cC$ if and only if 
it is closed under finite direct sums and satisfies one of the conditions 
above. 
}
\end{para}

\begin{proof}
Notice that for an object $x$ and a morphism $f\colon x\to y$ in $\cC$, 
$T^{-1}x=\hFib(x\to 0)$, $Tx=\Cone(0\to x)$ and  
$\Cone f\isoto T\hFib f$ and $\hFib f\isoto T^{-1}\Cone  f$ 
by Lemma~\ref{lem:Cone =THfib}.
\end{proof}

\begin{para}
\label{df:comp exact func between prenull cass}
{\bf Definition.}\ \ 
By mimicking the definition of complicial exact functors between 
complicial exact categories, 
we can define the notion to prenull classes. 
That is, let $\cC$ and $\cC'$ be complicial exact categories 
and let $\cN\subset \cC$ and $\cN'\subset \cC'$ be 
prenull classes of $\cC$ and $\cC'$ respectively. 
A {\it complicial exact functor} $(f,c)\colon \cN\to \cN'$ 
from $\cN$ to $\cN'$ 
is a pair $(f,c)$ consisting of a functor $\cN\to \cN'$ 
which is additive and preserves Frobenius admissible exact sequences 
and a natural equivalence $c\colon Cf\isoto fC$ satisfies the equality
\begin{equation}
c\cdot (\iota^{\cN'}\ast f)=f\ast \iota^{\cN}.
\end{equation}
Then we can show the similar statements of 
\ref{lemdf:cF}, \ref{lem:comp func and alpha beta} and 
\ref{lemdf:comp functor induce on C_h^[1]}.
\end{para}

\begin{para}
\label{df:pull-back by functors}
{\bf Definition (Pull-back by functors).}\ \ 
Let $f\colon \cX\to \cY$ be a functor 
between categories $\cX$ and $\cY$ 
and let 
$\cZ$ be a full subcategory of $\cY$ and 
let $w$ be a class of morphisms in $\cY$. 
Then we denote the full subcategory 
of $\cX$ consisting of those objects $x$ such that 
$f(x)$ is isomorphic to some object in $\cZ$ by $f^{-1}\cZ$ and 
we call it the {\it pull-back of $\cZ$ by $f$}. 
We also write $f^{-1}(w)$ for the class of morphisms in $\cX$ 
consisting of those morphisms $u\colon x\to y$ such that 
$f(u)$ is in $w$ and call it the 
{\it pull-back of $w$ by $f$}. 
\end{para}

\begin{para}
\label{lem:pull-back of semi-null classes}
{\bf Lemma (Pull-back of semi-null classes).}\ \ 
{\it
Let $(f,c)\colon\cC\to \cC'$ be a complicial exact functor 
between complicial exact categories $\cC$ and $\cC'$ and let 
$\cN$ be a prenull {\rm(}resp. semi-null, null{\rm)} class of $\cC'$.  
Then $f^{-1}\cN$ is also a prenull 
{\rm(}resp. semi-null, null{\rm)} class of $\cC$.
}
\end{para}

\begin{proof}
Let $u\colon x\to y$ be a morphism in $f^{-1}\cN$. 
Then by definition, there exists a pair of objects $x'$ and $y'$ 
in $\cN$ and 
a pair of isomorphisms 
$a\colon fx\isoto x'$ and $b\colon fy\isoto y'$ in $\cC'$. 
We set $u':=b\cdot fu\cdot a^{-1}\colon x'\to y'$. 
Then by \ref{lemdf:cF} and \ref{lemdf:comp functor induce on C_h^[1]}, 
$f(x\oplus y)$, $f(T^{\pm 1}x)$ 
and $f(\Cone u)$ are isomorphic to 
$x'\oplus y'$, $T^{\pm 1}x'$ and $\Cone u'$ respectively. 
Thus if $\cN$ is a prenull (resp. semi-null) class, then 
$x\oplus y$, $\Cone u$ (and $T^{-1}x$) are in $f^{-1}\cN$ 
and $f^{-1}\cN$ is a prenull (resp. semi-null) class. 

Next assume that $\cN$ is a null class of $\cC'$. 
Then since $f$ sends a $C$-contractible object in $\cC$ 
to a $C'$-contractible object in $\cC'$, $f^{-1}\cN$ 
contains all $C$-contractible objects. 
For an admissible exact sequences $x\rinf y\rdef z$ in $\cC$, 
if two of $f(x)$, $f(y)$ and $f(z)$ are in $\cN$, 
then the third one is also in $\cN$. 
Thus $f^{-1}\cN$ is a null class. 
\end{proof}

\begin{para}
\label{df:direct sum closure}
{\bf Definition (Direct sum closure).}\ \ 
Let $\cA$ be an additive category and let $\cO$ be 
a full subcategory of $\cA$. 
We write ${\langle \cO\rangle}_{\oplus}$ 
for the full subcategory consisting of those objects which are 
isomorphic to finite direct sums of some objects in $\cO$. 
In particular, ${\langle \cO\rangle}_{\oplus}$ 
contains all zero objects in $\cA$ and closed under finite 
direct sum operation and it is the smallest full subcategory which 
contains $\cO$ and closed under finite direct sum operation. 
We call it the {\it direct sum closure of $\cO$}.
\end{para}

\begin{para}
\label{df:T-closure}
{\bf Definition ($T$-closure).}\ \ 
Let $\cC$ be a complicial exact category and let 
$\cA$ be a full subcategory of $\cC$. 
We write ${\langle \cA\rangle}_T$ for the full subcategory of 
$\cC$ consisting of those objects which are isomorphic to 
$T^ny$ for some integers $n$ and some objects $y$ in $\cA$. 
We call ${\langle \cA\rangle}_T$ the {\it $T$-closure of $\cA$} in $\cC$.
\end{para}

\begin{para}
\label{lem:T-closure of add cat}
{\bf Lemma.}\ \ 
{\it
Let $\cC$ be a complicial exact category and 
let $\cA$ be a full subcategory of $\cC$. 
Assume that $\cA$ is closed under finite direct sums and 
the operation $T$. 
Namely for any objects $x$ and $y$ in 
$\cA$, $x\oplus y$ and $Tx$ are in $\cA$. 
Then ${\langle \cA\rangle}_T$ is 
also closed under finite direct sums.
}
\end{para}

\begin{proof}
Let $x$ and $y$ be objects in ${\langle \cA\rangle}_T$. 
Then there exists objects $x'$ and $y'$ in $\cA$ and isomorphisms 
$x\isoto T^mx'$ and $y\isoto T^ny'$ for suitable integers $m$ and $n$. 
We set $k:=\min\{m,n \}$. 
Then $x\oplus y\isoto T^k(T^{m-k}x'\oplus T^{n-k}y')$. 
Since $m-k$, $n-k\geq 0$ and $\cA$ is closed under the operation $T$ 
and finite direct sums, 
$T^{m-k}x'\oplus T^{n-k}y'$ is in $\cA$. 
Thus $x\oplus y$ is in ${\langle \cA\rangle}_T$. 
\end{proof}

\begin{para}
\label{lemdf:prenull closure}
{\bf Lemma-Definition (Prenull closure).}\ \  
Let $\cC$ be a complicial exact category and let 
$\calD$ be a full additive subcategory of $\cC$. 
For a non-negative integer $n$, 
we will inductively define 
$P_n(\calD)$ to be a full subcategory of $\cC$. 
First we set $\Ob P_0(\calD)=\Ob\calD$ and we set 
$\Ob P_{n+1}(\calD):=\{\Cone f \}_{f\in\Mor P_n(\calD)}$. 
For an object $x$ in $P_n(\calD)$, 
since we have an equality 
$x=\Cone (0\to x)$, $P_{n+1}(\calD)$ contains $P_n(\calD)$. 
Finally we set 
$\displaystyle{{\langle\calD \rangle}_{\prenull}:=\bigcup_{n\geq 0}P_n(\calD)}$. 
We claim that ${\langle\calD \rangle}_{\prenull}$ 
is the smallest full subcategory of $\cC$ which contains 
$\calD$ and we call it the {\it prenull closure of $\calD$}. 

For a full subcategory $\cS$ of $\cC$, 
We write ${\langle \cS\rangle}_{\prenull}$ for 
${\langle {\langle  \cS \rangle}_{\oplus} \rangle}_{\prenull} $ and also call it {\it prenull closure of $\cS$}. 
\end{para}

\begin{proof}
First we will show that for each non-negative integer $n$, 
$P_n(\calD)$ is closed under 
finite direct sum operations by induction on $n$. 
For $n=0$, it is clear by assumption. 
Assume that assertion is true for $n=k$ and we will prove for $n=k+1$. 
Let $x$ and $y$ be a pair of objects in $P_{k+1}(\calD)$. 
Then there exists a pair of morphisms 
$f\colon x'\to x''$ and $g\colon y'\to y''$ in $P_k(\calD)$ 
such that $x=\Cone f$ and $y=\Cone g$. 
Then $x\oplus y=\Cone \begin{pmatrix}f& 0\\ 0& g \end{pmatrix}$ 
is in $P_{k+1}(\calD)$. 

Next let $f\colon z\to w$ be a morphism in 
${\langle \calD\rangle}_{\prenull}$. 
Then there exists a non-negative integer $n$ such that 
$f$ is in $P_n(\calD)$. Then $\Cone f$ is in 
$P_{n+1}(\calD)\rinc {\langle \calD\rangle}_{\prenull}$. 
Thus ${\langle \calD\rangle}_{\prenull}$ is a prenull class in $\cC$. 

Finally let $\cP$ be  prenull class in $\cC$ which contains $\calD$. 
Then we can show for non-negative integer $n$, 
$P_n(\calD)$ is in $\cP$ 
by proceeding induction on $n$. 
Thus it turns out that ${\langle \calD\rangle}_{\prenull}$ 
is the smallest prenull class which contains $\calD$. 
\end{proof}

\begin{para}
\label{lemdf:seminull closure}
{\bf Definition (Semi-null closure).}\ \ 
Let $\cP$ be a prenull class of a complicial exact category $\cC$. 
Then we claim that 
${\langle \cP\rangle}_T$ is the smallest semi-null class which contains 
$\cP$. 
Therefore we write ${\langle \cP\rangle}_{\seminull}$ 
for ${\langle \cP\rangle}_T$ and 
we call it the {\it semi-null closure of $\cP$}.

For a full additive subcategory $\cS$ of $\cC$, we simply write 
${\langle \cS \rangle}_{\seminull}$ for 
${\langle{\langle \cS \rangle}_{\prenull}\rangle}_{\seminull} $ 
and also call it the {\it semi-null closure of $\cS$}.
\end{para}

\begin{proof}
Let $f\colon x\to y$ be a morphism in ${\langle \cP\rangle}_{\seminull}$. 
Then there exists a pair of integers $n$ and $m$ and a pair of 
objects $x'$ and $y'$ in $\cP$ and a pair of 
isomorphisms $a\colon T^nx'\isoto x$ and $b\colon T^my'\isoto y$. 
We set $k:=\min\{n,m\}$ and $f':=b^{-1}fa$. 
Then the morphism $T^{-k}f'\colon T^{n-k}x'\to T^{m-k}y'$ 
is in $\cP$ and 
$T^k\Cone f' $ and $T^{k-1}T^{n-k}x'$ are isomorphic to $\Cone f$ and $T^{-1}x$ respectively. 
Thus $\Cone f$ and $T^{-1}x$ are in 
${\langle \cP\rangle}_{T}$. 
Moreover $x\oplus y$ is in ${\langle \cP\rangle}_{T}$ 
by Lemma~\ref{lem:T-closure of add cat}. 

Let $\cS$ be a semi-null class of $\cC$ which contains $\cP$. 
Then for any object $x$ in $\cP$ and any integer $n$, 
$T^n x$ is in $\cS$. 
Thus $\cS$ contains ${\langle \cP\rangle}_{T}$. 
\end{proof}

\begin{para}
\label{lem:null class}
{\bf Lemma.}\ \ 
{\it
Let $\cC$ be a complicial exact category and let 
$\cN$ be a full subcategory of $\cC$. Then\\
$\mathrm{(1)}$ 
If $\cN$ is a null class of $\cC$, then\\
$\mathrm{(i)}$ 
$\cN$ is closed under $C$-homotopy equivalences. 
Namely an object in $\cC$ which is $C$-homotopy equivalent 
to an object in $\cN$ is also in $\cN$.\\
$\mathrm{(ii)}$ 
For any object $x$ in $\cN$, 
the objects $Cx$, $Tx$ and $T^{-1}x$ are also in $\cN$. 
In particular $\cN$ is a complicial exact subcategory of $\cC$.\\
$\mathrm{(2)}$ 
If $\cN$ is a semi-null class of $\cC$, 
then $C$-homotopy closure of $\cN$ in $\cC$ is 
the smallest null class of $\cC_{\frob}$ which contains $\cN$. 
}
\end{para}

\begin{proof}
A proof of $\mathrm{(1)}$ $\mathrm{(i)}$ is given in \cite[5.4]{Moc10}. 

\sn
$\mathrm{(ii)}$ 
For any object $z$ in $\cC$, $Cz$ is $C$-contractible and 
thus it is in $\cN$. 
For any object $x$ in $\cN$, there are Frobenius admissible exact sequences 
$x\rinf Cx \rdef Tx $ and $T^{-1}x \rinf CT^{-1}x \rdef x$ and 
$Cx$, $CT^{-1}x$ are in $\cN$. 
Thus $Tx$ and $T^{-1}x$ are also in $\cN$. 

\sn
$\mathrm{(2)}$ 
First notice that for any object $x$ in $\cN$, 
$Cx=\Cone \id_x$ is in $\cN$ and for any morphism 
$x\onto{f} y$ in $\cN$, $\Cyl f=y\oplus Cx$ is in $\cN$. 

Next we we will show that for any morphism 
$f\colon x\to y$ in $\cN_{\heq}$, 
the objects $\Cone f$, $T^{-1}x$ and $x\oplus y$ are in $\cN_{\heq}$. 
Then by definition, there exists objects $x'$, $y'$ in $\cN$ and 
$C$-homotopy equivalences $x'\onto{a} x$, $x\onto{a'} x'$ and 
$y\onto{b}y'$ and a $C$-homotopy $H\colon aa'\Rightarrow_C \id_x$. 
Then we can check that $(a',b,bfH)\colon [x\onto{f}y] \to [x'\onto{bfa}y']$ 
is a homotopy commutative square. 
Thus the morphism $\Cyl(a',b,bfH)\colon \Cyl f \to \Cyl bfa$ 
is a $C$-homotopy equivalence by 
Lemma~\ref{lem:homotopy commutative square} $\mathrm{(2)}$ $\mathrm{(ii)}$ and 
$\Cyl bfa$ is in $\cN$. 
Hence it turns out that $\Cone f$ is in $\cN_{\heq}$. 
By \cite[2.13 $\mathrm{(5)}$]{Moc10}, 
the morphism 
$\displaystyle{\begin{pmatrix}a' &0\\0& b \end{pmatrix}
\colon x\oplus y\to x'\oplus y'}$ 
is a $C$-homotopy equivalence and $x'\oplus y'$ are in $\cN$. 
Thus $x\oplus y$ is in $\cN_{\heq}$. 
Since $T^{-1}$ preserves $C$-homotopy equivalences 
by \cite[2.12 $\mathrm{(1)}$, 2.52]{Moc10}, 
$T^{-1}x \onto{T^{-1}a'} T^{-1}x'$ 
is a $C$-homotopy equivalence and $T^{-1}x' $ is in $\cN$. 
Thus $T^{-1}x$ is in $\cN_{\heq}$ and 
$\cN_{\heq}$ is a semi-null class in $\cC$. 
Hence by \ref{ex:heq} and 
\ref{prop:properties of weak equivalences} $\mathrm{(4)}$, 
it is a null class in $\cC_{\frob}$. 
Notice that there is no circular reasoning in our argument.

Let $\cN'$ be a null class of $\cC_{\frob}$ which contains $\cN$. 
Then $\cN_{\heq}\rinc {\cN'}_{\heq}=\cN'$. 
Thus $\cN_{\heq}$ is the smallest null class of $\cC_{\frob}$ 
which contains $\cN$. 
\end{proof}

\begin{para}
\label{df:null closure}
{\bf Definition (Null closure).}\ \ 
Let $\cC$ be a complicial exact category. 
Since the class of all null classes of $\cC$ is closed under intersection, 
for any full subcategory $\cO$ of $\cC$, 
there exists the smallest null class of $\cC$ which contains $\cO$. 
We denote it by ${\langle \cO\rangle}_{\nullclass}$ and call it 
the {\it null closure of $\cO$}. 

We write ${\langle \cO\rangle}_{\nullclass,\frob}$ for 
${({\langle \cO\rangle}_{\seminull})}_{\heq}$ and 
by the previous lemma~\ref{lem:null class}, 
it is the smallest null class in $\cC_{\frob}$ which contains $\cO$. 
We call it the {\it Frobenius null closure of $\cO$}. 
\end{para}

Recall the definition of coproduct of complicial exact categories 
from \ref{ex:coproduct of complicial exact categories}. 

\begin{para}
\label{ex:null closure of coproducts}
{\bf Example.}\ \ 
Let $\{\cC_i \}_{i\in\cI}$ be a family of 
complicial exact categories indexed by a set $\cI$ and let 
$\cO_i$ be a full subcategory of $\cC_i$ for all $i\in\cI$. 
Then we have equality
\begin{equation}
\label{eq:null closure of coproducts}
{\left\langle\bigvee_{i\in\cI}\cO_i \right\rangle}_{\nullclass,\frob,\bigvee_{i\in\cI}\cC_i}=\bigvee_{i\in\cI}{\langle \cO_i\rangle}_{\nullclass,\frob,\cC_i}.
\end{equation}
\end{para}

\begin{para}
\label{df:thi closure}
{\bf Definition (Thick closure).}\ \ 
Let $\cC$ be a category and let $\cF$ be 
a full subcategory of $\cC$. 
We write $\cF_{\thi,\cC}$ or simply $\cF_{\thi}$ 
for the full subcategory of $\cC$ consisting of 
those objects $x$ such that $x$ are retraction of some objects of $\cF$. 
We call it {\it thick closure of $\cF$} ({\it in $\cC$}) and 
then $\cF_{\thi}$ is closed under retractions. 
\end{para}

\begin{para}
\label{lem:functoriality of closures}
{\bf Lemma (Functoriality of closures).}\ \ 
{\it
Let $f\colon \cC \to \cC'$ be an additive functor between 
additive categories and let $\cS$ and $\cS'$ be a pair 
of full subcategories of $\cC$ and $\cC'$ respectively 
such that $f(\Ob\cS)\subset \cS'$ and let $w$ be a class of morphisms in 
$\cC$. Then

\sn
$\mathrm{(1)}$ 
We have the inclusions 
$f(\Ob{\langle\cS \rangle}_{\oplus})\subset 
\Ob{\langle\cS' \rangle}_{\oplus}$ and 
$f(\Ob{\cS}_{\thi})\subset \Ob{\cS'}_{\thi}$ and 
$f(\Ob{\cS}_w)\subset \Ob{\cS'}_{f(w)}$. 

\sn
$\mathrm{(2)}$ 
Moreover assume that $\cC$ and $\cC'$ be a complicial exact categories 
and there exists a natural equivalence $c\colon C'f\isoto fC$ such that the pair $(f,c)$ is a complicial exact functor $\cC\to \cC'$. 
Then we have the inclusion $f(\Ob{\langle\cS \rangle}_{\#})\subset \Ob{\langle\cS' \rangle}_{\#}$ where $\#\in \{\prenull,\seminull,\nullclass \}$.
}
\end{para}

\begin{proof}
$\mathrm{(2)}$ 
Let $\cN$ be a prenull (resp. semi-null, null) class in $\cC'$ 
which contains $\cS'$. 
Then for $\#=\prenull$ (resp. $\seminull$, $\nullclass$), 
we have the inclusions $\Ob\cS\subset \Ob{\langle \cS\rangle}_{\#} \subset \Ob f^{-1}\cN$ and we obtain the inclusion $f(\Ob{\langle \cS\rangle}_{\#})\subset \Ob\cN$. 
Thus we get the desired inclusion 
$f(\Ob{\langle \cS\rangle}_{\#})\subset \bigcap_{\cS'\subset \cN}\Ob\cN=\Ob{\langle\cS' \rangle}_{\#}$ where intersection taking all prenull (resp. semi-null, null) class which contains $\cS'$. 

\sn
A proof of $\mathrm{(1)}$ is similar. 
\end{proof}

\begin{para}
\label{df:Homotopy retraction}
{\bf Definition (Homotopy retraction).}\ \ 
Let $\cC$ be a complicial exact category and let $x$ and $y$ 
be a pair of objects in $\cC$. 
We say that $x$ is a ($C$-){\it homotopy retraction of $y$} 
if there exists a pair morphisms $i\colon x\to y$ and $p\colon y\to x$ 
and $C$-homotopy $H\colon \id_x\Rightarrow pi$. 
We say that a full subcategory $\cF$ in $\cC$ is 
{\it closed under $C$-homotopy retractions} if 
for a pair of objects $x$ and $y$ in $\cC$ such that $y$ is in $\cF$ and 
$x$ is a $C$-homotopy retraction of $y$, 
then $x$ is also in $\cF$. 
\end{para}

\begin{para}
\label{prop:fundamental properties of thick closure}
{\bf Proposition.}\ \ 
{\it
Let $\cC$ be an additive category and 
let $\cF$ be an additive full subcategory. Then

\sn
$\mathrm{(1)}$ 
We have an equalities
\begin{equation}
\label{eq:thick closure in idempotent completion}
{({\widehat{\cF}})}_{\thi,\widehat{\cC}}=
{({\widehat{\cF}})}_{\Isom,\widehat{\cC}},
\end{equation}
\begin{equation}
\label{eq:thick clousre in terms of idempotent completions}
{\cF}_{\thi,\cC}=
{({\widehat{\cF}})}_{\Isom,\widehat{\cC}}\cap\cC,
\end{equation}
where we denote the idempotent completion functor by $\widehat{(-)}$.

\sn
Moreover assume $\cC$ is a complicial exact category. Then

\sn
$\mathrm{(2)}$ 
If $\cF$ is a prenull {\rm(}resp. semi-null{\rm)} 
class of $\cC$, then 
${({\widehat{\cF}})}_{\Isom,\widehat{\cC}}$ 
is also a prenull 
{\rm(}resp. semi-null{\rm)} class of 
$\widehat{\cC}$. 

\sn
$\mathrm{(3)}$ 
If $\cF$ is closed under $C$-contractible objects and finite 
direct sum and retraction, then $\cF$ is closed under 
$C$-homotopy retractions. 
In particular $\cF$ is closed under $C$-homotopy equivalence. 

\sn
$\mathrm{(4)}$ 
If $\cF$ is closed under finite direct sum and contains 
all objects of the form $Cx$ for some object $x$ in $\cC$, 
then for an object $y$ in $\cF_{\thi}$, there exists an object $y'$ 
in $\cC$ such that $y\oplus y'$ is in $\cF$. 

\sn
$\mathrm{(5)}$ 
If $\cF$ is a prenull (resp. semi-null, null) class of $\cC$, 
then $\cF_{\thi}$ is also. 
}
\end{para}

\begin{proof}
$\mathrm{(1)}$ 
First we will show the equality 
$\mathrm{(\ref{eq:thick clousre in terms of idempotent completions})}$. 
What we need to prove are the following assertions:

\sn
$\mathrm{(a)}$ 
${(\widehat{\cF})}_{\Isom,\widehat{\cC}}\cap\cC$ 
is closed under retractions and contains $\cF$. 

\sn
$\mathrm{(b)}$ 
For a full subcategory $\cM$ of $\cC$, 
if $\cM$ contains $\cF$ and is closed under retractions, 
then $\cM$ contains ${(\widehat{\cF})}_{\Isom,\widehat{\cC}}\cap\cC$.

\begin{proof}[Proof of $\mathrm{(a)}$]
We regard $\cF$ as a full subcategory of $\widehat{\cF}$ 
by the canonical functor $\cF\to \widehat{\cF}$, $x\mapsto (x,\id_x)$. 
Thus $\cF\subset {(\widehat{\cF})}_{\Isom,\widehat{\cC}}\cap\cC$. 

Next 
let $x\onto{i}y\onto{p}x$ be a pair of morphisms such that $pi=\id_x$ 
and $y=(y,\id_y)$ is in ${(\widehat{\cF})}_{\Isom,\widehat{\cC}}$. 
Then by definition there exists an object $(z,e)$ in $\widehat{\cF}$ 
and a pair of morphisms $y\onto{a}z\onto{b}y$ 
in $\cC$ such that $ba=\id_y$ and $ab=e$. 
We set $e':=aipb$. Then $e'^2=aipbaipb=aipipb=aipb=e'$. 
Notice that $ai\colon (x,\id_x)\to (z,e')$ and 
$pb\colon (z,e')\to (x,\id_x)$ give isomorphisms 
between $(z,e')$ and $x$. 
Since $z$ is in $\cF$, $(z,e')$ is in $\widehat{\cF}$ and 
$x=(x,\id_x)$ 
is in ${(\widehat{\cF})}_{\Isom,\widehat{\cC}}$. 

\end{proof}

\begin{proof}[Proof of $\mathrm{(b)}$]
Let $(x,\id_x)$ be an object 
${(\widehat{\cF})}_{\Isom,\widehat{\cC}}\cap\cC$. 
By definition 
there exists an object $y$ in $\cF$ and a pair of morphisms 
$x\onto{a}y\onto x$ in $\cC$ with $ba=\id_x$. 
Since $y$ is in $\cF\subset \cM$ and $\cM$ is closed under 
retractions, $x$ is in $\cM$. 
\end{proof}

Next applying the equality 
$\mathrm{(\ref{eq:thick clousre in terms of idempotent completions})}$ 
to $\widehat{\cC}$ and $\widehat{\cF}$, we obtain the equalities 
${(\widehat{\cF})}_{\thi,\widehat{\cC}}= 
{(\widehat{\widehat{\cF}})}_{\Isom,\widehat{\widehat{\cC}}}\cap \widehat{\cC}= 
{(\widehat{\cF})}_{\Isom,\widehat{\cC}}$ 
where we regard $\widehat{\widehat{\cC}}=\widehat{\cC}$ and 
$\widehat{\widehat{\cF}}=\widehat{\cF}$. 

\sn
$\mathrm{(2)}$ 
Let $f\colon (x,e)\to (y,e')$ be a morphism in 
${(\widehat{\cF})}_{\Isom,\widehat{\cC}}$. 
By assumption there exists a pair of objects $x'$ and $y'$ and 
morphisms $a\colon x\to x'$, $b\colon x'\to x$, $c\colon y\to y'$ 
and $d\colon y'\to y$ such that 
$ae=aba=a$, $eb=bab=b$, $ba=e$, $ce'=cdc=c$, $e'd=dcd=d$, $dc=e'$. 
Then we can show that 
$(x'\oplus y',\begin{pmatrix}ab &0\\0& cd\end{pmatrix})$, 
$(T^{\pm 1}x',T^{\pm 1}ab)$ and $(\Cone (cfb),\Cone(ab,cd))$ 
are isomorphic to 
$(x\oplus y,\begin{pmatrix}e &0\\0& e'\end{pmatrix}) $, 
$(T^{\pm 1}x,T^{\pm 1}e)$ and 
$(\Cone f,\Cone(e,e'))$ respectively and 
it turns out that if $\cF$ is a prenull (resp. semi-null) class in $\cC$, 
then ${(\widehat{\cF})}_{\Isom,\widehat{\cC}}$ is also 
a prenull (resp. semi-null) class in $\widehat{\cC}$. 

\sn
$\mathrm{(3)}$ 
Let $x$ be a $C$-homotopy retraction of an object $y$ in $\cF$. 
Namely there is a triple of morphisms $i\colon x\to y$, $p\colon y\to x$ 
and $H\colon Cx\to y$ in $\cC$ 
such that $H\iota_x=\id_x-pi$. 
Then since we have the equality 
$\begin{pmatrix}H &p\end{pmatrix}
\begin{pmatrix}\iota\\ i\end{pmatrix}=\id_x$, 
$x$ is a retraction of $Cx\oplus y$. 
By assumptions, $Cx\oplus y$ is in $\cF$ and by assumption again, 
$x$ is also in $\cF$. 

\sn
$\mathrm{(4)}$ 
Let $x$ be an object in $\cF_{\thi}$. Then there is 
a retraction $x\onto{i}y\onto{p}x$ with $y\in \Ob\cF$ and $pi=\id_x$. 
Then by assumption $y\oplus Cx$ is in $\cF$ and it is isomorphic to 
$x\oplus \Cone i$ by the commutative diagram below. 
$$\xymatrix{
x\ar@{>->}[r]^{\tiny{\begin{pmatrix}i\\ \iota_x\end{pmatrix}}} \ar@{=}[d] &
y\oplus Cx \ar@{->>}[r]^{\tiny{\begin{pmatrix}\kappa_i & -\xi_i\end{pmatrix}}} \ar[d]_{\wr}^{\tiny{\begin{pmatrix}p & 0\\ \kappa_i & -\xi_i\end{pmatrix}}} &\Cone i \ar@{=}[d]\\
x\ar@{>->}[r]_{\tiny{\begin{pmatrix}\id_x\\ 0 \end{pmatrix}}} & 
x\oplus \Cone i \ar@{->>}[r]_{\tiny{\begin{pmatrix}0 & \id_{\Cone i}\end{pmatrix}}} & \Cone i.
}$$
Thus $x\oplus \Cone i$ is also in $\cF$. 

\sn
$\mathrm{(5)}$ 
Assume $\cF$ is a null class in $\cC$. 
We will show that $\cF_{\thi}$ is also a null class of $\cC$. 
What we need to prove are the following assertions:

\sn
$\mathrm{(a)}$ 
$\cF_{\thi}$ contains all $C$-contractible objects. 

\sn
$\mathrm{(b)}$ 
Let $x\overset{i}{\rinf} y \overset{p}{\rdef} z$ 
be an admissible exact sequence in $\cC$. 
If two of $x$, $y$ and $z$ are in $\cF_{\thi}$, 
then third one is also in $\cF_{\thi}$. 

\begin{proof}[Proof of $\mathrm{(a)}$]
Since $\cF$ contains all $C$-contractible objects and $\cF_{\thi}$ 
contains $\cF$, 
$\cF_{\thi}$ also contains all $C$-contractible objects. 
\end{proof}

\begin{proof}[Proof of $\mathrm{(b)}$]
Since $\cF$ is closed under finite direct sums and 
contains all $C$-contractible objects, 
we can use assertion $\mathrm{(4)}$. 
Assume there exists a pair of objects $x'$ and $z'$ 
(resp. $x'$ and $y'$, $y'$ and $z'$) such that 
both $x\oplus x'$ and $z\oplus z'$ (resp. $x\oplus x'$ and $y\oplus y'$, 
$y\oplus y'$ and $z\oplus z'$) are in $\cF$. 
Then by considering an admissible exact sequence 
$x\oplus x' \overset{\tiny{\begin{pmatrix}i\\ \id_{x'}\\ 0\end{pmatrix}}}{\rinf} y\oplus x'\oplus z' 
\overset{\tiny{\begin{pmatrix}p & 0 & \id_{z'}\end{pmatrix}}}{\rdef} 
z\oplus z'$ 
(resp. 
$x\oplus x' \overset{\tiny{\begin{pmatrix}i & 0\\ 0 & 0\\ 0 & \id_{x'}\\ 0 & 0\end{pmatrix}}}{\rinf} y\oplus  y'\oplus x'\oplus x 
\overset{\tiny{\begin{pmatrix}p & 0 & 0 & 0\\ 0 & \id_{y'} & 0 & 0\\ 0 & 0 & 0 & \id_x\end{pmatrix}}}{\rdef} 
z\oplus y'\oplus x$, 
$x\oplus y'\oplus z 
\overset{\tiny{\begin{pmatrix}i & 0 & 0\\ 0 & \id_{y'} & 0\\ 0 & 0 &\id_z\\ 
0 & 0 & 0\end{pmatrix}}}{\rinf} y\oplus y'\oplus  z\oplus z' 
\overset{\tiny{\begin{pmatrix}p & 0 & 0 & 0\\ 0 & 0 & 0 & \id_{z'}\end{pmatrix}}}{\rdef} 
z\oplus z' $), 
it turns out that 
$y\oplus x'\oplus z'$ 
(resp. $z\oplus y'\oplus x$, $x\oplus y'\oplus z$) is in $\cF$. 
Thus $y$ (resp. $z$, $x$) 
is in $\cF_{\thi}$. 

\end{proof}

Next let us assume that $\cF$ is a prenull (resp. semi-null) class in 
$\cC$. 
Then by $\mathrm{(2)}$, ${(\widehat{\cF})}_{\Isom,\widehat{\cC}}$ 
is also prenull (resp. semi-null) in $\widehat{\cC}$. 
Then the pull-back of ${(\widehat{\cF})}_{\Isom,\widehat{\cC}}$ 
by the strictly normal complicial exact functor 
$\cC\to\widehat{\cC}$, $x\mapsto (x,e)$ is also prenull (resp. semi-null) 
in $\cC$ 
by \ref{lem:pull-back of semi-null classes}. 
Finally by the equality 
$\mathrm{(\ref{eq:thick clousre in terms of idempotent completions})}$, 
it is just ${\cF}_{\thi}$ a thick closure of $\cF$. 
\end{proof}

\subsection{Relative complicial exact categories}
\label{subsec:rel compl exact cat}

In this subsection we study a special class of 
Waldhausen categories which we call relative complicial exact categories. 
First we recalling the notion of relative categories 
in the sense of Barwick and Kan \cite{BK12} and 
we recall the notion of relative exact categories 
which are relative categories 
whose underlying categories are Quillen exact categories. 
Next we recall the notion of relative complicial exact categories 
from \cite{Sch11}, \cite{Moc10} and \cite{Moc13b} 
(in the references, we call them complicial exact categories with weak equivalences or bicomplicial pairs).

\begin{para}
\label{df:rel cat}
{\bf (Relative categories).}\ \ 
A notion of {\it relative categories} 
introduced by Barwick and Kan in \cite{BK12} 
is a model of homotopy theory of homotopy theories. 
We briefly review the these notions from \cite{BK12}. 
A {\it relative category} $\bX$ is a pair $(\cX,v)$ consisting of 
a category $\cX$ and a class of morphisms $v$ in $\cX$ 
which is closed under finite compositions. 
Namely for any object $x$ in $\cX$, 
the identity morphism $\id_x$ is in $v$ and 
for any pair of composable morphisms $x\onto{f}y\onto{g}z$ in $v$, 
$gf$ is also in $v$. 
Thus we can regard $v$ as the subcategory of $\cX$. 
We call $\cX$ the {\it underlying category} of $\bX$ and denote it 
by $\cC_{\bX}$ and call $v$ the {\it class of weak equivalences} of $\bX$ 
and write $w_{\bX}$ for $v$. 
We say that a relative category $\bX$ is {\it small} 
if the underlying category $\cC_{\bX}$ of $\bX$ is small.

A {\it relative functor} between relative categories 
is a functor which preserves weak equivalences. 
Let $f$, $g\colon \bX\to\bY$ be relative functors 
between relative categories. 
A {\it natural weak equivalences} $\Theta$ from $f$ to $g$ is a natural transformation $\Theta\colon f\to g$ such that for any object $x$ in $\cC_{\bX}$, 
$\Theta_x$ is in $w_{\bY}$. 
We say that {\it the pair $f$ and $g$ are categorical homotopic} 
or {\it $f$ is categorical homotopic to $g$} 
if there exists a zig-zag sequence of natural weak equivalences 
which connects $f$ and $g$. 

Let $f\colon \bX\to \bY$ be a relative functor between 
relative categories $\bX$ and $\bY$. 
We say that $f$ is a {\it categorical homotopy equivalence} 
if there exists a relative functor $g\colon \bY\to \bX$ such that 
$fg$ and $gf$ are categorical homotopic to $\id_{\cC_{\bY}}$ and 
$\id_{\cC_{\bX}}$ respectively. 
We denote the category of small relative categories and relative functors 
by $\RelCat$. 
Let $\bX$ and $\bY$ be small relative categories. 
We define $\bX\times \bY$ to be a relative category by setting 
$\cC_{\bX\times\bY}:=\cC_{\bX}\times\cC_{\bY}$ and 
$w_{\bX\times \bY}:=w_{\bX}\times w_{\bY}$ and call it 
the {\it product of $\bX$ and $\bY$}. 
The relative category 
$\bX\times\bY$ is a categorical product of $\bX$ and $\bY$ in $\RelCat$. 
Next we write  $\bX^{\bY}$ for the relative category whose 
underlying category is the category of relative functors 
from $\bY$ to $\bX$ and natural transformations and whose 
weak equivalences are natural weak equivalences. 
Then $\bX^{\bY}$ is the internal hom object in $\RelCat$, 
namely for any object $\bZ$ in $\RelCat$, 
the functor $-\times \bZ\colon \RelCat\to\RelCat$ admits a right adjoint 
${(-)}^{\bZ}$. 
Thus $\RelCat$ is Cartesian closed. 

Let $\cS$ be a full subcategory of $\cC_{\bX}$ the underlying 
category of a relative category $\bX$. 
Then the pair $(\cS,w_{\bX}\cap\Mor\cS)$ 
is a relative category and call it the 
{\it restriction of $\bX$ to $\cS$} and denote it by $\bX|_{\cS}$. 
We often write $w|_{\cS}$ for $w_{\bX|_{\cS}}$. 
\end{para}

\begin{para}
\label{df:level weak equivalence}
{\bf Definition (Level weak equivalences).}\ \ 
Let $\bC=(\cC,w)$ be a relative category and let 
$\cI$ be a small category. 
We say that a morphism 
$f\colon x\to y$ in $\cC^{\cI}$ the 
category of $\cI$-diagrams 
is a {\it level weak equivalence} if 
for any object $i$, $f(i)\colon x(i)\to y(i)$ is in $w$. 
We denote the class of all level weak equivalences in $\cC^{\cI}$ 
by ${lw}_{\cC^{\cI}}$ or shortly $lw$. 

Assume if $\cC$ admits a zero object. 
Then we can consider $\Ch(\cC)$ the category of chain complexes on $\cC$ 
and $F_b\cC$ the category of bounded filtered objects on $\cC$. 
We can regard both of them as the full subcategory of $\cC^{\bbZ}$ 
the category of $\bbZ$-diagrams and we denote the 
restriction of ${lw}_{\cC^{\bbZ}}$ to 
$\Ch(\cC)$ and $F_b\cC$ by ${lw}_{\Ch(\cC)}$ and 
${lw}_{F_b\cC}$ or simply $lw$. 
We also denote the restriction of $lw$ to 
$\Ch_{\#}(\cC)$ by ${lw}_{\#}$ or shortly $lw$ 
for $\#\in\{\pm,b \}$. 

Moreover assume $(\cC,w)$ be a relative complicial exact category. 
Then we similarly say that a homotopy commutative diagram 
$(f,H)\colon x\to y$ in $F_{b,h}\cC$ is a 
{\it level weak equivalence} if 
for any integer $n$, $f_n\colon x_n\to y_n$ 
is in $w$ and we write ${lw}_{F_{b,h}\cC}$ or simply $lw$ 
for the class of all level weak equivalences in $F_{b,h}\cC$. 
\end{para}

\begin{para}
\label{df:rel exact cat}
{\bf Definition (Relative exact categories).}\ \ 
A {\it relative exact category} $\bE=(\cE,w)$ is a 
relative category whose underlying category is an exact category. 
A {\it relative exact functor} between 
relative exact categories $f\colon\bE=(\cE,w)\to\bF=(\cF,v)$ 
is a relative and exact functor. 
We denote the category of small relative exact categories and 
relative exact functors by $\RelEx$. 
We write $\underline{\RelEx}$ for the $2$-category of small relative 
exact categories and relative exact functors and 
natural weak equivalences. 
\end{para}

\begin{para}
\label{para:colim of rel ex cat}
{\bf (Colimit of relative exact categories).}\ \ 
Let $\cI$ be a small category and let 
$\bE\colon\cI\to\RelEx$ be an $\cI$-diagram of relative exact categories. 
Then $\displaystyle{\underset{\cI}{\colim} \bE=\underset{i\in\Ob\cI}{\colim}\bE_i}$ 
is defined as follows. 
The underlying category $\cC_{\underset{\cI}{\colim}\bE}$ of 
$\displaystyle{\underset{\cI}{\colim}\bE}$ is the colimit of 
the $\cI$-diagram 
$\displaystyle{\cI\onto{\bE}\RelEx\onto{\text{forget}}\ExCat}$ 
of small exact categories and the class of weak equivalences 
$w_{\underset{\cI}{\colim}\bE}$ of $\displaystyle{\underset{\cI}{\colim}\bE}$ 
is the subset of the class of morphisms in $\cC_{\underset{\cI}{\colim}\bE}$ 
consisting of those morphisms which are represented by 
morphisms in $w_i$ for some $i$ in $\Ob\cI$. 
\end{para}

\begin{para}
\label{ex:exact cat as rel exact cat}
{\bf Example (Exact categories).}\ \ 
Let $\cE$ be an exact category. 
We denote the class of all isomorphisms in $\cE$ by $i_{\cE}$ 
or simply $i$. Then the pair $(\cE,i_{\cE})$ is 
a relative category. 
In this way, we regard exact categories as relative exact categories. 
There is a functor from $\ExCat$ the category of exact categories to 
$\RelEx$ 
the category of relative exact categories which sends an exact category 
$\cE$ to the relative exact category $(\cE,i_{\cE})$. 
We often simply denote the relative category $(\cE,i_{\cE})$ by $\cE$. 
\end{para}

\begin{para}
\label{df:quasi isom}
{\bf Example (Quasi-isomorphisms).}\ \ 
We recall the notion of quasi-isomorphisms on $\Ch(\cE)$ 
the category of chain complexes on an exact category $\cE$. 
We say that a chain complex $x$ is {\it strictly acyclic} if 
for any integer $n$, the boundary morphism $d^x_n\colon x_n\to x_{n-1}$ 
factors as 
$$\xymatrix{
x_n \ar[rr]^{d^x_n} \ar@{->>}[rd]_{p_n} & & x_{n-1}\\
& z_{n} \ar@{>->}[ru]_{i_n} & 
}$$
where $p_n$ is a cokernel for $d_{n+1}^x$ and an admissible epimorphism 
and $i_n$ is a kernel for $d_{n-1}^x$ and an admissible monomorphism. 
We say that a chain complex $x$ is {\it acyclic} if it it chain homotopy 
equivalent to a strictly acyclic complex. 
We can show that if $\cE$ is idempotent complete, then 
acyclic complexes on $\cE$ are strictly acyclic complexes 
(see \cite[1.8]{Nee90}, \cite[11.2]{Kel96}).

We say that a chain morphism $f\colon x\to y$ is a {\it quasi-isomorphism} 
if $\Cone f$ is acyclic. 
We can show that if $\cE$ is essentially small and idempotent complete, then 
a morphism $f\colon x\to y$ in $\Ch(\cE)$ is a quasi-isomorphism 
if and only if it induces isomorphisms on $H_k(x)\isoto H_k(y)$ 
on homologies for any integer $k$ where homologies are 
taken in $\Lex\cE$ the category of left exact functor from 
$\cE^{\op}$ to the category of abelian groups via Yoneda embedding
$\cE\rinc \Lex\cE$ 
(see \cite[1.11.8]{TT90}).

We denote the class of all quasi-isomorphisms in $\Ch(\cE)$ 
by $\qis_{\Ch(\cE)}$ or simply $\qis$. 
Similarly on $\Ch_b(\cE)$, $\Ch_{+}(\cE)$ and $\Ch_{-}(\cE)$ 
the categories of bounded chain complexes, bounded below chain complexes 
and bounded above chain complexes on $\cE$ respectively, 
we define the class $\qis$ of quasi-isomorphisms. 
The pair $(\Ch_{\#}(\cE),\qis)$ ($\#\in\{b,\pm,\text{nothing}\}$) 
are typical examples of relative 
exact categories.

\end{para}

\begin{para}
\label{df:comp weak equiv}
{\bf Definition (Complicial weak equivalences).}\ \ 
Let $\cC$ be a complicial exact category. 
We say that a class of morphisms $w$ in $\cC$ is 
a {\it class of complicial weak equivalences} 
if it contains all $C$-homotopy equivalences and 
satisfies the extension and 
the saturation axioms in \cite[p.327]{Wal85}. 
In this case, we call the class of morphisms in $w$ 
a {\it complicial weak equivalences}. 

We say that a class of complicial weak equivalences on $\cC$ is {\it thick} 
if it closed under retractions in the morphisms category $\cC^{[1]}$. 

A {\it relative complicial exact category} $\bC=(\cC,w)$ 
is a relative category whose underlying category is a complicial 
exact category and whose class of morphisms is a class of 
complicial weak equivalences. 
A {\it relative complicial exact functor} between 
relative complicial categories $\bC=(\cC,w)\to \bC'=(\cC',w')$ 
is a complicial exact functor 
$(f,c)\colon\cC\to \cC'$ such that $f\colon\cC\to \cC'$ 
is a relative functor. 

A complicial natural weak equivalence between relative 
complicial exact functors $(f,c)$, 
$(f',c')\colon\bC=(\cC,w)\to\bC'=(\cC',w')$ 
is a complicial natural transformation $\theta\colon (f,c)\to(f',c')$ 
such that for any object $x$ in $\cC$, $\theta(x)$ is in $w'$. 
We denote the category of small relative complicial exact categories and 
relative complicial exact functors (resp. 
small thick normal ordinary relative complicial exact categories and 
strictly normal relative complicial exact functors) 
by $\RelComp$ (resp. $\RelComp_{\stnor}$) and 
we write $\underline{\RelComp}$ (resp. $\underline{\RelComp_{\stnor}}$) 
for the $2$-category of 
small relative complicial exact categories, relative complicial exact 
functors (resp. 
small thick normal ordinary relative complicial exact categories, 
strictly normal relative complicial exact functors) 
and complicial natural weak equivalences. 
\end{para}

\begin{para}
\label{lem:comp weak eq is Wald weak eq}
{\bf Lemma.}\ 
(\cf \cite[5.18]{Moc10}.)\ \ 
{\it
Let $(\cC,w)$ be a relative complicial exact category. 
Then the pair $(\cC,w)$ 
satisfies the functorial factorization axiom, the functorial 
cofactorization axiom, the gluing axiom and the cogluing axiom. 
In particular $(\cC,w)$ and $(\cC^{\op},w^{\op})$ are Waldhausen categories. 
Moreover $w$ is proper. 
Namely $w$ is stable under co-base change along any admissible monomorphisms 
and stable under base change along any admissible epimorphisms.
}
\qed
\end{para}

\begin{para}
\label{cor:homotopy proper}
{\bf Corollary.}\ \ 
{\it
Let $(\cC,w)$ be a relative complicial exact category. 
Then $w$ is homotopy proper in the following sense:

\sn 
$w$ is stable under canonical homotopy push-out along any morphisms and 
stable under canonical homotopy pull-back along any morphisms. 
}
\end{para}

\begin{proof}
Let $x\overset{f}{\leftarrow} y \onto{g}z$ be a pair of morphisms in $\cC$ 
such that $f$ is in $w$. 
Then in the homotopy push-out diagram 
$\mathrm{(\ref{eq:homotopy pushout diagram})}$ 
the morphism $i_f$ is just a composition of 
the $C$-homotopy equivalence 
$\displaystyle{\footnotesize{\begin{pmatrix}0\\ \id_z\end{pmatrix}}\colon z \to C(y)\oplus z}$ 
and $\displaystyle{f'\colon C(y)\oplus z \to }$ which is 
the co-base change of $f$ by the admissible monomorphism 
$\displaystyle{\footnotesize{\begin{pmatrix}\iota_y\\ -g\end{pmatrix}}\colon y\to Cy\oplus z}$. 
Thus by \ref{lem:comp weak eq is Wald weak eq}, 
$i_f$ is in $w$. 
A proof of the homotopy pull-back case is similar. 
\end{proof}

\begin{para}
\label{cor:zig-zag}
{\bf Corollary.}\ \ 
{\it
Let $(\cC,w)$ be a relative complicial exact category and 
let $x$ and $z$ be a pair of objects in $\cC$. 
Then the following two conditions are equivalent.

\sn
$\mathrm{(1)}$ 
There exists a pair of morphisms 
$x\overset{f}{\leftarrow}y\onto{g}z$ in $w$. 

\sn
$\mathrm{(2)}$ 
There exists a pair of morphisms 
$x\onto{f}y\overset{g}{\leftarrow}z$ in $w$. 
}
\end{para}

\begin{proof}
We assume the condition $\mathrm{(1)}$. 
Then in the canonical homotopy push-out diagram 
$\mathrm{(\ref{eq:homotopy pushout diagram})}$, 
$i_f$ is in $w$ by \ref{cor:homotopy proper} and it turns out that 
$i_g$ is also in $w$ by $2$ out of $3$ property of $w$. 
Hence we obtain the condition $\mathrm{(2)}$. 
A proof of the converse implication is similar. 
\end{proof}

\begin{para}
\label{para:derivable Waldhausen category}
{\bf (Derivable Waldhausen category).}\ \ 
In \cite[1.1]{Cis10b}, we say that a Waldhausen category 
$\bC=(\cC,w_{\bC})$ is {\it derivable} if 
it satisfies the saturation axiom in \cite[p.327]{Wal85} and 
(non-functorial) factorization axiom in \cite[1.1]{Cis10b}. 
Thus a relative complicial exact category is a derivable 
Waldhausen category by Lemma~\ref{lem:comp weak eq is Wald weak eq}. 
\end{para}

\begin{para}
\label{ex:generator of comp weak eq}
{\bf Example.}\ \ 
Let $\cC$ be a complicial exact category. 
For a family of classes of complicial (resp. thick complicial) 
weak equivalences $\{w_{\lambda}\}_{\lambda\in\Lambda}$ in $\cC$, 
$\displaystyle{\bigcap_{\lambda\in\Lambda}w_{\lambda}}$ 
is also a class of complicial (resp. thick complicial) weak 
equivalences. 
For any class of morphisms $v$ in $\cC$, 
we write $\langle v\rangle $ (resp. ${\langle v\rangle}_{\operatorname{thick}}$) for $\bigcap w$ where $w$ runs through all 
classes of complicial (resp. thick complicial) weak equivalences which 
contains $v$. 
We call $\langle v\rangle$ (resp ${\langle v\rangle}_{\operatorname{thick}}$) 
the {\it class of complicial} 
(resp. {\it thick complicial}) {\it weak equivalences spanned by $v$}.
\end{para}

\begin{para}
\label{ex:heq}
{\bf Example (Class of $C$-homotopy equivalences).}\ \ 
Let $\cC$ be a complicial exact category. 
Recall we denote the class of all $C$-homotopy equivalences in $\cC$ 
by $\heq_{\cC}$ or simply $\heq$. 
The class $\heq$ 
is the class of complicial weak equivalences in $\cC_{\frob}$ 
(see \cite[3.15 $\mathrm{(4)}$]{Moc13b}).
\end{para}

\begin{para}
\label{ex:quasi weak equivalences}
{\bf Example (Quasi-weak equivalences).}\ \ 
Let $\bE=(\cE,w)$ be a relative exact category. 
We denote the class of thick complicial weak equivalences in $\Ch_{\#}(\cE)$ 
spanned by $\qis$ and $lw$ by $qw_{\Ch_{\#}(\cE)}$ or simply $qw$ and 
we call a morphism in $qw$ a {\it quasi-weak equivalence}. 
In particular if $w=i_{\cE}$ the class of all isomorphisms in $\cE$, 
then the class $qw$ is 
just the class of all quasi-isomorphisms in $\Ch_{\#}(\cE)$. 
We denote the pair $(\Ch_{\#}(\cE),qw)$ by $\Ch_{\#}(\bE)$. 
Then for any relative exact functor $f\colon\bE=(\cE,w)\to \bF=(\cF,v)$, 
there is a strictly normal complicial exact functor 
$\Ch_{\#}(f)\colon\Ch_{\#}(\bE)\to\Ch_{\#}(\bF)$. 
The association 
$\Ch_{\#}\colon\underline{\RelEx}\to\underline{\RelComp_{\stnor}}$ 
is a $2$-functor. 
\end{para}

\begin{para}
\label{nt:Verdier correspondence}
{\bf (Correspondence between null classes and complicial weak equivalences).}\ 
(\cf \cite[3.5]{Moc13b}.)\ \ 
Let $\cC$ be a complicial small exact category. 
We denote the class of all null classes 
(resp. thick null classes)
in $\cC$ by 
$\NC(\cC)$ (resp. $\NC_{\thi}(\cC)$ and 
we write $\CW(\cC)$ (resp. $\CW_{\thi}(\cC)$) for 
the class of all classes of complicial weak equivalences 
(resp. thick complicial weak equivalences) in $\cC$.

For a null class in $\cC$, we define $w_{\cN}$ 
to be a class of those morphisms $f\colon x\to y$ in $\cC$ 
such that $\Cone f$ is in $\cN$. 
Then we can show that $w_{\cN}$ is a class of complicial weak equivalences 
in $\cC$ and moreover if $\cN$ is thick, then $w_{\cN}$ is also thick. 
For a class of complicial weak equivalences $u$ in $\cC$, 
$\cC^u$ the class of $u$-trivial objects in $\cC$ 
(See \ref{df:w-closure, w-trivial}) 
is a null class and if $u$ is thick, 
then $\cC^u$ is also thick. 

We can also show the equalities
\begin{equation}
\label{eq:Verdier correspondence 1}
\cC^{w_{\cN}}=\cN,\ \text{and}
\end{equation}
\begin{equation}
\label{eq:Verdier correspondence 2}
w_{\cC^u}=u.
\end{equation}
Thus the associations 
$\NC(\cC)\to \CW(\cC)$, $\cN\mapsto w_{\cN}$ and 
$\CW(\cC)\to \NC(\cC)$, $u\mapsto \cC^u$ give 
order preserving bijections between $\NC(\cC)$ and $\CW(\cC)$, 
and $\NC_{\thi}(\cC)$ and $\CW_{\thi}(\cC)$. 
\end{para}

\begin{para}
\label{cor:Cones invariance by co-base change}
{\bf Corollary.}\ \ 
{\it Let $\bC=(\cC,w)$ be a relative complicial exact category and 
let $i\colon x\to y$ be an admissible monomorphism 
{\rm(}resp. a Frobenius admissible monomorphism{\rm)} 
in $\cC$. 
Then

\sn
$\mathrm{(1)}$ 
The canonical morphism $\Cone i\to y/x$ is in $w$ 
{\rm(}resp. a $C$-homotopy equivalence{\rm)}. 

\sn
$\mathrm{(2)}$ 
Assume that $y/x$ is in $\cC^w$ 
{\rm(}resp. $C$-contractible{\rm)}, 
then $i$ is in $w$ 
{\rm(}resp. a $C$-homotopy equivalence{\rm)}.

\sn
$\mathrm{(3)}$ 
Let $f\colon x\to x'$ be a morphism and let 
$$
\footnotesize{
\xymatrix{
x \ar[r]^f \ar@{>->}[d]_i & x' \ar@{>->}[d]^{i'}\\
y \ar[r]_{f'} & y'
}}
$$
be a co-base change of $f$ along $i$. 
Then the morphism 
$\Cone(i,i')\colon \Cone f \to \Cone f'$ is 
in $w$ 
{\rm(}resp. a $C$-homotopy equivalence{\rm)}.
}
\end{para}

\begin{proof}
$\mathrm{(1)}$ 
By considering the following commutative diagram of admissible 
exact sequences, 
\begin{equation}
\label{eq:Cone ito y/x}
\footnotesize{
\xymatrix{
x \ar@{>->}[r]^{{\xi_1}_i} \ar@{=}[d] & \Cyl i \ar@{->>}[r]  \ar[d]_{\upsilon_i} & \Cone i \ar[d]\\ 
x \ar@{>->}[r]_i & y \ar@{->>}[r] & y/x,
}}
\end{equation}
it turns out that 
the canonical morphism $\Cone i\to y/x$ is in $w$ 
by the gluing axiom. 

\sn
$\mathrm{(2)}$ 
By $\mathrm{(1)}$, $\Cone i$ is in $\cC^w$. 
Thus $i$ is in $w$ by by \ref{nt:Verdier correspondence}.

\sn
$\mathrm{(3)}$ 
Notice that the induced morphism $f'/f\colon y/x\to y'/x'$ is an isomorphism and there is an admissible exact sequence 
$\Cone f\overset{\Cone(i,i')}{\rinf}\Cone f'\rdef \Cone f'/f$. 
Since $\Cone f'/f$ is $C$-contractible, we obtain the result by $\mathrm{(2)}$. 
\sn
Finally by applying assertions from $\mathrm{(1)}$ to $\mathrm{(3)}$ to the pair $(\cC_{\frob},\heq)$, we obtain assertions for `respectively' parts. 
\end{proof}

\begin{para}
\label{prop:properties of weak equivalences}
{\bf Proposition.}\ \ 
{\it
Let $\bC=(\cC,w)$ be a relative complicial exact category and let 
$\cN$ be a semi-null class of $\cC$. Then

\sn
$\mathrm{(1)}$ 
Let $u\colon x\to y$ be a morphism in $\cC$. Then

\sn
$\mathrm{(i)}$ 
$Cu\colon Cx\to Cy$ is in $w$. 

\sn
$\mathrm{(ii)}$ 
If $u$ is in $w$, then $T^{-1}u$ is in $w$. 

\sn
$\mathrm{(2)}$ 
Let $[f'\colon x'\to y']\onto{(u,v,H)}[f\colon x\to y]$ 
be a homotopy commutative square with $v\in w$. 
Then 

\sn
$\mathrm{(i)}$ 
$\Cyl(u,v,H)\colon\Cyl f'\to \Cyl f$ is in $w$.

\sn
$\mathrm{(ii)}$ 
Moreover if $u$ is in $w$, then 
$\Cone(u,v,H)\colon \Cone f'\to \Cone f$ is in $w$. 

\sn
$\mathrm{(3)}$ 
$w$ is closed under $C$-homotopic relations. 
Namely for a pair of morphisms 
$u$, $v\colon x\to y$ in $\cC$, if $u$ is in $w$ 
and there exists a $C$-homotopy $H\colon u\Rightarrow_C v$, 
then $v$ is also in $w$.

\sn
$\mathrm{(4)}$ 
If $\cN$ is closed under $w$-weak equivalences, 
then $\cN$ is a null class in $\cC$. 
}
\end{para}

\begin{proof}
$\mathrm{(1)}$ $\mathrm{(i)}$ 
Since $Cx$ and $Cy$ are $C$-contractible, the morphisms 
$Cx\to 0$ and $Cy\to 0$ are in $w$ and by $2$ out of $3$, 
$Cu$ is also in $w$. 

\sn
$\mathrm{(ii)}$ 
By considering cogluing axiom in the commutative diagram below, 
since $Pu$ and $u$ are in $w$, 
$T^{-1}u$ is also in $w$. 
$$
\xymatrix{
T^{-1}x \ar@{>->}[r]^{\iota_{T^{-1}x}} \ar[d]_{T^{-1}u} & Px \ar[d]_{Pu} 
\ar@{->>}[r]^{q_x} & x 
\ar[d]^{u}\\
T^{-1}y \ar@{>->}[r]_{\iota_{T^{-1}y}} & Py \ar@{->>}[r]_{q_y} & y. 
}
$$

\sn
$\mathrm{(2)}$ $\mathrm{(i)}$ 
Since $w$ is closed under extensions and $v$ and $Cu$ are in $w$, 
by considering the commutative diagram below, it turns out that 
$\Cyl(u,v,H)$ is in $w$. 
$$
\xymatrix{
y' \ar@{>->}[r]^{\tiny{\begin{pmatrix}0\\ \id_{y'} \end{pmatrix}}} \ar[d]_v & Cx'\oplus y' 
\ar@{->>}[r]^{\tiny{\begin{pmatrix}\id_{Cx'}& 0\end{pmatrix}}} 
\ar[d]_{\tiny{\begin{pmatrix}Cu & 0\\ -H & v\end{pmatrix}}}
& Cx' \ar[d]^{Cu}\\
y \ar@{>->}[r]_{\tiny{\begin{pmatrix}0\\ \id_{y} \end{pmatrix}}} & Cx\oplus y \ar@{->>}[r]_{\tiny{\begin{pmatrix}\id_{Cx}& 0\end{pmatrix}}} & Cx.
}
$$

\sn
$\mathrm{(ii)}$ 
In the commutative diagram below, we obtain the result by 
the gluing axiom. 
$$
\xymatrix{
x' \ar@{>->}[r]^{\xi_{1,f'}} \ar[d]_u & 
\Cyl f' \ar@{->>}[r]^{\eta_{f'}} \ar[d]_{\Cyl(u,v,H)} & 
\Cone f' \ar[d]^{\Cone (u,v,H)}\\
x \ar@{>->}[r]_{\xi_{1,f}} & \Cyl f \ar@{->>}[r]_{\eta_f}  & \Cone f.
}
$$

\sn
$\mathrm{(3)}$ 
If $u$ is in $w$, then $\Cone u$ is in $\cC^w$ 
by \ref{nt:Verdier correspondence}. 
Since there exists an isomorphism $\Cone(\id_x,\id_y,H)\colon 
\Cone u\isoto \Cone v$, $\Cone v$ is also in $\cC^w$ and 
it turns out that $v$ is in $w$ by \ref{nt:Verdier correspondence} again. 

\sn
$\mathrm{(4)}$ 
Since $w$ contains all $C$-homotopy equivalences, $\cN$ contains 
all $C$-contractible objects. 
Let $x\overset{i}{\rinf}y\overset{p}{\rdef}z$ be an admissible exact sequence in $\cC$. 
Since $\footnotesize{\begin{pmatrix}0&\id_y\end{pmatrix}}\colon Cx\oplus y\to y $ is a $C$-homotopy equivalence, 
by applying the gluing axiom to 
the commutative diagram below, 
it turns out that $\Cone(0,p)$ is also in $w$. 
$$
\xymatrix{
x \ar[d]_{\id_x} \ar@{>->}[r]^{\iota_i} & 
\Cyl i \ar[d]_{\tiny{\begin{pmatrix}0 &\id_y \end{pmatrix}}} \ar@{->>}[r]^{\eta_i} & 
\Cone i \ar[d]^{\Cone(0,p)}\\
x \ar@{>->}[r]_i & y \ar@{->>}[r]_p & z.
}
$$
Thus $z$ is in $\cN$ if and only if $\Cone i$ is in $\cN$. 
If $x$ and $y$ are in $\cN$, then $\Cone i$ and $z$ are also in $\cN$. 

If $y$ (resp. $x$) and $z$ are in $\cN$, then by 
applying argument above to 
an admissible exact sequence $y\rinf \Cone i \rdef Tx$ 
(resp. $\Cone i\rinf Tx\oplus Cy \rdef Ty$), 
$Tx$ and $x$ 
(resp. $Ty$ and $y$) are in $\cN$. 
\end{proof}

\begin{para}
\label{cor:T,P,C is relative functors}
{\bf Corollary.}\ \ 
{\it
Let $\bC=(\cC,w)$ be a relative complicial exact category. 
Then for $F\in\{T^{\pm 1},P,C \}$, $F$ preserves weak equivalences, 
namely, the pair $(F,\sigma_F^{-1})\colon\bC\to\bC$ 
is a relative complicial exact functor.
}
\end{para}

\begin{proof}
Let $f\colon x\to y$ be a morphism in $w$. 
Then $Cf$ and $Pf=CT^{-1}f$ are in 
$w$ by \ref{prop:properties of weak equivalences} 
$\mathrm{(1)}$ $\mathrm{(i)}$. 
Thus $Tf$ and $T^{-1}f$ are also in $w$ by 
the left and right commutative diagrams below and 
the gluing and cogluing axioms respectively.
$$
\footnotesize{
\xymatrix{
x \ar@{>->}[r]^{\iota_x} \ar[d]_f & Cx \ar@{->>}[r]^{\pi_x} \ar[d]_{Cf} & Tx \ar[d]^{Tf}\\
y \ar@{>->}[r]_{\iota_y} & Cy \ar@{->>}[r]_{\pi_y} & Ty,
}\ \ \ \ \ 
\xymatrix{
T^{-1}x \ar@{>->}[r]^{j_x} \ar[d]_{Pf} & Px \ar@{->>}[r]^{q_x} \ar[d]_{Pf} & x \ar[d]^{f}\\
T^{-1}y \ar@{>->}[r]_{j_y} & Py \ar@{->>}[r]_{q_y} & y.
}}
$$
\end{proof}

\begin{para}
\label{lem:closed condition of w-closure}
{\bf Lemma.}\ \ 
{\it
Let $\bC=(\cC,w_{\bC})$ be a relative complicial exact category and let 
$\cA$ be a full subcategory of $\cC$. 
If $\cA$ is closed under finite direct sums {\rm(}resp. 
the operation $T$, resp. the operation $T^{-1}${\rm)}, 
then $\cA_w$ the $w$-closure of $\cA$ is also.
}
\end{para}

\begin{proof}
Let $x$ and $y$ be objects in $\cA_w$. 
Then there exists a zig-zag sequences of morphisms 
$x\overset{a}{\leftarrow}x'\onto{a'}x''$ and 
$y\overset{b}{\leftarrow}y'\onto{b'}y''$ in $w_{\bC}$ such that 
$x''$ and $y''$ are in $\cA$. 

If $\cA$ is closed under finite direct sums, then 
$x''\oplus y''$ is in $\cA$ and the morphisms 
$\displaystyle{\begin{pmatrix}a & 0 \\ 0 & b\end{pmatrix}}$ and 
$\displaystyle{\begin{pmatrix}a' & 0 \\ 0 & b'\end{pmatrix}}$ 
are in $w_{\bC}$. 
Thus $x\oplus y$ is in $\cA_{w_{\bC}}$.   

If $\cA$ is closed under the operation $T$ (resp. $T^{-1}$), 
then $Tx''$ (resp. $T^{-1}x''$) is in $\cA$ and 
morphisms $Ta$ and $Ta'$ (resp. $T^{-1}a$ and $T^{-1}a'$) are in $w_{\bC}$. 
Thus $Tx$ (resp. $T^{-1}x$) is in $\cA_{w_{\bC}}$. 
\end{proof}

\begin{para}
\label{lem:w_bC-closed null class}
{\bf Lemma.}\ \ 
{\it
Let $\bC=(\cC,w_{\bC})$ be a relative complicial exact category and 
let $\cA$ be a null class in $\cC$. 
Then $\cA$ is $w_{\bC}$-closed if and only if 
it contains $\cC^{w_{\bC}}$. 
}
\end{para}

\begin{proof}
If $\cA$ is $w_{\bC}$-closed, then it contains all objects 
which is $w$-equivalent to $0$, 
in other words, $\cA$ contains $\cC^{w_{\bC}}$. 
Conversely assume that $\cA$ contains $\cC^{w_{\bC}}$. 
We have an equality $\cA=\cC^{w_{\cA}}$ 
by $\mathrm{(\ref{eq:Verdier correspondence 1})}$ and 
$w_{\cA}$ contains $w_{\bC}$ by assumption. 
Thus $\cA$ is $w_{\cA}$-closed and in particular $w_{\bC}$-closed.
\end{proof}

\begin{para}
\label{lemdf:w-closed Frobenius null closure}
{\bf Lemma-Definition ($w$-closed null closure).}\ \ 
Let $\bC=(\cC,w_{\bC})$ be a relative complicial exact category and let 
$\cA$ be a full subcategory of $\cC$. 
We set ${\langle \cA\rangle}_{\nullclass,w_{\bC}}:={\langle\cA\bigcup\cC^{w_{\bC}} \rangle}_{\nullclass,\frob}$. 
${\langle \cA\rangle}_{\nullclass,w_{\bC}}$ 
is the smallest $w_{\bC}$-closed null class which contains 
$\cA$ in $\cC$. 
We call ${\langle \cA\rangle}_{\nullclass,w_{\bC}}$ 
the {\it $w_{\bC}$-closed null closure of $\cA$}.
\end{para}

\begin{proof}
Since ${\langle \cA\rangle}_{\nullclass,w_{\bC}}$ 
contains $\cC^{w_{\bC}}$ and it is a null closure in $\cC_{\frob}$, 
${\langle \cA\rangle}_{\nullclass,w_{\bC}}$ 
is $w_{\bC}$-closed by Lemma~\ref{lem:w_bC-closed null class}. 
Then it turns out that 
${\langle \cA\rangle}_{\nullclass,w_{\bC}}$ 
is a null class in $\cC$ by 
Proposition~\ref{prop:properties of weak equivalences} $\mathrm{(4)}$. 
If $\cB$ is a $w_{\bC}$-closed null class which contains $\cA$, 
then $\cB$ contains $\cC^{w_{\bC}}$ and therefore contains 
${\langle \cA\rangle}_{\nullclass,w_{\bC}}$. 
\end{proof}

By \ref{lem:pull-back of semi-null classes} and 
\ref{nt:Verdier correspondence}, 
we obtain the following result.

\begin{para}
\label{lem:pull-back of complicial weak equivalences}
{\bf Lemma.}\ \ 
{\it
Let $(f,c)\colon \cC\to \cC'$ be a complicial exact functor 
between complicial exact categories and let $w$ be 
a class of complicial weak equivalences in $\cC'$. 
Then $f^{-1}(w)$ is a class of complicial weak equivalences in $\cC$. 
Moreover if $w$ is thick, then $f^{-1}(w)$ is also thick. 
}
\qed
\end{para}

\begin{para}
\label{df:stable weq}
{\bf Definition (Stable weak equivalences).}\ \ 
Let $\cX$ be 
a category with a specific zero object $0$ 
and $w$ a class of morphisms in $\cX$. 
We say that a morphism $f\colon x\to y$ in $F_b\cX$ is 
a {\it stable weak equivalence} 
if $f_{\infty}\colon x_{\infty}\to y_{\infty}$ 
is in $w$. 
We denote the class of all stable weak equivalences in $F_b\cX$ 
by ${w_{\st}}_{F_{b}\cX}$ or shortly $w_{\st}$. 

Assume $\cX$ is a complicial exact category. 
Then we similarly say that a homotopy commutative diagram 
$(f,H)\colon x\to y$ in $F_{b,h}\cX$ is a {\it stable weak equivalence} if 
$f_{\infty}\colon x_{\infty}\to y_{\infty}$ is in $w$. 
We similarly denote the class of all stable weak equivalences in $F_{b,h}\cX$ 
by ${w_{\st}}_{F_{b,h}\cX}$ or shortly $w_{\st}$. 

For any integer $k$, the degree shift functor 
${(-)}[k]\colon F_{b,h}\cC\to F_{b,h}\cC$ preserves $lw$ and $w_{\st}$.  

Moreover assume that the pair 
$\bX:=(\cX,w)$ is a normal ordinary relative complicial 
exact category. 
Then since the functor 
${(-)}_{\infty}\colon F_{b,h}\cX\to \cX$ 
is a strictly normal complicial exact functor 
and $w_{\st}$ is just a pull-back of $w$ in $\cC$ by this functor, 
$w_{\st}$ is a class of complicial weak equivalences in $F_{b,h}\cX$ 
by \ref{lem:pull-back of complicial weak equivalences}. 
We denote the normal ordinary relative complicial exact category 
Moreover if $w$ is thick, then $w_{\st}$ is also thick. 
We denote the normal ordinary relative complicial exact category 
$(F_{b,h}\cX,w_{\st})$ by $F_{b,h}\bX$. 
\end{para}

Recall the definition of categorical homotopic of 
relative functors from \ref{df:rel cat}. 

\begin{para}
\label{lem:degree shift in wst}
{\bf Lemma.}\ \ 
{\it
Let $(\cX,w)$ be a relative category and let 
$k$ be an integer. Then the degree shift functor 
$(-)[k]\colon F_{b}\cX\to F_{b}\cX$ is categorical homotopic to 
$\id_{F_{b}\cX}$.
}
\end{para}

\begin{proof}
There exists a natural weak equivalence 
${(-)}[k]\onto{\Theta^k}{(-)}[k+1]$ on $(F_b\cX,w_{\st})$. 
Namely for an object $x$ in $F_b\cX$ and an integer $n$, 
${\Theta^k(x)}_n:=i^x_{n+k}\colon x_{n+k}\to x_{n+k+1}$. 
Therefore there exists a sequence of natural 
weak equivalences from $\id_{F_b\cX}=(-)[0]$ to $(-)[k]$.

Notice that $\Theta^k$ is natural weak equivalence 
with respect to $w_{\st}$, not to $lw$ and if 
$\cX$ is a complicial exact category, then 
$\Theta^k$ is functorial on $F_b\cX$, but not on $F_{b,h}\cX$. 
\end{proof}

\subsection{The category of bounded chain complexes over relative complicial exact categories}
\label{subsec:chain complexes on cec}

Throughout this subsection, 
let $\bC=(\cC,w)$ be a relative complicial exact category. 
In this subsection, we will study compatibilities of 
several complicial structures and 
the class of complicial weak equivalences on $\Ch_b(\cC)$ 
the category of bounded chain complexes on $\cC$. 

\begin{para}
\label{nt:comp structure on ChbC}
{\bf (Complicial structures on $\Ch_b(\cC))$.}\ \ 
There are two kind of natural complicial structures on 
$\Ch_b(\cC)$. 
The first one is the standard complicial structure 
illustrated in \ref{ex:standard complicial structure}. 
We can regard $\Ch_b(\cC)$ as the full subcategory of 
$F\cC=\cC^{\bbZ}$ the category of functors from the totally ordered set 
$\bbZ$ to $\cC$ and the second one is the 
level complicial structure 
(see \ref{df:level complicial structure}) 
induced from $\cC^{\bbZ}$. 
We denote these two complicial structures by 
$(C^{\standard},\iota^{\standard},r^{\standard},\sigma^{\standard})$ and 
$(C^{\lv},\iota^{\lv},r^{\lv},\sigma^{\lv})$ 
respectively. 
We denote the associated suspension functor, path functor, cone functor and so on by $T^{\#}$, $P^{\#}$, $\Cone^{\#}$ and so on for 
$\#\in\{\standard,\lv\}$. 

We denote the forgetful functor $\CompEx\to\ExCat$ 
by $F$. 
For a small exact category $\cE$, we define 
$j_{\cE}\colon\cE\to \Ch_b(\cE)$ to be an exact functor 
by setting ${j_{\cE}(x)}_k$ is $x$ if $k=0$ and $0$ if $k\neq 0$. 
Then $j_{\cE}$ is natural on $\cE$, namely 
$j$ gives a natural transformation $j\colon\id_{\ExCat}\to F\Ch_b$. 
If $\cE$ is a complicial exact category, then 
the functor 
$j_{\cE}\colon \cE\to \Ch_b(\cE)$ 
is a strictly normal complicial exact functor 
with respect to $(C^{\lv},\iota^{\lv},r^{\lv},\sigma^{\lv})$. 
\end{para}

\begin{para}
\label{conv:chain complexes}
{\bf (Conventions of chain complexes).}\ \ 
Let $x$ be a chain complex on $\cC$ and let $f\colon x\to y$ 
be a chain morphism on $\cC$ and let $n$ be an integer. 
We define $\sigma_{\geq n}x$ and $\sigma_{\geq n}f$ to be 
a subcomplex of $x$ and a chain morphism $\sigma_{\geq n}x\to\sigma_{\geq n}y$ by setting 
\begin{equation}
\label{eq:sigmgeqn df}
{(\sigma_{\geq n}x)}_k=
\begin{cases}
x_k & \text{if $k\geq n$}\\
0 & \text{otherwise}
\end{cases},\ \ 
{(\sigma_{\geq n}f)}_k=
\begin{cases}
f_k & \text{if $k\geq n$}\\
0 & \text{otherwise}
\end{cases}
\end{equation}
Similarly we define $\sigma_{\leq n}x$ and 
$\sigma_{\leq n}f\colon \sigma_{\leq n}x\to \sigma_{\leq n}y$. 
We simply write $x[n]$ and $f[n]$ for ${(T^{\standard})}^n(x)$ and 
${(T^{\standard})}^n(f)$ respectively. 

For an integer $n\leq m$, 
we write $\Ch_{b,[n,m]}(\cC)$, $\Ch_{b,\geq n}(\cC)$ 
and $\Ch_{b,m\geq}(\cC)$
 for the 
full subcategories of $\Ch_b(\cC)$ consisting of those complexes 
$x$ such that $x_k=0$ if $k\notin [n,m]$, $k<n$ and $m<k$ respectively. 

For a complex $x$ in $\Ch_b(\cC)$, we set 
\begin{equation}
\label{eq:min df}
\min x:=\min\{n\in\bbZ;x_n\neq 0\},
\end{equation}
\begin{equation}
\label{eq:max df}
\max x:=\max\{n\in\bbZ;x_n\neq 0\},
\end{equation}
\begin{equation}
\label{eq:length df}
\length x:=\min\{b-a;x\in\Ob\Ch_{[a,b]}(\cC)\}
\end{equation}
and call it the {\it length of $x$}. 
\end{para}

\begin{para}
\label{lemdf:total functor}
{\bf Lemma-Definition (Exact functor $\Tot$).}\ \ 
There exists an exact functor 
$\Tot\colon\Ch_b(\cC)\to\cC$ and 
the natural equivalences 
$d^{\standard}\colon C\Tot\isoto\Tot C^{\standard}$ 
and 
$d^{\lv}\colon C\Tot\isoto \Tot C^{\lv}$ such that the pairs 
$(\Tot,d^{\standard})$ and 
$(\Tot,d^{\lv})$ are complicial exact functors with resepct to 
the standard complicial structure and the level complicial 
structure on $\Ch_b(\cC)$ respectively. 

For an integer $n$, we have the equality
\begin{equation}
\label{eq:Totj=id}
\Tot j_{\cC}(-)[n]=T^n. 
\end{equation}
Moreover we can regard $\Tot$ is the natural transformation 
$\Tot\colon\Ch_bF\to\id_{\CompEx}$ such that 
the pair $(j,\Tot)$ gives an adjunction of the pair 
$(\Ch_b,F)$ where $F\colon\CompEx\to \ExCat$ is the 
forgetful functor. 
\end{para}

\begin{proof}
We recall the construction of $\Tot$ from \cite[4.1, 4.9]{Moc10}. 
First we will define 
$\Tot_{\geq 0}\colon \Ch_{b,\geq 0}(\cC)\to \cC$, 
$d^{\standard}\colon C\Tot_{\geq 0}\isoto\Tot_{\geq 0} C^{\standard}$ 
and 
$d^{\lv}\colon C\Tot_{\geq 0}\isoto \Tot_{\geq 0} C^{\lv}$ 
to be exact functor and natural equivalences 
respectively by proceeding induction 
on the length of chain complexes. 
For an object $x$ and a morphism $f\colon x\to y$ in $\cC$, 
we set $\Tot_{\geq 0}\cdot j_{\cC}(x):=x$ and $\Tot_{\geq 0}j_{\cC}(f):=f$. 
Notice that we have the equality $Cx=\Cone\id_x=\Tot(x\onto{id_x}x)=\Tot C^{\standard}x$. 
We set $d^{\lv}_{j_{\cC}(x)}=d^{\standard}_{j_{\cC}(x)}=\id_{Cx}
\colon C\Tot_{\geq 0}j_{\cC}(x)=Cx=\Tot C^{\lv}j_{\cC}(x)=
\Tot C^{\standard}j_{\cC}(x)$. 

For a chain complex $x$ in $\Ch_{b,\geq 0}({\cC})$ and 
a chain morphism $f\colon x\to y$, we define 
\begin{equation}
\label{eq:Tot eq df 1}
\Tot_{\geq 0}(x):=\Cone(\Tot{\geq 0}((\sigma_{\geq 1}x)[-1])\to x_0),
\end{equation}
\begin{equation}
\label{eq:Tot eq df 2}
\Tot_{\geq 0}(f):=\Cone(\Tot_{\geq 0}((\sigma_{\geq 1}f)[-1]),f_0).
\end{equation}
Moreover for $\#\in\{\standard,\lv\}$, we 
define $d^{\#}_x\colon C\Tot x\isoto \Tot C^{\#}x$ by the compositions 
$${\scriptstyle{
C\Tot x=C\Cone(\Tot(\sigma_{\geq 1}x[-1])\onto{d_1^x}x_0)
\overset{c^{\Cone}}{\isoto}
\Cone(C\Tot(\sigma_{\geq 1}x[-1])\onto{Cd_1^x}Cx_0)
\overset{\Cone(d^{\#}_{\sigma_{\geq 1}x[-1]},\id_{Cx})}{\isoto}\Cone(\Tot(C^{\#}\sigma_{\geq 1}x[-1])\onto{Cd_1^x}Cx_0)
=\Tot C^{\#}x.}}
$$
Next 
for any positive integer $n$, 
We define 
$\Tot_{\geq -n}\colon \Ch_{b,\geq -n}(\cC)\to \cC$ 
and $X_{-n}\colon
\Tot_{\geq -(n+1)}|_{\Ch_{b,\geq -n}(\cC)}\isoto\Tot_{\geq -n}$ 
to be an exact functor and a natural equivalence 
by setting 
\begin{equation}
\label{eq:Totgeq-n df}
\Tot_{\geq -n}:={(T^{\cC})}^{-n}\cdot\Tot_{\geq 0}\cdot{(T^{\standard})}^n
\end{equation}
 and compositions of natural equivalences 
\begin{equation}
\label{eq:Xn df}
{\scriptstyle{
\Tot_{\geq -(n+1)}=
{(T^{\cC})}^{-(n+1)}\cdot\Tot_{\geq 0}\cdot{(T^{\standard})}^{n+1}\isoto 
{(T^{\cC})}^{-n}\cdot\Tot_{\geq 0}\cdot{(T^{\standard})}^{-1}\cdot{(T^{\standard})}^{n+1}
\isoto {(T^{\cC})}^{-n}\cdot\Tot_{\geq 0}\cdot{(T^{\standard})}^n=\Tot_{\geq -n}}}
\end{equation}
respectively. 
We will define $\Tot\colon\Ch_b(\cC)\to \cC$ to be 
an exact functor by patching the family of functors and 
natural equivalences 
$\{\Tot_{\geq -n},X_{-n} \}_{n>0}$. 
To do so, for any pair of integers $n>m$, 
we define $X_{n,m}\colon \Tot_{\geq n}\isoto \Tot_{\geq m}$ 
to be a natural equivalence by 
setting $X_{n,m}:=X_{m}X_{m+1}\cdots X_{n-1}$ and 
for any complex $x$ in $\Ch_b(\cC)$ 
we set
\begin{equation}
\label{eq:m(x) df}
m(x):=\min\{n;n\geq 0,x_{-k}=0\text{ for any $k\geq n$}\}.
\end{equation}
For a complex $x$ in $\Ch_b(\cC)$, we set 
$\Tot (x):=\Tot_{\geq -m(x)}(x)$. 
Next for a morphism $f\colon x\to y$, we take an integer 
$n\geq \max\{m(x),m(y)\}$ and we set 
$\Tot(f):=X_{n,m(y)}(y)\cdot\Tot_{\geq n}(f)\cdot{(X_{n,m(x)})(x)}^{-1}$. 
We can easily check that the definition of $\Tot(f)$ 
does not depend upon a choice of $n$ and 
that $\Tot$ is an exact functor on $\Ch_b(\cC)$. 

For $\#\in\{\standard,\lv\}$, 
we will define $d^{\#}\colon C^{\cC}\Tot\isoto \Tot C^{\#}$ 
to be a natural equivalence in the following way. 
For a complex $x$ in $\Ch_b(\cC)$, 
let us notice that we have equalities $m(x)=m(C^{\#}x)$ 
and $T^{\standard}C^{\#}=C^{\#}T^{\standard}$ 
for $\#\in\{\standard,\lv\}$. 
We define $d_x^{\#}\colon C^{\cC}\Tot(x)\isoto \Tot C^{\#}(x)$ 
to make the diagram below commutative
$$\xymatrix{
\Tot (x) \ar@{=}[r] \ar[d]_{\Tot\ast\iota_x^{\standard}} & 
{\scriptstyle{{(T^{\cC})}^{-m(x)}\cdot \Tot_{\geq 0}\cdot {(T^{\standard})}^{m(x)}(x)}} 
\ar@{=}[r] 
\ar[d]_{{\tiny{({(T^{\cC})}^{-m(x)}\cdot \Tot_{\geq 0}\cdot {(T^{\standard})}^{m(x)})\ast \iota_x}}} &
\Tot(x) \ar[d]^{\iota_{\Tot(x)}}\\
\Tot\cdot C^{\#}(x) \ar@{=}[r] &  
{\scriptstyle{{(T^{\cC})}^{-m(x)}\cdot \Tot_{\geq 0}\cdot {(T^{\standard})}^{m(x)}\cdot C^{\#}(x)}} 
\ar[r]_{\ \ \ \ \ \ \ \ \ \ \ \ \ A}^{\ \ \ \ \ \ \ \ \ \ \ \ \ \sim} & 
C^{\cC}\cdot\Tot(x)}$$
where the morphism $A$ is compositions of natural equivalences
$${\scriptstyle{ 
{(T^{\cC})}^{-m(x)}\cdot \Tot \cdot {(T^{\standard})}^{-m(x)}\cdot 
C^{\#} (x) 
\isoto {(T^{\cC})}^{-m(x)}\cdot C^{\cC}\cdot \Tot \cdot 
{(T^{\standard})}^{-m(x)}(x) 
\isoto C^{\cC}\cdot {(T^{\cC})}^{-m(x)}\cdot 
\Tot \cdot {(T^{\standard})}^{-m(x)}(x). 
}}$$
For $\#\in\{\standard,\lv\}$, by induction on the length of 
chain complexes, we can show the equality 
$d^{\#}\cdot(\iota^{\#}\ast\Tot)=\Tot\ast\iota$. 
\end{proof}

\begin{para}
\label{dfcor:Tot associativity}
{\bf Definition-Corollary.}\ \ 
Let $a\leq i\leq b$ be a triple of integers 
such that $b-a\geq 2$ and 
let $f\colon x\to y$ be a chain morphism in $\Ch_{[a,b]}(\cC)$. 
Then there is a canonical chain morphism 
$\sigma_{\geq i}x[-1]\to\sigma_{\leq i-1}x$ 
induced from the $i$th boundary morphism $d_i^x\colon x_i \to x_{i-1}$ 
and we denote it by $\fd_i^x\colon(\sigma_{\geq i}x)[-1]\to\sigma_{\leq i-1}x$. We set 
$$\Cone\Tot\fd_i^f:=\Cone((\sigma_{\geq i}f)[-1],\sigma_{\leq i-1}f)
\colon\Cone\Tot\fd_i^x\to\Cone\Tot\fd_i^y.
$$
Then the association $\Cone\Tot\fd_i\colon\Ch_{[a,b]}\cC\to\cC$ is 
an exact functor and there is a natural equivalence 
$\fe\colon\Tot\isoto\Cone\Tot\fd_i$. 
\end{para}

\begin{proof}
We will construct a natural equivalence $\fe$ by proceeding induction 
on the length of $\sigma_{\leq i-1}x$. 
By replacing $x$ with $x[N]$ for a suitable integer $N$, 
we shall assume $a=0$. 
Notice that if $b=2$ and $i=1$, then $\fe$ is defined in 
\ref{cor:ConeCone}. 
If $\length \sigma_{\leq i-1}x =1$, then 
by definition, we have the equality $\Tot x=\Cone\Tot \fd_i^x$. 
Assume $\length \sigma_{\leq i-1}x >1$, then 
there are functorial isomorphisms below by inductive hypothesis
\begin{eqnarray*}
\Cone\Tot \fd_i^x & = & 
\Cone(\Tot((\sigma_{\geq i}x)[-1])\to\Cone(\Tot((\sigma_{\geq 1}\sigma_{\leq i-1}x)[-1])\to x_0)\\
& \isoto & 
\Cone(\Cone(\Tot(\sigma_{\geq i}x)\to\Tot((\sigma_{\geq 1}\sigma_{\leq i-1}x)[-1]))\to x_0 )\\
& \isoto & 
\Cone(\Tot((\sigma_{\geq 1}x)[-1])\to x_0)\\
& = & \Tot x.
\end{eqnarray*}
\end{proof}

\ref{lem:functoriality of closures} and 
\ref{lemdf:total functor} imply the following corollary. 

\begin{para}
\label{cor:Tot restriction}
{\bf Corollary.}\ \ 
{\it
Let $\cF$ be an additive full subcategory of $\cC$. 
We regard $\cF$ as a full subcategory of $\Ch_b(\cC)$ by the 
functors $\cF\rinc\cC\onto{j_{\cC}} \Ch_b(\cC)$. 
Then the prenull closure and semi-null closure of $\cF$ 
in $\Ch_b(\cC)$ with respect to the standard complicial structure 
are $\Ch_{b,\geq 0}(\cF)$ and $\Ch_b(\cF)$ respectively. 
The functor $\Tot$ induces the 
complicial functors 
$\Tot_{\geq 0}\colon\Ch_{b,\geq 0}(\cF)\to{\langle\cF \rangle}_{\prenull}$ 
and 
$\Tot\colon\Ch_b(\cF)\to{\langle\cF \rangle}_{\seminull}$. 
}\qed
\end{para}

\begin{para}
\label{cor:ChbBthi}
{\bf Corollary.}\ \ 
{\it
Let $\cA$ be an additive category and let $\cB$ be an additive 
full subcategory of $\cA$. 
Then we have the equality
\begin{equation}
\label{eq:ChbBthi}
\Ch_b(\cB_{\thi,\cA})={(\Ch_b(\cB))}_{\thi,\Ch_b(\cA)}.
\end{equation}
}
\end{para}

\begin{proof}
We regard $\cA$ as an exact category with split exact sequences. 
Then $\Ch_b(\cA)$ is a complicial exact category with 
the standard complicial structure. 
We regard $\cB$ and $\cB_{\thi,\cA}$ as full subcategories of 
$\Ch_b(\cA)$ by the functor $j_{\cA}\colon\cA\to\Ch_b(\cA)$ and 
inclusions $\cB,\ \cB_{\thi,\cA}\rinc \cA$. 

We will show that $\Ch_b(\cB_{\thi,\cA})$ is thick in $\Ch_b(\cA)$. 
Let $y$ be a chain complex in $\Ch_b(\cB_{\thi,\cA})$ and let 
$x\onto{i}y\onto{p}x$ be a pair of chain morphisms in $\Ch_b(\cA)$ 
such that $pi=\id_x$. 
Then for each integer $n$, 
we have an equality 
$p_ni_n=\id_{x_n}$ and 
it turns out that 
$x_n$ is in $\cB_{\thi,\cA}$. 
Hence $x$ is in $\Ch_b(\cB_{\thi,\cA})$ and $\Ch_b(\cB_{\thi,\cA})$ 
is thick in $\Ch_b(\cA)$. 
In particular we have the inclusion 
${(\Ch_b(\cB))}_{\thi,\Ch_b(\cA)} \rinc \Ch_b(\cB_{\thi,\cB})$.

Next the inclusion $\cB\rinc \Ch_b(\cB)$ implies the inclusion 
$\cB_{\thi,\cA}\rinc {(\Ch_b(\cB))}_{\thi,\Ch_b(\cA)} $. 
Notice that by \ref{prop:fundamental properties of thick closure} 
$\mathrm{(5)}$, 
${(\Ch_b(\cB))}_{\thi,\Ch_b(\cA)}$ is a semi-null class in $\Ch_b(\cA)$. 
Thus we obtain the inclusion $\Ch_b(\cB_{\thi,\cA})={\langle\cB_{\thi,\cA} \rangle}_{\seminull,\Ch_b(\cA)}\rinc {(\Ch_b(\cB))}_{\thi,\Ch_b(\cA)}$. 
\end{proof}

\begin{para}
\label{lemdf:tw=qw}
{\bf Lemma-Definition.}\ \ 
We say that a morphism $f\colon x\to y$ in $\Ch(\cC)$ is 
a {\it Frobenius quasi-isomorphism} if 
it is a quasi-isomorphism with respect to the Frobenius admissible 
exact structure of $\cC$. 
We denote the class of all Frobenius quasi-isomorphisms in $\Ch_{\#}(\cC)$ 
by $\qfis_{\Ch_{\#}(\cC)}$ or simply $\qfis$. 
We write ${q^fw}_{\Ch_{\#}(\cC)}$ or simply $q^fw$ for 
the class of thick complicial weak equivalences 
spanned by $\qfis$ and $lw$ and call it 
the {\it class of Frobenius quasi-weak equivalences} 
({\it in $\Ch_{\#}(\cC)$}). 
Assume that $w$ is thick. 
We denote the class of all morphisms $f\colon x\to y$ in $\Ch_b(\cC)$ 
such that $\Tot f$ is in $w$ by $tw$ and call it 
the {\it class of total quasi-weak equivalences}. 
Then there are equalities
\begin{equation}
\label{eq:tw=qw}
q^fw=tw=qw.
\end{equation}
\end{para}

\begin{proof}
Since the class $tw$ is the pull-back of $w$ by the complicial exact functor 
$\Tot$, 
it is a class of thick complicial weak equivalences in $\cC$ 
by \ref{lem:pull-back of complicial weak equivalences}. 
We will show the inclusions 
$$q^fw\underset{\textbf{I}}{\subset}qw\underset{\textbf{II}}{\subset}tw\underset{\textbf{III}}{\subset} q^fw.$$
Since a Frobenius admissible exact sequence in $\cC$ is an 
admissible exact sequence in $\cC$, 
the inclusion $\textbf{I}$ is obvious. 
To show inclusion $\textbf{II}$, what we need to prove is the 
following two inclusions:

\sn
\textbf{IV} $\qis\subset tw$.

\sn
\textbf{V} $lw\subset tw$.

\sn
The inclusion $\textbf{IV}$ follows form \cite[4.5]{Moc13b}. 
We will show the inclusion $\textbf{V}$. 
Let $f\colon x\to y$ be a chain morphism in $\Ch_b(\cC)$. 
By considering $f[N]\colon x[N]\to y[N]$ for sufficiently large 
integer $N$, without loss of generality, 
we shall assume that $f$ is in $\Ch_{b,\geq 0}(\cC)$. 
Assume that $f$ is in $lw$. 
Since by definition, we have the equality 
$\Tot(f)=\Cone((\sigma_{\geq 1}f)[-1],f_0)$, 
we obtain the result by \ref{prop:properties of weak equivalences} 
$\mathrm{(2)}$ $\mathrm{(ii)}$ and induction on 
$\max\{\length(x),\length(y) \}$. 

Next we will show the inclusion $\textbf{III}$. 
What we need to show is the inclusion 
${(\Ch_b(\cC))}^{tw}\subset {(\Ch_b(\cC^w))}_{\qfis}$ and 
it follows from \ref{lemdf:uk and vk} below.
\end{proof}

\begin{para}
\label{lemdf:uk and vk}
{\bf Lemma-Definition.}\ (\cf \cite[4.10]{Moc13b}.)\ \  
Let $k$ be an integer. 
We write $\fraki\colon\Ch_{b,\leq k-1}(\cC)\to\Ch_{b,\leq k}(\cC)$ 
for the canonical inclusion functor. 
We define $\fQ_k\colon\Ch_{b,\leq k}(\cC)\to\Ch_{b,\leq k}(\cC)$ 
and $\fL_k\colon\Ch_{b,\leq k}(\cC)\to \Ch_{b,\leq k-1}(\cC)$ 
and $\fu_k\colon\id_{\Ch_{b,\leq k}(\cC)}\to \fQ_k$ 
and $\fv_k\colon \fraki\fL_k\to\fQ_k$ 
be a pair of functors and a pair of natural transformations 
respectively by setting 
for a chain morphism $f\colon x\to y$ in $\Ch_{b,\leq k}(\cC)$,
\begin{equation}
\label{eq:fQk df}
{\scriptscriptstyle{
{\fQ_k(x)}_n:=
\begin{cases}
x_n & \text{if $n\geq k-2$}\\
\Cone d_k^x & \text{if $n=k-1$}\\
Cx_k & \text{if $n=k$}\\
0 & \text{if $n>k$}, 
\end{cases}\ \ 
{d_n^{\fQ_k(x)}}:=
\begin{cases}
d^x_n & \text{if $n\geq k-2$}\\
\Cone(0, d_{k-1}^x) & \text{if $n=k-1$}\\
\mu_{d^x_k} & \text{if $n=k$}\\
0 & \text{if $n>k$}, 
\end{cases}\ \ 
{\fQ_k(f)}_n:=
\begin{cases}
f_n & \text{if $n\geq k-2$}\\
\Cone(f_k,f_{k-1}) & \text{if $n=k-1$}\\
Cf_k & \text{if $n=k$}\\
0 & \text{if $n>k$}, 
\end{cases}}}
\end{equation}
\begin{equation}
\label{eq:fLk df}
{\scriptscriptstyle{
{\fL_k(x)}_n:=
\begin{cases}
x_n & \text{if $n\geq k-2$}\\
\Cone d_k^x & \text{if $n=k-1$}\\
0 & \text{if $n>k-1$}, 
\end{cases}\ \ 
d_n^{\fL_k(x)}:=
\begin{cases}
d^x_n & \text{if $n\geq k-2$}\\
\Cone(0,d_{k-1}^x) & \text{if $n=k-1$}\\
0 & \text{if $n>k-1$}, 
\end{cases}\ \ 
{\fL_k(f)}_n:=
\begin{cases}
f_n & \text{if $n\geq k-2$}\\
\Cone(f_k,f_{k-1}) & \text{if $n=k-1$}\\
0 & \text{if $n>k-1$}. 
\end{cases}}}
\end{equation}
\begin{equation}
\label{eq:fu fv df}
{\fu_k(x)}_n:=
\begin{cases}
\id_{x_n} & \text{if $n\geq k-2$}\\
\kappa_{d_k^x} & \text{if $n=k-1$}\\
\iota_{x_k} & \text{if $n=k$}\\
0 & \text{if $n>k$}, 
\end{cases}\ \ \ \ \ \ \ \ \ 
{\fv_k(x)}_n:=
\begin{cases}
\id_{x_n} & \text{if $n\geq k-2$}\\
\id_{\Cone d_k^x} & \text{if $n=k-1$}\\
0 & \text{if $n>k-1$}.
\end{cases}
\end{equation}
Then for each complex $x$ in $\Ch_{b,\leq k}(\cC)$, 
$\fu_k(x)$ is a Frobenius quasi-isomorphism by 
\ref{lem:kel90} and $\fv_k(x)$ is a level weak equivalence. 

For a complex $x$ in $\Ch_{b,\leq k}(\cC)$, 
we denote the zig-zag sequence of morphisms 
$x\onto{\fu_k(x)}\fQ_k(x)\overset{\fv_k(x)}{\leftarrow}\fL_k(x)$ 
by $\cU_x^k$. If for an integer $n$, 
$x[n]$ is also in $\Ch_{b,\leq k}(\cC)$, 
then we have the equality
\begin{equation}
\label{eq:calU degree shift}
\cU_{x[n]}^k=\cU^{k-n}_x[n].
\end{equation}
For a pair of integers $a<b$, we write $\fL^{[a,b]}$ 
for the compositions 
$$\fL_{a+1}\fL_{a+2}\cdots \fL_b\colon\Ch_{b,\leq b}(\cC)\to \Ch_{b,\leq a}(\cC)$$
and for a complex $x$ in $\Ch_{b,\leq b}(\cC)$, we write $\cU^{[a,b]}_x$ for the composition of 
zig-zag sequence of morphisms $\cU^b_x\cU^{b-1}_x\cdots \cU_x^{a+1}$ 
which connects $x$ and $\fL^{[a,b]}x$. 
Moreover we have the equality
\begin{equation}
\label{eq:L[a,b]x}
\fL^{[a,b]}(x)=j_{\cC}(T^m\Tot x)[-m]
\end{equation}
where $m=\min x-a$. 
Thus for a complex $x$ in $\Ch_b(\cC)$, there exists 
a zig-zag sequence of morphisms in $lw\cup\qfis$ which connects 
$x$ and $j_{\cC}(T^m\Tot(x))[-m]$ for a suitable integer $m$. 
\qed
\end{para}

\begin{para}
\label{cor:qw=qwthi}
{\bf Corollary.}\ \ 
{\it
For a relative complicial exact category $\bC=(\cC,w)$, 
we have equalities
\begin{equation}
\label{eq:lwthi}
l({\langle w \rangle}_{\thi,\cC})={\langle lw \rangle}_{\thi,\Ch_b(\cC)},
\end{equation}
\begin{equation}
\label{eq:qw=qwthi}
qw=q({\langle w \rangle}_{\thi,\cC}).
\end{equation}
}
\end{para}

\begin{proof}
By correspondence \ref{nt:Verdier correspondence}, the equality 
$\mathrm{(\ref{eq:lwthi})}$ is equivalent to the equality 
$${(\Ch_b(\cC))}_{\thi,\Ch_b(\cC)}=\Ch_b({(\cC^w)}_{\thi,\cC})$$
and it follows from the equality $\mathrm{(\ref{eq:ChbBthi})}$. 

The inclusion $w\subset {\langle w\rangle}_{\thi,\cC}$ 
implies the inclusion $qw\subset q({\langle w\rangle}_{\thi,\cC})$. 
On the other hand, since the class $qw$ is thick, the inclusion 
$lw\subset qw$ implies 
$l({\langle w\rangle}_{\thi,\cC})={\langle lw\rangle}_{\thi,\Ch_b(\cC)}
\subset qw$. 
Thus we obtain the converse inclusion 
$q({\langle w\rangle}_{\thi,\cC})\subset qw$.
\end{proof}

\begin{para}
\label{cor:adjunction j and Tot}
{\bf Corollary.}\ \ 
{\it
We denote the forgetful functor $\RelComp\to\RelEx$ by $F$. 
Then the pair $(j,\Tot)$ can be regarded as
natural transformations 
$j\colon \id_{\RelEx}\to F\Ch_b$ and 
$\Tot\colon \Ch_b F\to \id_{\RelComp}$ 
between functors $\Ch_b\colon\RelEx\to\RelComp$ and 
$F\colon \RelComp\to\RelEx$ and it gives an adjunction of the pair 
$(\Ch_b,F)$.  
}\qed
\end{para}

\begin{para}
\label{df:canonical factorization}
{\bf Lemma-Definition.}\ \ 
Let $\bC=(\cC,w_{\bC})$ be a relative complicial exact category and 
let $n$ be a non-negative integer and let 
$x$ be a complex in $\Ch_{[0,n]}(\cC)$. 
We define $\bar{x}$ to 
be a complex in $\Ch_{[0,n]}(\cC)$ in the following way. 
First by convention, we set
\begin{equation}
\label{eq:{d'}_k^x df 1}
{d'}_k^x:=\begin{cases}
0\colon 0\to x_n & \text{if $k=n+1$}\\
d_n^x\colon x_n\to x_{n-1} & \text{if $k=n$}.
\end{cases}
\end{equation}
For $k\geq 1$, we inductively set 
\begin{equation}
\label{eq:{d'}_k^x df 2}
{d'}_{n-k}^x:=\Cone(0,d_{n-k}^x)\colon \Cone {d'}_{n-k+1}^x \to x_{n-k-2}.
\end{equation}
$$
\xymatrix{
x_n \ar[r] \ar[d]_{{d'}_n^x} & 0 \ar[d]\\
x_{n-1} \ar[r]_{d_{n-1}^x} & x_{n-1},
}\ \ 
\xymatrix{
\Cone {d'}_{n-k+2}^x \ar[r] \ar[d]_{{d'}_{n-k+1}^x} & 0 \ar[d]\\
x_{n-k} \ar[r]_{d_{n-k}^x} & x_{n-k-1}.
}
$$
Well-definedness of ${d'}_{n-k}^x$ follows from the equality
\begin{equation}
\label{eq:dn-kd'n-k+1=0}
d^x_{n-k}{d'}^x_{n-k+1}=0.
\end{equation}
Then we set
\begin{equation}
\label{eq:xbar df}
\bar{x}_k:= \Cone {d'}_{k+1}^x\ \ \  \text{and}\ \ \  d^{\bar{x}}_k=\kappa_{{d'}_k^x}{d'}_k^x.
\end{equation}
We define $\fu_x\colon x\to \bar{x}$ to be a chain morphism by setting
\begin{equation}
\label{eq:fuxk df}
{\fu_x}_k:=\kappa_{{d'}^x_{k+1}}\colon x_k \to \bar{x}_k.
\end{equation}
Notice that we can identify $\bar{x}_0=\Tot x$.

Next let $y$ be an object in $\cC$ and 
let $u\colon \Tot x \to y$ be a morphism in $\cC$. 
Then we define $z_u$ and $g_u\colon \bar{x}\to z_u$ 
to be a complex and 
a chain morphism in $\Ch_{[0,n]}(\cC)$ by setting
\begin{equation}
\label{eq:zu df}
{z_u}_k:=\begin{cases}
C(\bar{x}_k) & \text{if $k\geq 1$}\\
y & \text{if $k=0$}
\end{cases}\ \ \ 
\text{and}\ \ \ 
d_k^{z_u}:=\begin{cases}
C(d_k^{\bar{x}}) & \text{if $k\geq 2$}\\
u\mu_{{d'}_1^x} & \text{if $k=1$}
\end{cases}
\ \ \ \text{and}
\end{equation}
\begin{equation}
\label{eq:g_u df}
{g_u}_k:=\begin{cases}
\iota_{x_k} & \text{if $k\geq 1$}\\
\id_{\bar{x}_0} & \text{if $k=0$}.
\end{cases}
\end{equation}
For a complex $s$ in $\Ch_{[0,n]}(\cC)$, we 
denote the canonical morphism $\Tot \sigma_{\geq 1}s \to s_0$ induced from 
$d_0^s$ by $\underline{d}_0^s$. 
Then we have an isomorphism 
$c^{\Tot}_{\sigma_{\geq 1}\bar{x}}\colon C(\Tot\sigma_{\geq 1}\bar{x})\isoto \Tot(C^{\lv}(\sigma_{\geq 1}\bar{x}))=\Tot\sigma_{\geq 1}z_u$ and 
there is a homotopy commutative square 
$(0,\id_y,-\underline{d}_0^{z_u}r_{\Tot\sigma_{\geq 1}\bar{x}}{(c^{\Tot}_{\sigma_{\geq 1}\bar{x}})}^{-1})\colon[\Tot\sigma_{\geq 1}z_u \to y] \to [0\to y]$. 
We set
\begin{equation}
\label{eq:a_u}
a_u:=\Cone(0,\id_y,-\underline{d}_0^{z_u}r_{\Tot\sigma_{\geq 1}\bar{x}}{(c^{\Tot}_{\sigma_{\geq 1}\bar{x}})}^{-1})\colon\Tot z_u \to y. 
\end{equation}
Then $a_u$ is in $w_{\bC}$ and we have an equality
\begin{equation}
\label{eq:u=auTotguux}
u=a_u\Tot(g_u\fu_x).
\end{equation}
\end{para}

\begin{proof}
The equality $\mathrm{(\ref{eq:dn-kd'n-k+1=0})}$ follows from
$$d_{n-k}^x{d'}_{n-k+1}^x=\Cone(0,d^x_{n-k}d^x_{n-k+1})=\Cone(0,0)=0.$$
To show that $\bar{x}$ and $z_u$ are complexes, 
we check equalities
$$d_{k-1}^{\bar{x}}d_k^{\bar{x}}=\kappa_{{d'}^x_{k-1}}{d'}^x_{k-1}\kappa_{{d'}^x_k}{d'}_k^x=\kappa_{{d'}_{k-1}^x}d_{k-1}^xd_k^x=0, $$
$$d_{k-1}^{z_u}d_{k}^{z_u}=C(d_{k-1}^{\bar{x}}d_{k}^{\bar{x}})=0\ \ \ (k\geq 3), $$
\begin{multline*}
d_1^{z_u}d_2^{z_u}=u\mu_{{d'}_1^x}C(d_2^{\bar{x}})=
u\mu_{{d'}_1^x}C(\kappa_{{d'}_2^x}{d'}_2^x)=
u\kappa_{{d'}_1^x}{d'}_1^x\mu_{{d'}_2^x}\\
=u\kappa_{{d'}_1^x}\Cone(0,d_1^x)\Cone(\id_{\Cone {d'}_3^x},{d'}_2^x)=
u\kappa_{{d'}_1^x}\Cone(0,d_1^x{d'}_2^x)=0.
\end{multline*}
$a_u$ is in $w_{\bC}$ 
by Proposition~\ref{prop:properties of weak equivalences} 
$\mathrm{(2)}$ $\mathrm{(ii)}$ and 
applying Lemma~\ref{lem:Cone(0,g,-gmufrx)Cone(iotax,kappaf)} below to 
homotopy commutative squares 
$${\scriptstyle{[\Tot\sigma_{\geq 1}x \onto{\underline{d}_0^x} x_0] 
\onto{(\Tot\sigma_{\geq 1}\fu_x,\kappa_{\underline{d}_0^x})} 
[\Tot\sigma_{\geq 1}\bar{x}\onto{\underline{d}_0^{\bar{x}}} \Tot x]
\onto{(c^{\Tot}_{\sigma_{\geq 1}\bar{x}}\iota_{\Tot\sigma_{\geq 1}\bar{x}},u)} 
[\Tot\sigma_{\geq 1}z_u\onto{\underline{d}_0^{z_u}}y] 
\onto{(0,\id_y,-\underline{d}_0^{z_u}r_{\Tot\sigma_{\geq 1}\bar{x}}{(c^{\Tot}_{\sigma_{\geq 1}\bar{x}})}^{-1})} 
[0\to y],}}$$
we obtain the equality $\mathrm{(\ref{eq:u=auTotguux})}$. 
\end{proof}

\begin{para}
\label{lem:Cone(0,g,-gmufrx)Cone(iotax,kappaf)}
{\bf Lemma.}\ \ 
{\it
Let $\cC$ be a complicial exact category and let 
$f\colon x\to y$ and $g\colon \Cone f\to z$ be morphisms in $\cC$. 
Then the homotopy commutative squares 
$[f\colon x\to y]\onto{(\iota_x,\kappa_f)}[Cx\onto{\mu_f}\Cone f]\onto{(0,g,-g\mu_fr_x)}[0\to z]$ induces an equality
\begin{equation}
\label{eq:Cone(0,g,-gmufrx)Cone(iotax,kappaf)}
g=\Cone(0,g,-g\mu_fr_x)\Cone(\iota_x,\kappa_f).
\end{equation}
}
\end{para}

\begin{proof}
It follows from the following commutative diagram of exact sequences.
$$
\xymatrix{
x \ar@{>->}[r]^{\tiny{\begin{pmatrix}-f \\ \iota_x\end{pmatrix}}} \ar[d] & 
y\oplus Cx \ar@{->>}[r]^{\tiny{\begin{pmatrix}\kappa_f & \mu_f \end{pmatrix}}} 
\ar[d]_{\tiny{\begin{pmatrix}g\kappa_f & g\mu_f\end{pmatrix}}} & 
\Cone f \ar[d]^g\\
0 \ar[r] & z \ar[r]_{\id_z} & z.
}
$$
\end{proof}

\subsection{Homotopy category of relative categories}
\label{subsec:Homotopy category of relative categories}

In this subsection, we review the notion of 
homotopy categories of relative categories and 
that of derived categories of relative exact categories 
from \cite{Sch11} and \cite{Moc13b} with slightly different conventions. 
In this subsection except \ref{cor:complification}, we assume that 
underlying categories of 
relative categories are small. 

\begin{para}
\label{df:homotopy category of rel cat}
{\bf (Homotopy category of relative category).}\ 
(\cf \cite[Chapter one 1.1, 1.2]{GZ67}.)\ \ 
Let $\bC=(\cC,w)$ be a relative category. 
A {\it Gabriel-Zisman localization of $\cC$ with respect to $w$} 
is a pair $(w^{-1}\cC,Q_w)$ consisting of a category 
$w^{-1}\cC$ and a functor $Q_w\colon\cC\to w^{-1}\cC$ 
which satisfies the following two conditions:

\begin{itemize}
\item
The functor $Q_w$ sends a morphism in $w$ 
to an isomorphism in $w^{-1}\cC$.

\item
For any small category $\calD$, $Q_w$ induces an isomorphism of 
categories $\calD^{Q_w}\colon\calD^{w^{-1}\cC}\isoto \calD^{\cC}_{w-\operatorname{inv}}$ where 
$\calD^{\cC}_{w-\operatorname{inv}}$ 
is the full subcategory of $\calD^{\cC}$ 
the category of functors from $\cC$ to $\calD$ 
consisting of those functors $f\colon\cC\to\calD$ 
such that it sends all morphisms in $w$ to 
isomorphisms in $\calD$. 
\end{itemize}
We denote $w^{-1}\cC$ 
by $\Ho(\bC)$ or $\Ho(\cC;w)$ and call it the {\it homotopy category of 
a relative category $\bC$}. 
\end{para}

\begin{para}
\label{df:strongly saturation axiom}
{\bf Definition (Strongly saturation axiom).}\ \ 
Let $\bC=(\cC,w_{\bC})$ be a relative category. 
We say that $w_{\bC}$ satisfies the {\it strongly saturation axiom} 
or $\bC$ is {\it strongly saturated} 
if for any morphism $f\colon x\to y$ in $\cC$, 
$f$ is in $w_{\bC}$ if and only if $f$ is an isomorphism in $\Ho(\bC)$.
\end{para}

\begin{para}
\label{df:tricat}
{\bf (Triangulated categories).}\ \ 
We recall the conventions of triangulated categories. 
Basically we follows the notations in \cite{Kel96} and \cite{Nee01}. 
We denote a triangulated category by $(\cT,\Sigma,\Delta)$ or simply 
$(\cT,\Sigma)$ or $\cT$ 
where $\cT$ is an additive category, $\Sigma$ is an additive self category equivalence on $\cT$ 
which is said to be the {\it suspension functor} and 
$\Delta$ is a class of sequences in $\cT$ of the form 
\begin{equation}
\label{eq:sigma triangle}
x \onto{u} y \onto{v} z \onto{w} \Sigma x
\end{equation}
such that $vu=0$ and $wv=0$ which we denote by $(u,v,w)$ and call it 
a {\it distinguished triangle} and they satisfies the usual Verdier axioms 
from {\bf (TR $1$)} to {\bf (TR $4$)}. 
In the sequence $\mathrm{(\ref{eq:sigma triangle})}$, 
we sometimes write $\Cone u$ for the object $z$. 

A {\it triangle functor} between triangulated categories 
from $(\cT,\Sigma)$ 
to $(\cT',\Sigma')$ 
is a pair $(f,a)$ consisting of an additive functor 
$f\colon\cT \to \cT'$ and 
a natural equivalence 
$a\colon f\Sigma \to \Sigma'f$ such that 
they preserves distinguished triangles. 
Namely for a distinguished triangle 
$\mathrm{(\ref{eq:sigma triangle})}$ in $\cT$, 
the sequence 
$f(x)\onto{f(u)}f(y)\onto{f(v)}f(z)\onto{af(w)} \Sigma'f(x)$ 
is a distinguished triangle in $\cT'$. 
We say that a triangle functor $(f,a)$ is {\it strictly normal} 
if $f\Sigma=\Sigma' f$ and $a=\id_{f\Sigma}$. 

A {\it triangle natural transformation} 
$\theta\colon (f,a)\to (g,b)$ between triangulated functors 
$(f,a)$, $(g,b)\colon (\cT,\Sigma)\to (\cT',\Sigma')$ is 
a natural transformation $\theta\colon f\to g$ 
which satisfies the equality 
$(\Sigma'\ast\theta)\cdot a=b\cdot (\theta\ast\Sigma)$. 

We denote the category of small triangulated categories and 
triangle functors (resp. strictly normal triangle functors) by 
$\TriCat$ (resp. $\TriCat_{\stnor}$) and 
we write $\underline{\TriCat}$ (resp. $\underline{\TriCat_{\stnor}}$) 
for the $2$-category of 
small triangulated categories, triangulated functors 
(resp. strictly normal triangle functors) and 
triangle natural transformations. 

Let $(\cT,\Sigma)$ be a triangulated category. 
We say that 
a full subcategory $\calD$ of $\cT$ is 
a {\it quasi-triangulated subcategory} ({\it of $\cT$}) 
if $(\calD,\Sigma)$ is a triangulated category and the inclusion functor 
$(\iota,\id_{\Sigma}):\calD \to \cT$ is a triangulated functor. 
We say that a quasi-triangulated subcategory $\calD$ of $\cT$ is 
a {\it triangulated subcategory} ({\it of $\cT$}) if 
it is closed under isomorphisms. 
Namely an additive full subcategory $\calD$ of $\cT$ is 
a triangulated subcategory if 
$\Sigma^{\pm}(\Ob\calD)\subset \Ob\calD$ and if 
for any distinguished triangle $x \to y \to z \to \Sigma x$ 
in $\cT$, 
if $x$ and $y$ are in $\calD$, then $z$ is also in $\calD$. 
Assuming the condition $\Sigma^{\pm}(\Ob\calD)\subset \Ob\calD$, 
the last condition is equivalent to the condition that 
if two of $x$, $y$ and $z$ are in $\calD$, 
then the other one is also in $\calD$. 

We say that 
a triangulated subcategory $\calD$ of $\cT$ is {\it thick} 
if $\calD$ is closed under direct summand. 
Namely for any objects $x$ and $y$ in $\cT$, 
if $x\oplus y$ is in $\calD$, 
then both $x$ and $y$ are also in $\calD$. 
\end{para}

\begin{para}
\label{df:stable category of complicial eact categories}
{\bf (Stable category of complicial exact category).}\ \ 
Let $\cC$ be a complicial exact category. 
We write $\pi_0(\cC)$ for the quotient category of 
$\cC$ by $C$-homotopic relations. 
Namely we set $\Ob\pi_0(\cC)=\Ob\cC$ and for 
any pair of objects $x$ and $y$ in $\Ob\cC$, 
we set 
\begin{equation}
\label{eq:Hompi0C df}
\Hom_{\pi_0(\cC)}(x,y):=\Hom_{\cC}(x,y)/(\text{$C$-homotopic relation}).
\end{equation}
Since composition of morphisms in $\cC$ 
compatible with $C$-homotopic relation 
(see \cite[2.13]{Moc10}), 
it induces the compositions of morphisms in $\pi_0(\cC)$. 
Moreover since the suspension functor on $\cC$ preserves 
$C$-homotopic relations, it induces an endofunctor on 
$\pi_0(\cC)$ and we denote it by the same letter $T$. 
We call $\pi_0(\cC)$ the {\it stable category associated with $\cC$}. 
Then we can make the pair $(\pi_0(\cC),T)$ into a triangulated category by 
declaring that a sequence 
$x'\to y'\to z'\to Tx'$ is a distinguished triangle if 
it is isomorphic to a sequence of the following form
\begin{equation}
\label{eq:standard dist tri}
x\onto{f} y \onto{\kappa_f} \Cone f \onto{\psi_f} Tx.
\end{equation}

We define 
the association from $\NC(\cC_{\frob})$ 
the set of all null classes of $\cC_{\frob}$ 
to $\Tri(\pi_0(\cC))$ 
the set of all triangulated subcategories of $\pi_0(\cC)$ 
by sending a null class $\cN$ to $\pi_0(\cN)$. 
It give rises to bijective correspondences between 
$\NC(\cC_{\frob})$ and $\Tri(\pi_0(\cC))$, and 
$\NC_{\thi}(\cC)$ and $\Tri_{\thi}(\pi_0(\cC))$ the set of 
all thick subcategories of $\pi_0(\cC)$ 
(see \cite[3.15]{Moc13b}). 
\end{para}

\begin{para}
\label{ex:pi0 of FbhC and underline FBhC}
{\bf Example.}\ \ 
Let $\cC$ be a normal ordinary complicial exact category and 
let $(f,H) $, $(f,H')\colon x\to y$ be a pair of 
homotopy commutative diagrams. Assume that 
$(f,H)$ and $(f,H')$ are $C$-homotopic as morphisms in $F_{b,h}\cC$. 
Then there exists a $C$-homotopy $(0,S)\colon (f,H)\Rightarrow_C (f,H')$. 
Then for each integer $n$, $S_n$ can be regarded as a $CC$-homotopy 
from $H_n$ to $H'_n$. 
Namely the pair $(f,H)$ and $(f,H')$ are $CC$-homotopic. 
Thus the canonical functor $F_{b,h}\cC\to\underline{F_{b,h}\cC}$ 
induces an equivalence 
of triangulated categories $\pi_0(F_{b,h}\cC)\isoto \pi_0(\underline{F_{b,h}\cC})$. 
For a pair of integers $a<b$, similar statements hold 
for $F_{[a,b],h}\cC$ and so on. 
\end{para}

\begin{para}
\label{prop:hom cat of complicial relative exact category}
{\bf Proposition (Homotopy category of relative complicial exact category).}\ 
(\cf \cite[3.29]{Moc10}.)\ \ 
{\it Let $\bC=(\cC,w)$ be a relative complicial exact category. 
Then the identity functor on $\cC$ induces an equivalence of categories 
$\Ho(\bC)$ the homotopy category of $\bC$ and $\pi_0(\cC)/\pi_0(\cC^w)$ 
the Verdier quotient of $\pi_0(\cC)$ with respect to 
thick closure of the subcategory $\pi_0(\cC^w)$. 
In particular 
we can make the pair $(\Ho(\cC),T)$ into a triangulated category 
to make the equivalence above an equivalence of triangulated categories. 
The association $\Ho\colon\underline{\RelComp}\to\underline{\TriCat}$ is 
a $2$-functor and it induces a $2$-functor 
$\underline{\RelComp_{\stnor}}\to\underline{\TriCat_{\stnor}}$.
}\qed
\end{para}

\begin{para}
\label{df:Q_C df}
{\bf Definition.}\ \ 
For a relative complicial exact category $\bC=(\cC,w)$, 
we denote the canonical functor $\cC\to \Ho(\bC)$ by $\cQ_{\bC}$.
\end{para}

Recall the definition of thick \ref{df:comp weak equiv} and 
strongly saturated \ref{df:strongly saturation axiom} 
relative complicial exact categories.

\begin{para}
\label{lem:thick is strongly saturated}
{\bf Lemma.}\ \ (\cf \cite[3.2.18]{Sch11}.)\ \ 
{\it
Let $\bC=(\cC,w_{\bC})$ is a relative complicial exact category.
If $\bC$ is thick, then $\bC$ is strongly saturated. 
}
\end{para}

\begin{proof}
Let $f\colon x\to y$ be a morphism in $\cC$ such that 
its image in $\Ho(\bC)$ is an isomorphism. 
What we need to prove is that $f$ is in $w_{\bC}$. 
Then $\Cone f$ is trivial in $\Ho(\bC)$ and 
$\pi_0(\cC^{w_{\bC}})$ is thick by assumption and 
\ref{nt:Verdier correspondence}. 
Moreover $\cC^{w_{\bC}}$ is closed under 
$C$-homotopy equivalences by \ref{lem:w_bC-closed null class}. 
Thus $\Cone f$ is in $\cC^{w_{\bC}}$ and 
it means $f$ is in $w_{\bC}$ by 
$\mathrm{(\ref{eq:Verdier correspondence 2})}$. 
Thus we obtain the result. 
\end{proof}

\begin{para}
\label{lem:HS85}
{\bf Lemma.}\ (\cf \cite{HS85}, \cite[1.9.2]{TT90}, \cite[3.16]{Moc13b}.)\ \ 
{\it
Let $\bC=(\cC,w)$ be a relative complicial exact category. 
Then $\id_{\cC}$ the identity functor of $\cC$ 
induces the following equivalences of triangulated categories.
\begin{equation}
\label{eq:HS85}
\Ho(\cC_{\frob};w)\isoto\Ho(\cC_{\frob};{\langle w\rangle}_{\thi})\isoto\Ho(\cC;w)\isoto\Ho(\cC;{\langle w\rangle}_{\thi}).
\end{equation}
}\qed
\end{para}

\begin{para}
\label{df:relative homotopy category}
{\bf Definition (Relative homotopy category).}\ \ 
Let $\bC=(\cC,w)$ be a relative complicial exact category 
and let $\cF$ be a full subcategory of $\cC$. 
We regard ${\langle \cF\rangle}_{\nullclass,\frob}$ as a complicial exact subcategory of $\cC_{\frob}$ and 
write $\Ho(\inn{\cF}{\bC})$ for the homotopy category of 
the relative exact category 
$({\langle \cF\rangle}_{\nullclass,\frob},
w|_{{\langle \cF\rangle}_{\nullclass,\frob}})$ 
and call 
it a {\it relative homotopy category of $\cF$} ({\it in $\bC$.}) 

Assume that $w$ is thick. 
We say that a {\it null class $\cN$ of $\cC_{\frob}$ compatible with $w$} 
if for any morphism $u\colon x\to y$ from an object $x$ in $\cN$ 
to an object $y$ in $\cC^w$, there exists a pair of 
morphisms $u'\colon x\to z$ and $u''\colon z\to y$ with 
$z\in\Ob(\cN\cap\cC^w)$ and a $C$-homotopy $H\colon u\Rightarrow_C u''u'$. 
We say that {\it a full subcategory $\cF$ is compatible with $w$} 
if the null class ${\langle \cF\rangle}_{\nullclass,\frob} $ of $\cC_{\frob}$ is compatible with $w$. 
In this case the inclusion functor 
${\langle \cF\rangle}_{\nullclass,\frob}\to \cC$ 
induces a fully faithful 
functor $\Ho(\inn{\cF}{\bC})\to\Ho(\bC)$ and the inclusion functor 
${\langle\cF\rangle}_{\nullclass,\frob}\rinc{\langle\cF \rangle}_{\nullclass,w}$ induces an equivalence of triangulated categories
\begin{equation}
\label{eq:w-envelope invariance}
\Ho(\inn{\cF}{\bC})\isoto\Ho(\inn{{\langle\cF\rangle}_{\nullclass,w}}{\bC})
\end{equation}
by \cite[10.3]{Kel96}.
\end{para}

\begin{para}
\label{df:derived cat}
{\bf (Derived category of relative exact category).}\ \ 
Let $\bE=(\cE,w)$ be a relative exact category. 
Then we write $\calD_b(\bE)$ or $\calD_b(\cE;w)$ for 
the homotopy category of 
$\Ch_b(\bE)$ the relative category of bounded chain complexes on $\cE$ 
and we call it the ({\it bounded}) {\it derived category of $\bE$}. 
By virtue of 
Proposition~\ref{prop:hom cat of complicial relative exact category}, 
$\calD_b(\bE)$ naturally becomes a triangulated category. 
The association $\calD_b\colon\underline{\RelEx}\to\underline{\TriCat_{\stnor}}$ 
is a $2$-functor. 
\end{para}

\begin{para}
\label{df:derived equivalence}
{\bf Definition (Derived equivalence).}\ \ 
We say that a relative exact functor $f\colon\bE \to \bE'$ 
between relative exact categories $\bE$ and $\bE'$ is 
a {\it derived equivalence} if it induces an 
equivalence of triangulated categories $\calD_b(f)\colon\calD_b(\bE)\isoto 
\calD_b(\bE')$ 
on bounded derived categories. 
We denote the class of derived equivalences in $\RelEx$ by 
$\deq_{\RelEx}$ or simply $\deq$. 
Since the association $\calD_b$ is $2$-functorial and 
it sends a natural weak equivalence to a triangle natural equivalence, 
it turns out that a categorical homotopy equivalence in $\underline{\RelEx}$ 
is a derived equivalence. 

By virtue of Proposition~\ref{prop:comp hom and der} below, for a 
relative complicial exact functor $(f,c)\colon\bC\to\bC'$, the following 
two conditions are equivalent. 
\begin{itemize}
\item
$(f,c)$ induces an equivalence of triangulated categories 
$\Ho(\bC)\isoto\Ho(\bC')$.
\item
The relative functor $f\colon\bC\to\bC'$ between relative exact categories 
is a derived equivalence in $\RelEx$. 
\end{itemize}
In this case, we say that $(f,c)$ is {\it derived equivalence} and 
we denote the class of derived equivalences in $\RelComp$ by 
$\deq_{\RelComp}$ or simply $\deq$.
\end{para}

\begin{para}
\label{rem:derived equivalence}
{\bf Remark.}\ \ 
Let $(f,c)\colon\bC=(\cC,w)\to \bC'=(\cC',w')$ be a 
relative complicial exact functor between relative 
complicial exact categories and let $\cF$ be 
a full subcategory of $\cC$ and let $\bbD$ be a $T$-system of $\cC$ 
(for definition of $T$-system, see \ref{df:properties of family}). 
Then
\begin{enumerate}
\enumidef
\item
By Theorem~1.5 in \cite{BM11}, 
if the following two conditions hold, then 
$f$ is a derived equivalence.

\begin{itemize}
\item[\bf (App 1).]
For a morphism $u\colon x\to y$ in $\cC$, 
$u$ is in $w$ if and only if $f(u)$ is in $w'$.

\item[\bf (App 2).]
For an object $x$ in $\cC$ and 
an object $y$ in $\cC'$ and a morphism 
$u\colon f(x)\to y$, 
there exists an object $z$ in $\cC$ and a morphism 
$s\colon x\to z$ in $\cC$ and a morphism 
$v\colon f(z)\to y$ in $w'$ such that $u=vf(s)$. 
\end{itemize}
In {\bf (App 2)}, we call the triple $(z,s,v)$ 
a {\it factorization of $u$} ({\it along $f$}).

\item
In the following cases, 
the condition {\bf (App 1)} is automatically verified. 
Here $\cF$ is a strict exact subcategory of $\cC_{\frob}$ 
or a prenull class of $\cC$. 

\begin{enumerate}
\item
The relative complicial exact functor 
$\Tot\colon (\Ch_b(\cF),tw)\to ({\langle\cF \rangle}_{\nullclass,\frob},w|_{{\langle\cF \rangle}_{\nullclass,\frob}})$. 
(See \ref{lemdf:total functor}, \ref{lemdf:tw=qw}.)

\item
The relative complicial exact functor 
${(-)}_{\infty}\colon(F_{b,h}^{\bbD|_{\cF}}{\langle\cF \rangle}_{\nullclass,\frob},w_{\st})\to ({\langle\cF \rangle}_{\nullclass,\frob},w|_{{\langle\cF \rangle}_{\nullclass,\frob}})$ 
(see \ref{df:stable weq} and for the definition of the category 
$F_{b,h}^{\bbD|_{\cF}}{\langle\cF \rangle}_{\nullclass,\frob}$ 
of d\'evissage filtrations, see \ref{df:devissage filtrations}).
\end{enumerate}

\item
In condition {\bf (App 2)}, 
by replacing $s$ with ${\xi_1}_s$, $z$ with $\Cyl s$ 
and $u$ with $uf(\upsilon_s)$, 
we shall assume that $s$ is an admissible monomorphism.

\item 
Assume that $\bC$ and $\bC'$ are thick. 
Then by Lemma~\ref{lem:thick is strongly saturated}, 
they are strongly saturated 
(see \ref{df:strongly saturation axiom}). 
Since relative complicial exact categories are derivable 
Waldhausen categories by \ref{para:derivable Waldhausen category}, 
we can replace condition {\bf (App 2)} with the more 
weaker condition {\bf (App 2)'} below by 
\cite[Th{\'e}or{\`e}me 2.9.]{Cis10b}.

\sn
{\bf (App 2)'} 
For an object $x$ in $\cC$ and an object $y$ in $\cC'$ and 
a morphism $u\colon f(x)\to y$, 
there exists a morphism $b\colon y\to y'$ in $w'$ such that 
$bu\colon f(x)\to y'$ admits a factorization $(z,s,v)$.

\item
We can replace condition {\bf (App 2)} with the following 
more weaker condition {\bf (App 2)''}.

\sn
{\bf (App 2)''} 
For an object $x$ in $\cC$ and an object $y$ in $\cC'$ and 
a morphism $u\colon f(x)\to y$, 
there exists a 
{\it homotopy factorization of $u$}. 
Namely there exists an object $z$ in $\cC$ and a morphism 
$s\colon x\to z$ in $\cC$ and a morphism 
$v\colon f(z)\to y$ in $w'$ and a $C$-homotopy 
$H\colon u\Rightarrow_C vf(s)$. 

\sn
In {\bf (App 2)''}, 
we call the quadruple $(z,s,v,H)$ 
a {\it homotopy factorization of $u$} ({\it along $f$}).
\end{enumerate}
\end{para}

\begin{proof}
If we have a homotopy factorization 
$s\colon x\to z$, $v\colon f(z)\to y$ 
and $H\colon u\Rightarrow_C vf(s)$ of $u\colon f(x) \to y$, 
then by setting $z':=z\oplus Cx$ and 
$\displaystyle{s':=\begin{pmatrix}s\\ -\iota_x\end{pmatrix}}$ 
and $v':=\begin{pmatrix}a & -Hc_x^{-1}\end{pmatrix}$ where 
$c_x^{-1}\colon fCx\isoto Cfx$ is the canonical isomorphism, 
we obtain a factorization $f(x)\onto{f(s')} f(z')\onto{v'}y$ of $u$.
\end{proof}

\begin{para}
\label{lem:App2 test}
{\bf Lemma.}\ \ 
{\it
Let $(f,c)\colon\bE=(\cE,w_{\bE})\to \bC=(\cC,w_{\bC})$ 
be a relative complicial exact functor between 
relative complicial exact categories and let 
$\cF$ be an additive full subcategory of $\cC$. 
Assume that
\begin{enumerate}
\enumidef
\item
$\cC$ is strictly ordinary and
\item
for any object $x$ in $\cE$, $f(x)$ is in ${\langle\cF\rangle}_{\nullclass,\frob}$ and
\item
for any object $x$ in $\cE$ and any object $y$ in $\cF$, 
a morphism $u\colon f(x)\to y$ admits a factorization 
$f(x)\onto{f(v)}f(z)\onto{a}y$ with $a\in w_{\bC}$. 
\end{enumerate}
Then $(f,c)\colon\bE\to ({\langle\cF\rangle}_{\nullclass,\frob},w_{\bC}|_{{\langle\cF\rangle}_{\nullclass,\frob}})$ satisfies} {\bf (App 2)''}.
\end{para}

\begin{proof}
Let $x$ be an object in $\cE$ and let $y$ be an object 
in ${\langle\cF\rangle}_{\nullclass,\frob}$ 
and let $u\colon f(x)\to y$ be a morphism in $\cC$. 
We shall show that $u$ admits a homotopy factorization. 
Since ${\langle\cF\rangle}_{\nullclass,\frob}$ 
is a $C$-homotopy equivalences closure of 
${\langle\cF\rangle}_{\seminull}$ 
by Lemma~\ref{lem:null class} $\mathrm{(2)}$, 
there exists an object $y'$ in ${\langle\cF\rangle}_{\seminull}$ and 
a $C$-homotopy equivalences $a'\colon y'\to y$ and $b'\colon y'\to y$ 
and a $C$-homotopy $H'\colon\id_{y'}\Rightarrow_C b'a'$. 
If the composition $a'u\colon f(x)\to y'$ 
admits a homotopy factorization $(z,c,a,H)$, 
then the quadruple $(z,c,b'a,b'H+H'Ca)$ 
gives a homotopy factorization of $u$ along $f$. 
Therefore we shall assume that $y$ is in ${\langle\cF\rangle}_{\seminull}$. 
Moreover we shall assume that $y$ is in ${\langle\cF\rangle}_{\prenull}$ 
by replacing $f(x)\onto{u}y$ with 
$f(T^nx)\overset{c^{-1}_{T^n}}{\isoto} T^nf(x)\onto{T^n u}T^ny$ 
for suitable positive integer $n$. 
Here $c_{T^n}\colon T^n f\isoto fT^n$ is compositions of 
$T^nf\overset{T^{n-1}\ast c_T}{\isoto}T^{n-1}fT\overset{T^{n-2}\ast c_T\ast T}{\isoto}\cdots\overset{T\ast c_T\ast T^{n-2}}{\isoto}TfT^{n-1}\overset{c_T\ast T^{n-1}}{\isoto}fT^n$.
Recall the definition of $\cP_m(\cF)$ from 
Lemma-Definition~\ref{lemdf:prenull closure}. 
There exists a non-negative integer $m$ such that $y$ is in $\cP_m(\cF)$. 
If $m=0$, then $f(x)\to y$ admits a factorization by 
assumption $\mathrm{(3)}$. 
Assume that $m\geq 1$ and we proceed by induction on $m$. 
Namely, assuming that if $y$ is in $\cP_{m-1}(\cF)$, then $u$ admits 
a factorization, we will prove that if $y$ is in $\cP_m(\cF)$, 
then $u$ admits a factorization. 
By definition of $\cP_m(\cF)$, there exists a morphism 
$d\colon y_1\to y_0$ in $\cP_{m-1}(\cF)$ such that $y=\Cone d$. 
Applying an inductive hypothesis to the composition 
$f(x)\onto{u} y\onto{\psi_d} Ty_1$, 
we obtain a factorization $(s,v',a')$ of $\psi_du\colon f(x)\to Ty_1$. 
We apply Lemma~\ref{lem:replacement of factorizations} $\mathrm{(2)}$ 
below to the composition 
$f(\Cone v')\onto{\Cone(u,a')}\Cone\psi_d
\onto{{\tiny{\begin{pmatrix}\id_{Ty_0}& 0\end{pmatrix}}}}Ty_0$, 
we obtain a factorization $(t,v'',b')$ of $f(\Cone v')\to Ty_0$ 
by inductive hypothesis. 
We set $z:=\Cone v''\kappa_{v'}$ and 
$a:=\sigma_T^{\Cone}\cdot\Cone(a',b')\cdot d^{\Cone}_{v''\kappa_{v'}}\colon f(z)\to Ty$ 
and $v:=\Cone(v'',0,v''\mu_{v'})\colon Tx\to z$. 
Then $a$ is in $w_{\bC}$ by 
Proposition~\ref{prop:properties of weak equivalences} 
$\mathrm{(2)}$ $\mathrm{(ii)}$ and 
$a\cdot f(v)=Tu\cdot {(c^{-1}_T)}_x$ by 
Lemma~\ref{lem:Cone(psif,0,Thetaf)} and 
Lemma~\ref{lem:ThetafCg} below. 
Thus we obtain a factorization of $u$ by 
Lemma~\ref{lem:replacement of factorizations} $\mathrm{(1)}$ below. 
We complete the proof.
\end{proof}

\begin{para}
\label{lem:ThetafCg}
{\bf Lemma.}\ \ 
{\it
Let $\cC$ be a strictly ordinary complicial exact category and let 
$f\colon x\to y$, $g\colon z\to \Cone f$, $h\colon z\to w$ 
and $k\colon w\to Tx$ be morphisms in $\cC$ such that 
$\psi_f g=kh$. 
Then we have an equality
\begin{equation}
\label{eq:=ThetafCg}
-\begin{pmatrix}\id_{Ty} & 0 \end{pmatrix}\Cone(g,k)\mu_h=\Theta_f Cg.
\end{equation}
Namely we have a commutative diagram of homotopy squares
\begin{equation}
\label{eq:comm diagram of homotopy squares}
\xymatrix{
[z\to 0] \ar[r]^{(h,0,\mu_h)} \ar[d]_{(g,0,0)} & 
[w\onto{\kappa_h}\Cone h] 
\ar[d]^{(k,-{\tiny{\begin{pmatrix}\id_{Ty}& 0\end{pmatrix}}}\Cone(g,k),0)}\\
[\Cone f \to 0] \ar[r]_{(\psi_f,0,\Theta_f)} & [Tx\onto{Tf}Ty].
}
\end{equation}
In particular we have an equality
\begin{equation}
\label{eq:Cone(psif,0,Thetaf)Tg=}
\Cone(\psi_f,0,\Theta_f)Tg=\Cone(k,-\begin{pmatrix}\id_{Ty}& 0\end{pmatrix}\Cone(g,k))\Cone(h,0,\mu_h).
\end{equation}
}
\end{para}

\begin{proof}
By definition of $\Cone(g,k)$, the following diagram is commutative and 
it implies the equality $\mathrm{(\ref{eq:=ThetafCg})}$. 
$$
\xymatrix{
z \ar@{>->}[r]^{{\tiny{\begin{pmatrix}-h\\ \iota_z\end{pmatrix}}}} \ar[d]_g & 
w\oplus Cz \ar[d]^{{\tiny{\begin{pmatrix}k & 0\\ 0 & Cg\end{pmatrix}}}} 
\ar@{->>}[r]^{{\tiny{\begin{pmatrix}\kappa_h & \mu_h\end{pmatrix}}}} & 
\Cone h \ar[d]^{\Cone(g,k)}\\
\Cone f \ar@{>->}[r]_{{\tiny{\begin{pmatrix}-\psi_f\\ \iota_{\Cone f}\end{pmatrix}}}} & Tx\oplus C\Cone f 
\ar@{->>}[r]_{{\tiny{\begin{pmatrix}-Tf & -\Theta_f\\ \iota_{Tx} & C\psi_f
\end{pmatrix}}}} & Ty\oplus CTx.
}
$$
\end{proof}

\begin{para}
\label{lem:replacement of factorizations}
{\bf Lemma.}\ \ 
{\it
Let $(f,c)\colon\bC=(\cC,w_{\bC})\to\bC'=(\cC',w_{\bC'})$ 
be a relative complicial exact functor and let $x$ and $y$ be 
objects in $\cC$ and $\cC'$ respectively. Then

\sn
$\mathrm{(1)}$ 
For a morphism $u\colon f(x)\to Ty$, if 
the composition $F(T^{-1}x)\overset{c^{-1}_{T^{-1}}}{\isoto}T^{-1}f(x)\onto{T^{-1}u}T^{-1}Ty\overset{\beta^{-1}_y}{\isoto}y$ admits a factorization $(z,s,a)$, 
then the triple $(Tz,Ts\cdot\alpha^{-1}_x,Ta\cdot{(c^{-1}_T)}_z)$ 
is a factorization of $u$. 

\sn
$\mathrm{(2)}$ 
For a morphism $v\colon f(x)\to y$, if 
the composition $f(Tx)\overset{{(c_T^{-1})}_x}{\isoto} Tf(x)\onto{Tu} Ty$ 
admits a factorization $(z,s,a)$, 
then the triple $(T^{-1}z,T^{-1}s\cdot\beta_y,
\beta_y^{-1}\cdot T^{-1}a\cdot{(c^{-1}_{T^{-1}})}_z)$ 
is a factorization of $v$.
}
\end{para}

\begin{proof}
Assertion $\mathrm{(1)}$ (resp. $\mathrm{(2)}$) follows from the 
left (resp. right) commutative diagram below 
by Lemma~\ref{lem:comp func and alpha beta}.
$$
\tiny{
\xymatrix{
& & f(x) \ar[lld]^{\sim}_{f\ast\alpha^{-1}_x} \ar[d]_{\wr}^{\alpha^{-1}\ast f(x)}  \ar[r]^u & Ty \ar[d]^{\alpha^{-1}\ast Ty}_{\wr}\\
f(TT^{-1}x) \ar[r]^{\sim}_{{(c^{-1}_{T})}_{T^{-1}x}} \ar[d]_{f(Ts)} & 
Tf(T^{-1}x) \ar[r]^{\sim}_{T{(c^{-1}_{T^{-1}})}_x} \ar[rd]_{Tf(s)} &
T^{-1}f(x) \ar[r]^{TT^{-1}u} & TT^{-1}Ty \ar[d]_{\wr}^{T\ast\beta_y^{-1}}\\
f(Tz) \ar[rr]^{\sim}_{{(c^{-1}_T)}_z} & & Tf(z) \ar[r]_{Ta} & Ty,
}
\xymatrix{
& & f(x) \ar[r]^u \ar[lld]^{\sim}_{f\ast\beta_y} \ar[d]_{\wr}^{\beta\ast f(x)} &
y \ar[d]_{\wr}^{\beta_y}\\
f(T^{-1}Tx) \ar[r]^{\sim}_{c^{-1}_T\ast Tx} \ar[d]_{f(T^{-1}s)} & 
T^{-1}f(Tx) \ar[r]^{\sim}_{T^{-1}{(c^{-1}_T)}_x} \ar[dr]_{T^{-1}f(s)} & 
T^{-1}Tf(x) \ar[r]_{T^{-1}Tu} & T^{-1}Ty\\
f(T^{-1}z) \ar[rr]^{\sim}_{{(c^{-1}_{T^{-1}})}_z} & & T^{-1}f(z). \ar[ru]_{T^{-1}a} 
}}
$$
\end{proof}

Recall the definition of 
the zig-zag sequence of morphisms $\cU_x^{[a,b]}$ and 
the functor $\fL^{[a,b]}$ 
from \ref{lemdf:uk and vk}.

\begin{para}
\label{lem:TotcU[a,b]x}
{\bf Lemma.}\ \ 
{\it
Let $\bC=(\cC,w)$ be a relative complicial exact category and 
let $\cP$ be a prenull class of $\cC$ and let 
$a<b$ be a pair of integers and let 
$x$ be a complex in $\Ch_{b,\leq b}(\cP)$. 
Then we can regard the zig-zag sequence $\cU_x^{[a,b]}$ 
as an isomorphism $x\isoto \fL^{[a,b]}x$ in $\calD_b(\cP,w|_{\cP})$. 
Moreover we have the equality
\begin{equation}
\label{eq:TotcU[a,b]x}
\Tot\cU_x^{[a,b]}=\id_{\Tot x}
\end{equation}
as a morphism in $\Ho(\inn{\cP}{\bC})$.
}
\end{para}

\begin{proof}
Since $\cU_x^{[a,b]}$ is the zig-zag sequence of morphisms 
in $lw\cup\qfis$, it can be regarded as an isomorphism in 
$\calD_b(\cP,w|_{\cP})$. 
To show the equality $\mathrm{(\ref{eq:TotcU[a,b]x})}$, 
by \ref{dfcor:Tot associativity} and 
the equality $\mathrm{(\ref{eq:calU degree shift})}$, 
by replacing $x$ with $(\sigma_{b-1}x)[-b+1]$, 
we shall assume that $a=0$ and $b=1$ and $x$ is in $\Ch_{[0,1]}(\cC)$. 
In this case, the equality follows form the equality 
$\mathrm{(\ref{eq:fundamental triangle of cone})}$. 
$$
\xymatrix{
\Tot x \ar[rrr]^{\!\!\!\!\!\!\!\!\!\!\!\!\!\!\!\!\!\!\!\!\Tot\fu_1(x)=\Cone(\iota_{x_1},\kappa_{d_1^x})} \ar[rrrd]_{\id_{\Tot x}} & & & \Tot\fQ_1(x)=\Cone\mu_{d_1^x} 
\ar[d]^{\!\!\!\!\!\scriptscriptstyle{\Cone(0,\id_{\Cone d_1^x},-\mu_{d_1^x}r_{x_1})}} 
& & &
\Tot\fL_1(x)=\Tot x \ar[llld]^{\id_{\Tot x}}
\ar[lll]_{\ \ \ \ \ \ \ \ \ \Tot \fv_1(x)=\Cone(0,\id_{\Cone f})}
\\
& & & \Tot x.
}$$
\end{proof}

\begin{para}
\label{prop:comp hom and der}
{\bf Proposition (Comparison of homotopy and derived categories).}\ 
(\cf \cite[4.15]{Moc13b}.)\ \ 
{\it
Let $\bC=(\cC,w)$ be a relative complicial exact category and 
let $\cP$ be a full additive subcategory of $\cC$. Then

\sn
$\mathrm{(1)}$ 
If $\cP$ is a prenull class, then 
the triangulated structures on $\Ho(\Ch_b(\cP);q(w|_{\cP}))$ 
induced from the standard and the level complicial 
structures on $\Ch_b(\cP)$ are equivalence.

\sn
$\mathrm{(2)}$
Assume that either condition $\mathrm{(i)}$ or $\mathrm{(ii)}$ below

\sn
$\mathrm{(i)}$ 
$\cP$ is a prenull class or 

\sn
$\mathrm{(ii)}$ 
$\cC$ is strictly ordinary and $\cP$ contains all $C$-contractible objects 
and $\cP$ is closed under the operations $T^{\pm1}$.

\sn
Then the functor 
$\Tot\colon\Ch_b(\cP)\to{\langle \cP\rangle}_{\nullclass,\frob}$ 
induces an equivalence of 
triangulated categories 
\begin{equation}
\label{eq:comparison of Ho}
\Ho(\Ch_b(\cP);q(w|_{\cP}))\isoto\Ho(\inn{\cP}{\bC}).
\end{equation}

\sn
$\mathrm{(3)}$ 
The canonical functor $j_{\cC}\colon\bC \to \Ch_b(\bC)$ 
which sends an object $x$ to a complex $j_{\bC}(x)$ 
such that ${j_{\bC}(x)}_k$ is $x$ if $k=0$ and $0$ if $k\neq 0$ 
induces an equivalence of triangulated categories 
$\Ho(\bC)\isoto \calD_b(\bC)$. 

\sn
$\mathrm{(4)}$ 
If $\cP$ is a strict exact subcategory of $\cC$, then 
the identity functor of $\Ch_b\cF$ 
induces an equivalence of triangulated categories 
$\calD_b(\cF;w|_{\cF})(=\Ho(\Ch_b(\cF);q(w|_{\cF})))\isoto
\Ho(\inn{\Ch_b\cF}{\Ch_b\bC})$.
}
\end{para}

\begin{proof}
$\mathrm{(1)}$ 
By \ref{lemdf:uk and vk} and the equality $\mathrm{(\ref{eq:L[a,b]x})}$, 
the suspension functors on $\Ho(\Ch_b(\cP);q(w|_{\cP}))$ induced from 
the standard and the level complicial 
structures on $\Ch_b(\cP)$ 
are equivalence of functors. 
By \ref{lem:HS85}, 
the class of distinguished triangles on $\Ho(\Ch_b(\cP);q(w|_{\cP}))$ 
only depends upon the Frobenius exact structure on $\Ch_b(\cP)$. 
Thus we obtain the result.

\sn
$\mathrm{(2)}$ 
First we assume that $\cP$ is a prenull class. 
By \ref{cor:qw=qwthi} and \ref{lem:HS85}, replacing $w$ with ${\langle w\rangle}_{\thi,\cC}$, we shall assume that $w$ is thick. 
We will show that the functor $\Ho(\Tot)$ is essentially surjective. 
Let $x$ be an object in ${\langle \cP \rangle}_{\seminull,\cC}$, 
then there exists an integer $m$ and an object $y$ in $\cP$ such that 
$T^my$ is isomorphic to $x$. 
Then by the equality $\mathrm{(\ref{eq:Totj=id})}$, 
$\Tot j_{\cP}(y)[m]=T^my$ is isomorphic to $x$. 
Since the inclusion functor $\pi_0({\langle \cP \rangle}_{\seminull,\cC})\to 
\pi({\langle \cP \rangle}_{\nullclass,\frob})$ 
is an equivalence of triangulated categories, 
we obtain the result. 

Next for a pair of chain complexes $x$ and $y$ in $\Ch_b(\cP)$, 
we will show that the map 
$$\Tot\colon\Hom_{\Ho(\Ch_b(\cP);q(w|_{\cP}))}(x,y)\to
\Hom_{\Ho(\inn{\cP}{\bC})}(\Tot x,\Tot y)$$
is surjective. 
Since the suspension functors on both $\Ho(\Ch_b(\cP);q(w|_{\cP}))$ and 
$\Ho(\inn{\cP}{\bC})$ are equivalences, 
by replacing $x$ and $y$ with $x[N]$ and $y[N]$ for sufficiently large 
integer $N$ respectively, we shall assume that $x$ and $y$ are 
in $\Ch_{b,\geq 0}(\cP)$. 
Next by \ref{lemdf:uk and vk}, replacing $x$ and $y$ with 
$j_{\cP}(\Tot x)[m]$ and 
$j_{\cP}(\Tot y)[m']$ for 
a suitable integers $m$ and $m'$, we shall assume that 
$x$ and $y$ are in $\cP$. 
In this case, surjectivity comes from commutative of the diagram below and 
the commutativity follows from \ref{lem:TotcU[a,b]x}. 
$$
\xymatrix{
\Hom_{\Ho(\Ch_b(\cP);q(w|_{\cP}))}(j_{\cP}(T^mx),j_{\cP}(T^{m'}y)) 
\ar[d]_{\wr} & 
\Hom_{\Ho(\inn{\cP}{\bC})}(T^mx,T^{m'}y) 
\ar@{=}[d] \ar[l]_{\ \ \ \ \ \ \ \ \ \ \ \ \ \ \Ho(j_{\cP})}\\
\Hom_{\Ho(\Ch_b(\cP);q(w|_{\cP}))}(j_{\cP}(x)[m],j_{\cP}(y)[m']) 
\ar[r]_{\ \ \ \ \ \ \ \ \ \ \ \ \ \ \Ho(\Tot)} & 
\Hom_{\Ho(\inn{\cP}{\bC})}(T^mx,T^{m'}y).
}
$$
By the equality $\mathrm{(\ref{eq:tw=qw})}$, 
$\Ho(\Tot)$ has the trivial kernel and 
thus $\Ho(\Tot)$ is equivalence of triangulated categories by 
\cite[3.18]{Bal07}. 

Next we assume that $\cC$ is strictly ordinary and $\cP$ contains 
all $C$-contractible objects and $\cP$ is closed under the operations 
$T^{\pm 1}$. 
Let $x$ be an object in $\Ch_b(\cP)$ and let $y$ be an object $y$ in $\cP$ 
and let $u\colon\Tot x\to y$ be a morphism in $\cP$. 
We prove that $u$ admits a factorization of $u$ along $\Tot$. 
By replacing $u\colon\Tot x\to y$ with $T^mu\colon \Tot x[m] \to T^my$ 
for a suitable integer $m$, 
we shall assume that $x$ is in $\Ch_{[0,n]}(\cP)$ for 
some non-negative integer $n$. 
Then by Lemma-Definition~\ref{df:canonical factorization}, 
the triple $(z_u,g_u\fu_x,a_u)$ gives a factorization of $u$ along $\Tot$. 
Notice that by construction and assumption $z_u$ is in $\Ch_{[0,n]}(\cP)$. 
Thus by Lemma~\ref{lem:App2 test}, 
$\Tot\colon(\Ch_b\cP,q(w|_{\cP}))\to 
({\langle\cP \rangle}_{\nullclass,\frob},w|_{{\langle\cP \rangle}_{\nullclass,\frob}})$ 
satisfies {\bf (App 2)'} and we obtain the equivalence 
$\mathrm{(\ref{eq:comparison of Ho})}$ by 
Remark~\ref{rem:derived equivalence} $\mathrm{(1)}$.

\sn
$\mathrm{(3)}$ 
By the equality $\mathrm{(\ref{eq:Totj=id})}$, 
it turns out that $\calD_b(j_{\cC})$ is the right inverse functor of 
$\calD_b(\Tot)$ which is an equivalence of triangulated categories 
$\calD_b(\bC)\isoto\Ho(\bC)$ by $\mathrm{(2)}$. 
Thus $\calD_b(j_{\cC})$ is also an equivalence of triangulated categories.

\sn
$\mathrm{(4)}$ 
Since $\Ch_b\cF$ is a prenull class in $\Ch_b\cC$, relative complicial 
exact functors 
$$(\Ch_b\cF,q(w|_{\cF}))\onto{j_{\Ch_b\cF}}(\Ch_b\Ch_b\cF,qq(w|_{\cF}))\onto{\Tot}({\langle\Ch_b\cF \rangle}_{\nullclass,\frob},q(w|_{{\langle\Ch_b\cF\rangle}_{\nullclass,\frob}}))$$
induce equivalences of triangulated categories 
$$
\Ho(\Ch_b\cF;q(w|_{\cF}))\isoto\Ho(\Ch_b\Ch_b\cF;qq(w|_{\cF}))\isoto\Ho({\langle\Ch_b\cF \rangle}_{\nullclass,\frob};q(w|_{{\langle\Ch_b\cF\rangle}_{\nullclass,\frob}}))
$$
by $\mathrm{(2)}$ and $\mathrm{(3)}$.
\end{proof}

\begin{para}
\label{cor:complification}
{\bf Corollary.}\ \ 
{\it
For the (large) relative categories $(\RelEx,\deq_{\RelEx})$ and 
$(\RelComp,\deq_{\RelComp})$, 
the adjoint pair $(\RelEx,\deq_{\RelEx})\rlto{\Ch_b}{F}(\RelComp,\deq_{\RelComp})$ with the adjunction $(j,\Tot)$ gives 
categorical homotopy equivalences. 
}
\qed
\end{para}

\begin{para}
\label{cor:w-closure envelope}
{\bf Corollary.}\ \ 
{\it
Let $\bC=(\cC,w)$ be a relative complicial exact category 
and let $\cA$ be an additive full subcategory of $\cC$. Then

\sn
$\mathrm{(1)}$ 
Assume that $\cA$ contains all $C$-contractible objects. 
Then $\Ch_b\cA$ is compatible with $qw$ 
{\rm(see \ref{df:relative homotopy category})}.

\sn
$\mathrm{(2)}$ 
Assume that $\cA$ is a strict exact subcategory of $\cC$. 
Then the inclusion functor
$${\langle\Ch_b\cA\rangle}_{\nullclass,\frob,\Ch_b\cC}\rinc{\langle\Ch_b\cA\rangle}_{\nullclass,qw,\Ch_b\cC}$$
induces an equivalence of categories
\begin{equation}
\label{eq:ChbA qw-envelope invariance}
\Ho(\inn{\Ch_b\cA}{\Ch_b\bC})\isoto\Ho(\inn{{\langle\Ch_b\cA\rangle}_{\nullclass,qw,\Ch_b\cC}}{\Ch_b\bC}).
\end{equation}
}
\end{para}

\begin{proof}
$\mathrm{(1)}$ 
First notice that $\cA$ has enough objects to resolve in $\cC$ in 
the following sense:

\sn
For an object $y$ in $\cC$, there exists an admissible epimorphism 
$P(y)\overset{q_y}{\rdef}y$ with $P(y)\in\Ob\cA$.

\sn
Thus by \cite[1.9.5]{TT90}, for any morphism $f\colon x\to y$ with 
$x\in\Ob\Ch_b\cA$ and $y\in\Ob{(\Ch_b\cC)}^{qw}$, 
there exists an object $z\in\Ob\Ch_b\cA$ and an admissible monomorphism 
$x\overset{i}{\rinf}z$ and a quasi-isomorphism $z\onto{a}y$ 
such that $f=ai$. 
Then since ${(\Ch_b\cC)}^{qw}$ is 
$\qis$-closed by Lemma~\ref{lem:w_bC-closed null class}, 
$z$ is in ${(\Ch_b\cC)}^{qw}$. 
Thus $\Ch_b\cA$ is compatible with $qw$. 

\sn
$\mathrm{(2)}$ 
First notice that ${\langle\Ch_b\cA\rangle}_{\nullclass,\frob}={\langle{(\Ch_b\cA)}_{\heq}\rangle}_{\nullclass,\frob}$. 
Thus the inclusion functor $\Ch_b\cA\rinc{(\Ch_b\cA)}_{\heq}$ 
induces an equivalence of triangulated categories 
$\Ho(\inn{\Ch_b\cA}{\Ch_b\bC})\isoto\Ho(\inn{{(\Ch_b\cA)}_{\heq}}{\Ch_b\bC})$. 
In the commutative diagram of triangulated categories
$${\tiny{
\xymatrix{
\Ho(\inn{{(\Ch_b\cA)}_{\heq}}{\Ch_b\bC}) \ar[r]^{\!\!\!\!\!\!\!\!\!\!\!\!\!\!\!\!\!\!\!\!\textbf{II}} \ar[d]_{\textbf{I}} &
\Ho(\inn{\Ch_b{(\Ch_b\cA)}_{\heq}}{\Ch_b\Ch_b\bC}) \ar[r]^{\!\!\!\!\!\!\!\!\!\!\!\!\!\!\!\!\!\!\!\!\textbf{III}} &
\Ho(\inn{{\langle\Ch_b{(\Ch_b\cA)}_{\heq}\rangle}_{\nullclass,qqw}}{\Ch_b\Ch_b\bC})\ar@{=}[d]\\
\Ho(\inn{{\langle\Ch_b\cA\rangle}_{\nullclass,qw}}{\Ch_b\bC}) \ar[r]_{\!\!\!\!\!\!\!\!\!\!\!\!\!\!\!\!\!\!\!\!\textbf{II}} &
\Ho(\inn{\Ch_b{\langle\Ch_b\cA\rangle}_{\nullclass,qw}}{\Ch_b\Ch_b\bC}) \ar[r]_{\!\!\!\!\!\!\!\!\!\!\!\!\!\!\!\!\!\!\!\!\textbf{III}} &
\Ho(\inn{{\langle\Ch_b{\langle\Ch_b\cA\rangle}_{\nullclass,qw}\rangle}_{\nullclass,qqw}}{\Ch_b\Ch_b\bC}),
}}}
$$
the triangulated functors $\textbf{II}$ and $\textbf{III}$ are equivalences 
of triangulated categories by Proposition~\ref{prop:comp hom and der} $\mathrm{(4)}$ and $\mathrm{(1)}$ respectively. 
Thus the triangulated functor $\textbf{I}$ above is also an equivalence of 
triangulated categories and we complete the proof. 
\end{proof}

\subsection{Homology theory on complicial exact categories}
\label{subsec:homology theory on complicial exact categories}

\begin{para}
\label{df:category of exact categories}
{\bf (Category of admissible exact sequences).}\ \ 
Let $\cE$ be an exact category and let 
$\cA$ and $\cB$ be a pair of full subcategory of $\cE$. 
We write $E(\cA,\cE,\cB)$ for the full subcategory of $F_{[0,2]}\cE$ 
consisting of those objects $x$ such that 
$x_0$ is in $\cA$ and 
$x_2$ is in $\cB$ and 
the sequence $x_0\onto{i^x_0}x_1\onto{i^x_1}x_2$ is an admissible exact sequence in $\cE$. 
if $\cA=\cB=\cE$, then we denote $E(\cA,\cE,\cB)$ by $E(\cE)$.  
\end{para}

\begin{para}
\label{lemdf:homology theory of complicial exact category}
{\bf Lemma-Definition (Homology theory of complicial exact categories).}\ 
(\cf \cite[1.7]{Wal85}.)\ \  
Let $\cE$ be an exact category and let $\cB$ be an abelian category. 
A {\it homology theory on $\cE$} ({\it with values in $\cB$}) 
is a family $\calH=\{H_n\colon\cC\to \cB,\delta_n\colon H_n({(-)}_2)\to 
H_{n-1}({(-)}_0) \}_{n\in\bbZ}$ 
indexed by the set of all integers 
consisting of 
additive functors $H_n$ and 
natural transformations $\delta^n$ between the functors from $E(\cE)$ to $\cB$ 
where ${(-)}_0$ and ${(-)}_2 $ are functors $ $ which sends an object 
$x\rinf y\rdef z$ in $E(\cE)$ to $x$ and $z$ respectively which subjects to 
the following condition:

\sn
For any admissible exact sequence $x\overset{i}{\rinf}y\overset{p}{\rdef}z$ 
in $\cE$, the long sequence 
\begin{equation}
\label{eq:les ass to xyz}
\cdots \onto{\delta_{n+1}^{(i,p)}} H_{n+1}(x) \onto{H_{n+1}(i)} H_{n+1}(y) \onto{H_{n+1}(p)} H_{n+1}(z)
\onto{\delta_n^{(i,p)}}H_{n}(x) \onto{H_{n}(i)} H_{n}(y) \onto{H_{n}(p)} 
H_{n}(z)\onto{\delta_{n-1}^{(i,p)}} \cdots
\end{equation}
is exact in $\cE$. 
For a homology theory $\calH=\{H_n,\delta_n\}_{n\in\bbZ}$ on $\cE$, 
we say that a morphism $f\colon x\to y$ is 
a {\it $\calH$-quasi-isomorphism} if 
$H_n(f)$ is an isomorphism for all integer $n$. 
We denote the class of all $\calH$-quasi-isomorphisms in $\cE$ by 
$w_{\calH,cE}$ or simply $w_{\calH}$.

Moreover, we assume that $\cE$ is a complicial exact category. 
Then for a homology theory $\calH=\{H_n,\delta_n\}_{n\in\bbZ}$ on $\cE$, 
the following conditions are equivalent:

\sn
$\mathrm{(1)}$ 
For an object $x$ in $\cE$ and an integer $n$, $H_n(r_x\colon CCx\to Cx)=0$. 

\sn
$\mathrm{(2)}$ 
For an object $x$ in $\cE$ and an integer $n$, $H_n(Cx)=0$. 

\sn
$\mathrm{(3)}$ 
For a pair of morphisms $f$, $g\colon x\to y$ in $\cE$, 
if $f$ and $g$ are $C$-homotopic, then $H_n(f)=H_n(g)$ for any 
integer $n$.

If $\calH$ satisfies the conditions above, we say that {\it $\calH$ is 
$C$-homotopy invariant}. 
For a $C$-homotopy invariant homology theory $\calH$ on $\cE$, 
we can show that 
the class $w_{\calH}$ of $\calH$-quasi-isomorphisms in $\cE$ 
is a class of complicial weak equivalences in $\cE$. 
\end{para}

\begin{proof}
Assume condition $\mathrm{(1)}$. 
Then for an object $x$ in $\cE$, $\id_{H_n(Cx)}=H_n(r_x\iota_x)=0$ for any 
integers. 
Thus condition $\mathrm{(2)}$ holds. 
Next assume condition $\mathrm{(2)}$ and let $f$, $g\colon x\to y$ 
be a pair of morphisms in $\cE$ and let $H\colon f\Rightarrow_C g$ be a 
$C$-homotopy from $f$ to $g$. Then 
$H_n(f)-H_n(g)=H_n(H\iota_x)=0$ for any integer $n$.  
Thus condition $\mathrm{(3)}$ holds. 
Finally assume condition $\mathrm{(3)}$. 
Then for an object $x$ in $\cE$, $r_xr_{Cx}\colon r_x\Rightarrow_C0$ is 
a $C$-homotopy from $r_x$ to $0$. 
Hence $H_n(r_x)=0$ for any integer $n$. Thus condition $\mathrm{(1)}$ holds. 
\end{proof}

\begin{para}
\label{ex:homoogy theory associated with t-structure}
{\bf Example (Homology theory associated with $t$-structures).}\ \ 
A typical example of a homology theory on a relative exact category $\bC$ 
comes from a $t$-structure on $\calD_b(\bC)$ the derived category of 
$\bC$. 
In this article, we use {\it homological $t$-structure} as 
in \cite[\S 1.2]{Lur17}. 
We will briefly recall 
the notions and fundamental results about (homological) $t$-structures 
from \cite{BBD82} 
(but conventions are written by cohomological notions in Ibid).

Let $(\cT,\Sigma)$ be a triangulated category and let 
$\calD$ be a full subcategory of $\cT$. 
For any integer $n$, we write 
$\calD[n]$ for the full subcategory of $\cT$ consisting of 
those objects $x[n]$ for some object $x$ in $\calD$ and 
we denote the full subcategory of $\cC$ consisting of those 
objects $x$ such that $\Hom_{\cT}(x,y)=0$ (resp. $\Hom_{\cT}(y,x)=0$) 
for any object $y$ in $\calD$ by ${}^{\perp}\calD$ (resp. $\calD^{\perp}$). 
We say that $\calD$ is a ({\it homological}) {\it $t$-structure of $\cT$} 
if $\calD$ satisfies the following two conditions. 

\sn
{\bf (Suspension and isomorphisms closed condition).}\ \ 
$\calD$ is closed under isomorphisms and closed under the suspension functor, 
that is $\calD[1]\subset\calD$.  

\sn
{\bf (Decomposition condition).}\ \ 
For any object $x$ in $\calD$, there exists a pair of 
object $\tau_{\geq 0}x$ in 
$\calD$ and $\tau_{\leq -1}x$ in $\calD^{\perp}$ and 
a distinguished triangle 
\begin{equation}
\label{eq:decomposition condition}
\tau_{\geq 0}x\to x\to \tau_{\leq -1}x\to \Sigma\tau_{\geq 0}x.
\end{equation}
In this case, we write 
$\calD_{\geq n}$ and $\calD_{\leq n}$ for $\calD[n]$ and 
${(\calD[n+1])}^{\perp}$ respectively. 
We can show that for an integer $n$, 
there exists a left (resp. right) adjoint functor of 
the inclusion functor $\calD_{\leq n}\rinc\calD$ 
(resp. $\calD_{\geq n}\rinc\calD$) which we denote by 
$\tau_{\leq n}\colon\calD\to\calD_{\leq n}$ 
(resp. $\tau_{\geq n}\colon\calD\to\calD_{\geq n}$) 
and 
$\calD^{\heartsuit}:=\calD_{\geq 0}\bigcap\calD_{\leq 0}$ 
is an abelian category and call it the {\it heart of $\calD$}. 
Moreover, the associations $x\mapsto \tau_{\geq 0}\tau_{\leq 0} x$ 
gives homological functor 
$H^{\calD}_0\colon \cT\to\calD^{\heartsuit}$ which we call 
the {\it homological functor associated with the $t$-structure $\calD$}. 

Thus for a relative complicial exact category $\bC=(\cC,w)$ and a homological 
$t$-structure $\calD$ 
on $\calD_b(\bC)$ the bounded derived category of $\bC$, 
the composition with the canonical functor $\cC\to\calD_b(\bC)$ with 
the homological functors $H_0^{\calD}(T^n(-))\colon\calD_b(\bC)\to\calD^{\heartsuit}$ 
gives a homotopy invariant homology theory on $\bC$. 
\end{para}

\begin{para}
\label{lemdf:Puppe exact sequence}
{\bf Lemma-Definition (Puppe exact sequences).}\ \ 
Let $\cC$ be a complicial exact category and let $\cB$ be an abelian category 
and let $\calH=\{H_n,\delta_n\}_{n\in\bbZ}$ be 
a $C$-homotopy invariant homology theory on $\cC$ with values in $\cB$. 
Then there exists a family of natural transformations 
$\{\delta_n'\colon H_n(\Cone )\to H_{n-1}(\dom)\}_{n\in\bbZ}$ 
between the functors from $\cC^{[1]}$ to $\cB$ indexed 
by the set of integers such that the long sequence 
\begin{equation}
\label{eq:long exact seq}
\cdots \onto{{\delta'}_{n+1}^f} H_{n+1}(x) \onto{H_{n+1}(f)} H_{n+1}(y) \onto{H_{n+1}(\kappa_f)} H_{n+1}(\Cone f)
\onto{{\delta'}_n^f}H_{n}(x) \onto{H_{n}(f)} H_{n}(y) \onto{H_{n}(\kappa_f)} H_{n}(\Cone f)\onto{{\delta'}_{n-1}^f} \cdots
\end{equation}
induced from a morphism $f\colon x\to y$ in $\cC$ 
is exact in $\cB$. 
We call a sequence $\mathrm{(\ref{eq:long exact seq})}$ 
the {\it Puppe exact sequence} ({\it associated 
with the morphism $f\colon x\to y$}). 
\end{para}

\begin{proof}
Let $f\colon x\to y$ be a morphism in $\cC$. 
Notice that $\Cyl f\onto{\upsilon_f} y$ is a $C$-homotopy equivalence and 
thus by considering the admissible exact sequence in 
the top line in the commutative diagram below,
$$
\footnotesize{
\xymatrix{
x \ar@{>->}[r]^{{\xi_1}_f} \ar[rd]_f & 
\Cyl f \ar[d]_{\upsilon_f} \ar@{->>}[r]^{\eta_f} & 
\Cone f\\
& y \ar[ru]_{\kappa_f}, &
}}
$$
we obtain the long exact sequence $\mathrm{(\ref{eq:long exact seq})}$ 
where we set 
${\delta'}^f_n:=\delta_n^{({\xi_1}_f,\eta_f)}$ for all integer $n$. 
\end{proof}

\section{D\'evissage spaces}
\label{sec:devissage spaces}

The main theme of this section is to define 
d\'evissage spaces associated with 
cell structures and by using this phraseology, 
we will discuss a derived version of d\'evissage condition. 
In the first subsection \ref{subsec:cell str}, 
we will introduce an idea of cell structures 
on relative complicial exact categories 
which is a resemblance of $t$-structure 
or weight structure 
on triangulated categories 
and the next subsection \ref{subsec:dev filt}, 
we will renew the conception of d\'evissage filtrations with respect to 
a cell structure. 
In subsection \ref{subsec:derived devissage condition}, 
we will innovate a derived version of 
d\'evissage conditions and discuss relationship this conditions 
with Quillen's, Raptis' and Waldhausen's one. 
The last two subsections are devoted to categorified calculation of 
d\'evissage spaces. 
In subsection \ref{subsec:derived quasi-split sequence}, 
we will study a derived quasi-split sequence of 
relative exact categories. 
In the context of quasicategory theory, analogous concept is studied 
in \cite{FP07}. 
The final section \ref{subsec:flag}, 
we provide natural derived flag structures on the category of 
d\'evissage filtrations by utilizing results in the previous subsections.

\subsection{Cell structures}
\label{subsec:cell str}

\begin{para}
\label{df:properties of family}
{\bf Definition.}\ \ 
Let $\cC$ be a complicial exact category and let $w$ 
be a class of morphisms in $\cC$ and let 
$\bbD=\{\calD_n\}_{n\in\bbZ}$ be a family of full subcategories of 
$\cC$ indexed by the set of all integers $\bbZ$. 
We say that $\bbD$ is {\it topologizing} 
(resp. {\it semi-Serre}, {\it Serre}, {\it prenull}, {\it semi-null}, 
{\it closed under finite direct sums}, {\it thick}, {\it $w$-closed}) 
if for all integer $n$, $\calD_n$ is a topologizing 
(resp. semi-Serre and so on) 
subcategory of $\cC$. 
We say that $\bbD$ is a {\it $T$-system} 
if for all integer $n$, $T(\Ob \calD_n)\subset \Ob \calD_{n+1}$. 

Let $\bbD=\{\calD_n\}_{n\in\bbZ}$ and $\bbD'=\{\calD'_n\}_{n\in\bbZ}$ be 
a pair of $T$-systems of relative complicial exact categories 
$\bC=(\cC,w)$ and $\bC'=(\cC',w')$ respectively. 
Assume that $\bbD'$ are isomorphisms closed. 
Namely for each integer $n$, $\calD'_n$ is closed under isomorphisms. 
A {\it morphism of $T$-system} $(f,d)\colon\bbD\to\bbD'$ 
is a relative complicial 
exact functor $\bC\to\bC'$ such that for each integer $n$, 
$f(\Ob\calD_n)$ is contained in $\Ob\calD'_{n}$. 
\end{para}

\begin{para}
\label{df:cell str}
{\bf Definition (Cell structure).}\ \ 
Let $\cC$ be a complicial exact category. 
A {\it cell structure $\bbD=\{\calD_{\leq m},\calD_{\geq n} \}_{n,m\in\bbZ}$ 
of $\cC$} 
is a family of class of full subcategories in $\cC$ satisfying 
the following conditions:
\begin{itemize}
\item
For any integer $n$, $\calD_{\leq n}\subset \calD_{\leq n+1}$ and 
$\calD_{\geq n+1}\subset \calD_{\geq n}$. 

\item
For any integer $n$, 
$\calD_{\geq n}$ and $\calD_{\leq n}$ are closed under $C$-homotopy 
equivalences.

\item
For any integer $n$ and an object $x$ in $\cC$, 
\begin{itemize}
\item
if $x$ is in $\calD_{\geq n}$, then $Tx$ is in $\calD_{\geq n+1}$ and 
\item
if $x$ is in $\calD_{\leq n}$, then $T^{-1}x$ is in $\calD_{\leq n-1}$.
\end{itemize}

\item
For any morphism $f\colon x\to y$ in $\cC$, 
\begin{itemize}
\item
if $y$ is in $\calD_{\geq n-1}$ and $\Cone f$ is in $\calD_{\geq n}$, 
then $x$ is in $\calD_{\geq n-1}$, and

\item
if $x$ is in $\calD_{\leq n}$ and $\Cone f$ is in $\calD_{\leq n}$, 
then $y$ is in $\calD_{\leq n}$. 
\end{itemize}
\end{itemize}
\item
We set $\calD_n:=\calD_{\leq n}\cap\calD_{\geq n}$ for any integer $n$ 
and set $\bbD^{\heartsuit}:=\{\calD_n\}_{n\in\bbZ}$ the family of full subcategory of 
indexed by the set of all integers. 
Then 
we will show that 
$\bbD^{\heartsuit}$ is a $T$-system of $\cC$ 
in \ref{lem:cell structure fundamental result}. 
We call $\bbD^{\heartsuit}$ the {\it heart of $\bbD$}. 
Moreover let $w$ be a class of complicial weak equivalences of $\cC$. 
We say that {\it $\bbD$ is compatible with $w$} or 
$\bbD$ is a {\it cell structure of $(\cC,w)$}, 
if $\displaystyle{\cC^w\subset\underset{n\in\bbZ}{\bigcap}\calD_n}$.

We say that a cell structure $\bbD$ is {\it topologizing} 
(resp. {\it semi-Serre}, {\it Serre}, {\it prenull}, {\it semi-null}, 
{\it closed under finite direct sums}, {\it thick}, {\it $w$-closed}) 
if the heart $\bbD^{\heartsuit}$ of $\bbD$ 
is a topologizing (resp. semi-Serre and so on). 
We say that a cell structure $\bbD$ is {\it strongly topologizing} 
(resp. {\it strongly semi-Serre}, {\it strongly Serre}, {\it strongly prenull}, {\it strongly semi-null}, 
{\it strongly closed under finite direct sums}, {\it strongly thick}, 
{\it strongly $w$-closed}) if 
for any integer $n$, 
both $\calD_{\geq n}$ and $\calD_{\leq n}$ 
are topologizing (resp. semi-null and so on). 
\end{para}

\begin{para}
\label{df:restriction of cell structures}
{\bf Definition (Restriction of cell structures).}\ \ 
Let $\bC=(\cC,w)$ be a relative complicial exact category and let 
$\bbD=\{\calD_{\geq n},\calD_{\leq m} \}_{n,m\in\bbZ}$ 
be a cell structure of $\bC$ and 
let $\cF$ be a full subcategory of $\cC$. 
Then for each integer $n$, we set 
${(\calD|_{\cF})}_{\geq n}:=\calD_{\geq n}\bigcap {\langle\cF\rangle}_{\nullclass,\frob}$ and ${(\calD|_{\cF})}_{\leq n}:=\calD_{\leq n}\bigcap {\langle\cF\rangle}_{\nullclass,\frob}$ and 
$\bbD|_{\cF}:=\{{(\calD|_{\cF})}_{\geq n},
{(\calD|_{\cF})}_{\leq m}\}_{n,m\in\bbZ}$. 
Then $\bbD|_{\cF}$ is a cell structure of $({\langle\cF\rangle}_{\nullclass,\frob},w|_{{\langle\cF\rangle}_{\nullclass,\frob}})$ 
and call it the {\it restriction 
of $\bbD$} ({\it to $\cF$}). 
\end{para}

Recall the definition of homology theory on a complicial exact categories 
from \ref{lemdf:homology theory of complicial exact category}. 

\begin{para}
\label{ex:cell str associated with homology theory}
{\bf Example (Cell structure associated with homology theory).}\ 
(\cf \cite[1.7]{Wal85})\ \ 
Let $\cC$ be a complicial exact category and 
let $\calH=\{H_n,\delta_n\}_{n\in\bbZ}$ be a homology theory on $\cC$ 
with values in an abelian category $\cB$. 
Then we define a family 
$\bbD_{\calH}=\{\calD_{\leq n},\calD_{\geq m}\}_{n,m\in\bbZ}$ 
of full subcategories in $\cC$ by the formula 
\begin{equation}
\label{eq:cell st ass to homo th 1}
\Ob\calD_{\leq n}:=\{x\in\Ob \cX;\text{$H_k(x)=0$ for $k>n$}\}\ \text{ and}
\end{equation}
\begin{equation}
\label{eq:cell st ass to homo th 2}
\Ob\calD_{\geq m}:=\{x\in\Ob \cX;\text{$H_k(x)=0$ for $k<m$}\}.
\end{equation}
Then the family $\bbD_{\calH}$ is a cell structure on $\cC$ 
by Puppe exact sequence \ref{lemdf:Puppe exact sequence}. 
\end{para}

\begin{para}
\label{ex:cell structure ass with full subcat}
{\bf Example (Cell structure associated with a full subcategory).}\ \ 
Let $\bC=(\cC,w)$ be a relative complicial exact category and let 
$\cA$ be a full subcategory of $\cC$. 
Then we set 
$\calD^{\cA}_{\geq n}=\calD^{\cA}_{\leq m}={\langle\cA\rangle}_{\nullclass,w}$ 
for any integers $n$ and $m$. 
The family $\bbD_{\cA}:=\{\calD^{\cA}_{\geq n},\calD^{\cA}_{\leq m} \}_{n,m\in\bbZ} $ forms a cell structure of $\bC$. 
We call $\bbD_{\cA}$ the 
{\it cell structure associated with a full subcategory $\cA$}. 
\end{para}

\begin{para}
\label{df:m-connected morphism}
{\bf Definition ($m$-connected morphisms).}\ \ 
Let $\cC$ be a complical exact category and $\bbD=\{\calD_{\geq n},\calD_{\leq m}\}_{n,m\in\bbZ}$ a cell structure of $\cC$ and let 
$f\colon x\to y$ be a morphism in $\cC$ and let $m$ be an integer. 
We say that $f$ is {\it $m$-connected} ({\it with respect to $\bbD$}) or 
{\it $(m,\bbD)$-connected} if 
$\Cone f$ is in $\calD_{\geq m+1}$. 
\end{para}

\begin{para}
\label{lem:cell structure fundamental result}
{\bf Lemma.}\ \ 
{\it
Let $\cC$ be a complicial exact category and 
let $\bbD=\{\calD_{\leq n},\calD_{\geq m}\}_{n,m\in\bbZ}$ be a 
cell structure of $\cC$ and let $f\colon x\to y$ 
and $g\colon y\to z$ be a pair of morphisms in $\cC$ 
and let $n$ be an integer. 
Then

\sn
$\mathrm{(1)}$ 
If $x$ is in $\calD_{\geq n}$ and 
if $\Cone f$ is in $\calD_{\geq n}$, 
then $y$ is in $\calD_{\geq n}$.

\sn
$\mathrm{(2)}$ 
If $x$ is in $\calD_{\geq n-1}$ and 
if $y$ is in $\calD_{\geq n}$, 
then $\Cone f$ is in $\calD_{\geq n}$. 

\sn
$\mathrm{(3)}$ 
If $y$ is in $\calD_{\leq n}$ and if $x$ is in $\calD_{\leq n-1}$, 
then $\Cone f$ is in $\calD_{\leq n}$. 

\sn
$\mathrm{(4)}$ 
If $y$ is in $\calD_{\leq n-1}$ and if $\Cone f$ is in $\calD_{\leq n}$, 
then $x$ is in $\calD_{\leq n-1}$.

\sn
$\mathrm{(5)}$ 
If $x$ is in $\calD_{\leq n}$, 
then $Tx$ is in $\calD_{\leq n+1}$. 

\sn
$\mathrm{(6)}$ 
If $x$ is in $\calD_{\geq n}$, 
then $T^{-1}x$ is in $\calD_{\geq n-1}$. 

\sn
$\mathrm{(7)}$ 
If $\Cone f$ is in $\calD_{\geq n}$ and 
if $\Cone g$ is in $\calD_{\geq n}$, 
then $\Cone gf$ is in $\calD_{\geq n}$.

\sn
$\mathrm{(8)}$ 
If $\Cone f$ is in $\calD_{\geq n-1}$ and 
if $\Cone gf$ is in $\calD_{\geq n}$, 
then $\Cone g$ is in $\calD_{\geq n}$.

\sn
$\mathrm{(9)}$ 
If $\Cone g$ is in $\calD_{\geq n}$ and 
if $\Cone gf$ is in $\calD_{\geq n-1}$, 
then $\Cone f$ is in $\calD_{\geq n-1}$.

\sn
$\mathrm{(10)}$ 
If $\Cone f$ is in $\calD_{\leq n}$ and 
if $\Cone g$ is in $\calD_{\leq n}$, 
then $\Cone gf$ is in $\calD_{\leq n}$.

\sn
$\mathrm{(11)}$ 
If $\Cone f$ is in $\calD_{\leq n-1}$ and 
if $\Cone gf$ is in $\calD_{\leq n}$, 
then $\Cone g$ is in $\calD_{\leq n}$.

\sn
$\mathrm{(12)}$ 
If $\Cone g$ is in $\calD_{\leq n}$ and 
if $\Cone gf$ is in $\calD_{\leq n-1}$, 
then $\Cone f$ is in $\calD_{\leq n-1}$.

\sn
$\mathrm{(13)}$ 
$T\colon \cC\to\cC$ induces equivalences of categories 
$T\colon\calD_{\leq n}\isoto \calD_{\leq n+1}$, 
$T\colon\calD_{\geq n}\isoto\calD_{\geq n+1}$ and 
$T\colon\calD_n\isoto\calD_{n+1}$. 
In particular, $\bbD^{\heartsuit}$ is a $T$-system.
}
\end{para}

\begin{proof}
By \ref{ex:Conekappaf}, 
the canonical morphism $\Cone \kappa_f\to Tx$ is 
a $C$-homotopy equivalence. 
Thus by applying the axiom of cell structure 
to the morphism $\kappa_f\colon y\to\Cone f$, 
we obtain assertions from $\mathrm{(1)}$ to $\mathrm{(4)}$. 
Next Applying $\mathrm{(3)}$ to the morphism $x\to 0$ and $\mathrm{(2)}$ 
to the morphism $T^{-1}x \to 0$, we obtain assertions 
$\mathrm{(5)}$ and $\mathrm{(6)}$. 
Next by \ref{ex:cone of composition}, 
the canonical morphism 
$\Cone(\kappa_f,\kappa_{gf})\colon\Cone g\to \Cone\Cone(\id_x,g)$ 
is a $C$-homotopy equivalence. 
Thus by applying the axiom of cell structures and assertions from 
$\mathrm{(1)}$ to $\mathrm{(4)}$ 
to the morphism $\Cone(\kappa_f,\kappa_{gf})$, 
we obtain assertions from $\mathrm{(7)}$ to $\mathrm{(12)}$. 
Finally $\mathrm{(13)}$ follows from the axioms of cell structures and 
$\mathrm{(5)}$ and $\mathrm{(6)}$.
\end{proof}

\begin{para}
\label{rem:omittion of bbD}
{\bf Remark.}\ \ 
In \ref{lem:cell structure fundamental result}, 
$\bbD$ is essentially determined from the pair $\calD_{\geq 0}$ and 
$\calD_{\leq 0}$ by $\mathrm{(13)}$. 
Thus we often abbreviate $\bbD$ as $\{\calD_{\geq 0},\calD_{\leq 0}\}$. 
\end{para}

\begin{para}
\label{df:bounded cell structure}
{\bf Definition (Bounded cell structure).}\ \ 
Let $\cC$ be a complicial exact category and let 
$\bbD=\{\calD_{\geq 0},\calD_{\leq 0}\}$ be a cell structure on $\cC$. 
For a pair of integers $n$ and $m$, we set 
$\calD[n,m]:=\calD_{\geq n}\cap\calD_{\leq m}$. 
We say that $\bbD$ is {\it bounded} if 
\begin{equation}
\label{eq:bounded cell structure df}
\cC=\bigcup_{\substack{(n,m)\in\bbZ\\ n<m}}\calD[n,m].
\end{equation}
\end{para}

\begin{para}
\label{df:ordinary cell structure}
{\bf Definition (Ordinary cell structure).}\ \ 
Let $\bC=(\cC,w)$ be a relative complicial exact category and 
let $\bbD=\{\calD_{\geq 0},\calD_{\leq 0}\}$ 
be a cell structure on $\bC$. 
We say that $\bbD$ is {\it ordinary} if for any pair of integers 
$n<m$, $\calD[m,n]=\cC^w$. 
\end{para}

\begin{para}
\label{lem:invariance connectedness under co-base change}
{\bf Lemma.}\ \ 
{\it
Let $\bC=(\cC,w)$ be a relative complicial exact category and let 
$\bbD=\{\calD_{\geq 0},\calD_{\leq 0}\}$ be a cell structure of 
$\bC$ and let $i\colon x\rinf y$ be an admissible monomorphism in $\cC$ 
and let $f\colon x\to z$ be a morphism in $\cC$ and let $m$ be an integer. 
Assume that $i$ is a Frobenius admissible monomorphism or $\bbD$ 
is strongly $w$-closed. 
Then if $\Cone i$ is in $\calD_{\geq m}$ {\rm(}resp. $\calD_{\leq m}${\rm)}, 
then $\Cone(i'\colon z\rinf z\sqcup_x y)$ is also in 
$\calD_{\geq m}$ {\rm(}resp. $\calD_{\leq m}${\rm)}.
}
\end{para}

\begin{proof}
Assume that $i$ is a Frobenius admissible monomorphism 
(resp. $\bbD$ is strongly $w$-closed). 
Then the canonical morphism $\Cone i\to \Cone i'$ 
is a $C$-homotopy equivalence 
(resp. is in $w$) by \ref{cor:Cones invariance by co-base change}. 
Thus if $\Cone i$ is in $\calD_{\geq m}$ (resp. $\calD_{\leq m}$), 
then $\Cone i'$ is also in $\calD_{\geq m}$ (resp. $\calD_{\leq m}$).
\end{proof}

\subsection{D\'evissage filtrations}
\label{subsec:dev filt}

\begin{para}
\label{df:devissage filtrations}
{\bf Definition (D\'evissage filtrations).}\ \ 
Let $\cC$ be a complicial exact category and 
let $\bbD=\{\calD_n \}_{n\in\bbZ}$ 
be a family of full subcategories of $\cC$ indexed by 
the set of all integers $\bbZ$. 
A {\it d\'evissage filtration} ({\it with respect to the family 
$\bbD$ in $\cC$}) is a bounded filtered objects (see Definition~\ref{df:filtered objects}) $x$ on $\cC$ such that 
for any integer $n$, $\Cone i_n^x$ is in $\calD_{n+1}$. 
We write $F_{b}^{\bbD}\cC$ (resp. $F_{b,h}^{\bbD}\cC$) 
for the full subcategory of $F_b\cC$ (resp. $F_{b,h}\cC$) 
consisting of d\'evissage filtrations with respect to the family $\bbD$. 
We also set for any pair of integers $a\leq b$, 
$F_{[a,b]}^{\bbD}\cC:=F_{[a,b]}\cC\bigcap F_b^{\bbD}\cC$ and so on. 
\end{para}

\begin{para}
\label{ex:C[1]hD0D1}
{\bf Example.}\ \ 
Let $\cC$ be a complicial exact category and let $\calD_0$ 
and $\calD_1$ be a pair of full subcategories of $\cC$. 
Then we set 
$\cC^{[1]}_h(\calD_0,\calD_1):=F_{[0,1],h}^{\{\calD_0,\calD_1\}}\cC$. 
Namely $\cC^{[1]}_h(\calD_0,\calD_1)$ 
is a full subcategory of $\cC^{[1]}_h$ consisting 
of those objects $[f\colon x\to y]$ such that $x$ is in $\calD_0$ 
and $\Cone f$ is in $\calD_1$. 
\end{para}

\begin{para}
\label{prop:herditary properties}
{\bf Proposition.}\ \ 
{\it
Let $\cC$ be a normal ordinary complicial exact category and let 
$\bbD=\{\calD_n\}_{n\in \bbZ}$ be a family of full subcategories 
of $\cC$ indexed by the set of all integers $\bbZ$. Then

\sn
$\mathrm{(1)}$ 
If $\bbD$ is topologizing {\rm(}resp. semi-Serre, Serre, 
closed under finite direct sums, prenull, semi-null, thick{\rm)} in $\cC$, 
then $F_{b,h}^{\bbD}\cC$ is also in $F_{b,h}\cC$. 

\sn
$\mathrm{(2)}$ 
Moreover let $w$ be a class of complicial weak equivalences of $\cC$. 
If $\bbD$ is $w$-closed in $\cC$, 
then $F_{b,h}^{\bbD}\cC$ is $lw$-closed in 
$F_{b,h}\cC$.
}
\end{para}

\begin{proof}
$\mathrm{(1)}$ 
For a pair of objects $x$ and $y$ in $F_{b,h}\cC$ and an integer $n$, 
we have the canonical isomorphism $\Cone i_n^{x\oplus y}\isoto \Cone i_n^x\oplus \Cone i_n^y$. 
Thus if $\bbD$ is closed under finite direct sum in $\cC$, 
then $F_{b,h}^{\bbD}\cC$ 
is also closed under finite direct sum in $F_{b,h}\cC$. 

For a level Frobenius admissible exact sequence 
$\mathrm{(\ref{eq:Frob exact in Fbh})}$ in $F_{b,h}\cC$ and 
an integer $n$, the sequence 
$\Cone i_n^x \rinf \Cone i_n^y \rdef \Cone i_n^z$ is 
a Frobenius admissible sequence in $\cC$ by 
\ref{cor:exactness of Cyl and Cone}. 
Thus it turns out that if $\bbD$ is topologizing (resp. semi-Serre, Serre), 
then $F_{b,h}^{\bbD}\cC$ is also in $F_{b,h}\cC$. 

Let $(f,H)\colon x\to y$ be a homotopy commutative diagram in $F_{b,h}\cC$. 
Then by \ref{cor:ConeCone}, 
for each integer $n$, there is a natural isomorphisms 
$\Cone i_n^{\Cone(f,H)}\isoto \Cone \Cone(f_n,f_{n+1},H_n)$. 
Thus if $\calD_n$ is closed under taking the functor 
$\Cone$ for all $n$, 
then $F_{b,h}^{\bbD}\cC$ is 
also closed under taking $\Cone$ in $F_{b,h}\cC$. 
If there exists a homotopy commutative diagram $(g,K)\colon y\to x$ 
such that $(g,K)(f,H)=\id_x$, then for any integer $n$, we have 
the equality 
$\Cone(g_n,g_{n+1},K_n)\Cone(f_n,f_{n+1},H_n)=\id_{\Cone i_n^x}$. 
Thus if $y$ is in $F_{b,h}^{\bbD}\cC$ and $\bbD$ is thick, then 
$x$ is also in $F_{b,h}^{\bbD}\cC$. 

\sn
$\mathrm{(2)}$ 
Let $(f,H)\colon x\to y$ be 
a homotopy commutative diagram in 
$F_{b,h}\cC$ and assume that $(f,H)$ is in $lw$ and $x$ (resp. $y$) 
is in $F_{b,h}^{\bbD}\cC$. 
Then for each integer $n$, $\Cone(f_n,f_{n+1},H_n)\colon 
\Cone i_n^x\to \Cone i_n^y$ is in $w$ 
by \ref{prop:properties of weak equivalences} 
$\mathrm{(2)}$ $\mathrm{(ii)}$. 
Thus $\Cone i_n^y$ (resp. $\Cone i_n^x$) is in $\calD_{n+1}$ 
and it turns out that $y$ (resp. $x$) is in $F_{b,h}^{\bbD}\cC$.
\end{proof}

Recall the functors $\fj$ and $(-)[k]$ from 
\ref{df:filtered objects} and \ref{df:degree shift}. 

\begin{para}
\label{lemdf:devissage filtration with level weak equivalences}
{\bf Lemma-Definition.}\ \ 
Let $(\cC,w)$ be a relative complicial exact category. 
For an integer $k$ and a full subcategory $\cS$ in $\cC$, 
we denote the full subcategory of $F_{b,\leq k,h}\cC$ consisting 
of those objects $x$ such that $x_k$ is in $\cS$ and 
$x_i$ is $C$-contractible for any $i<k$ by 
$F_{b,h}(\inn{\cS[k]}{\cC})$. 
We define $\fj(-)[k]\colon\cS\to F_{b,h}(\inn{\cS[k]}{\cC})$ and 
$(-)_k\colon F_{b,h}(\inn{\cS[k]}{\cC})\to \cS$ to be a pair of functors 
by sending an object $x$ in $\cS$ to $\fj(x)[k]$ and 
an object $z$ in $F_{b,h}(\inn{\cS[k]}{\cC})$ to $z_k$ respectively. 
We have the equality ${(\fj(-)[k])}_k=\id_{\cS}$ and there exists 
a natural weak equivalence 
$\fj({(-)}_k)[k]\to\id_{F_{b,h}(\inn{\cS[k]}{\cC})}$ 
with respect to the class of level weak equivalences of $F_{b,h}(\inn{\cS[k]}{\cC})$. 
Thus the pair of relative categories $(\cS,w|_{\cS})$ and 
$(F_{b,h}(\inn{\cS[k]}{\cC}),lw)$ are categorical homotopy equivalent. 

Next assume that $(\cC,w)$ is thick normal ordinary and 
let $a<b$ be a pair of integers. 
We set $\cF:=F_{[a,b],h}\cC$ and $\cG:=\cF^{[1]}_h(F_{[a,b-1],h}\cC,F_{b,h}(\inn{\cC[b]}{\cC}))$. 
Recall the definition of the category $\cF^{[1]}_h(-,-)$ 
from \ref{ex:C[1]hD0D1}. 
We define 
$\fG\colon\cG\to \cF$ 
and $\fH\colon\cF\to \cG$ to be strictly normal complicial functors 
by sending an object $[f\colon x\to y]$ in $\cG$ to $y$ in $\cF$ and 
$y$ in $\cF$ to $[\sigma_{\leq b-1}y\to y]$ in $\cG$. 
We have the equality $\fG\fH=\id_{\cF}$ and 
there exist a natural weak equivalence 
$\id_{\cG}\to \fH\fG$ with respect to the class $llw$ of level level weak equivalences of $\cG$. 
Thus the pair of relative categories 
$(\cG,llw)$ and $(\cF,lw)$ are categorical homotopy equivalent each other.

Moreover let $\bbD=\{\calD_n\}_{n\in\bbZ}$ be a family of 
full subcategories of $\cC$ indexed by the set of integers $\bbZ$. 
Assume that for any integer $n$, $\calD_n$ contains 
all $C$-contractible objects in $\cC$. 
Then we set $\cF':=F_{[a,b],h}^{\bbD}\cC$ and 
$\cG':={\cF'}^{[1]}_h(F_{[a,b-1],h}^{\bbD}\cC,F_{b,h}(\calD_b[b],\cC))$. 
Then the restrictions of $\fG$ and $\fH$ to $\cG'$ and $\cF'$ induce 
the categorical homotopy equivalences 
$(\cG',llw)\underset{\fH}{\overset{\fG}{\rightleftarrows}} (\cF',lw)$ and 
$({\langle\cG'\rangle}_{\nullclass,\frob},llw)\underset{\fH}{\overset{\fG}{\rightleftarrows}}({\langle\cF'\rangle}_{\nullclass,\frob},lw)$ 
by \ref{lem:functoriality of closures}. 
\qed
\end{para}

\begin{para}
\label{lem:dev fil}
{\bf Lemma.}\ \ 
{\it 
Let $(\cC,w)$ be a relative complicial exact category and 
$\bbD=\{\calD_{\leq n},\calD_{\geq m} \}_{n,\ m\in\bbZ}$ 
be a cell structure of $(\cC,w)$ and 
$x$ a d\'evissage filtration with respect to $\bbD$ 
having amplitude contained in $[a,b]$. 
Then\\
$\mathrm{(1)}$ 
$x_k$ is in $\calD_{\leq k}$ for any $k$.\\ 
$\mathrm{(2)}$ 
If $x_{b}$ is in $\cC^w$, then $x_k$ is in 
$\calD_{k}$ for any $k$.
}
\end{para}

\begin{proof}
$\mathrm{(1)}$ 
We proceed by induction on $k$. 
Since $x_k=0$ for $k<a$, $x_k$ is in $\calD_k$. 
Assume that $x_k$ is in $\calD_{\leq k}$. 
Then by considering an admissible sequence
\begin{equation}
\label{eq:ixk}
x_k \onto{i^x_k} x_{k+1} \onto{\kappa_{i^x_k}}\Cone i^x_k
\end{equation}
with 
$x_k \in \calD_{\leq k}$ and $\Cone i_k^x  \in \calD_{k+1}$, 
it turns out that $x_{k+1}$ is in $\calD_{\leq k+1}$. 

\sn
$\mathrm{(2)}$ 
We proceed by descending induction on $k$. 
Since $x_k=x_b$ for $k\geq b$, $x_k$ is in $\calD_k$. 
Assume that $x_{k+1}$ is in $\calD_{k+1}$. 
Then by considering a sequence $\mathrm{(\ref{eq:ixk})}$ with 
$x_{k+1} \in \calD_{k+1}$ and $\Cone i_k^x  \in \calD_{k+1}$, 
it turns out that $x_{k}$ is in $\calD_{\geq k}$. 
Thus $x_k$ is in $\calD_k$.
\end{proof}

Recall the functors $\fj$, $(-)[k]$, $\fc_k$ and $\fs_k$ from 
\ref{df:filtered objects}, \ref{df:degree shift}, 
\ref{df:Cone functor} and \ref{df:skip functor}. 

\begin{para}
\label{lemdf:acyclic devissage filtration with level weak equivalences}
{\bf Lemma-Definition.}\ \ 
Let $(\cC,w)$ be a thick normal ordinary relative 
complicial exact category and let $\cS$ be a full subcategory and let 
$k$ be an integer. 
Then there are a pair of functors
$\cS\rlto{\fj(-)[k]}{{(-)}_k}
\underline{F_{[k,k],h}\cS}$ and they gives an isomorphism of categories. 

For a triple of integers $a\leq k\leq b$, 
we define $\cE$ to be a full subcategory of 
$\underline{F_{[a,b],h}\cC}$ consisting of those objects $x$ 
such that $x_{b-1}$ is $C$-contractible and 
let $\cE'$ be a subcategory of $\cE$ such that 
$\Ob\cE':=\Ob\cE$ and 
$\Mor\cE':=\{(f,H)\colon x\to y\in\Mor\cE;i_{b-1}^yf_{b-1}=f_bi_{b-1}^x,\ 
H_{b-1}=0 \}$. 
Then $\cE'$ is a complicial exact category such that 
the inclusion functor $\cE'\rinc\cE$ is a 
strictly normal complicial functor and reflects 
exactness. 
For an object $x$ and a homotopy commutative 
diagram $(f,H)\colon x\to y$ in $\cE$, we define 
$\fr x$ and $\fr (f,H)\colon \fr x\to \fr y$ 
to be an object and a morphism in $\cE'$ by setting
$$
\footnotesize{
{(\fr x)}_n:=
\begin{cases}
x_n & \text{if $n\leq b-1$}\\
x_b\oplus Cx_{b-1} & \text{if $n\geq b$},
\end{cases}\ \ \ 
i_n^{\fr x}:=
\begin{cases}
i_n^x & \text{if $n\leq b-2$}\\
\begin{pmatrix}
i^x_{b-1}\\
-\iota_{x_{b-1}}
\end{pmatrix} & \text{if $n=b-1$}\\
id_{x_b\oplus Cx_{b-1}} & \text{if $n\geq b$},
\end{cases}\ \ \ 
{(\fr (f,H))}_n:=
\begin{cases}
(f_n,H_n) & \text{if $n\leq b-2$}\\
(f_{b-1},0) & \text{if $n=b-1$}\\
(\Cyl(f_{b-1},f_b,H_b),0) & \text{if $n\geq b$}.
\end{cases}
}$$
Then the association $\fr\colon\cE\to\cE'$ is an exact functor. 
Moreover the natural equivalence $\sigma\colon CC\isoto CC$ 
induces a natural equivalence 
$d^{\fr}\colon C^{\lv}\fr\isoto \fr C^{\lv}$ and the pair 
$(\fr,d^{\fr})\colon (\cE,lv)\to (\cE',lv)$ is a relative 
complicial exact functor. 
We denote the inclusion functor by $\fraki\colon\cE'\rinc \cE$. 
Then by \ref{lem:CC fundamental results} $\mathrm{(6)}$, 
there are natural weak equivalences 
$\fraki\fr\to\id_{\cE}$ and $\fr\fraki\to\id_{\cE'}$ with 
respect to level weak equivalences. 
Restrictions of $\fs_{b-1}$ and $\fc_b$ induces a pair of 
complicial exact functors 
$\cE^{w_{\st}}\rlto{\fs_{b-1}}{\fc_b}{\left(\underline{F_{[a,b-1],h}\cC}\right)}^{w_{\st}}$ and there are natural weak equivalences 
$\id_{\cE'}\to \fc_b\fs_{b-1}|_{{\cE'}^{w_{\st}}}$ and 
$\id_{{\left(\underline{F_{[a,b-1],h}\cC}\right)}^{w_{\st}}}
\to\fs_{b-1}\fc_{b}$. 
Thus $(\cE^{w_{\st}},lw)$, $({\cE'}^{w_{\st}},lw)$ 
and $({\left(\underline{F_{[a,b-1],h}\cC}\right)}^{w_{\st}},lw)$ 
are categorical homotopy equivalences each other.

We set $\cF:={\left (\underline{F_{[a,b],h}\cC}\right)}^{w_{\st}}$ and 
$\cG:=\cF^{[1]}_h
(\underline{F_{[b-1,b-1],h}\cC},\cF,\cE)$. 
Then there exists a pair of 
strictly normal complicial exact functors 
$\cF\rlto{\fL}{\fM}\cG$ defined by sending an object $x$ in $\cF$ 
to $[x_{b-1}[b-1]\to x]$ in $\cG$ and an object $[x\to y]$ in $\cG$ to $y$ in $\cF$. 
We have the equality $\fM\fL=\id_{\cF}$ 
and there exists a natural transformation 
$\id_{\cG}\to \fL\fM$.  
Thus the pair of relative categories $(\cF,lw)\rlto{\fL}{\fM}(\cG,llw)$ are 
categorical homotopy equivalence.

Moreover let $\bbD=\{\calD_{\leq n},\calD_{\geq m} \}_{n,m\in\bbZ}$ 
be a cell structure of $(\cC,w)$ such that it is 
closed under $C$-homotopy equivalences. 
We set 
$\cE'':=\cE\cap\underline{F_{[a,b],h}^{\bbD^{\heartsuit}}\cC}$, 
$\cF':={\left(F^{\bbD^{\heartsuit}}_{[a,b],h}\cC \right)}^{w_{\st}}$ 
and $\cG':={\cF'}^{[1]}_h(\underline{F_{[b-1,b-1],h}\calD_{b-1}},\cF',\cE'')$. 
Then by \ref{lem:functoriality of closures} and 
\ref{lem:dev fil}, 
restriction of $\fL$ and $\fM$ induce 
a pair of strictly normal complicial exact functors 
$\cF'\rlto{\fL'}{\fM'}\cG'$ and it induces a categorical homotopy equivalences 
$(\cF',lw)\rlto{\fL'}{\fM'}(\cG',llw)$ and 
$({\langle\cF'\rangle}_{\nullclass,\frob},lw)\rlto{\fL'}{\fM'}
({\langle\cG'\rangle}_{\nullclass,\frob},llw)$. 
\end{para}

\subsection{Derived d\'evissage condition}
\label{subsec:derived devissage condition}

\begin{para}
\label{df:Derived devissage condition}
{\bf Definition (Derived d\'evissage condition).}\ \ 
Let $\bC=(\cC,w)$ be a thick normal ordinary relative 
complicial exact category and let $\bbD$ be a 
cell structure of $\bC$. 
Then we write 
$F_{b,h}^{\bbD^{\heartsuit}}\bC$ 
for the relative category 
$(F_{b,h}^{\bbD^{\heartsuit}}\cC,w_{\st})$ and call 
it the {\it d\'evissage space of $\bC$ with respect to 
the cell structure $\bbD$}. 

We say that {\it $\bC$ satisfies derived d\'evissage condition} 
({\it with respect to the cell structure $\bbD$}) 
if the relative complicial exact functor 
${(-)}_{\infty}\colon F_{b,h}^{\bbD^{\heartsuit}}\cC\to\bC$ 
induces an equivalence 
$\Ho(\inn{F_{b,h}^{\bbD^{\heartsuit}}\cC}{F_{b,h}\bC})\isoto\Ho(\bC)$
of triangulated categories. 

Moreover let $\cF$ be a full subcategory of $\cC$. 
Then we write $F_{b,h}^{\bbD^{\heartsuit}}\bC|_{\cF}$ for 
the relative complicial exact category 
$(F_{b,h}^{\bbD^{\heartsuit}}{\langle\cF \rangle}_{\nullclass,\frob},
w_{\st}|_{F_{b,h}^{\bbD^{\heartsuit}}{\langle\cF \rangle}_{\nullclass,\frob}})$ and 
call it the {\it relative d\'evissage space of $\cF$ in 
$\cC$ with respect to the cell structure $\bbD|_{\cF}$}. 
We say that $\cF$ satisfies 
{\it derived d\'evissage condition with respect to the cell structure $\bbD$ 
in $\bC$} if the relative complicial exact functor 
${(-)}_{\infty}\colon F_{b,h}^{{(\bbD|_{\cF})}^{\heartsuit}}{\langle\cF \rangle}_{\nullclass,\frob}\to {\langle\cF \rangle}_{\nullclass,\frob}$ induces an equivalence $\Ho(\inn{F_{b,h}^{{(\bbD|_{\cF})}^{\heartsuit}}{\langle\cF \rangle}_{\nullclass,\frob}} F_{b,h}\bC)\isoto \Ho(\inn{\cF}\bC)$ of triangulated categories. 
\end{para}

\begin{para}
\label{df:devissage decomposition}
{\bf Definition (D\'evissage decomposition).}\ \ 
Let $\cC$ be a category with cofibration and let $\cA$ be a full subcategory of $\cC$ which is closed under isomorphisms 
and let $f\colon x\rinf y$ be a cofibration and let 
$m$ be a non-negative integer. 
We write $[m]$ for the linear ordered set $[m]=\{k\in\bbZ;0\leq k\leq m\}$ with the usual ordering. 
A {\it $(m,\cA)$-d\'evissage decomposition of $f$} or 
simply {\it $\cA$-d\'evissage decomposition of $f$} 
is a functor $z\colon [m]\to \cC$ such that 
\begin{enumerate}
\enumidef
\item
$z(0)=x$ and $z(m)=y$, 
\item
for all $0\leq i\leq m-1$, 
$z(i\geq i+1)\colon z(i)\to z(i+1)$ 
is a cofibration, 
\item
$z(i+1)/z(i)$ is in $\cA$ for all $0\leq i\leq m-1$ and 
\item
the compositions 
$z(m-1\leq n)z(m-2\leq m-1)\cdots z(0\leq 1)$ is equal to $f$. 
\end{enumerate}
\end{para}

\begin{para}
\label{df:classical Devissage condition}
{\bf (D\'evissage conditions).}\ \ 
Let $\cC$ be a category with cofibration and let $\cA$ a 
full subcategory of $\cC$ which is closed under isomorphisms. 

\begin{enumerate}
\enumidef
\item[]
{\bf (Quillen's d\'evissage condition).}\ \ 
We say that 
{\it $\cA$ satisfies Quillen's d\'evissage condition} ({\it in $\cC$}) 
if for any object $x$ in $\cC$, 
the canonical morphism $0\to x$ admits a $\cA$-d\'evissage decomposition.

Recall the definition of Serre radical ${}^S\!\!\!\sqrt{\ }$ from \ref{df:Serre subcategory} and assume that $\cC$ is a noetherian abelian category and 
$\cA$ is a topologizing subcategory 
(see \ref{df:toplogizing subcat}) of $\cC$. 
Then we can show that 
$\cA$ satisfies Quillen's d\'evissage condition in $\cC$ 
if and only if $\cC={}^S\!\!\!\sqrt{\cA}$. 
(\cf \cite[3.1]{Her97}, \cite[2.2]{Gar09}.) 

\item[]
{\bf (Raptis' d\'evissage condition).}\ \ (\cf \cite[5.1]{Rap18}.)\ \ 
We say that 
{\it $\cA$ satisfies Raptis' d\'evissage condition} ({\it in $\cC$}) 
if for any cofibration $f\colon x\rinf y$ in $\cC$, 
there exists $\cA$-d\'evissage decomposition of $f$.
\end{enumerate}
\end{para}

Relationship between Quillen's and Raptis' d\'evissage conditions for 
abelian categories 
is summed up with the following lemma in \cite[5.12, 5.13]{Rap18}.

\begin{para}
\label{lem:Raptis characterization}
{\bf Lemma.}\ \ 
{\it
Let $\cA$ be an abelian category and let $\cB$ be a full subcategory 
closed under quotients. 
Namely for an object $x$ in $\cA$, if there exists an epimorphism 
$y\rdef x$ with $y\in\Ob\cB$, then $x$ is also in $\cB$. 
Then the following conditions are equivalent.
\begin{enumerate}
\enumidef
\item
$\cB$ satisfies Quillen's d\'evissage condition in $\cA$.
\item
$\cB$ satisfies Raptis' d\'evissage condition in $\cA$.
\item
$\Ch_b(\cB)$ satisfies Raptis' d\'evissage condition in $\Ch_b(\cA)$.
\end{enumerate}
}
\qed
\end{para}

Recall the definition of cell structures 
from Definition~\ref{df:cell str}.

\begin{para}
\label{df:homotopy devissage decomposition}
{\bf Definition (Homotopy d\'evissage decomposition).}\ \ 
Let $\bC=(\cC,w_{\bC})$ be a relative complicial exact category and 
let $\bbD$ be a cell structure of $\bC$ and let 
$f\colon x\rinf y$ be an admissible monomorphism in $\cC$ and 
let $m>n$ be integers. 
A {\it $([m,n],\bbD)$-homotopy d\'evissage decomposition of $f$} or 
simply {\it $\bbD$-homotopy d\'evissage decomposition of $f$} 
is a pair $(z,a)$ consisting of an object $z$ in 
$F_{[m,n]}^{\bbD}\cC$ and a morphism $a\colon z_n=z_{\infty}\to y$ in 
$w_{\bC}$ such that 
\begin{enumerate}
\enumidef
\item
$x=z_m$ and 
\item
the composition 
\begin{equation}
\label{eq:Waldhausen devissage factorization}
x=z_{m}\onto{i^z_{m}}z_{m+1}\onto{i^z_{m+1}}z_{m+2}\onto{i^z_{m+2}}\cdots 
\onto{i^z_{n-1}}z_n=z_{\infty}\onto{a}y
\end{equation}
is equal to $f$. 
\end{enumerate}
\end{para}

Recall the definition of cell structures associated with full subcategories 
from \ref{ex:cell structure ass with full subcat}. 

\begin{para}
\label{df:homotopy devissage condition}
{\bf Definition (Homotopy d\'evissage condition).}\ \ 
Let $\bC=(\cC,w_{\bC})$ be a relative complicial exact category and 
let $\cA$ be a full subcategory of $\cC$. 
We say that $\cA$ satisfies the {\it homotopy d\'evissage condition} 
(in $\bC$) if for any morphism $f\colon x\to y$ in $\cC$, 
there exists a $\bbD_{\cA}$-homotopy d\'evissage decomposition of $f$. 
\end{para}

Relationship between Raptis' and homotopy d\'evissage conditions 
is summed up with the following lemma.

\begin{para}
\label{lem:Raptis and homotopy devissage}
{\bf Lemma.}\ \ 
{\it
Let $\bC=(\cC,w)$ be a relative complicial exact category and let 
$\cA$ be a full subcategory closed under isomorphisms and let 
$f\colon x\to y$ be a morphism in $\cC$. Then

\sn
$\mathrm{(1)}$ 
If $z$ is a $(m,\cA)$-d\'evissage decomposition of $f$, then 
the pair $(z,\id_z)$ is a 
$([0,m],\bbD_{\cA})$-d\'evissage decomposition of $f$. 

\sn
$\mathrm{(2)}$ 
If $\cA$ satisfies Raptis' d\'evissage condition in $\cC$, 
then $\cA$ satisfies homotopy d\'evissage condition in $\bC$. 
}
\end{para}

\begin{proof}
$\mathrm{(1)}$ 
For each $0\leq i\leq m-1$, the canonical morphism 
$\Cone z(i\leq i+1)\to z(i+1)/z(i)$ is in $w$ by 
Corollary~\ref{cor:Cones invariance by co-base change} $\mathrm{(1)}$. 
Thus $\Cone z(i\leq i+1)$ is in ${\langle \cA\rangle}_{\nullclass,w}$. 

\sn
$\mathrm{(2)}$ 
It follows from $\mathrm{(2)}$. 
\end{proof}

\begin{para}
\label{df:Waldhausen devissage conditions}
{\bf Definition (Waldhausen's d\'evissage condition).}\ 
(\cf \cite[\S 7 Hypothesis]{Wal85}) 
Let $\bC=(\cC,w)$ be a 
relative complicial exact category and let 
$\bbD=\{\calD_{\leq n},\calD_{\geq m} \}_{n,m\in\bbZ}$ 
be a cell structure of $\bC$. 
We say that {\it the pair $\bC=(\cC,w)$ satisfies 
Waldhausen's d\'evissage condition} 
({\it with respect to a cell structure $\bbD$}) 
if the following two conditions hold.
\begin{itemize}
\item
$\displaystyle{\cC=\bigcup_{m\in\bbZ}\calD_{\geq m}}$.
\item
For any $(m,\bbD)$-connected morphism $f\colon x\to y$ in $\cC$, 
there exits a $([m,n],\bbD)$-homotopy d\'evissage decomposition of $f$ for 
a suitable integer $n\geq m$.  
\end{itemize}
\end{para}

\begin{para}
\label{rem:Waldhausen devissage condition}
{\bf Remark.}\ \ 
$\mathrm{(1)}$ 
In the second assertion in Waldhausen's d\'evissage conditions, 
by replacing $z$ with $\fr_{[m,n]}z$ 
(for definition of $\fr_{[m,n]}$, see \ref{df:rab}) 
and $a\colon z_{\infty}\to y$ with the composition 
of $a$ with the projection ${\fr_{[m,n]}z}_{\infty} \to z_{\infty}$, 
we shall assume that 
$i_j^z$ is a Frobenius admissible monomorphism for $m\leq j\leq n-1$. 

\sn
$\mathrm{(2)}$ 
If $\bbD$ is a bounded (\ref{df:bounded cell structure}) 
and ordinary (\ref{df:ordinary cell structure}) cell structure, then 
we can replace the second assertion in Waldhausen's d\'evissage 
conditions with the following weaker assertion:

\sn
For any $(m,\bbD)$-connected morphism $f\colon x\to y$, 
there exists a factorization $x\onto{i}z\onto{f'}y$ of $f$ such that $f'$ is 
$(m+1,\bbD)$-connected and $\Cone i$ is in $\calD_{m+1}$. 

\sn
$\mathrm{(3)}$ 
The following condition implies assertion in $\mathrm{(2)}$.

\sn
For any integer $m$ and any object $x$ in $\calD_{\geq m}$, 
there exists an object $y$ in $\calD_{\geq m+1}$ and 
a morphism $g\colon x\to y$ in $\cC$ such that $\Cone g$ is in $\calD_{m+1}$.
\end{para}

\begin{proof}
$\mathrm{(2)}$ 
Let $f\colon x\to y$ be a morphism in $\cC$. 
By bounded assumption of $\bbD$, we shall assume that 
$\Cone f$ is in $\calD[m+1,n]$ for some pair of integers $m+1<n$. 
Then by hypothesis, there exists a factorization 
$z_{m}=x \onto{i^z_{m}}z_{m+1} \onto{f'}y$ with 
$\Cone i^z_{m}$ is in $\calD_{m+1}$ and 
$\Cone f'$ is $(m+1,\bbD)$-connected. 
Then by \ref{lem:cell structure fundamental result} $\mathrm{(11)}$, 
$\Cone f'$ is in 
$\calD^{[m+2,n]}$. 
Now by proceeding by induction, we will finally obtain a factorization 
$\mathrm{(\ref{eq:Waldhausen devissage factorization})}$ 
of $f$ such that $a$ is in $w$ by ordinarily of $\bbD$.

\sn
$\mathrm{(3)}$ 
Let $f\colon x\to y$ be a morphism in $\cC$ with 
$\Cone f$ is in $\calD_{\geq m+1}$. 
Then applying assumption to $\Cone f$, 
there exists an object $z$ in $\calD_{\geq m+2}$ and 
a morphism $g\colon \Cone f\to z$ with $\Cone g$ is in 
$\calD_{m+2}$. 
Then there exists a commutative diagram of distinguished triangles 
in $\pi_0\cC$ below
$$
\footnotesize{
\xymatrix{
x \ar[r]^f \ar[d]_{i'} & y \ar[r]^{\kappa_f} \ar[d]_{\id_y} & \Cone f 
\ar[r] \ar[d]_g & Tx \ar[d]^{Ti'}\\
u \ar[r]_{f'} \ar[d]_{\kappa_{i'}} & y \ar[r]_{g\kappa_f} \ar[d] 
& z \ar[d]_{\kappa_g} \ar[r] & Tu \ar[d]^{T\kappa_{i'}}\\
\Cone i' \ar[r] & 0 \ar[r] & \Cone g \ar[r] & T\Cone i'.
}}
$$
Thus there exists a $C$-homotopy $H\colon f\Rightarrow_C f'i'$. 
We shall assume that $f=f'i'$ and $\Cone f'$ is 
$C$-homotopy equivalent to $z$ 
by replacing $i'$ with $\xi_{1,i'}\colon x\to u\oplus Cx$ and $f'$ with 
$\displaystyle{\tiny{\begin{pmatrix}f'\\ -H\end{pmatrix}}\colon 
u\oplus Cx\to y}$. 
Hence $f'$ is $(m+1,\bbD)$-connected. 
By $3\times 3$-lemma $\Cone g$ is $C$-homotopy equivalent to $T\Cone i'$. 
Therefore $\Cone i'$ is in $\calD_{m+1}$ 
by \ref{lem:cell structure fundamental result} $\mathrm{(13)}$. 
\end{proof}

\begin{para}
\label{prop:sufficient conditions of derived devissage}
{\bf Proposition.}\ \ 
{\it
Let $\bC=(\cC,w)$ be a normal ordinary relative complicial exact category 
and let $\cA$ be a full subcategory of $\cC$ and 
$\bbD$ be a cell structure of $\bC$. Then
\begin{enumerate}
\enumidef
\item
{\rm(\cf \cite[1.7.2]{Wal85}.)}\ \ 
Waldhausen's d\'evissage condition 
with respect to $\bbD$ implies 
derived d\'evissage condition with respect to $\bbD$.

\item
Homotopy d\'evissage condition with respect to $\cA$ 
implies derived d\'evissage condition with respect to $\bbD_{\cA}$. 
\end{enumerate}
}
\end{para}

\begin{proof}
We only give a proof of assertion $\mathrm{(1)}$. 
A proof of assertion $\mathrm{(2)}$ is similar. 
We will show that ${(-)}_{\infty}\colon F_{b,h}^{\bbD}\cC\to \cC$ 
satisfies \textbf{App 2} in 
\ref{rem:derived equivalence}. Then by \ref{rem:derived equivalence} 
$\mathrm{(2)}$ $\mathrm{(b)}$, 
it turns out that $\bC$ satisfies 
derived d\'evissage condition with respect to $\bbD$. 
We proceed by induction on $\dim x$ 
(for definition of $\dim x$, see \ref{df:filtered objects}). 
If $\dim x=-1$, then applying Waldhausen's d\'evissage condition 
to the morphism $0\to y$, we obtain the result. 

For $m\geq 0$, 
by applying the inductive hypothesis 
to the composition 
$\sigma_{\leq m-1}x \to x\to y$, 
we obtain a level admissible monomorphism 
$i'\colon\sigma_{\leq m-1}x \rinf z''$ and a morphism 
$a''\colon z''_{\infty}\to y$ in $w$ 
which makes diagram below commutative. 
$$
\footnotesize{
\xymatrix{
{(\sigma_{\leq m-1}x)}_{\infty} \ar[r] \ar@{>->}[d]_{i'_{\infty}} & x_{\infty} \ar[d]^f\\
z''_{\infty} \ar[r]_{a''} & y
}}
$$
Then we will show the following assertions.

\sn
$\mathrm{(a)}$ 
For each integer $i$, the compositions 
${z''}_i\onto{i_i^{z''}}{z''}_{i+1}\onto{i_{i+1}^{z''}}\cdots \to{z''}_{\infty}\onto{a''}y$ is $(i-1,\bbD)$-connected.

\sn
$\mathrm{(b)}$ 
The canonical morphism $x_m\sqcup_{x_{m-1}}{z''}_{m-1}\to y$ is $(m-1,\bbD)$-connected. 

\begin{proof}[Proof of $\mathrm{(a)}$]
We proceed by descending induction on $i$. 
For $n=\dim z''$, 
the morphism ${z''}_n ={z''}_{\infty}\onto{a''} y$ is in $w$, 
in particular $(n-1,\bbD)$-connected. 
In general, applying \ref{lem:cell structure fundamental result} $\mathrm{(7)}$ to the pair of 
composable morphisms ${z''}_i\onto{i^{z''}_i}{z''}_{i+1}$ and 
${z''}_{i+1}\to y$, 
it turns out that the composition ${z''}_i\to y$ 
is $(i-1,\bbD)$-connected by 
inductive hypothesis. 
\end{proof}

\begin{proof}[Proof of $\mathrm{(b)}$]
By \ref{lem:invariance connectedness under co-base change}, 
$\Cone(z''_{m-1}\rinf x_m\sqcup_{x_{m-1}}z''_{m-1})$ is in $\calD_m$. Now by applying \ref{lem:cell structure fundamental result} $\mathrm{(8)}$ 
to the pair of 
composable monomorphisms $\displaystyle{z''_{m-1}\rinf x_m\sqcup_{x_{m-1}}{z''}_{m-1}\to y}$, 
it turns out that the canonical morphism 
$\displaystyle{x_m\sqcup_{x_{m-1}}{z''}_{m-1}\to y}$ is 
$(m-1,\bbD)$-connected. 
\end{proof}

By applying Waldhausen's d\'evissage condition to the morphism 
$x_m\sqcup_{x_{m-1}}{z''}_{m-1}\to y$, 
there exists an object $z'$ in $F_{b,h,\geq m-1}^{\bbD}\cC$ 
and a morphism $a'\colon {z'}_{\infty}\to y$ in $w$ 
such that 
$\displaystyle{{z'}_{m-1}=x_m\sqcup_{x_{m-1}}{z''}_{m-1} }$ 
and the composition 
${z'}_{m-1}\onto{i^{z'}_{m-1}}{z'}_{m}\onto{i^{z'}_{m}}\cdots \to 
z'_{\infty}\onto{a'} y$ 
is equal to the morphism 
$\displaystyle{x_m\sqcup_{x_{m-1}}{z''}_{m-1}\to y}$. 
We write $I'$ for the composition of 
$\displaystyle{{z''}_{m-1}\rinf x_m\sqcup_{x_{m-1}}{z''}_{m-1}\rinf {z'}_m}$. 
By \ref{lem:invariance connectedness under co-base change}, 
$\displaystyle{\Cone ({z''}_{m-1}\rinf x_m\sqcup_{x_{m-1}}{z''}_{m-1})}$ 
is in $\calD_m$. 
Then it turns out that $\Cone I'$ is also in $\calD_m$ 
by applying \ref{lem:cell structure fundamental result} 
$\mathrm{(7)}$ and $\mathrm{(10)}$ to the pair of composable morphisms 
$\displaystyle{{z''}_{m-1}\rinf x_m\sqcup_{x_{m-1}}{z''}_{m-1}\rinf {z'}_m}$. 

We define $z$, $I\colon x\rinf z$ and $a\colon z_{\infty}\to y$ 
to be an object in $F_{b,h}^{\bbD}\cC$ and a level 
admissible monomorphism and 
a morphism in $w$ 
such that the composition 
$x_{\infty}\overset{I}{\rinf} z_{\infty}\onto{a} y$ is equal to $f$ 
by setting
$$z_k:=
\begin{cases}
{z''}_k & \text{if $k\leq m-1$}\\
{z'}_k & \text{if $k\geq m$},
\end{cases}\ \ \ 
i^z_k:=
\begin{cases}
i^{z''}_k & \text{if $k\leq m-2$}\\
I' & 
\text{if $k=m-1$}\\
i^{z'}_k & \text{if $k\geq m$}.
\end{cases}
$$
$$
I_k:=
\begin{cases}
i'_k & \text{if $k\leq m-1$}\\
\displaystyle{x_m\to x_m\sqcup_{x_{m-1}}z''_{m-1}\to z_k} & \text{if $k\geq m$}
\end{cases}\ \ \ 
\text{and $a:=a'$}.
$$
Now we complete the proof. 
\end{proof}

\begin{para}
\label{df:weight structure}
{\bf (Weight structure).}\ \ 
The notion of weight structures on triangulated categories 
is introduced by Bondarko \cite{Bon10} and \cite{Bon12} and 
in this article we will use homological weight structure 
as in \cite{Fon18}. 
Let $\cT$ be a triangulated category. 
A {\it homological weight structure} $\cW$ of $\cT$ is 
an additive full subcategory of $\cT$ closed under retract in $\cT$ 
which satisfies the following properties:\\
We set $\cW_{\leq n}:=\Sigma^n(\cW)$ and 
$\cW_{\geq n}:=\Sigma^{n-1}({}^{\perp}\cW)$.

\begin{itemize}
\item[]
{\bf (Suspension closed).}\ \ 
$\cW_{\geq 1}\subset\cW_{\geq 0}$ and $\cW_{\leq -1}\subset\cW_{\leq 0}$.
\item[]
{\bf (Weight decomposition).}\ \ 
For an object $x$ in $\cT$, there exists a distinguished triangle 
$b\to x\to a \to \Sigma b$ 
with $b\in\cW_{\leq 0}$ and $a\in\cW_{\geq 1}$.
\end{itemize}

Let $\bC=(\cC,w)$ be a relative complicial exact category and let 
$\cW$ be a homological weight structure on $\Ho(\bC)$. 
Then we define $\calD_{\geq n}^{\cW}$ and $\calD_{\leq m}^{\cW}$ 
to be full subcategories of $\Ho(\bC)$ by setting
\begin{equation}
\label{eq:calDgeqnW df}
\Ob\calD_{\geq n}^{\cW}:=\{x\in\Ob\cC;\cQ_{\bC}(x)\text{ is in $\cW_{\geq n}$}\},
\end{equation}
\begin{equation}
\label{eq:calDleqmW df}
\Ob\calD_{\leq n}^{\cW}:=\{x\in\Ob\cC;\cQ_{\bC}(x)\text{ is in $\cW_{\leq n}$}\}.
\end{equation}
Here $\cQ_{\bC}\colon\cC\to\Ho(\bC)$ is the canonical functor. 
Then it turns out that a family 
$\bbD^{\cW}:=\{\calD^{\cW}_{\geq n},\calD^{\cW}_{\leq m}\}_{n,m}$ is 
a cell structure of $\bC$ by \cite[1.3.3]{Bon10}. 
We call $\bbD^{\cW}$ 
the {\it cell structure of $\bC$ associated with a homological weight structure $\cW$ on $\Ho(\bC)$}. 
In Lemma~\ref{lem:decomposition cond and devissage} below, we see a relationship between weight decomposition condition and Waldhausen's d\'evissage condition with respect to a cell structure 
associated with a homological weight structure.
\end{para}

\begin{para}
\label{lem:decomposition cond and devissage}
{\bf Lemma.}\ \ 
{\it
Let $\bC=(\cC,w)$ be a relative complicial exact category and 
let $\cW$ be a homological weight structure of $\Ho(\bC)$ and 
let $\bbD^{\cW}$ be an associated cell structure of $\cW$ in $\bC$. 
If $\cW$ is bounded below, then $\bC$ satisfies 
Waldhausen's d\'evissage condition with respect to the cell 
structure $\bbD^{\cW}$. 
}
\qed
\end{para}

\begin{proof}
We shall check condition $\mathrm{(3)}$ 
in Remark~\ref{rem:Waldhausen devissage condition}. 
Let $x$ be an object in $\calD^{\cW}_{\geq m}$. 
Then there exists an object $y$ in $\calD^{\cW}_{\geq m+1}$ 
and a morphism $g\colon x\to y$ such that the image 
of $T^{-1}\Cone g$ in $\Ho(\bC)$ is in $\cW_{\leq m}$. 
Then $\Cone g$ is in $\calD_{m+1}^{\cW}$ by 
Lemma~\ref{lem:cell structure fundamental result} 
$\mathrm{(2)}$ and $\mathrm{(5)}$. 
\end{proof}

\subsection{Derived quasi-split sequences}
\label{subsec:derived quasi-split sequence}

In this subsection, 
we give a typical example of a derived quasi-split sequence. 
We start by recalling the definition of derived exact sequences 
of relative exact categories from \cite{Moc13b}. 

\begin{para}
\label{df:derived exaxt seq}
{\bf (Derived exact sequences).}\ \ 
We say that a sequence of triangulated categories 
and triangle functors 
\begin{equation}
\label{eq:seq of tri cat}
(\cT,\Sigma)\onto{i}(\cT',\Sigma')\onto{p}(\cT'',\Sigma'')
\end{equation}
is {\it exact} if $pi$ is the zero functor and $i$ and $p$ induce 
the equivalences of triangulated categories 
$\cT\isoto\Ker p$ and $\cT'/\Ker p\isoto \cT''$. 
We say that the sequence $\mathrm{(\ref{eq:seq of tri cat})}$ is 
{\it weakly exact} 
if the functor $i$ is fully faithful and if the composition $pi$ is 
the zero morphism and if the induced morphism 
$\cT'/\cT \to \cT''$ is {\it cofinal}. 
The last condition means that it is fully faithful and for any object 
$x$ in $\cT''$, there exists a pair of objects $y$ in $\cT''$ and 
$z$ in $\cT'$ such that $x\oplus y$ is isomorphic to $z$. 

We say that a sequence of relative exact categories 
and relative exact functors 
\begin{equation}
\label{eq:derived spltting sequence}
\bE\onto{i}\bE'\onto{p}\bE''
\end{equation} 
is {\it derived exact} (resp. {\it derived weakly exact}) if 
the induced sequence 
\begin{equation}
\label{eq:associated quasi-split seq}
\calD_b(\bE)\onto{\calD_b(i)}\calD_b(\bE')\onto{\calD_b(p)}\calD_b(\bE'')
\end{equation}
of triangulated categories and triangle functors 
is exact (resp. weakly exact). 
\end{para}

We can show the following.

\begin{para}
\label{lem:five lemma for derived equiv}
{\bf Lemma.}\ \ 
{\it
Let 
$$
\xymatrix{
\bE \ar[r]^i \ar[d]_a & \bE' \ar[r]^p \ar[d]^b & \bE'' \ar[d]^c\\
\bF \ar[r]_j & \bF' \ar[r]_q & \bF''
}
$$ 
be a commutative diagram of derived exact sequences of 
relative exact categories. 
If $a$ and $b$ {\rm(}resp. $b$ and $c${\rm)} are 
derived equivalences, then $c$ {\rm(}resp. $a${\rm)} 
is also a derived equivalence.
}\qed
\end{para}

By \ref{prop:comp hom and der}, 
we obtain the following result.

\begin{para}
\label{lem:Chb is exact}
{\bf Lemma.}\ \ 
{\it
The functor 
$\Ch_b\colon\RelEx\to\RelComp_{\stnor}$ preserves 
derived exact sequences and derived weakly exact sequences. 
Namely if the sequence $\mathrm{(\ref{eq:derived spltting sequence})}$ 
of relative exact categories and 
relative exact functors 
is derived exact 
{\rm(}resp. derived weakly exact{\rm)}, 
then the sequence of relative complicial exact categories 
and relative complicial exact functors 
\begin{equation}
\label{eq:ass Chb exact seq}
\Ch_b(\bE)\onto{\Ch_b(i)}\Ch_b(\bE')\onto{\Ch_b(p)}\Ch_b(\bE'') 
\end{equation}
is also derived exact {\rm(}resp. derived weakly exact{\rm)}.
}\qed
\end{para}

\begin{para}
\label{df:triangle adjoint}
{\bf (Triangle adjoint).}\ (\cf \cite[\S 8]{Kel96}.)\ \ 
Let $(R,\rho)\colon(\cS,\Sigma) \to (\cT,\Sigma')$ 
and $(L,\lambda)\colon\cT \to \cS$ be two triangle functors 
between triangulated categories $(\cS,\Sigma)$ and $(\cT,\Sigma')$ 
such that 
$L$ is left adjoint to $R$. 
Let $A\colon\id_{\cT} \to RL$ and $B\colon LR\to \id_{\cS}$ 
be adjunction morphisms. 
For any objects $x$ in $\cT$ and $y$ in $\cS$, 
we write $\mu(x,y)$ for the bijection 
$\Hom_{\cS}(Lx,y) \to \Hom_{\cT}(x,Ry)$, 
$f \mapsto (Rf)\cdot A_x$. 
We say that $(L,\lambda)$ (resp. $(R,\rho)$) is 
{\it left} (resp. {\it right}) 
{\it triangle adjoint} to $(R,\rho)$ (resp. $(L,\lambda)$) 
if for any objects $x$ in $\cT$ and $y$ in $\cS$, 
the diagram below is commutative.
$$
{\footnotesize{\xymatrix{
\Hom_{\cS}(Lx,y) \ar[r]^{\!\!\!\!\!\Sigma} \ar[d]_{\mu(x,y)} & 
\Hom_{\cS}(\Sigma Lx,\Sigma y) \ar[r]^{\Hom(\lambda x,\Sigma y)} & 
\Hom_{\cS}(L\Sigma' x,\Sigma y) \ar[d]^{\mu(\Sigma' x,\Sigma y)}\\
\Hom_{\cT}(x,Ry) \ar[r]_{\!\!\Sigma'} & 
\Hom_{\cT}(\Sigma x,\Sigma' Ry) \ar[r]_{\Hom(\Sigma' x,R\Sigma y)} & 
\Hom_{\cT}(\Sigma' x,R\Sigma y).
}}}$$
$\mu(\Sigma',\Sigma)\Hom(\lambda,\Sigma)\Sigma=
\Hom(\Sigma',\rho^{-1})\Sigma'\mu$. 

Let 
$(g,\rho')\colon\cT' \to \cT$ 
be a triangle functor, 
and let $f$ be a left (resp. right) adjoint of $g$ 
with $\Phi\colon fg \to \id_{\cT}$ and $\Psi\colon \id_{\cT'}\to gf$ 
(resp. $\Phi\colon\id_{\cT} \to fg$ and $\Psi\colon gf\to \id_{\cT'}$) 
adjunction morphisms 
and $\lambda':=(\Phi\ast\Sigma' f)
\cdot(f\ast{\rho'}^{-1}\ast f)\cdot(f\Sigma\ast \Psi)$ 
(resp. $\lambda':={\{(f\Sigma\ast \Psi)\cdot(f\ast\rho'\ast f)
\cdot(\Phi\ast\Sigma' f)\}}^{-1}$). 
Then the pair $(f,\lambda')$ is a triangle functor and 
it is a left (resp. right) 
triangle adjoint of $(g,\rho')$ (\cf \cite[8.3]{Kel96}).
\end{para}

\begin{para}
\label{df:quasi split sequence}
{\bf (Quasi-split sequences).}\ (\cf \cite[6.3]{Moc13b}.)\ \ 
We say that a sequence $\mathrm{(\ref{eq:seq of tri cat})}$ 
is {\it right} (resp. {\it left}) {\it quasi-split} 
if both $i$ and $p$ admit right (resp. left) adjoint functors 
$q\colon \cT'\to \cT$ and $j\colon \cT''\to \cT'$ 
with adjunction morphisms 
$iq\onto{A}\id_{\cT'}$, $\id_{\cT}\onto{B}qi$, 
$\id_{\cT'}\onto{C}jp$ and $pj\onto{D}\id_{\cT''}$ 
(resp. $\id_{\cT'}\onto{A}iq$, $qi\onto{B}\id_{\cT}$, 
$jp\onto{C}\id_{\cT'}$ and $\id_{\cT''}\onto{D}pj$) 
respectively such that $B$ and $D$ are natural equivalences 
and if there exists a triangulated natural transformation 
$jp\onto{E}\Sigma'iq$ (resp. $iq\onto{E}\Sigma' jp$) 
such that a triangle $(A,C,E)$ (resp. $(C,A,E)$) 
is a distinguished triangle in $\cT'$. 
We call a system $(j,q,A,B,C,D,E)$ (resp. $(j,q,C,D,A,B,E)$) 
or shortly $(j,q,E)$ a 
{\it right} (resp. {\it left}) {\it splitting of a sequence $(i,p)$}.

We say that a sequence 
$\mathrm{(\ref{eq:derived spltting sequence})}$ 
of relative exact categories and 
relative exact functors 
is {\it derived right} (resp. {\it left}) {\it quasi-split} 
if the induced sequence 
$\mathrm{(\ref{eq:associated quasi-split seq})}$ of 
bounded derived categories 
is a right (resp. left) quasi-split sequence of triangulated categories. 
\end{para}

\begin{para}
\label{rem:derived quasi-split}
{\bf Remark.}\ \ 
If the sequence $\mathrm{(\ref{eq:seq of tri cat})}$ 
of triangulated categories and triangle functors is 
a right quasi-split sequence with $(j,q,E)$ a right quasi-splitting, 
then the sequence $(j,q,E)$ is a left quasi-split sequence with 
a left quasi-splitting $(i,p,E)$. 

Thus if the sequence $\mathrm{(\ref{eq:derived spltting sequence})}$ of 
relative exact categories and relative exact functors 
is a derived right quasi-split sequence and if for a pair of 
relative exact functors $\bE'\onto{q}\bE$ and $\bE''\onto{j}\bE'$, 
induced triangle functors $\calD_b(\bE')\onto{\calD_b(q)}\calD_b(\bE)$ and 
$\calD_b(\bE'')\onto{\calD_b(j)}\calD_b(\bE')$ give a 
part of a right quasi-splitting of the sequence 
$\mathrm{(\ref{eq:associated quasi-split seq})}$,  
then the sequence $\bE''\onto{j}\bE'\onto{q}\bE$ 
is a derived left quasi-split sequence. 

For simplicity, 
we will only mention the cases for derived right quasi-split 
sequences and similar statements also hold for 
the left case. 
\end{para}

By \ref{prop:comp hom and der}, we can show the following. 

\begin{para}
\label{lem:Chb preserves quasi-split seq}
{\bf Lemma.}\ \ 
{\it
The functor $\Ch_b\colon\RelEx\to\RelComp_{\stnor}$ 
preserves right quasi-split sequences. 
Namely if the sequence of relative exact categories and relative exact 
functors $\mathrm{(\ref{eq:derived spltting sequence})}$ 
is a derived right quasi-split, then 
the sequence $\mathrm{(\ref{eq:ass Chb exact seq})}$ 
of relative complicial exact categories 
is also a derived right quasi-split sequence. 
Moreover if there exists $q\colon\bE'\to\bE$ 
{\rm(}resp. $j\colon\bE''\to\bE'${\rm)} such that 
$\calD_b(q)$ (resp. $\calD_b(j)$) 
gives a part of a right quasi-splitting of the 
sequence $\mathrm{(\ref{eq:associated quasi-split seq})}$ 
of triangulated categories, 
then $\calD_b(\Ch_b(q)) $ {\rm(}resp. $\calD_b(\Ch_b(j))${\rm)} 
also gives a part of right quasi-splitting of the sequence 
$\calD_b(\Ch_b(\bE))\onto{\calD_b(\Ch_b(i))}\calD_b(\Ch_b(\bE'))\onto{\calD_b(\Ch_b(p))}\calD_b(\Ch_b(\bE''))$.
}\qed
\end{para}

\begin{para}
\label{lem:quasi split is exact}
{\bf Lemma.}\ (\cf \cite[6.6]{Moc13b}.)\ \ 
{\it
$\mathrm{(1)}$ 
A right quasi-split sequence of triangulated categories 
and triangle functors is exact. 

\sn
$\mathrm{(2)}$ 
A derived right quasi-split sequence of relative exact categories and 
relative exact functors is derived exact.
}\qed
\end{para}

\begin{para}
\label{nt:setting of additivity}
{\bf (Conventions).}\ \ 
Let $\bC=(\cC,w)$ be a normal ordinary relative complicial exact category 
and let $\cS$ and $\cT$ be a pair of full subcategories of $\cC$. 
Assume that $\cS$ contains the zero objects. 
We set $\cC':=\cC^{[1]}_h(\cS,\cT_{\heq})$. 
We define $i\colon\cC\to \cC^{[1]}_h$, 
$q\colon \cC^{[1]}_h\to \cC$, 
$j\colon \cC\to\cC^{[1]}_h$ to be a strictly normal exact functors 
by sending an object $x$ in $\cC$ to $[x\onto{\id_x}x]$, 
an object $ $ in $\cC^{[1]}_h$ to $ $ and an object $y$ in 
$\cC$ to $[0\to y]$ in $\cC^{[1]}_h$ respectively. 
For simplicity we write $(p,c)$ for the complicial exact functor 
$(\Cone,\sigma^{\Cone})\colon\cC^{[1]}_h\to \cC$. 
We also denote the restriction of $i$, $q$ to full subcategories 
of $\cC$ and $\cC^{[1]}_h$ 
etc 
and the induced functors of $i$, $q$ on homotopy categories 
by the same letters $i$ and $q$ respectively. 
Moreover we denote the restriction of $w$ and $lw$ to 
full subcategories of $\cC$ and $\cC^{[1]}_h$ by the same letters 
$w$ and $lw$. 
\end{para}

\begin{para}
\label{prop:typical derived quasi-split sequence}
{\bf Proposition (Typical derived quasi-split sequence).}\ \ 
{\it
In the convention \ref{nt:setting of additivity}, the sequence 
\begin{equation}
\label{eq:typical dqss}
({\langle\cS \rangle}_{\nullclass,\frob},w)\onto{i}
({\langle \cC'\rangle}_{\nullclass,\frob},lw)
\onto{p}
({\langle\cT_{\heq} \rangle}_{\nullclass,\frob},w)
\end{equation}
is a right derived quasi-split sequence and the sequence
\begin{equation}
\label{eq:typical left dqss}
({\langle\cT_{\heq} \rangle}_{\nullclass,\frob},w)\onto{j}
({\langle \cC'\rangle}_{\nullclass,\frob},lw)\onto{q}
({\langle\cS \rangle}_{\nullclass,\frob},w)
\end{equation}
is a left derived quasi-split sequence. 
}
\end{para}

\begin{proof}
We define $B\colon \id_{\cC^{[1]}_h}\to jp$ and 
$C\colon jp\to Tiq$ to be a pair of 
complicial natural transformations by 
setting for an object $x\onto{f}y$ in $\cC^{[1]}_h$, 
$B(x\onto{f}y):=(0,\kappa_f,-\mu_f)$ and 
$C(x\onto{f}y):=(0,\psi_f,0)$. 
For simplicity we denote $\Ho(B)$ and $\Ho(C)$ by 
the same letter $B$ and $C$ respectively. 
Moreover for an object $x\onto{f}y$ in $\cC^{[1]}_h$, 
we set $A(x\onto{f}y):=(\id_x,f,0)\colon [x\onto{\id_x}x]\to [x\onto{f}y]$ and for a morphism 
$(a,b,H)\colon [x\onto{f}y]\to [x'\onto{f'}y']$ 
in $\cC^{[1]}_h$, we set 
$A(a,b,H):=(0,-H,H\cdot c_x)\colon (a,b,H)(\id_x,f,0)\Rightarrow_C 
(\id_x',f',0)(a,a,0)$. 
Then $A$ induces a triangle natural transformation 
$iq\to \id_{\Ho(\cC^{[1]}_h)}$ on homotopy 
categories and we denote it by the same letter $A$. 
We can show that there exists a functorial distinguished 
triangle $iq\onto{A}\id_{\Ho(\cC^{[1]}_h)}\onto{B}jp\onto{C}Tiq$ and the equalities $qi=\id_{\cC}$, $pj=\id_{\cC}$, $B\ast j=\id_{j}$, 
$p\ast B=\id_p$, $A\ast i=\id_i$ 
and $q\ast A=\id_q$. 
Thus it turns out that the sequence $\mathrm{(\ref{eq:typical dqss})}$ 
is a derived right quasi-split sequence. 
A proof of the sequence $\mathrm{(\ref{eq:typical left dqss})}$ is similar. 
\end{proof}

\begin{para}
\label{rem:typical derived quasi-split sequence}
{\bf Remark.}\ \ 
In \ref{prop:typical derived quasi-split sequence}, 
replacing $\cC^{[1]}_h(\cS,\cT_{\heq})$ with $\cC^{[1]}(\cS,\cT_{\heq})$, 
similar statement also holds. 
The same proof works fine. 
\end{para}

\begin{para}
\label{ex:typical derived quasi-split sequence}
{\bf Example.}\ \ 
Let $\bC=(\cC,w)$ be a normal ordinary relative complicial exact 
category and let $a<b$ be a pair of integers and 
let $\bbD=\{\calD_{\leq n},\calD_{\geq m} \}_{n,m\in\bbZ}$ 
be a $C$-homotopy closed cell structure of $\bC$. 
Then the inclusion functor $F_{[a,b-1],h}\cC\rinc F_{[a,b],h}\cC$ and 
the functors $\Cone i_b^{(-)}$, ${(-)}_{b-1}\colon F_{[a,b],h}\cC\to \cC$ 
which sends an object $x$ to 
$\Cone i^x_{b-1}$ and $x_{b-1}$ respectively 
induce a right and a left derived flags below respectively by 
\ref{ex:pi0 of FbhC and underline FBhC}, 
\ref{lemdf:devissage filtration with level weak equivalences} and 
\ref{lemdf:acyclic devissage filtration with level weak equivalences}. 
\begin{equation}
\label{eq:devissage split seq right}
({\langle F_{[a,b-1],h}^{\bbD^{\heartsuit}}\cC \rangle}_{\nullclass,\frob},lw)\to 
({\langle F_{[a,b],h}^{\bbD^{\heartsuit}}\cC \rangle}_{\nullclass,\frob},lw)\to 
({\langle\calD_b \rangle}_{\nullclass,\frob},lw),
\end{equation}
\begin{equation}
\label{eq:devissage split seq left}
({\langle {(F_{[a,b-1],h}^{\bbD^{\heartsuit}}\cC \rangle}_{\nullclass,\frob})}^{w_{\st}},lw)\to 
({\langle {(F_{[a,b],h}^{\bbD^{\heartsuit}}\cC \rangle}_{\nullclass,\frob})}^{w_{\st}},lw)\to 
({\langle\calD_{b-1} \rangle}_{\nullclass,\frob},lw).
\end{equation}
\end{para}

\subsection{Derived flag of relative exact categories}
\label{subsec:flag}

In this subsection, we give examples of right flags on 
the category of d\'evissage filtrations with level weak equivalences 
which enable us to calculate the $\bbK$-theory of d\'evissage spaces in 
\S \ref{subsec:dev thm}. 
We start by recalling the notion of right derived flag of 
relative exact categories from \cite[6.11]{Moc13b}. 

\begin{para}
\label{df:derived flag}
{\bf (Derived flag).}\ \ 
Let $\cT$ be a triangulated category and let $\bE$ be a relative 
exact category. 
A {\it right} (resp. {\it left}) {\it flag} 
of $\cT$ is a finite sequence of fully faithful 
triangle functors 
\begin{equation}
\label{eq:right flag df}
\{0\}=\cT_0\onto{k_0}\cT_1\onto{k_1}\cdots \onto{k_{n-1}}\cT_n=\cT
\end{equation}
such that for each $0\leq i\leq n-1$, the canonical sequence 
$\cT_i\onto{k_i}\cT_{i+1}\to\cT_{i+1}/\cT_{i}$ of triangulated categories 
is a right (resp. left) quasi-split exact sequence  
(See \ref{df:quasi split sequence}) of triangulated categories. 

A {\it derived right} (resp. {\it left}) {\it flag} of $\bE$ 
is a finite sequence of 
fully faithful relative exact functors 
\begin{equation}
\label{eq:right derived flag df}
\{0\}=\bE_0\onto{k_0}\bE_1\onto{k_1}\cdots \onto{k_{n-1}}\bE_n=\bE
\end{equation}
such that it induces a right (resp. left) flag of $\calD_b(\bE)$
$$\{0\}=\calD_b(\bE_0)\onto{\calD_b(k_0)}\calD_b(\bE_1)\onto{\calD_b(k_1)}\cdots \onto{\calD_b(k_{n-1})}\calD_b(\bE_n)=\calD_b(\bE)$$
on bounded derived categories. 

A {\it derived right} (resp. {\it left}) {\it retractions} of 
the derived right (resp. left) flag 
$\mathrm{(\ref{eq:right derived flag df})}$ 
of $\bE$ is a family of 
relative exact functors 
$\calR=\{u_i\colon\bE_i\to\bE_{i-1} \}_{1\leq i\leq n}$ such that 
for each $1\leq i\leq n$, 
$\calD_b(u_i)\colon \calD_b(\bE_i)\to\calD_b(\bE_{i-1})$ 
gives a part of a right (resp. left) quasi-splitting of the sequence 
$\calD_b(\bE_{i-1})\onto{\calD_b(k_i)}\calD_b(\bE_i)\to\calD_b(\bE_i)/\calD_b(\bE_{i-1})$. 

An {\it associated graded relative exact categories} of the 
derived right (resp. left) flag $\mathrm{(\ref{eq:right derived flag df})}$ 
is a family of relative exact categories and relative exact functors 
$\cG=\{\bG_i,\ t_i\colon\bE_i\to\bG_i \}_{1\leq i\leq n}$ such that for each $1\leq i\leq n$, 
the sequence $\bE_{i-1}\onto{k_{i-1}}\bE_i\onto{t_i}\bG_i$ 
is a derived right (resp. left) quasi-split sequence. 
\end{para}

By \ref{lem:Chb preserves quasi-split seq}, 
we obtain the following lemma. 
We will only mention the case for right flags and 
similar statements for left flags also hold. 

\begin{para}
\label{lem:Ch_b preserves flag}
{\bf Lemma.}\ \ 
{\it
Let $\bE$ be a relative exact category and let 
$\mathrm{(\ref{eq:right derived flag df})}$ 
be a derived right flag of $\bE$. Then

\sn
$\mathrm{(1)}$ 
The sequence
\begin{equation}
\label{eq:Ch_b(flag)}
\{0\}=\Ch_b(\bE_0)\onto{\Ch_b(k_0)}\Ch_b(\bE_1)\onto{\Ch_b(k_1)}\cdots 
\onto{\Ch_b(k_{n-1})}\Ch_b(\bE_n)=\Ch_b(\bE)
\end{equation}
is a derived right flag of $\Ch_b(\bE)$.

\sn
$\mathrm{(2)}$ 
If a family $\calR=\{u_i\colon\bE_i\to\bE_{i-1} \}_{1\leq i\leq n}$ 
is a derived right retractions of 
the derived right flag $\mathrm{(\ref{eq:right derived flag df})}$ 
of $\bE$, 
then a family 
$\Ch_b(\calR)=\{\Ch_b(u_i)\colon\Ch_b(\bE_i)\to
\Ch_b(\bE_{i-1}) \}_{1\leq i\leq n}$ 
is a derived right retractions of 
the derived right flag $\mathrm{(\ref{eq:Ch_b(flag)})}$ of $\Ch_b(\bE)$.

\sn
$\mathrm{(3)}$ 
If a family $\cG=\{\bF_i,\ t_i\colon\bE_i\to\bF_i \}_{0\leq i\leq n}$ 
is an associated graded relative exact categories 
of a derived right flag $\mathrm{(\ref{eq:right derived flag df})}$, 
then the family 
$\Ch_b(\cG)=\{\Ch_b(\bF_i),\ \Ch_b(t_i)\colon
\Ch_b(\bE_i)\to\Ch_b(\bF_i) \}_{0\leq i\leq n}$ 
is an associated graded relative exact categories 
of a derived right flag $\mathrm{(\ref{eq:Ch_b(flag)})}$ of $\Ch_b(\bE)$. 
}\qed
\end{para}

\begin{para}
\label{ex:canonical flag}
{\bf Example.}\ \ 
Let $\bC=(\cC,w)$ be a normal ordinary relative complicial exact 
category and let $\bbD=\{\calD_{\leq n},\calD_{\geq m}\}_{n,\ m\in\bbZ}$ 
be a cell structure of $\bC$ and let $a<b$ be a pair of integers. 
For each integer $1\leq k\leq b-a+1$, 
we set 
$\cE_k:={\langle F_{[a,a+k-1],h}^{\bbD^{\heartsuit}}\cC \rangle}_{\nullclass,\frob}$, 
$\cG_k:={\langle \calD_{a+k-1}\rangle}_{\nullclass,\frob}$, 
$\bE_k=(\cE_k,lw|_{\cE_k})$, $\bG_k=(\cG_k,w|_{\cG_k})$ and 
we define $k_i\colon\bE_i\to\bE_{i+1}$, 
$u_i\colon\bE_i\to\bE_{i-1}$ and $t_i\colon\bE_i\to\bG_i$ 
to be restrictions of relative complicial exact functors 
$F_{[a,a+i-1],h}\cC\rinc F_{[a,a+i],h}\cC$, 
$\sigma_{\leq a+i-1}\colon F_{[a,a+i],h}\cC\to F_{[a,a+i-1]}\cC$ and 
$F_{[a,a+i-1],h}\cC\to \cC$ where 
the last functor is defined by  
sending an object $x$ in $F_{[a,a+i-1],h}\cC$ to 
$\Cone i^x_{a+i-1}$ in $\cC$. 
Then the sequence 
$\{0\}=\bE_0\to\bE_1\onto{k_1}\cdots \onto{k_{b-a}}\bE_{b-a+1}$ 
is a right derived flag of $\bE_{b-a+1}$ and the families 
$\{u_i\colon\bE_i\to\bE_{i-1} \}_{1\leq i\leq b-a+1}$ and 
$\{\bG_i,t_i\colon \bE_i\to\bG_i\}_{1\leq i\leq b-a+1}$ 
are a right derived retractions 
and an associated graded relative exact categories of this flag 
by \ref{ex:typical derived quasi-split sequence}. 
By \ref{rem:typical derived quasi-split sequence}, 
there exists a similar right derived flag on 
$({\langle F_{[a,b],h}^{\bbD^{\heartsuit}}\cC \rangle}_{\nullclass,\frob},lw)$. 

Next for each integer $1\leq k\leq b-a+1$, 
we set 
${cE'}_k:={\langle{\left (\underline{F_{[a,a+k-1],h}^{\bbD^{\heartsuit}}\cC} 
\right )}^{w_{\st}}\rangle}_{\nullclass,\frob}$, 
${\bE'}_k:=({\cE'}_k,lw|_{{\cE'}_k})$ and we define 
${k'}_i\colon{\bE'}_i\to{\bE'}_{i+1}$, 
${u'}_i\colon{\bE'}_i\to{\bE'}_{i+1}$ and 
${t'}_i\colon{\bE'}_i\to{\bG}_{i-1}$ to be restrictions of 
complicial exact functors 
$\underline{F_{[a,a+i-1],h}\cC}\rinc\underline{F_{[a,a+i],h}\cC}$, 
$\fs_{a+i-1}\colon\underline{F_{[a,a+i],h}\cC}\to
\underline{F_{[a,a+i-1],h}\cC}$ and 
${(-)}_{a+i-2}\colon \underline{F_{[a,a+i-1],h}\cC}\to \cC$ for 
$2\leq i\leq b-a+1$. 
where by convention we set $\bG_0:=\{0\}$. 
Then the sequence 
$\{0\}={\bE'}_0\onto{{k'}_0}{\bE'}_1\onto{{k'}_1}
\cdots \onto{{k'}_{b-a}}{\bE'}_{b-a+1}$ 
is a left derived flag of ${\bE'}_{b-a+1}$ 
and the families 
$\{{u'}_i\colon {\bE'}_i\to {\bE'}_{i-1} \}_{1\leq i\leq b-a+1}$ and 
$\{\bG_{i-1},{t'}_i\colon{\bE'}_i\to\bG_{i-1} \}_{1\leq i\leq b-a+1}$ 
are a left derived retractions and 
an associated graded relative exact categories 
of this flag respectively 
by \ref{ex:typical derived quasi-split sequence} again. 

Since by \ref{ex:pi0 of FbhC and underline FBhC}, for each integer $k$, 
the canonical morphism 
${(F_{[a,a+k],h}^{\bbD^{\heartsuit}}\cC)}^{w_{\st}}\to  
{\left(\underline{F_{[a,a+k],h}^{\bbD^{\heartsuit}}\cC}\right)}^{w_{\st}}
$ 
is a derived equivalence, 
there exists a similar left derived flag on 
$({\langle{(F_{[a,b],h}^{\bbD^{\heartsuit}}\cC)}^{w_{\st}}\rangle}_{\nullclass,\frob},lw)$ and $({\langle{(F_{[a,b]}^{\bbD^{\heartsuit}}\cC)}^{w_{\st}}\rangle}_{\nullclass,\frob},lw)$ by \ref{rem:typical derived quasi-split sequence}.
\end{para}

\section{Non-connective $K$-theory}
\label{sec:nonconnective}

A goal of this section is a d\'evissage theorem of 
non-connective $K$-theory \ref{thm:divissage for complicial exact categories}. 
In the first subsection \ref{subsec:def of nonconnective} and 
the second subsection \ref{subsec:fund prop of nonconnective}, 
we will recall 
the definition and fundamental properties of 
non-connective $K$-theory 
of relative exact categories. 
In subsection \ref{subsec:homotopy additivity theorem}, 
we will study a homotopy version of additivity 
theorem \ref{cor:homotopy additivity theorem} which enable us to 
compute $K$-theory of d\'evissage filtrations 
with level equivalences \ref{cor:homotopy invariance of C and D} 
by utilizing the canonical flag structures 
introduced in \ref{ex:canonical flag}. 
In the final subsection \ref{subsec:dev thm}, we will discuss 
a d\'evissage theorem of non-connective $K$-theory 
and as its special cases, 
we will reprove 
Waldhausen's cell 
filtration theorem~\ref{cor:Waldhausen's cell filtration theorem}, 
theorem of heart~\ref{cor:Theorem of heart}, 
Raptis' d\'evissage theorem~\ref{cor:Raptis' devissage theorem} and 
Quillen's d\'evissage theorem~\ref{cor:Quillen's devissage theorem}. 
In this section, we assume that 
underlying categories of relative categories are essentially small. 

\subsection{Definition of non-connective $K$-theory}
\label{subsec:def of nonconnective}

\begin{para}
\label{df:Schlichting nonconnective}
{\bf (Schlichting non-connective $K$-theory).}\ \ 
In the series of papers \cite{Sch04}, \cite{Sch06} and \cite{Sch11}, 
Marco Schlichting define and establish non-connective $\bbK$-theory 
for exact categories, Frobenius pairs and complicial exact 
categories with weak equivalences respectively. 
We write $\bbK^{S,\operatorname{exact}}(\cE)$ for 
the non-connective $K$-theory of an exact category $\cE$ 
and denote the non-connective $K$-theory of 
a relative complicial exact category $\bC=(\cC,w)$ 
by $\bbK^{S,\operatorname{comp}}(\bC)$ or $\bbK^{S,\operatorname{comp}}(\cC;w)$. 
\end{para}

\begin{para}
\label{df:nonconnective}
{\bf Definition (Non-connective $K$-theory).}\ \ 
For a relative exact category $\bE=(\cE,w)$, 
we set 
\begin{equation}
\label{eq:bbK df}
\bbK(\bE)(=\bbK(\cE;w)):=\bbK^{S,\operatorname{comp}}(\Ch_b(\bE)). 
\end{equation}
Moreover let $\bC=(\cC,w)$ be a relative complicial category and let 
$\cF$ be a full subcategory of $\cC$ and let 
$\bbD=\{\calD_n\}_{n\in\bbZ}$ be a $T$-system (for definition see \ref{df:properties of family}) 
of $\cC$. Then we set 
\begin{equation}
\label{df:rel bbK df}
\bbK(\inn{\cF}{\bC}):=\bbK^{S,\operatorname{comp}}
({\langle\cF \rangle}_{\nullclass,\frob};w|_{{\langle\cF \rangle}_{\nullclass,\frob}}),\ \ \text{and}
\end{equation}
\begin{equation}
\bbK(\inn{\bbD}{\bC}):=\underset{T}{\hocolim}\bbK(\inn{\calD_n}{\bC}).
\end{equation}
We call $\bbK(\inn{\cF}{\bC})$ and $\bbK(\inn{\bbD}{\bC})$ the {\it relative $\bbK$-theory of $\cF$ in $\bC$} and 
the {\it relative $\bbK$-theory of $\bbD$ in $\bC$}. 

Recall the definition of morphisms of $T$-systems from \ref{df:properties of family}. 
Let $\bbD$ and $\bbD'$ be a pair of $T$-systems of relative complicial 
exact categories $\bC=(\cC,w)$ and $\bC'=(\cC',w')$ respectively. 
Assume that $\bbD'$ is closed under isomorphisms and let 
$(f,d)\colon\bbD\to\bbD'$ be a morphism of $T$-system. 
Then $(f,d)$ induces a morphism of spectra 
$\bbK(\inn{\bbD}{\bC})\to\bbK(\inn{\bbD'}{\bC'})$.
\end{para}

\begin{para}
\label{df:Euler character}
{\bf Definition (Euler character).}\ \ 
Let $\bC=(\cC,w)$ be a relative complicial exact category and 
let $\bbD=\{\calD_n\}_{n\in\bbZ}$ be a $T$-system of $\bC$. 
If we regard $\cC$ as a $T$-system $\{\cC_n\}_{n\in\bbZ}$ by 
setting $\cC_n:=\cC$ for all integer $n$, 
then the family of inclusion functors 
$\calD_n\rinc \cC$ for all integers can be regarded as 
a morphism of $T$-systems and it induces a morphism of 
spectra 
\begin{equation}
\label{eq:Euler character df}
E_{\bbD}\colon \bbK(\inn{\bbD}{\bC})\to\underset{T}{\hocolim} \bbK(\bC)
\end{equation}
which we call the {\it Euler characteristic morphism} ({\it associated with 
the $T$-system $\bbD$}). 
In \ref{cor:homotopy additivity for functors}, 
we will show the equality $\mathrm{(\ref{eq:KT=-id})}$ and 
by this equality it turns out 
that $\displaystyle{\underset{T}{\hocolim} \bbK(\bC)}$ 
is homotopy equivalent to $\bbK(\bC)$. 
We denote composition of $E_{\bbD}$ with 
a non-canonical homotopy equivalence 
$\displaystyle{\underset{T}{\hocolim} \bbK(\bC)\to\bbK(\bC)}$ by the same letters $E_{\bbD}$ and call it 
{\it Euler characteristic morphism} too. 
For a full subcategory $\cF$ of $\cC$, we similarly 
construct a morphism of spectra
\begin{equation}
\label{eq:relative Euler character df}
E_{\bbD,\cF}\colon \bbK(\inn{\bbD|_{\cF}}{\bC|_{\cF}})\to\underset{T}{\hocolim}\ \bbK(\bC|_{\cF}).
\end{equation} 
\end{para}

\begin{para}
\label{prop:comp of non-connective K}
{\bf Proposition (Comparison theorem of non-connective $K$-theory).}\ 
(\cf \cite[0.3]{Moc13b})\ \ 
{\it 
Let $\bE=(\cE,w)$ be a relative exact category. Then

\sn
$\mathrm{(1)}$ 
If $w$ is the class of all isomorphisms in $\cE$ , 
then there is a natural homotopy equivalence of spectra
$$\bbK(\bE)\isoto\bbK^{S,\operatorname{exact}}(\cE).$$
This equivalence is functorial in the sense that it 
gives a natural equivalence between the functors from 
the category of (essentially small) exact categories and exact functors 
to the stable category of spectra. 

\sn
$\mathrm{(2)}$ 
If $\bE$ is a relative complicial exact category, 
then there is a natural homotopy equivalence of spectra
$$\bbK(\bE)\isoto \bbK^{S,\operatorname{comp}}(\bE).$$
This equivalence is functorial in the sense that it 
gives a natural equivalence between the functors from 
the category of (essentially small) relative complicial exact categories 
and relative complicial exact functors 
to the stable category of spectra.
}\qed
\end{para}

\subsection{Fundamental properties of non-connective $K$-theory}
\label{subsec:fund prop of nonconnective}

Recall the definition of categorical homotopic of relative 
(exact) functors from \ref{df:rel cat}. 

\begin{para}
\label{thm:categorical homotopy equivalence}
{\bf Theorem (Categorical homotopy invariance).}\ 
(\cf \cite[3.25]{Moc13b}.)\ \ 
{\it
Let $f$, $f'\colon \bE \to \bE'$ be a pair of relative 
exact functors between relative exact categories $\bE$ and $\bE'$. 
If $f$ and $f'$ are categorical homotopic in $\underline{\RelEx}$, then 
they induce same maps $\bbK(f)=\bbK(f')\colon\bbK(\bE)\to\bbK(\bE')$ on 
$\bbK$-theory. 
}\qed
\end{para}

\begin{para}
\label{thm:localization theorem}
{\bf Theorem (Localization).}\ (\cf \cite[0.2]{Moc13b}.)\ \ 
{\it
Let $\mathrm{(\ref{eq:derived spltting sequence})}$ 
be a derived weakly exact sequence 
{\rm(}\ref{df:derived exaxt seq}{\rm)} of 
relative exact categories and relative functors. Then 
the induced sequence of spectra 
$$\bbK(\bE)\onto{\bbK(i)}\bbK(\bE')\onto{\bbK(p)}\bbK(\bE'')$$
is a fibration sequence of spectra.
}\qed
\end{para}

\begin{para}
\label{rem:categorical homotopy equivalence}
{\bf Remark.}\ \ 
In the statements of \cite[0.2, 3.25]{Moc13b}, 
there is an assumption that 
we assume that 
$w$ is consistent. 
In this paper we change the definition of 
the class of quasi-weak equivalences 
as in Example~\ref{ex:quasi weak equivalences} 
and our proof in Ibid with the new definition 
of the class of quasi-weak equivalences still works fine in the setting of 
Theorem~\ref{thm:categorical homotopy equivalence}
and Theorem~\ref{thm:localization theorem}.
\end{para}

\begin{para}
\label{thm:approximation}
{\bf Corollary (Cofinality theorem).}\ \ 
{\it
Let $f\colon\bE\to\bF$ be a relative exact functor between 
relative exact categories. 
If $f$ induces a cofinal functor on their 
derived categories $\calD_b(\bE)\to\calD_b(\bF)$, 
then $f$ induces a homotopy equivalence 
of spectra $\bbK(\bE)\to\bbK(\bF)$ on $\bbK$-theory.
}\qed
\end{para}

\begin{para}
\label{cor:relative Gillet-Waldhausen theorem}
{\bf Corollary (Relative Gillet-Waldhausen theorem).}\ \ 
{\it
Let $\bC=(\cC,w)$ be a relative complicial exact category and let 
$\cA$ be a strict exact subcategory of $\cC$. Then
the inclusion functor $j_{\cF}\colon\cF\to\Ch_b\cF$ 
induces a homotopy equivalence of spectra 
$\bbK(\cF;w|_{\cF})\to\bbK(\inn{\Ch_b\cF}{\Ch_b\bC})$ 
on $\bbK$-theory.
}
\end{para}

\begin{proof}
It follows from Proposition~\ref{prop:comp hom and der} $\mathrm{(4)}$ and Corollary~\ref{thm:approximation}. 

\end{proof}

\begin{para}
\label{cor:w-envelope invariance}
{\bf Corollary ($w$-envelope invariance).}\ \ 
Let $\bC=(\cC,w)$ be a relative complicial exact category and 
let $\cF$ be a full subcategory of $\cC$. 
Assume that $\cF$ is compatible with $w$ 
{\rm(see \ref{df:relative homotopy category})}. 
Then the inclusion functor 
${\langle\cF \rangle}_{\nullclass,\frob}
\rinc{\langle\cF \rangle}_{\nullclass,w}$ 
induces a homotopy equivalence of spectra 
$\bbK(\inn{\cF}{\bC})\to\bbK(\inn{{\langle\cF \rangle}_{\nullclass,w}}{\bC})$ on $\bbK$-theory. 
\end{para}

\begin{proof}
It follows from the equivalence 
$\mathrm{(\ref{eq:w-envelope invariance})}$ and 
Corollary~\ref{thm:approximation}.
\end{proof}

\begin{para}
\label{cor:w-envelope invariance II}
{\bf Corollary ($w$-envelope invariance II).}\ \ 
{\it
Let $\bC=(\cC,w)$ be a relative complicial exact category and let 
$\cA$ be an additive full subcategory of $\cC$. 
Assume either 
condition $\mathrm{(i)}$ or condition $\mathrm{(ii)}$ below. 

\sn
$\mathrm{(i)}$ 
$\cA$ contains all $C$-contractible objects of $\cC$,

\sn
$\mathrm{(ii)}$ 
$\cA$ is a strict exact subcategory of $\cC$. 

\sn
Then the inclusion functor 
${\langle\Ch_b\cA\rangle}_{\nullclass,\frob,\Ch_b\cC}\rinc{\langle\Ch_b\cA\rangle}_{\nullclass,qw,\Ch_b\cC}$ 
induces a homotopy equivalence of spectra 
$\bbK(\inn{\Ch_b\cA}{\Ch_b\bC})\to\bbK(\inn{{\langle\Ch_b\cA\rangle}_{\nullclass,qw,\Ch_b\cC}}{\Ch_b\bC}) $ on $\bbK$-theory. 
}
\end{para}

\begin{proof}
It follows from 
Corollary~\ref{cor:w-closure envelope} and 
Corollary~\ref{cor:w-envelope invariance}.
\end{proof}

\begin{para}
\label{cor:purity theorem}
{\bf Corollary (Purity theorem).}\ \ 
{\it
Let $\bC=(\cC,w)$ be a thick relative complicial exact category 
and let $\cP$ 
be a full additive subcategory of $\cC$. 
Assume that either condition 
$\mathrm{(i)}$ or condition $\mathrm{(ii)}$ below.

\sn
$\mathrm{(i)}$ 
$\cP$ is a prenull class or 

\sn
$\mathrm{(ii)}$ $\cC$ is strictly ordinary and $\cP$ contains all $C$-contractible objects and $\cP$ is closed under operations $T^{\pm 1}$. 

\sn
Then the relative complicial exact functor $\Tot\colon \Ch_b(\bC|_{\cP})\to ({\langle\cP\rangle}_{\nullclass,\frob},w|_{{\langle\cP\rangle}_{\nullclass,\frob}})$ 
induces a homotopy equivalence of spectra 
$\bbK^{S,\operatorname{comp}}(\Ch_b(\bC|_{\cP}))\isoto \bbK(\inn{\cP}{\bC})$ 
on $\bbK$-theory.
}
\end{para}

\begin{proof}
It follows from \ref{prop:comp hom and der} and \ref{thm:approximation}. 
\end{proof}

\subsection{Homotopy additivity theorem}
\label{subsec:homotopy additivity theorem}

Recall the definition of derived quasi-split sequences of 
relative exact categories from \ref{df:quasi split sequence}. 

\begin{para}
\label{thm:additivity theorem}
{\bf Theorem (Additivity theorem).}\ (\cf \cite[7.10]{Moc13b}.)\  
{\it
In the sequence $\mathrm{(\ref{eq:derived spltting sequence})}$ 
of relative exact categories, 
if there exists a relative exact functor 
$q\colon\bE'\to \bE$ such that 
$\calD_b(q)$ gives a part of right splitting of 
a sequence $\mathrm{(\ref{eq:associated quasi-split seq})}$ 
of triangulated categories, 
then the projection functor 
$\bE'\onto{\footnotesize{\begin{pmatrix}q\\ p\end{pmatrix}}}\bE\times \bE'$ 
induces a homotopy equivalence of spectra
$\bbK(\bE')\to\bbK(\bE)\vee \bbK(\bE'')$ on $\bbK$-theory.
}\qed
\end{para}

By Proposition~\ref{prop:typical derived quasi-split sequence} and Theorem~\ref{thm:additivity theorem}, 
we obtain the following result.

\begin{cor}
\label{cor:homotopy additivity theorem}
{\bf Corollary (Homotopy additivity theorem).}\ \ 
{\it
In the convention~\ref{nt:setting of additivity}, 
the complicial exact functor 
$\footnotesize{\begin{pmatrix}q\\ p\end{pmatrix}}\colon 
{\langle\cC' \rangle}_{\nullclass,\frob}\to {\langle\cS \rangle}_{\nullclass,\frob}\times {\langle{\cT}_{\heq} \rangle}_{\nullclass,\frob} $ induces a homotopy equivalence of spectra 
$$\bbK(\inn{\cC'}{\bC^{[1]}_h})\to \bbK(\inn{\cS}{\bC})\times 
\bbK(\inn{{\cT}_{\heq}}{\bC})$$
on $\bbK$-theory. 
If we replace $\bC^{[1]}_h$ with $\bC^{[1]}$, 
then we also have similar result.
}
\qed
\end{cor}

\begin{para}
\label{df:cone of complicial natural transformation}
{\bf Definition (Cone of complicial natural transformation).}\ \ 
Let $(f,c)$, $(g,d)\colon\cC\to\cC'$ be a pair of 
complicial exact functors between complicial 
exact categories $\cC$ and $\cC'$ and let 
$\theta\colon (f,c)\to (g,d)$ be a complicial natural transformation. 
Then we define $\Cone\theta\colon\cC\to\cC'$ to be 
an exact functor by sending an object $x$ in $\cC$ to 
$\Cone \theta_x$ in $\cC'$ and a morphism $a\colon x\to y$ in $\cC$ to 
$\Cone(f(a),g(a))\colon\Cone \theta_x\to\Cone \theta_y$ in $\cC'$. 
There exists a natural equivalence $c^{\Cone\theta}\colon C'\Cone\theta\isoto \Cone\theta C$ which is characterized by the following equalities
\begin{equation}
\label{eq:conetheta char 1}
c^{\Cone\theta}\cdot(c'\ast\kappa\ast{\theta})=(\kappa\ast{\theta}\ast C)\cdot d,\ \text{and}
\end{equation}
\begin{equation}
\label{eq:conetheta char 2}
c^{\Cone\theta}\cdot(c'\ast\mu\ast\theta)=(\mu\ast C)\cdot(C'\ast c)\cdot(\sigma\ast f).
\end{equation}
Then the pair $(\Cone\theta,c^{\Cone\theta})$ 
is a complicial exact functor. 
There is a complicial natural transformation $\kappa\ast\theta\colon 
(g,d)\to (\Cone\theta,c^{\Cone\theta})$. 

If $w$ and $w'$ are classes of complicial weak equivalences of $\cC$ and 
$\cC'$ respectively and $f$ and $g$ are relative functors, 
then $\Cone\theta$ is also relative exact functor.
\end{para}

\begin{para}
\label{cor:homotopy additivity for functors}
{\bf Corollary.}\ \ 
{\it
Let $(f,c)$, $(g,d)\colon\bC\to\bC'$ be a relative complicial exact functors 
between relative complicial exact categories $\bC$ and $\bC'$ 
and let $\theta\colon (f,c)\to (g,d)$ 
be a complicial natural transformation. 
Then for the morphisms of spectra $\bbK(f)$, $\bbK(g)$, 
$\bbK(\Cone\theta)\colon \bbK(\bC)\to\bbK(\bC')$, 
we have the equality
\begin{equation}
\label{eq:homotopy additivity for functors}
\bbK(g)=\bbK(f)\vee\bbK(\Cone\theta).
\end{equation}
If $\theta$ is an admissible monomorphisms in the category of 
complicial functors, then 
we have the equality
\begin{equation}
\label{eq:cone=coker}
\bbK(\Cone \theta)=\bbK(\coker \theta).
\end{equation}
In particular for the suspension functor $T\colon\bC\to \bC$, 
we have the equality
\begin{equation}
\label{eq:KT=-id}
\bbK(T)=-\id_{\bbK(\bC)}.
\end{equation}
}
\end{para}

\begin{proof}
By virtue of \ref{prop:comp hom and der}, 
replacing $\bC$ and $\bC'$ with $\Ch_b\bC$ and $\Ch_b\bC'$, 
without loss of generality, we shall assume that 
$\bC=(\cC,w)$ and $\bC'=(\cC',w')$ are thick normal ordinary. 
We define $U$, $V\colon\bC\to({\cC'}^{[1]}_h,lw)$ and 
$W\colon({\cC'}^{[1]}_h,lw')\to (\cC',w')$ 
and $X\colon ({\cC'}^{[1]}_h,lw)\to (\cC',w')\times (\cC',w') $ 
to be relative complicial exact functors 
by sending an object $x$ in $\cC$ to $[f(x)\onto{\theta_x}g(x)] $ and 
$[f(x)\onto{\tiny{\begin{pmatrix}\id_{f(x)}\\ 0\end{pmatrix}}}f(x)\oplus\Cone\theta (x)]$ in ${\cC'}^{[1]}_h$ and 
an object $[x\onto{u} y]$ in $\cC'$ to $y$ in $\cC'$ 
and $(x,\Cone u)$ in $\cC'\times\cC'$ 
respectively. 
We have the equalities $\bbK(g)=\bbK(W)\bbK(U)$ and 
$\bbK(f)\vee\bbK(\Cone\theta)=\bbK(W)\bbK(V)$ and 
$\displaystyle{\bbK(X)\bbK(U)=\begin{pmatrix}\bbK(f)\\ \bbK(\Cone\theta)
\end{pmatrix}\underset{\textbf{I}}{=}
\begin{pmatrix}\bbK(f)\\ \bbK(C\oplus\Cone\theta)\end{pmatrix}
=\bbK(X)\bbK(V)}$
where the equality $\textbf{I}$ follows from 
categorical homotopy equivalence~\ref{thm:categorical homotopy equivalence}. 
By homotopy additivity theorem~\ref{cor:homotopy additivity theorem}, 
$\bbK(X)$ is a homotopy equivalence of spectra and thus 
we obtain the equality $\mathrm{(\ref{eq:homotopy additivity for functors})}$.
If $\theta$ is an admissible monomorphism in the category of 
relative complicial exact functors, then 
$\coker \theta$ and $\Cone \theta$ is a categorical homotopic. 
Thus the equality $\mathrm{(\ref{eq:cone=coker})}$ follows 
from \ref{thm:categorical homotopy equivalence}.

Finally the equality $\mathrm{(\ref{eq:KT=-id})}$ follows form the equality 
$T=\Cone(\id_{\cC}\to 0)$ and the equality 
$\mathrm{(\ref{eq:homotopy additivity for functors})}$.
\end{proof}

Recall the functors $\fj$, $\sigma_{\geq k}$, $(-)[k]$, $\fc_k$ from 
\ref{df:filtered objects}, 
\ref{df:truncations on FbhC}, 
\ref{df:degree shift} and 
\ref{df:Cone functor}. 

\begin{para}
\label{cor:bbKcn+1}
{\bf Corollary.}\ \ 
{\it
Let $\bC=(\cC,w)$ be a normal ordinary complicial relative 
exact category and let $n$ be an integer. 
Then we have the equality
\begin{equation}
\label{eq:bbKcn+1}
\bbK(\fc_{n+1})=\id_{\bbK(F_{b,h}\cC)}\vee\bbK(\sigma_{\geq n+1}T). 
\end{equation}
In particular we also have the equality.
\begin{equation}
\label{eq:bbKcn+1j}
\bbK(\fc_{n+1}(\fj(-)[n]))=\bbK(\fj(-)[n])\vee\bbK(\fj(T(-))[n]).
\end{equation}
}
\end{para}

\begin{proof}
Notice that there is a level Frobenius admissible 
exact sequence of relative complicial exact functors on $F_{b,h}\cC$.
$$\id_{F_{b,h}\cC}\rinf \fc_{n+1}\rdef \sigma_{\geq n+1}T.$$
$$
\xymatrix{
x_n \ar@{>->}[r] \ar[d]_{i_n^x} & x_n \ar@{->>}[r] \ar[d]_{Ci^x_{n}\iota_{x_n}} & 0 \ar[d]\\
x_{n+1} \ar@{>->}[r]^{\iota_{x_{n+1}}} \ar[d]_{i^{x}_{n+1}} & 
Cx_{n+1} \ar@{->>}[r]^{\pi_{x_{n+1}}} \ar[d]_{Ci^x_{n+1}} & 
Tx_{n+1} \ar[d]^{Ti^x_{n+1}}\\
\vdots & \vdots & \vdots.
}
$$
Thus we obtain the equalities 
$\mathrm{(\ref{eq:bbKcn+1})}$ and $\mathrm{(\ref{eq:bbKcn+1j})}$ 
by \ref{cor:homotopy additivity for functors}. 
\end{proof}

From \ref{thm:additivity theorem}, 
proceeding by induction, we obtain the result for $\bbK$-theory of 
derived flags below. 
We will only state the case for right derived flags and 
the case for left derived flags also holds.  

\begin{para}
\label{cor:K-theory of derived flag}
{\bf Corollary ($\bbK$-theory of derived flags).}\ \ 
{\it
Let $\bE$ be a relative exact category and let 
$\mathrm{(\ref{eq:right derived flag df})}$ 
be a right derived flag of $\bE$ with 
a derived right retractions 
$\calR=\{u_i\colon\bE_i\to\bE_{i-1} \}_{1\leq i\leq n}$ 
and an associated graded relative exact categories 
$\cG=\{\bG_i,\ t_i\colon\bE_i\to\bG_i \}_{1\leq i\leq n}$. 
Then the relative complicial exact functor 
$\bE\onto{\tiny{\begin{pmatrix}t_n\\ t_{n-1}u_n\\ t_{n-2}u_{n-1}u_n \\ \vdots\\ t_1u_2\cdots u_{n-1}u_n \end{pmatrix}}}
\bG_n\times\bG_{n-1}\times\cdots\times\bG_1$ 
induces a homotopy equivalence of spectra 
$\bbK(\bE)\isoto \bigvee_{i=1}^n\bbK(\bG_i)$.
}\qed
\end{para}

Recall the functors $\fj$, $(-)[k]$, $\fc_k$ from 
\ref{df:filtered objects}, \ref{df:degree shift} and 
\ref{df:Cone functor}. 

\begin{para}
\label{df:C and D df}
{\bf Definition.}\ 
Let $(\cC,w)$ be a thick normal ordinary 
relative complicial exact category and let 
$a<b$ be a pair of integers. 
For a positive integer $m$, we 
write $\cC^m$ for 
$\underbrace{\cC\times\cdots\times\cC}_{m}$ 
the $m$-times products of $\cC$. 
We define $A_{[a,b]}\colon \cC^{b-a+1}\to F_{[a,b],h}\cC$ 
and $B_{[a,b]}\colon \cC^{b-a}\to {(F_{[a,b],h}\cC)}^{w_{\st}}\cC$ 
to be complicial exact functors 
by setting $A_{[a,b]}:=\footnotesize{
\fj(-)[a]\oplus \fj(-)[a+1]\oplus \cdots \oplus \fj(-)[b]
}$ and 
$B_{[a,b]}:=\footnotesize{
\fc_{a+1}\fj(-)[a] \oplus \fc_{a+2}\fj(-)[a+1] \oplus 
\cdots \oplus \fc_b\fj(-)[b-1]
}$ 
and we define 
$C_{[a,b]}\colon F_{[a,b],h}\cC\to \cC^{b-a+1}$ and 
$D_{[a,b]}\colon{(F_{[a,b],h}\cC)}^{w_{\st}}\cC\to \cC^{b-a}$ to be 
complicial 
exact functors by sending an object $x$ in $F_{[a,b],h}\cC$ to 
$(x_a,\Cone i^x_a,\Cone i^x_{a+1},\cdots, \Cone i^x_{b-1})$ 
in $\cC^{b-a+1}$ 
and an object $x$ in ${(F_{[a,b],h}\cC)}^{w_{\st}}$ to 
$(x_a,x_{a+1},\cdots,x_{b-1})$ in $\cC^{b-a}$ 
and we also define $E_{[a,b]}=(e_{ij})\colon\cC^{b-a}\to\cC^{b-a+1}$ 
to be a complicial exact functor 
by setting $e_{ij}=\begin{cases}\id_{\cC} & \text{if $i=j$}\\ 
T & \text{if $i=j+1$}\\
0 & \text{otherwise}
\end{cases}$.
We also write the same letters $A_{[a,b]}$, $B_{[a,b]}$, $C_{[a,b]}$, $D_{[a,b]}$ and $E_{[a,b]}$ 
for the restrictions of $A_{[a,b]}$, $B_{[a,b]}$, 
$C_{[a,b]}$, $D_{[a,b]}$ and $E_{[a,b]}$ 
to subcategories of 
$\cC^{b-a+1}$, $\cC^{b-a}$, 
$F_{[a,b],h}\cC$, ${(F_{[a,b],h}\cC)}^{w_{\st}}$ and $\cC^{b-a}$. 
Since $\bbD$ is a cell structure, these functors are well-defined.
\end{para}

\begin{para}
\label{cor:homotopy invariance of C and D}
{\bf Corollary.}\ \ 
{\it
Let $\bC=(\cC,w)$ be a thick normal ordinary relative 
complicial exact category and let $a<b$ be a pair of 
integers and let $\bbD$ be a $C$-homotopy invariance closed 
cell structure of $\bC$. 
Then $C_{[a,b]}$ and $D_{[a,b]}$ induce the homotopy equivalence 
of spectra
\begin{equation}
\label{eq:C homotopy equivalence}
\footnotesize{
\bbK(\inn{F^{\bbD^{\heartsuit}}_b\cC}{(F_b\cC,lw)})\to 
\displaystyle{\prod_{k\in [a,b]}\bbK(\inn{\calD_k}{\bC})},\ \ 
\bbK(\inn{F^{\bbD^{\heartsuit}}_{b,h}\cC}{(F_{b,h}\cC,lw)})\to 
\displaystyle{\prod_{k\in [a,b]}\bbK(\inn{\calD_k}{\bC})}}
\end{equation}
\begin{equation}
\label{eq:D homotopy equivalence}
\footnotesize{
\bbK(\inn{{(F^{\bbD^{\heartsuit}}_b\cC)}^{w_{\st}}}{(F_b\cC,lw)})\to 
\displaystyle{\prod_{k\in [a,b-1]}\bbK(\inn{\calD_k}{\bC})},\ \text{and}\ 
\bbK(\inn{{(F^{\bbD^{\heartsuit}}_{b,h}\cC)}^{w_{\st}}}{(F_{b,h}\cC,lw)})\to
\displaystyle{\prod_{k\in [a,b-1]}\bbK(\inn{\calD_k}{\bC})}}
\end{equation}
and $A_{[a,b]}$ and $B_{[a,b]}$ induce the inverse morphisms 
of $\bbK(C_{[a,b]})$ and $\bbK(D_{[a,b]})$ respectively. 
Moreover the following diagram is commutative
\begin{equation}
\label{eq:ABE commutative diagram}
{\footnotesize{
\xymatrix{
{\displaystyle{\prod_{k\in [a,b-1]}\bbK(\inn{\calD_k}{\bC})}} 
\ar[r]^{\!\!\!\!\!\!\!\!\!\!\!\!\!\bbK(B_{[a,b]})} \ar[d]_{\bbK(E_{[a,b]})} & 
\bbK(\inn{{(F^{\bbD^{\heartsuit}}_{b,h}\cC)}^{w_{\st}}}{(F_{b,h}\cC,lw)})
\ar[d]\\
{\displaystyle{\prod_{k\in [a,b]}\bbK(\inn{\calD_k}{\bC})}} 
\ar[r]_{\!\!\!\!\!\!\!\!\!\!\!\!\!\bbK(A_{[a,b]})} & 
\bbK(\inn{F^{\bbD^{\heartsuit}}_{b,h}\cC}{(F_{b,h}\cC,lw)}).
}}}
\end{equation}
}
\end{para}

\begin{proof}
The first assertion follows from \ref{ex:pi0 of FbhC and underline FBhC}, \ref{ex:canonical flag}, 
\ref{thm:approximation} and 
\ref{cor:K-theory of derived flag}. 
We have an equality $C_{[a,b]}A_{[a,b]}=\id_{\cC^{b-a+1}}$ and 
there is a natural weak equivalence $D_{[a,b]}B_{[a,b]}\to\id_{\cC^{b-a}}$. 
Thus $\bbK(A_{[a,b]})$ and $\bbK(B_{[a,b]})$ are the inverse morphisms 
of $\bbK(C_{[a,b]})$ and $\bbK(D_{[a,b]})$ respectively. 
Commutativity of the diagram $\mathrm{(\ref{eq:ABE commutative diagram})}$ 
follows form the equality 
$\mathrm{(\ref{eq:bbKcn+1j})}$. 
\end{proof}

\begin{para}
\label{cor:K-theory of Fb and Fbh}
{\bf Corollary.}\ \ 
{\it
Let $(\cC,w)$ be a thick normal ordinary relative 
complicial exact category and let 
$\bbD$ be a $C$-homotopy equivalent cell structure of $(\cC,w)$. 
Then the inclusion functors induces homotopy equivalences 
of spectra
\begin{equation}
\label{eq:comparison Fb and Fbh 1}
\bbK(\inn{F^{\bbD^{\heartsuit}}_b\cC}{(F_b\cC,lw)})\to 
\bbK(\inn{F^{\bbD^{\heartsuit}}_{b,h}\cC}{(F_{b,h}\cC,lw)}), 
\end{equation}
\begin{equation}
\label{eq:comparison Fb and Fbh 2}
\bbK(\inn{({F^{\bbD^{\heartsuit}}_b\cC})^{w_{\st}}}{(F_b\cC,lw)})\to 
\bbK(\inn{F^{\bbD^{\heartsuit}}_{b,h}\cC}{(F_{b,h}\cC,lw)})\ 
\text{and}
\end{equation}
\begin{equation}
\label{eq:comparison Fb and Fbh 3}
\bbK(\inn{F^{\bbD^{\heartsuit}}_b\cC}{F_b\bC})\to 
\bbK(\inn{F^{\bbD^{\heartsuit}}_{b,h}\cC}{F_{b,h}\bC}).
\end{equation}
}
\end{para}

\begin{proof}
Let $a<b$ be a pair of integers. 
Then there are the following commutative diagrams of $\bbK$-theory. 
$$
\footnotesize{
\xymatrix{
\bbK(\inn{F^{\bbD^{\heartsuit}}_b\cC}{(F_b\cC,lw)})\ar[rr] \ar[rd]_{\bbK(C_{[a,b]})} & &
\bbK(\inn{F^{\bbD^{\heartsuit}}_{b,h}\cC}{(F_{b,h}\cC,lw)})\ar[ld]^{\bbK(C_{[a,b]})} \\
& \displaystyle{\prod_{k\in [a,b]}\bbK(\inn{\calD_k}{\bC})}
}},
$$
$$
\footnotesize{
\xymatrix{
\bbK(\inn{{(F^{\bbD^{\heartsuit}}_b\cC)}^{w_{\st}}}{(F_b\cC,lw)})\ar[rr] \ar[rd]_{\bbK(D_{[a,b]})} & &
\bbK(\inn{{(F^{\bbD^{\heartsuit}}_{b,h}\cC)}^{w_{\st}}}{(F_{b,h}\cC,lw)})\ar[ld]^{\bbK(D_{[a,b]})} \\
& \displaystyle{\prod_{k\in [a,b-1]}\bbK(\inn{\calD_k}{\bC})}.
}}
$$
By \ref{cor:homotopy invariance of C and D}, horizontal morphisms 
in the diagrams 
are homotopy equivalence of spectra. 
By taking colimits, we obtain the homotopy equivalences 
$\mathrm{(\ref{eq:comparison Fb and Fbh 1})}$ and 
$\mathrm{(\ref{eq:comparison Fb and Fbh 2})}$. 

The last homotopy equivalence of spectra 
$\mathrm{(\ref{eq:comparison Fb and Fbh 3})}$ 
follows from the diagram of fibration sequences below.
$$\footnotesize{\xymatrix{
\bbK(\inn{{(F^{\bbD^{\heartsuit}}_b\cC)}^{w_{\st}}}{(F_b\cC,lw)})\ar[r] \ar[d] &
\bbK(\inn{F^{\bbD^{\heartsuit}}_b\cC}{(F_b\cC,lw)})\ar[r] \ar[d] &
\bbK(\inn{F^{\bbD^{\heartsuit}}_b\cC}{F_b\bC}) \ar[d]\\
\bbK(\inn{{(F^{\bbD^{\heartsuit}}_{b,h}\cC)}^{w_{\st}}}{(F_{b,h}\cC,lw)})\ar[r] &
\bbK(\inn{F^{\bbD^{\heartsuit}}_{b,h}\cC}{(F_{b,h}\cC,lw)})\ar[r] &
\bbK(\inn{F^{\bbD^{\heartsuit}}_{b,h}\cC}{F_{b,h}\bC}).
}}
$$
\end{proof}

\begin{para}
\label{cor:degree shift is id}
{\bf Corollary.}\ \ 
{\it
Let $\bC=(\cC,w)$ be a normal ordinary relative 
complicial exact category and 
let $\bbD$ be a $C$-homotopy closed cell structure of $\bC$. 
Then the relative complicial exact functor 
$T(-)[1]\colon (F_{b,h}^{\bbD^{\heartsuit}}\cC,w_{\st})\to (F_{b,h}^{\bbD^{\heartsuit}}\cC,w_{\st})$ 
which sends an object $x$ to an object $T(x)[1]$ 
induces $-\id_{\bbK(F_{b,h}^{\bbD^{\heartsuit}}\bC)}$ 
on $\bbK$-theory. 
}
\end{para}

\begin{proof}
We will show $T(-)[1]\colon
(F_{[a,b]}^{\bbD^{\heartsuit}}\cC,w_{\st})\to
(F_{[a,b]}^{\bbD^{\heartsuit}}\cC,w_{\st})$ 
induces $-\id$ on $\bbK$-theory 
for any pair of integers $a<b$. 
Then by taking direct limit as $a$ goes to $-\infty$ and 
$b$ goes to $+\infty$ and by \ref{cor:K-theory of Fb and Fbh}, 
we will obtain the result. 
For $m>b$, since 
there are stable weak equivalences 
$\sigma_{\geq n}(-)[1]\to \id[1] \leftarrow \id$ and $\sigma_{\geq n}C\to 0$ 
and there exists an admissible exact sequences 
relative complicial exact functors 
$\sigma_{\geq n}(-)[1]\rinf 
\sigma_{\geq n}C(-)[1]\rdef \sigma_{\geq n}T(-)[1]$, 
we have equalities by \ref{thm:categorical homotopy equivalence} 
and \ref{cor:homotopy additivity for functors} 
$$\bbK(T[1])=\bbK(\sigma_{\geq n}T[1])=-\bbK(\sigma_{\geq n}(-)[1])=-\id.$$

\end{proof}

\begin{para}
\label{df:Euler character in devissage space}
{\bf Definition (Euler characteristic in 
$K$-theory of d\'evissage spaces).}\ \ 
Let $\bC=(\cC,w)$ be a normal ordinary 
relative complicial exact category 
and let $\bbD=\{\calD_{\geq n},\calD_{\geq m} \}_{n,m\in\bbZ}$ 
be a cell structure of $\bC$. 
Then for each integer $n$, the commutative diagram below 
$$
\footnotesize{
\xymatrix{
\calD_n \ar[r]^{\fj(-)[n]} \ar[d]_T & F_{b,h}^{\bbD^{\heartsuit}}\cC 
\ar[d]^{T(-)[1]}\\
\calD_{n+1} \ar[r]_{\fj(-)[n+1]} & F_{b,h}^{\bbD^{\heartsuit}}\cC
}}
$$
induces a morphism 
\begin{equation}
\label{eq:Euler characteristic in devissage spaces}
E_{\bbD}\colon\bbK(\inn{\bbD}{\bC})\to \underset{T(-)[1]}{\hocolim} \bbK(F_{b,h}^{\bbD^{\heartsuit}}\bC)
\end{equation}
which we call the {\it Euler characteristic morphism} 
({\it in $\bbK$-theory of 
d\'evissage space associated with $\bbD^{\heartsuit}$}). 
As a variant, for a full subcategory $\cF$ in $\cC$, 
we also define the {\it relative Euler characteristic morphism}
\begin{equation}
\label{eq:Euler characteristic in devissage spaces relative}
E_{\bbD,\cF}\colon
\bbK(\inn{{(\bbD|_{\cF})}}{\bC|_{\cF}})
\to\underset{T(-)[1]}{\hocolim} 
\bbK(\inn{F_{b,h}^{{(\bbD|_{\cF})}^{\heartsuit}}{\langle\cF\rangle}_{\nullclass,\frob}}{F_{b,h}^{{\bbD}^{\heartsuit}}}\bC).
\end{equation} 
For simplicity, 
we set $E:=\bbK(\inn{F_{b,h}^{{(\bbD|_{\cF})}^{\heartsuit}}{\langle\cF\rangle}_{\nullclass,\frob}}{F_{b,h}^{{\bbD}^{\heartsuit}}}\bC)$. 
Then by \ref{cor:degree shift is id}, 
$\displaystyle{\underset{T(-)[1]}{\hocolim}\ E}$ 
is homotopy equivalent to $E$. 
We denote the composition of $E_{\bbD,\cF}$ 
with a non-canonical homotopy 
equivalence 
between these two spectra 
by the same letters $E_{\bbD,\cF}$ and 
call it {\it relative Euler characteristic morphism} too. 

For each integer $i$, we set $\cC_i:=\cC$ and we regard $T\colon \cC\to \cC$ 
as the functor $\cC_i\to \cC_{i+1}$ for each integer $i$. 
Then it induces a relative complicial exact functor 
$\displaystyle{T\colon \bC^{\vee \bbZ}\to \bC^{\vee \bbZ}}$ 
where we denote the pair 
$\displaystyle{\left (\bigvee_{i\in\bbZ}\cC_i,lw \right )}$ 
by $\bC^{\vee \bbZ}$. 
For each integer $n$, we write $\calD_{[-\infty,n]}$ 
for the full subcategory of $\displaystyle{\bigvee_{i\in\bbZ}\calD_i}$ 
consisting of 
those objects ${(x_i)}_{i\in\bbZ} $ such that $x_i=0$ for $i>n$. 
There is a relative complicial exact functor 
$T\colon\calD_{[-\infty,n]}\to\calD_{[-\infty,n+1]}$ which 
sends an object ${(x_i)}_{i\in\bbZ} $ to ${(Tx_i)}_{i\in\bbZ}$ 
where we regard $Tx_i$ as an object $\calD_{i+1}$ 
for each integer $i$. 
Thus we can regard the family $\{\calD_{[-\infty,n]} \}_{n\in\bbZ}$ 
as a $T$-system of $\bC^{\vee\bbZ}$. 
Then we write $\bbK(\inn{\calD_{[-\infty,\infty]}|_{\cF}}{\bC^{\vee\bbZ}})$ 
for 
$\bbK(\inn{\{\calD_{[-\infty,n]} \}_{n\in\bbZ}}{\bC^{\vee \bbZ}})$ 
and for each integer $n$, 
the functors $A_{[a,b]}$ and $B_{[a,b]}$ in \ref{df:C and D df} induces the functors 
$A_{[-\infty,n]}\colon\calD_{[-\infty,n]}|_{\cF}\to{\left( F_{b,h,\leq n}^{{\bbD|_{\cF}}^{\heartsuit}}{\langle\cF \rangle}_{\nullclass,\frob}\right)}^{w_{\st}} $ and 
$B_{[-\infty,n]}\colon\calD_{[-\infty,n]}|_{\cF}\to F_{b,h,\leq n}^{{\bbD|_{\cF}}^{\heartsuit}}{\langle\cF \rangle}_{\nullclass,\frob}$. 
We define $\nabla_{[-\infty,n]}
\colon\calD_{[-\infty,n]}|_{\cF}\to \calD_n|_{\cF}$ 
to be a functor by sending an object 
${(x_i)}_{i\leq n}$ in $\calD_{[-\infty,n]}|_{\cF}$. 
Then they induce a triple of morphisms of spectra 
\begin{equation}
\label{eq:K(A[-infty,n])}
\bbK(A_{[-\infty,\infty]})\colon
\bbK(\inn{\calD_{[-\infty,\infty]}|_{\cF}}{\bC^{\vee\bbZ}})\to 
\underset{T(-)[1]}{\hocolim}\bbK(\inn{{(F_{b,h}^{{(\bbD|_{\cF})}^{\heartsuit}}{\langle\cF \rangle}_{\nullclass,\frob})}^{w_{\st}}}{(F_{b,h}^{\bbD^{\heartsuit}}\cC,lw)}),
\end{equation}
\begin{equation}
\label{eq:K(B[-infty,n])}
\bbK(B_{[-\infty,\infty]})\colon
\bbK(\inn{\calD_{[-\infty,\infty]}|_{\cF}}{\bC^{\vee\bbZ}})\to 
\underset{T(-)[1]}{\hocolim}\bbK(\inn{F_{b,h}^{{(\bbD|_{\cF})}^{\heartsuit}}{\langle\cF \rangle}_{\nullclass,\frob} }{(F_{b,h}^{\bbD^{\heartsuit}}\cC,lw)})\ \ \ \text{and}
\end{equation}
\begin{equation}
\label{eq:K(nabla)}
\bbK(\nabla_{[-\infty,\infty]})\colon
\bbK(\inn{\calD_{[-\infty,\infty]}|_{\cF}}{\bC^{\vee\bbZ}})\to
\bbK(\inn{\bbD^{\heartsuit}|_{\cF}}{\bC|_{\cF}})
\end{equation}
respectively. 
\end{para}

We can show the following lemma.

\begin{para}
\label{lem:commutativity of K-theory of devissage spaces}
{\bf Lemma.}\ \ 
{\it
In the convention \ref{df:Euler character in devissage space}, 
the following diagrams of spectra are commutative.
\begin{equation}
\label{eq:Euler character diagram}
\xymatrix{
\bbK(\inn{\bbD|_{\cF}}{\bC|_{\cF}})
\ar[r]^{\!\!\!\!\!\!\!\!\!\!\!\!\!\!\!\!\!\!\!\!\!\!\!\!\!\!\!\!\!\!\!\!\!\!\!\!\!\!\!\!\!\!\!\!\!\!\!\!\!\!\!\!\!\!\!\!\!\!\!\!E_{\bbD,\cF}} \ar[rd]_{E_{\bbD,\cF}} & 
\underset{T(-)[1]}{\hocolim} 
\bbK(\inn{F_{b,h}^{{(\bbD|_{\cF})}^{\heartsuit}}{\langle\cF\rangle}_{\nullclass,\frob}}{F_{b,h}^{{\bbD}^{\heartsuit}}}\bC)
\ar[d]^{{(-)}_{\infty}}\\
& \underset{T}{\hocolim}\bbK(\bC|_{\cF}),
}
\end{equation}
\begin{equation}
\label{eq:K(B)K(nabla)diagram}
\xymatrix{
\bbK(\inn{\calD_{[-\infty,\infty]}|_{\cF}}{\bC^{\vee\bbZ}})
\ar[r]^{\!\!\!\!\!\!\!\!\!\!\!\!\!\!\!\!\!\!\!\!\!\!\!\!\!\!\!\!\!\!\!\!\!\!\!\!\!\!\!\!\!\!\!\!\!\!\!\!\!\!\!\!\!\!\!\!\!\!\!\!\bbK(B_{[-\infty,n]})} \ar[d]_{\bbK(\nabla_{[-\infty,n]})} &
\underset{T(-)[1]}{\hocolim}\bbK(\inn{F_{b,h}^{{(\bbD|_{\cF})}^{\heartsuit}}{\langle\cF \rangle}_{\nullclass,\frob} }{(F_{b,h}^{\bbD^{\heartsuit}}\cC,lw)})
\ar[d]\\
\bbK(\inn{\bbD|_{\cF}}{\bC|_{\cF}})
\ar[r]_{\!\!\!\!\!\!\!\!\!\!\!\!\!\!\!\!\!\!\!\!\!\!\!\!\!\!\!\!\!\!\!\!\!\!\!\!\!\!\!\!\!\!\!\!\!\!\!\!\!\!\!\!\!\!\!\!\!\!\!\!E_{\bbD,\cF}} &
\underset{T(-)[1]}{\hocolim} 
\bbK(\inn{F_{b,h}^{{(\bbD|_{\cF})}^{\heartsuit}}{\langle\cF\rangle}_{\nullclass,\frob}}{F_{b,h}^{{\bbD}^{\heartsuit}}}\bC).
}
\end{equation}
}\qed
\end{para}

\subsection{D\'evissage theorem}
\label{subsec:dev thm}

The goal of this subsection is to give d\'evissage theorems 
for relative complicial exact categories. 
The specific feature in our proof is `motivic' in the sense that 
the properties of $K$-theory we will utilize to prove the theorems 
are only 
localization (additivity and derived invariance), 
categorical homotopy invariance and cocontinuity. 
Namely our proof works to any other localizing theory 
on relative complicial exact categories 
(for the precise statement, see Corollary~\ref{cor:devissage for localization theory}). 
Recall the definition of the Euler characteristic morphisms 
from \ref{df:Euler character in devissage space}.

\begin{para}
\label{thm:divissage for complicial exact categories}
{\bf Theorem ($\bbK$-theory of D\'evissage spaces).}\ \ 
{\it
Let $\bC=(\cC,w)$ be a normal ordinary 
relative complicial exact category and let 
$\bbD$ be a cell structure of $\bC$ and 
let $\cF$ be a full subcategory of $\cC$. 
Then the Euler characteristic morphism 
$E_{\bbD,\cF}\colon
\bbK(\inn{\bbD^{\heartsuit}|_{\cF}}{\bC|_{\cF}})\to
\underset{T(-)[1]}{\hocolim}\bbK(\inn{F_{b,h}^{{(\bbD|_{\cF})}^{\heartsuit}}{\langle\cF \rangle}_{\nullclass,\frob}}F_{b,h}^{{\bbD}^{\heartsuit}}\bC)$ 
is a homotopy equivalence of spectra. 
In particular if $\cF$ satisfies relative derived d\'evissage condition with respect to the cell structure $\bbD$ in $\bC$, then 
$E_{\bbD,\cF}\colon \bbK(\inn{\bbD|_{\cF}}\bC|_{\cF})\to \bbK(\inn{\cF}\bC)$ is a homotopy equivalence of spectra. 
}
\end{para}

\begin{proof}
By commutative diagrams $\mathrm{(\ref{eq:ABE commutative diagram})}$ and 
$\mathrm{(\ref{eq:K(B)K(nabla)diagram})}$, 
there exists a commutative diagram of fibration sequences of spectra
$$
\xymatrix{
\bbK(\inn{\calD_{[-\infty,\infty]}|_{\cF}}\bC^{\vee \bbZ}) 
\ar[r]^{\!\!\!\!\!\!\!\!\!\!\!\!\!\!\!\!\!\!\!\!\!\!\!\!\!\!\!\!\!\!\!\!\!\!\!\!\!\!\!\!\!\!\!\!\!\!\!\!\!\!\!\bbK(B_{[-\infty,\infty]})} 
\ar[d]_{\bbK(E_{[-\infty,\infty]})} & 
\underset{T(-)[1]}{\hocolim}\bbK(\inn{{(F_{b,h}^{{(\bbD|_{\cF})}^{\heartsuit}}{\langle\cF \rangle}_{\nullclass,\frob})}^{w_{\st}}}(F_{b,h}^{\bbD^{\heartsuit}}\cC,lw)) 
\ar[d]\\
\bbK(\inn{\calD_{[-\infty,\infty]}|_{\cF}}\bC^{\vee \bbZ}) 
\ar[r]^{\!\!\!\!\!\!\!\!\!\!\!\!\!\!\!\!\!\!\!\!\!\!\!\!\!\!\!\!\!\!\!\!\!\!\!\!\!\!\!\!\!\!\!\!\!\!\!\!\!\!\!\bbK(A_{[-\infty,\infty]})} 
\ar[d]_{\bbK(\nabla_{[-\infty,\infty]})} & 
\underset{T(-)[1]}{\hocolim}\bbK(\inn{{(F_{b,h}^{{(\bbD|_{\cF})}^{\heartsuit}}{\langle\cF \rangle}_{\nullclass,\frob})}}(F_{b,h}^{\bbD^{\heartsuit}}\cC,lw)) \ar[d]\\
\bbK(\inn{\bbD|_{\cF}}\bC|_{\cF}) 
\ar[r]_{\!\!\!\!\!\!\!\!\!\!\!\!\!\!\!\!\!\!\!\!\!\!\!\!\!\!\!\!\!\!\!\!\!\!\!\!\!\!\!\!\!\!\!\!\!\!\!\!\!\!\!E_{\bbD,\cF}} & 
\underset{T(-)[1]}{\hocolim}\bbK(\inn{{(F_{b,h}^{{(\bbD|_{\cF})}^{\heartsuit}}{\langle\cF \rangle}_{\nullclass,\frob})}}F_{b,h}^{\bbD^{\heartsuit}}\bC).
}
$$
By Corollary~\ref{cor:homotopy invariance of C and D}, 
$\bbK(A_{[-\infty,\infty]})$ and $\bbK(B_{[-\infty,\infty]})$ 
are homotopy equivalence of spectra. 
Thus we obtain the result.
\end{proof}

By Proposition~\ref{prop:sufficient conditions of derived devissage} and 
Theorem~\ref{thm:divissage for complicial exact categories}, 
we obtain the following 
Corollary~\ref{cor:Waldhausen's cell filtration theorem} and 
Corollary~\ref{cor:Homotopy devissage theorem}.

\begin{para}
\label{cor:Waldhausen's cell filtration theorem}
{\bf Corollary (Waldhausen's cell filtration theorem for non-connective $K$-theory).}\ 
(\cf \cite[1.7.1]{Wal85}.)\ \ 
{\it
Let $\bC=(\cC,w)$ be a thick normal ordinary 
relative complicial exact category and 
let $\bbD$ be a cell structure of $\bC$. 
Assume that $\bbD$ satisfies Waldhausen's d\'evissage condition. 
Then the Euler characteristic morphism 
$E_{\bbD}\colon\bbK(\inn{\bbD^{\heartsuit}}{\bC})\to\bbK(\bC)$ 
is a homotopy equivalence of spectra. 
}
\qed
\end{para}

\begin{para}
\label{cor:Homotopy devissage theorem}
{\bf Corollary (Homotopy d\'evissage theorem).}\ \ 
{\it
Let $\bC=(\cC,w)$ be a thick normal ordinary 
relative complicial exact category and 
let $\cA$ be a full subcategory of $\bC$. 
Assume that $\cA$ satisfies homotopy d\'evissage condition in $\bC$. 
Then the Euler characteristic morphism 
$E_{\bbD}\colon\bbK(\inn{{\langle\cA\rangle}_{\nullclass,w}}{\bC})\to\bbK(\bC)$ 
is a homotopy equivalence of spectra. 
}\qed
\end{para}

By Proposition~\ref{prop:sufficient conditions of derived devissage}, 
Lemma~\ref{lem:decomposition cond and devissage} and 
Corollary~\ref{cor:Waldhausen's cell filtration theorem}, 
we obtain the following. 

\begin{para}
\label{cor:Theorem of heart}
{\bf Corollary (Theorem of heart).}\ \ 
{\it
Let $\bC=(\cC,w)$ be a normal ordinary relative complicial exact 
category and let $\cW$ be a bounded below homological weight structure 
on $\Ho(\bC)$. 
Then the Euler characteristic morphism 
$E_{\bbD^{\cW}}\colon\bbK(\inn{{(\bbD^{\cW})}^{\heartsuit}}{\bC})\to 
\bbK(\bC)$ 
is a homotopy equivalence of spectra. 
}
\qed
\end{para}

By Lemma~\ref{lem:Raptis and homotopy devissage} and 
Corollary~\ref{cor:Homotopy devissage theorem}, we obtain the following.

\begin{para}
\label{cor:Raptis' devissage theorem}
{\bf Corollary (Raptis' d\'evissage theorem for non-connective $K$-theory).}\ \ 
{\it
Let $\bC=(\cC,w)$ be a thick normal ordinary 
relative complicial exact category and 
let $\cA$ be a full subcategory of $\bC$. 
Assume that $\cA$ satisfies Raptis' d\'evissage condition in $\bC$. 
Then the inclusion functor $\cA\rinc\cC$ induces a 
homotopy equivalence of spectra 
$\bbK(\inn{{\langle\cA\rangle}_{\nullclass,w}}{\bC})\to\bbK(\bC)$ on $\bbK$-theory. 
}
\qed
\end{para}

\begin{para}
\label{cor:Quillen's devissage theorem}
{\bf Corollary (Quillen's d\'evissage theorem 
for non-connective $K$-theory).}\ \ 
{\it
Let $\cA$ be an abelian category and 
let $\cB$ be a topologizing subcategory of $\cA$. 
Assume that $\cB$ satisfies Quillen's d\'evissage condition in $\cA$. 
Then the inclusion functor $\cB\rinc \cA$ induces a homotopy equivalence 
of spectra $\bbK(\cB)\to\bbK(\cA)$ on $\bbK$-theory. 
}
\end{para}

\begin{proof}
By Lemma~\ref{lem:Raptis characterization}, 
$\Ch_b(\cB)$ satisfies Raptis' d\'evissage condition in 
$\Ch_b\cA$. 
Thus by Corollary~\ref{cor:Raptis' devissage theorem}, 
the inclusion functor 
$\Ch_b\cB\rinc\Ch_b\cA$ 
induces a homotopy equivalence of spectra 
$$\bbK(\inn{{\langle\Ch_b\cB\rangle}_{\nullclass,\qis}}{(\Ch_b\cA,\qis)})\to\bbK(\Ch_b\cA;\qis)$$
on $\bbK$-theory. 
Thus by Corollary~\ref{cor:relative Gillet-Waldhausen theorem} and 
Corollary~\ref{cor:w-envelope invariance II}, 
we obtain the result. 
\end{proof}

\begin{para}
\label{df:Localization theory} 
{\bf Definition (Localization theory).}\ \ 
We define the full subcategory of $\RelComp^{[2]}$ consisting 
of those objects $x$ such that the sequence $x(0)\to x(1)\to x(2)$ 
are derived weakly exact sequences 
(\ref{df:derived exaxt seq}) 
of relative complicial 
exact categories and relative complicial exact functors by 
$\mathfrak{EX}(\RelComp)$. 
The functor $[0]\to [2]$, $0\mapsto k$ induces functor 
$q_k\colon\mathfrak{EX}(\RelComp)\to\RelComp$ for $0\leq k\leq 2$. 

Let $(\cT,\Sigma)$ be a triangulated category. 
A localization theory on $\RelComp$ to $\cT$ is a pair 
$(\calL,\partial)$ consisting of functor 
$\calL\colon\mathfrak{EX}(\RelComp)\to\cT$ and a natural transformation 
$\partial \colon \calL q_2\to\Sigma\calL q_0$ such that 
they satisfies the following conditions:
\begin{itemize}
\item[]
{\bf (Categorical homotopy invariance).}\ \ 
Let $f$, $f'\colon \bC \to \bC'$ be a pair of relative complicial 
exact functors between relative complicial exact categories 
$\bC$ and $\bC'$. 
If $f$ and $f'$ are categorical homotopic in $\underline{\RelComp}$, 
then they induce same maps 
$\calL(f)=\calL(f')\colon\calL(\bC)\to\calL(\bC')$. 

\item[]
{\bf (Localization).}\ \ 
Let $\mathrm{(\ref{eq:derived spltting sequence})}$ 
be a derived weakly exact sequence (\ref{df:derived exaxt seq}) of 
relative complicial exact categories and 
relative complicial exact functors. 
Then the induced sequence 
$$\calL(\bE)\onto{\calL(i)}\calL(\bE')\onto{\calL(p)}\calL(\bE'')$$
is a distinguished triangle in $\cT$. 

\item[]
{\bf (Cocontinuity).}\ \ 
$\calL$ sends filtered colimits to homotopy colimits in $\cT$. 
\end{itemize}
\end{para}

As alluded at the outset of this subsection, 
we obtain the following result. 

\begin{para}
\label{cor:devissage for localization theory}
{\bf Corollary (D\'evissage theorem for localization theory).}\ \ 
{\it
For a localization theory 
$(\calL,\partial)$ from $\RelComp$ 
to $(\cT,\Sigma)$, 
similar statements of 
\ref{thm:approximation}, 
\ref{cor:relative Gillet-Waldhausen theorem}, 
\ref{cor:w-envelope invariance}, 
\ref{cor:w-envelope invariance II}, 
\ref{cor:purity theorem}, 
\ref{thm:divissage for complicial exact categories}, 
\ref{cor:Waldhausen's cell filtration theorem}, 
\ref{cor:Homotopy devissage theorem}, 
\ref{cor:Theorem of heart}, 
\ref{cor:Raptis' devissage theorem} and 
\ref{cor:Quillen's devissage theorem} 
hold.
}
\qed
\end{para}

\appendix

\section{A d\'evissage theorem for modular Waldhausen exact categories}
\label{subsec:dev modular}

In this section, 
we recall Quillen's original d\'evissage theorem in \cite{Qui73} 
with an adequate generalization. 
Throughout this section, let $\cE$ be a small exact category. 

\begin{para}
\label{df:modular exact category}
{\bf (Modular exact categories, modular Waldhausen categories).}\ \ 
We say that $\cE$ is {\it modular} 
if it satisfies the following two conditions:\\ 
${\bf (Mod 1).}$ 
For any admissible monomorphisms 
$x \rinf z$, $y \rinf z$ in $\cE$, 
there exists a fiber product $x \times_z y$ and 
the canonical morphisms $x \times_z y \to x$, 
$x \times_z y \to y$ and 
$\displaystyle x \underset{x \times_z y}{\sqcup} y \to z$ are 
admissible monomorphism.

\sn
${\bf (Mod 2).}$ 
In the following commutative diagram of admissible monomorphisms,
$$\xymatrix{
x \ar@{>->}[r] \ar@{>->}[d]& y \ar@{>->}[d]\\
z \ar@{>->}[r] & w
}$$
the canonical induced morphism $x \to z \times_w y$ is also 
an admissible monomorphisms.

\sn
For a modular exact category $\cE$ and admissible monomorphisms 
$x\rinf z$ and $y\rinf z$, we set 
$x\cap_z y:=x\times_z y$ and $x\cup_z y:=x\underset{x\cap_z y}{\sqcup} y$. 
For simplicity we often write $x\cap  y$ and $x\cup y$ for 
$x\cap_z y$ and $x\cup_z y$ respectively. 
Then by \cite[p.406 step 1]{Kel90}, 
the commutative diagram of admissible monomorphisms 
$$ 
\xymatrix{
x\cap y \ar@{>->}[r] \ar@{>->}[d] & x \ar@{>->}[d]\\
y \ar@{>->}[r] & x\cup y.
}$$
is a biCartesian square. 
Hence we have an admissible exact sequence
\begin{equation}
\label{eq:lattice exact}
x\cap y \rinf x\oplus y \rdef x\cup y.
\end{equation}

Moreover for any sequence of admissible monomorphisms 
$x\rinf z \rinf a \linf u \linf y$ in $\cE$, 
we can produce the commutative diagram below: 
$$
\xymatrix{
x\cap y \ar@{>->}[r] \ar@{>->}[d] & x \oplus y \ar@{->>}[r] \ar@{>->}[d] & x \cup y\ar@{>->}[d]\\
z\cap w \ar@{>->}[r] \ar@{->>}[d] & z \oplus u \ar@{->>}[r] \ar@{->>}[d] & z \cup w \ar@{->>}[d]\\
z\cap w/x\cap y \ar[r] & z/x\oplus u/y \ar[r] & z\cup w/x\cup y.
}
$$
Here by the sequence $\mathrm{(\ref{eq:lattice exact})}$, 
the first and the second horizontal lines are admissible short exact sequences 
and we can check that the left, the middle and the right vertical lines are also 
admissible short exact sequences. 
Thus by $3\times 3$-lemma, the third horizontal line 
\begin{equation}
\label{eq:modular exact seq}
z\cap w/x\cap y \rinf  z/x\oplus u/y \rdef z\cup w/x\cup y
\end{equation}
is also an admissible short exact sequence.

\sn
We say that a pair $(\cE,w)$ is a {\it modular Waldhausen exact category} 
if the pair $(\cE,w)$ and 
the pair of the opposite categories $(\cE^{\op},w^{\op})$
are Waldhausen categories and $\cE$ is modular 
and moreover $w$ satisfies the following {\it strong cogluing axiom}:

\sn
In the commutative diagram below
$$\xymatrix{
x \ar@{>->}[r] \ar[d]_a & y \ar[d]_b & z \ar@{>->}[l] \ar[d]^c\\
x' \ar@{>->}[r] & y' & z' \ar@{>->}[l],
}$$
if $a$, $b$ and $c$ are in $w$, 
then the induced canonical morphism 
$a\times_b c\colon x \times_y z \to x' \times_{y'} z'$ 
is also in $w$.
\end{para}

In the rest of this section, 
we assume that the pair $(\cE,w)$ is a modular Waldhausen exact category.

\begin{para}
\label{rem:modularity}
{\bf Remark.}\ \ 
To explain where the naming of modularity comes from, 
we need to prepare some notations. 
Let $L$ be a lattice, namely a partially ordered set closed under 
$\sup$ and $\inf$ of any finite sets. 
For any elements $a$ and $b$ in $L$, 
we write $a\cup b$ and $a\cap b$ for $\sup\{a,b\}$ and $\inf\{a,b\}$ 
respectively. 
Recall that a lattice $L$ is {\it modular} 
if for any triple of elements $a\leq b$ and $c$ in $L$, the following 
equality called {\it modular law} holds:
\begin{equation}
\label{eq:modular law}
a\cup (b\cap c)= (a\cup b) \cap c.
\end{equation}
For any object $e$, we denote the isomorphism class of admissible 
subobjects of $e$ by $\cP(e)$. 
Then we will show that $\cP(e)$ is a modular lattice. 
Indeed, for any triple of admissible 
subobjects $x\rinf y$ and $z$ of $e$, 
by applying the admissible short exact sequence 
$\mathrm{(\ref{eq:modular exact seq})}$ to 
$x$, $y$, $z$ and $u=y$, 
it turns out that 
the conditions $x\cap y=z\cap y$ and $x\cup y=z\cup y$ imply $x\isoto z$. 
Thus by Dedekind's criterion of modularity 
(see for example \cite[1.5]{MY14}), 
$\cP(e)$ is a modular lattice.
\end{para}

\begin{para}
\label{df:E(n,w)}
For any non-negative integer $n$, 
we define $\cE(n,w)$ to be the full subcategory of 
$\cE^{[1]}$ the functor category 
from the ordinal $[n]=\{0,1,\cdots,n\}$ to $\cE$ consisting of those objects 
$x$ such that for any pair $i\leq j$ of elements of $[n]$, 
the morphism $x(i\leq j)$ from $x(i)$ to $x(j)$ is in $w$. 

We natural make $\cE(n,w)$ into an exact category 
where a sequence $x\onto{f}y\onto{g}z$ of morphisms in $\cE(n,w)$ is 
an admissible exact sequence if 
$x(i)\onto{f(i)}y(i)\onto{g(i)}z(i)$ 
is an admissible exact sequence 
in $\cE$ for any $i\in[n]$. 

The association $[n] \mapsto \cE(n,w)$, gives 
a simplicial exact category(= a simplicial object 
in the category of exact categories and exact functors). 
\end{para}

\begin{para}
\label{lem:her of modularity}
{\bf Lemma.}\ \ 
{\it
Let $(\cE,w)$ be a modular Waldhausen exact category. 
Then for any non-negative integer $n$, 
the exact category $\cE(n,w)$ is a modular exact category.
}
\end{para}

\begin{proof}
Recall that an admissible monomorphism in 
$\cE(n,w)$ is a term-wise admissible monomorphisms in $\cE$ 
(see \ref{df:E(n,w)}). 
Only non-trivial assertion is existence of a fiber product 
$x\times_z y$ for any pair of admissible monomorphisms 
$x\rinf z$ and $y\rinf z$ in 
$\cE(n,w)$. 
Since $\cE$ is modular it exists term-wisely, and 
therefore it exists in $\cE^{[n]}$ the functor category 
from $[n]$ to $\cE$. 
By virtue of strong cogluing axiom, 
it is actually in $\cE(n,w)$. 
\end{proof}

\begin{para}
\label{df:Devissage conditions} 
{\bf (Relative d\'evissage condition).}\ \ 
Let $\calD$ be a topologizing subcategory of $\cE$ 
(see \ref{df:toplogizing subcat}). 
Then $\calD$ naturally becomes a strict exact subcategory of 
$\cE$ 
(see the definition of strict exact subcategories 
for \ref{nt:exact categories} and 
see \cite[5.3]{Moc13a} for the proof of this fact). 

For a non-negative integer $n$, 
we can check that $\calD(n,w)$ is 
again a topologizing subcategory of $\cE(n,w)$ 
by utilizing strong gluing axiom. 
For simplicity we also write the same letter $w$ for $w\cap\Mor\calD$. 
We say that the inclusion functor 
$(\calD,w)\rinc (\cE,w)$ satisfies the {\it relative d\'evissage condition} 
if for any non-negative integer $n$, 
the inclusion functor $\calD(n,w)\rinc \cE(n,w)$ satisfies 
the d\'evissage condition 
(see \ref{df:classical Devissage condition}).
\end{para}

\begin{para}
\label{thm:devissage}
{\bf Theorem (D\'evissage).}\ \ 
{\it
If the inclusion functor 
$(\calD,w) \rinc (\cE,w)$ satisfies the relative d\'evissage condition, 
then it induces a homotopy equivalence 
$wS_{\cdot}\calD\to wS_{\cdot}\cE$ on $K$-theory.
}
\end{para}

\begin{proof}
We fix a non-negative integer $n$ and 
we apply Quillen's d\'evissage theorem in \cite{Qui73} 
to the inclusion functor 
$\calD(n,w)\rinc \cE(n,w)$. 
In \cite{Qui73}, we need to assume that $\cE(n,w)$ is abelian. 
But this hypothesis is only using to prove the fact that 
for any object $a$ in $\cE(n,w)$, the isomorphism class of 
subobjects of $a$ becomes a (modular) lattice and the fact that 
for any sequence of admissible monomorphisms 
$z\rinf x\rinf a \linf y \linf u$ in $\cE(n,w)$, 
there exists an admissible short exact sequence 
$\mathrm{(\ref{eq:modular exact seq})}$ 
in $\cE(n,w)$. 
These facts are proven in \ref{df:modular exact category}, 
\ref{lem:her of modularity} and \ref{df:Devissage conditions}. 
Thus the proof in \cite{Qui73} works fine for our situation. 
Hence it turns out that the inclusion functor 
$Q\calD(n,w)\rinc Q\cE(n,w)$ 
is a homotopy equivalence. 
Here the letter $Q$ means the Quillen's $Q$-construction. 
Then by virtue of Waldhausen's $Q=s$ theorem in \cite[1.9]{Wal85}, 
the inclusion functor $\calD\rinc \cE$ induces 
a homotopy equivalence 
$wS_n\calD \to wS_n\cE$
for any non-negative integer $n$. 
Finally by 
the realization lemma in \cite[Appendix A]{Seg74} or \cite[5.1]{Wal78}, 
we obtain the result.
\end{proof}

\sn
{\bf Acknowledgement.}\ \ 
The author wishes to express his deep gratitude to Masana Harada, 
Kei Hagihara, Sho Saito, Seidai Yasuda, Toshiro Hiranouchi and 
Gon{\c{c}}alo Tabuada for useful conversations in the early stage of this work.

\mn
SATOSHI MOCHIZUKI\\
{\it{DEPARTMENT OF MATHEMATICS,
CHUO UNIVERSITY,
BUNKYO-KU, TOKYO, JAPAN.}}\\
e-mail: {\tt{mochi@gug.math.chuo-u.ac.jp}}\\

\end{document}